# Explanation of *The Dynamics* involved in The 3N+1 Problem
## A proof for *The Collatz Conjecture*

By Mr. R. Bruun (BSc), Copenhagen, September 1st, 2025
**e-mail: r.collatz.bruun@gmail.com**

The scope of the present work is to explain **why** it is true that all
N ∊ ℕ have a distinct **position** in the below shown ***Collatz Tree***.

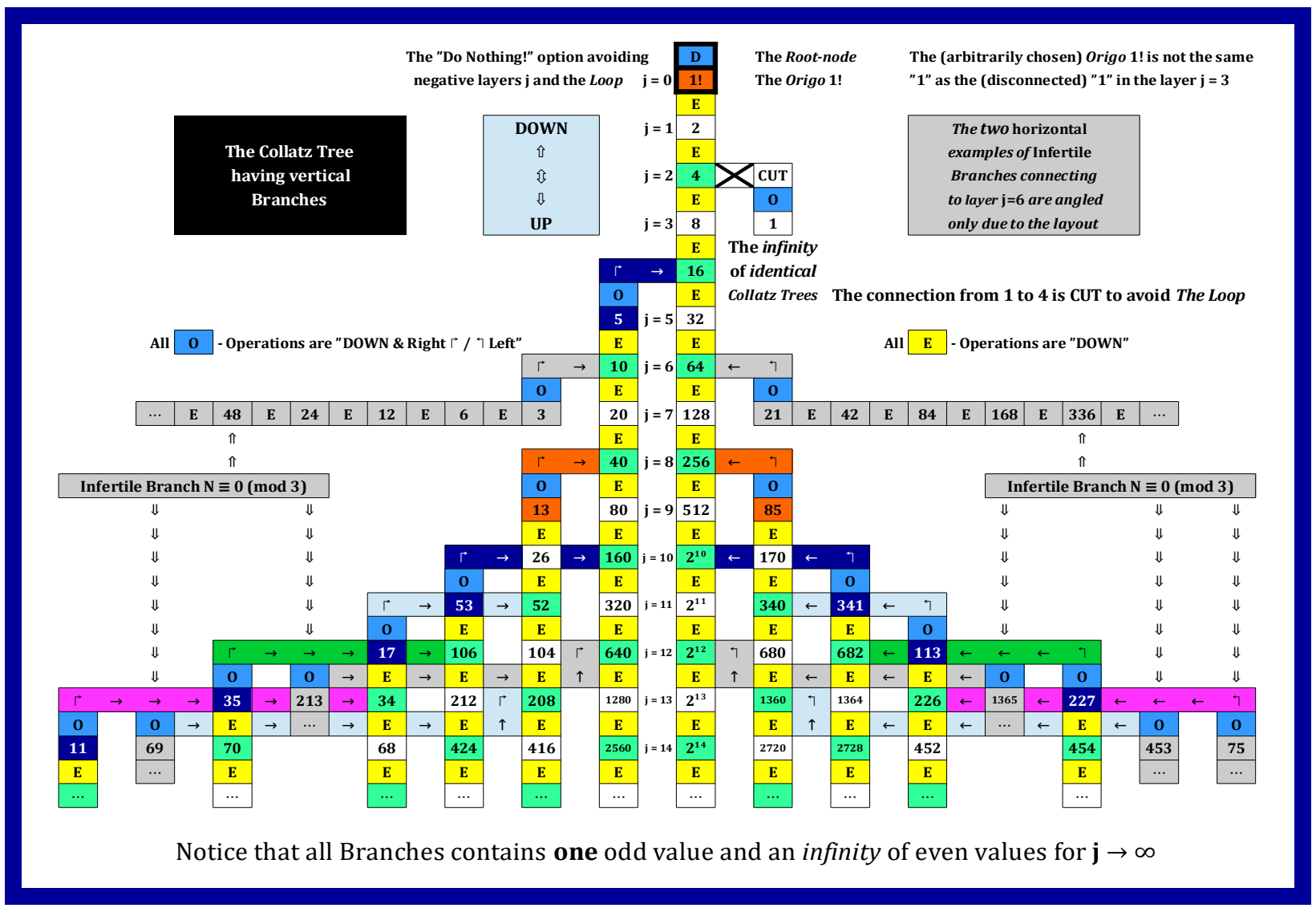

**Figure 1a.** ***The Collatz Tree*** extrapolated from ***The Origo*** **1!**

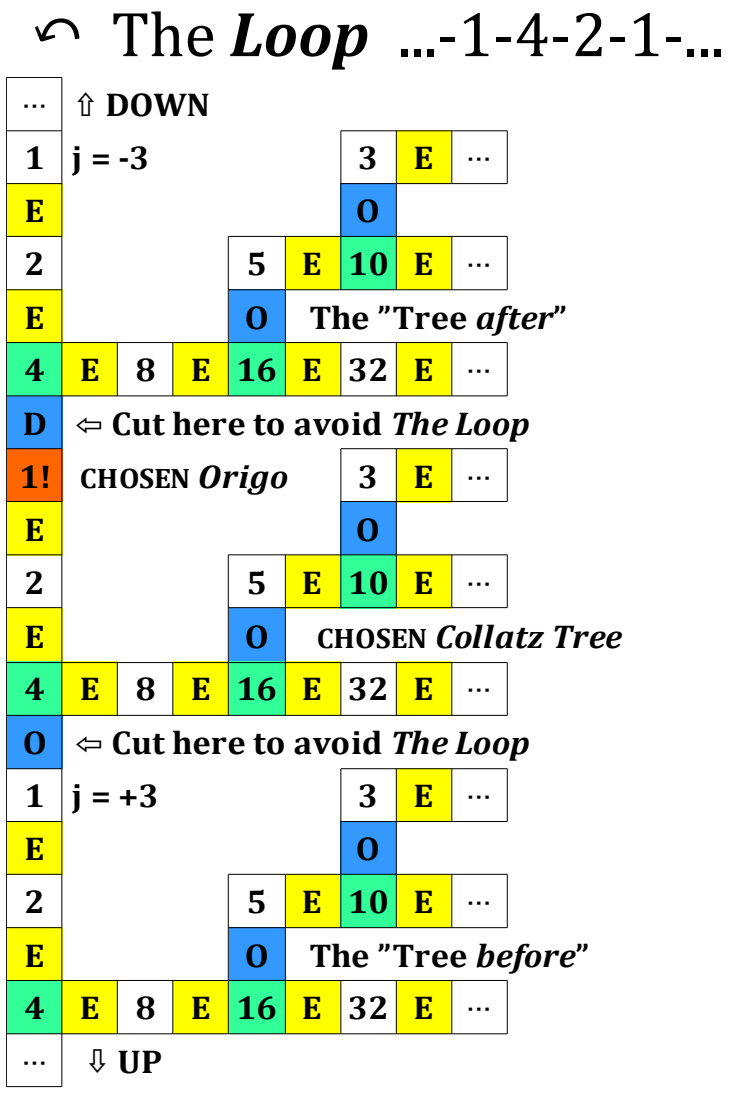

**Figure 1b.** *Other* **Trees.**

# 1. Introduction

**The 3N+1 Problem** concerns the behaviour of natural numbers N ∊ ℕ using these rules:

- If the the number N is **O**dd multiply N by **3** and add **1** (Type **Odd** operation)
- If the the number N is **E**ven divide N by **2** (Type **Even** operation)
- Repeat with the new value (i.e. with the **Result** of the iteration/operation)

For example for $N_0$ = 3 as the **I**nitial **V**alue: $3_0$ ↗ $10_1$ ↘ $5_2$ ↗ $16_3$ ↘ $8_4$ ↘ $4_5$ ↘ $*2_6$ ↘ $^\dagger 1_7$
The above number-sequence is informally called an iterationseries (for the **I**nitial **V**alue $N_0$)
while the *order of operations* **Odd** & **Even** is called the parity-sequence for the iterationseries.

***The Collatz Conjecture*** states, that **eventually** any iterationseries will reach the ***Loop***-value $^\dagger$1
(alternatively that iterationseries for all $N_0$ except {1,2,4,8} reaches the "*bottleneck-value*" 16).

**Definition 1**

- Let **j** ∊ $ℕ_0$ count the Total number of operations in an iterationseries (**j** is the *index* in $N_j$)
- Let **s** ∊ $ℕ_0$ count the number of type **Odd** operations in an iterationseries
- Let **r** ∊ $ℕ_0$ count the number of type **Even** operations in an iterationseries

It is observed that **s+r = j** i.e. **#Odd(s) + #Even(r) = #Total(j)** (where # is number of pcs.)

All the iterationseries which reach the **End** Value $N_j$ = $^\dagger 1_j$ are called **Convergent** iterationseries
and it is possible to identify the iterationseries with the corresponding **End** $^\dagger$parity-sequence.

All iterationseries which reach a **Reduced Value* $*N_j < N_0$ are **Converging* iterationseries and it
is possible to identify the iterationseries with the corresponding **reducing* **parity-sequence*.
**NB**! Everywhere in the present work a star/asterisk * represents that something is **reducing*.

All iterationseries that reach $N_j = {}^{\dagger}1_j$ *continues* in ***The Loop*** | ⋯ $1_j$ ↗ $4_{j+1}$ ↘ $2_{j+2}$ ↘ $1_{j+3}$ ⋯

Parity-sequences can be constructed from **O E-tiles** & **E-tiles** like Domino-tiles, because for every odd N = (2C-1), C ∊ ℕ, the **Result** after a type Odd operation is even 3(2C-1)+1 = (6C-2), so the next operation is a type Even (6C-2)/2 = (3C-1). The operations together is considered a **unit** named an **O E-tile** and the *other* **unit** for type Even operations following other type Even operations is named an **E-tile. OBS**! The parity of (3C-1) ≡**?**(mod 2) depends on C ∊ ℕ.

In the above example it is observed, that the **I**nitial **V**alue $3_0$ reach the **End** Value $^{\dagger}1_7$ in seven operations of which **s** = 2 are type Odd and **r** = 5 are type Even. The **End** $^{\dagger}$parity-sequence is:

$^{\dagger}$**R O E - O E - E - E - E - D** where **R** indicates a sta**R**ting point and **D** the En**D** point.

It is also observed, that the **I**nitial **V**alue $3_0$ reach the **Reduced Value* $^*2_6$ in six operations of which **s** = 2 are type Odd and ***r**$_s$ = $^*4_2$ are type Even. The **reducing* **parity-sequence* is:

***R O E - O E - E - E - T** where **T** indicates the **T**otal/**T**erminating **parity-sequence*.

The **E-tile** indicates the *last* type Even operation in the **reducing* **parity-sequence*.

**Comment**

It is emphasized that the **Reduced Value* $^*N_j$ is defined to be the **first** value (the lowest **j**) for which it is true that $^*N_j < N_0$, and that the **End** Value $^{\dagger}1_j$ is defined to be the **first** value (the lowest **j**) for which it is true that $N_j = {}^{\dagger}1_j$.

**A brief history of The Problem**

The **Q**N±D Problem, N ∊ ℕ, for **Q** and D odd natural numbers have been known for a long time. The version for **Q** = 3 & D = 1 have puzzled mathematicians and hopeful amateurs for a more or less equally long time. **The 3N+1 Problem** is these days known as ***The Collatz Conjecture*** and is named after Lothar Collatz (1910-1990†) a German mathematician who - according to Wikipedia[I] - worked on The Problem in the late 1930's and introduced The Problem to an increasing number of mathematicians during the 1950's, *slowly* bringing it into the public eye.

***The Collatz Conjecture*** is easy to state, but surprisingly difficult to prove. It has gained some notoriety amongst professional mathematicians for being a "rabbit hole" that leads to a **lot** of wasted time (which this author can confirm is absolutely true) but not to **useful** results.

Jeffrey C. Lagarias have made a comprehensive study of The Problem, collecting known literature in "The 3x+1 Problem: An Annotated Bibliography (1963-1999)[II]". Some progress have been made over the years and it is well known that The Conjecture holds for "Almost All" integers; Riho Terras[III] (1976) shows that the **Asymptotic Density** for the **Diverging** parity-sequences tends to zero when **modulo** tends to infinity in the **Union** of residue-classes **and** very importantly that the **densities** (the **density** of a residue-class is defined to be the inverse of the **modulo** i.e. **1/modulo**) can be **calculated** using the ***Binomial-distribution*** (**NB**! In the present work **Formulas** are shown, which produce results ***identical*** to results from R. Terras).

**Axiom 1**

If the conjecture is true, and all iterationseries for all N ∊ ℕ as the **I**nitial **V**alue $N_0$ **eventually** do reach ***The Loop***, then it is of course also true, that the series contains a **Reduced Value* $^*N_j$. Vice versa: **If** it is true for all **iterationseries* for all N ∊ ℕ as the **I**nitial **V**alue $N_0$ that the series contains a **Reduced Value* $^*N_j$ **then** it is also true, that all series **eventually** reach ***The Loop***. If **all** $N_0 > 2$ reach a **Reduced Value* $^*N_j < N_0$ then the $^*N_j$ is the new $N_0$ in a new **iterationseries*.

**Axiom 2**

If the conjecture is false, then it is true that there exist a **lowest** counter**e**xample $N_{ce} \in \mathbb{N}$ which have an iterationseries (with a parity-sequence) that **never** reach the **End** Value $N_j = {}^{\dagger}1_j$ and in **that** case it is also true, that the iterationseries for $N_{ce}$ as the **I**nitial **V**alue $N_0 = N_{ce}$ **never** reach a **Reduced Value* $^*N_j < N_{ce}$, because if it *does*, then $N_{ce}$ is **not** the **lowest** counterexample.

**Axiom 3**

For all $N \in \mathbb{N}$, $N \neq 2^c$, it is true that there exist a **c** $\in \mathbb{N}$ for which it is true that: $2^{c-1} < N < 2^c$

**Corollary 1**

From the above **Axioms** follows that the hypothetical $N_{ce}$ have a **finite size**: $2^{r_{ce}-1} < N_{ce} < 2^{r_{ce}}$

**Introduction to the method and model**

Most arguments are based on **modular arithmetic**, **graph-theory** and **combinatorics**.

Following the rules in **The 3N+1 Problem** the present work is an analysis of the behaviour of **residue-classes** $N \equiv \mathbf{P} \pmod{2^r}$, $\mathbf{P} \leq 2^r$, in a ***Covering System*** including all $\mathbf{P} \leq 2^r$, $\mathbf{P} \in \mathbb{N}$, $\mathbf{r} \in \mathbb{N}$.

Different notations can be used to describe the same ***modular set*** $N \equiv \mathbf{P} \pmod{A}$, $\mathbf{P} \leq A$. In the present work [A**X**-B], B < A, is used to define the ***set*** of all $N_X = (A\mathbf{X}-B) \Leftrightarrow N_X \equiv (\mathbf{A-B}) \pmod{A}$. It is demonstrated that in this **model** the **modulo** A can only take on certain values $A = 2^r3^s$. **OBS**! Notice that in a **class** [A**X**-B] the *independent variable* **X** $\in \mathbb{N}$ describes the **position** of $N_X$.

It is demonstrated, that it is possible to perform Collatz-operations on **congruence-classes** (residue-classes) resulting in **Class-series** for **I**nitial **V**alue classes $N \equiv \mathbf{P}_{IV} \pmod{2^r}$, $\mathbf{P} \leq 2^r$, which have a finite number of terms and **ends** at **T**erminal **V**alue classes $N \equiv \mathbf{P}_{TV} \pmod{3^s}$. For example the *$^*$Converging* Class-series $^*[2^{*2_1}\mathbf{X}\text{-}^*3_0]_{IV}$ $\rightleftarrows_O$ $[2^23^1\mathbf{X}\text{-}8_1]$ $\rightleftarrows_E$ $[2^13^1\mathbf{X}\text{-}4_2]$ $\rightleftarrows_E$ $[3^1\mathbf{X}\text{-}2_3]_{TV}$ & the $^{\Phi}$**Diverging** Class-series $^{\Phi}[2^2\mathbf{X}\text{-}1_0]_{IV}$ $\rightleftarrows_O$ $[2^23^1\mathbf{X}\text{-}2_1]$ $\rightleftarrows_E$ $[2^13^1\mathbf{X}\text{-}1_2]$ $\rightleftarrows_O$ $[2^13^2\mathbf{X}\text{-}2_3]$ $\rightleftarrows_E$ $[3^2\mathbf{X}\text{-}1_4]_{TV}$ where ***conclusions*** can be made from the difference in the size of the **modulo** $2^r > 3^s$ **or** $2^r < 3^s$. **NB**! Everywhere in the present work the symbol $^{\Phi}$ represents that something is $^{\Phi}$**Diverging**. Class-series (for **IV**-classes) can - like iterationseries for *values* $N_0$ - be identified with a *parity-sequence* with a finite number of **tiles** as all the elements in the class *share* the same *sequence*.

It is shown, that it is possible to perform Collatz-operations on values $N \in \mathbb{N}$ in **Branches** for $N = (2C\text{-}1)2^r$ i.e. **r** > 0 implies a Type Even operation and **r** = 0 implies a Type Odd operation. A system of *flowgraphs* (containing all odd N) is constructed, showing Collatz-operations on classes $[2^r3^s\mathbf{X}\text{-}B]$ and the development for changing (increasing) **modulo**. The *flowgraphs* are used to illustrate ***parity-sequences*** (for Class-series), which are analysed using *binary* graphs.

For parity-sequences in the Infinite Complete Binary Trees - a (rooted) **Tree** is a **graph** having a **Root**-node with a unique **route** (path) to any other node - the following **Axiom** is important:

**Axiom 4**

In an **Infinite Complete Binary Tree** all nodes in a layer/level **k** $\in \mathbb{N}_0$ have two *Child*-nodes in the next layer **k+1**, so all nodes in a layer **k** $\in \mathbb{N}$ have a **sub-tree** identical to the **sub-tree** for the Root-node in layer **k** = **0**. All nodes in layer **k** are *Parent*-nodes (i.e. *to* the two Left/Right-*Children* in the next layer **k+1**). All nodes except the Root-node are *Child*-nodes (i.e. *from* the *Parent*-node in layer **k-1**). In each layer **k** there exist $2^k$ unique nodes. There exist a *countable infinity* of layers for **k** $\rightarrow \infty$. This implies that the **Total** number of nodes (i.e. of unique **routes**) in the **Tree** is an **Uncountable Infinity**. When the *countable infinity* of layers **k** is included, it is only possible to index **one** node in each layer with **k** which leaves $2^k$-1 nodes *without* index. **NB**! It is well known[IV] that the **power-set** of a *countable infinity* is an **Uncountable Infinity**.

## Collatz-operations on congruence-classes

Suppose a **class** N ≡ **P** (mod A) is described by three variables; N(A,**X**,B) = [A**X**-B], B = A-**P**. When the **modulo** A is even then A**X** is even and this implies, that the parity of [A**X**-B] is then determined *only* by the parity of the constant B < A and is independent of the variable **X** ∊ ℕ. This implies, that when A is even it **is** possible to perform a **Collatz-operation** on the class:

- Even A & **odd** $B_j$ in [A**X**-$B_j$] ⇒ a type Odd operation: [A**X**-$B_j$] ⇄O [**3**A**X**-(**3**$B_j$-**1**)$_{j+1}$]
- Even A & **even** $B_j$ in [A**X**-$B_j$] ⇒ a type Even operation: [A**X**-$B_j$] ⇄E [(A/**2**)**X**-($B_j$/**2**)$_{j+1}$]
- Odd A in [A**X**-$B_j$] ⇒ No operation is possible - the parity of the class now depends on **X**

Notice that in a class [A**X**-$B_j$] the **modulo** A is the distance between neighbour-elements $N_X$ & $N_{X+1}$ = ($N_X$+A), and **P** = (A-$B_j$) = $N_1$ is the *lowest* value in the *infinite set*; $N_X$ ∊ {$N_1$,$N_2$,$N_3$, ... ,$N_X$, ...}

Initial Value classes $[2^r$**X**-$B_0]_{IV}$ are **defined** to have **modulo** $2^r$, **r** ∊ $ℕ_0$, and for the **modulo** A there are **two** possibilities for the new **modulo 3**A *or* A/**2**, so **all** the used **moduli** are A = $2^r3^s$.

| |
|---|
| **r** > 0 & **odd** $B_j$ in [$2^r3^s$**X**-$B_j$] ⇒ a type Odd operation: [$2^r3^s$**X**-$B_j$] ⇄O [$2^r3^{s+1}$**X**-(**3**$B_j$-**1**)$_{j+1}$] |
| **r** > 0 & **even** $B_j$ in [$2^r3^s$**X**-$B_j$] ⇒ a type Even operation: [$2^r3^s$**X**-$B_j$] ⇄E [$2^{r-1}3^s$**X**-($B_j$/**2**)$_{j+1}$] |
| **r** = 0 ⇒ The **class** is a Terminal Value class $[3^s$**X**-$B_{j=s+r}]_{TV}$ and the **Class-series** Terminates |

All classes in the used **model** have the form $N_X$ ∊ [$2^r3^s$**X**-$B_j$], **X** ∊ ℕ, 0 ≤ $B_j$ < $2^r3^s$, **r** ∊ $ℕ_0$, **s** ∊ $ℕ_0$.
$N_X$ ∊ $[2^03^0$**X**-$0_0]_{IV/TV}$, **X** ∊ ℕ ⇒ $N_X$ = **X** is the **Universe-class** (the **Entire Domain**) in the **model**.

It is not possible to perform an operation on the **Universe-class**, because the **modulo** is odd, but when N ≡ **1** (mod 1) is ***split*** into odd values N ≡ **1** (mod 2) and even values N ≡ **2** (mod 2) it **is** possible to perform operations on the two *Child*-classes $[2^1$**X**-$1_0]_{IV}$ ∪ $[2^1$**X**-$0_0]_{IV}$ = $[2^0$**X**-$0_0]_{IV}$:

| $N_0$ | | $N_1$ | | $N_2$ | |
|---|---|---|---|---|---|
| $^*[2^1$**X**-$^*0_0]_{IV}$ | ⇄E | $[3^0$**X**-$0_1]_{TV}$ | | *R E-T | Class-series for Even $N_0$ |
| $^ϕ[2^1$**X**-$1_0]_{IV}$ | ⇄O | $[2^13^1$**X**-$2_1]$ | ⇄E | $[3^1$**X**-$1_2]_{TV}$ | Class-series for Odd $N_0$ |

For all **X** ∊ ℕ it is true that $[$**X**$]_{TV}$ < $[2^1$**X**$]_{IV}$ ⇒ $^*[2^{*1_0}$**X**-$^*0_{Con}]_{IV}$ is a **Converging* Initial Value class.
For all **X** ∊ ℕ it is true that $[3^1$**X**-$1_2]_{TV}$ > $[2^1$**X**-$1_0]_{IV}$ ⇒ $^ϕ[2^1$**X**-$1_0]_{IV}$ is a $^ϕ$**Diverging IV**-class.

When the **modulo** $2^r$ is doubled to $2^{r+1}$ the above classes are *split* and operations are possible:

| $N_0$ | | $N_1$ | | $N_2$ | | $N_3$ | | $N_4$ |
|---|---|---|---|---|---|---|---|---|
| $^Δ[2^2$**X**-$0_0]_{IV}$ | ⇄E | $[2^1$**X**-$0_1]$ | ⇄E | $[3^0$**X**-$0_2]_{TV}$ | | The four Class-series for $A_{IV} = 2^2$ | | |
| $^ϕ[2^2$**X**-$1_0]_{IV}$ | ⇄O | $[2^23^1$**X**-$2_1]$ | ⇄E | $[2^13^1$**X**-$1_2]$ | ⇄O | $[2^13^2$**X**-$2_3]$ | ⇄E | $[3^2$**X**-$1_4]_{TV}$ |
| $^Δ[2^2$**X**-$2_0]_{IV}$ | ⇄E | $[2^1$**X**-$1_1]$ | ⇄O | $[2^13^1$**X**-$2_2]$ | ⇄E | $[3^1$**X**-$1_3]_{TV}$ | | |
| $^*[2^2$**X**-$^*3_0]_{IV}$ | ⇄O | $[2^23^1$**X**-$8_1]$ | ⇄E | $[2^13^1$**X**-$4_2]$ | ⇄E | $[3^1$**X**-$2_3]_{TV}$ | | *R O E-E-T |

The two *Children* from the **Converging* **IV**-class $^*[2^{*1_0}$**X**-$^*0_{Con}]_{IV}$ = $^Δ[2^2$**X**-$0_0]_{IV}$ ∪ $^Δ[2^2$**X**-$2_0]_{IV}$ are **proper subsets** of the **reducing Parent*-class $^*[2^{*1_0}$**X**-$^*0_{Con}]_{IV}$, i.e. all the elements in the two classes are **reduced* in one type **E** -operation. **Proper subsets** (with longer parity-sequences) of **Converging* **IV**-classes $^*[2^{*r_s}$**X**-$^*B_{Con}]_{IV}$ are named $^Δ$**Redundant IV**-classes $^Δ[2^r$**X**-$B_0]_{IV}$, **r** > $^*$**r**$_s$.

For **X** = **1** it is true that $[3^1$**X**-$2_3]_{TV}$ = $[2^2$**X**-$3_0]_{IV}$ i.e. ***The Loop*** ... $1_j$ ↗ $4_{j+1}$ ↘ $2_{j+2}$ ↘ $1_{j+3}$ ...
For all **X** > **1** it is true that $[3^1$**X**-$2_3]_{TV}$ < $[2^2$**X**-$3_0]_{IV}$ ⇒ $^*[2^{*2_1}$**X**-$^*3_{Con}]_{IV}$ is a **Converging* **IV**-class.
For all **X** ∊ ℕ it is true that $[3^2$**X**-$1_4]_{TV}$ > $[2^2$**X**-$1_0]_{IV}$ ⇒ $^ϕ[2^2$**X**-$1_0]_{IV}$ is a $^ϕ$**Diverging IV**-class.

For all $[2^r$**X**-$B_0]_{IV}$, $B_0$ < $2^r$, the Class-series can be constructed; $[2^r$**X**-$B_0]_{IV}$ → ... → $[3^s$**X**-$B_{j=r+s}]_{TV}$

**Covering systems**

In the context of the present work a ***Covering System*** is a ***set*** of **congruence-classes** [A**X**-B] for which it is true that the *Union* of all the classes *Covers* the ***set*** of natural numbers N ∊ ℕ. **OBS**! In a proper/genuine ***Covering System*** all the unique classes [A**X**-B] are *pairwise disjoint*.

**NB**! Notice that **X** ∊ ℕ but *if* **X** ∊ ℤ in the classes [A**X**-B], 0 ≤ B < A the ***System*** *Covers* all N ∊ ℤ.

Below are shown examples of different ***Covering Systems*** encountered in the **model**:

$^{\Phi}[2^0\mathbf{X}\text{-}0]_{IV/TV}$ = ℕ = $^*[2^{*1_0}\mathbf{X}\text{-}{}^*0_{Con}]$ ∪ $^{\Phi}[2^1\mathbf{X}\text{-}1]$ =

$^{\Delta}[2^2\mathbf{X}\text{-}0]$ ∪ $^{\Phi}[2^2\mathbf{X}\text{-}1]$ ∪ $^{\Delta}[2^2\mathbf{X}\text{-}2]$ ∪ $^*[2^{*2_1}\mathbf{X}\text{-}{}^*3_{Con}]$ =

$^{\Delta}[2^3\mathbf{X}\text{-}0]$ ∪ $^{\Phi}[2^3\mathbf{X}\text{-}1]$ ∪ ⋯ ∪ $^{\Delta}[2^3\mathbf{X}\text{-}6]$ ∪ $^{\Delta}[2^3\mathbf{X}\text{-}7]$ =

$^{\Delta}[2^4\mathbf{X}\text{-}0]$ ∪ $^{\Phi}[2^4\mathbf{X}\text{-}1]$ ∪ ⋯ ∪ $^{\Delta}[2^4\mathbf{X}\text{-}14]$ ∪ $^{\Delta}[2^4\mathbf{X}\text{-}15]$ =

| Meaning of the symbols |
|---|
| $^{\Phi}$ means $^{\Phi}$**Diverging** |
| * means **Converging* |
| $^{\Delta}$ means $^{\Delta}$**Redundant** |

⋯ **NB**! All the above shown $[2^r\mathbf{X}\text{-}B_0]_{IV}$ are **I**nitial **V**alue classes for **modulo** $A_{IV} = 2^r$.

$[2^13^1\mathbf{X}\text{-}0]$ ∪ $[2^13^1\mathbf{X}\text{-}1]$ ∪ $[2^13^1\mathbf{X}\text{-}2]$ ∪ $[2^13^1\mathbf{X}\text{-}3]$ ∪ $[2^13^1\mathbf{X}\text{-}4]$ ∪ $[2^13^1\mathbf{X}\text{-}5]$ =

$[2^23^1\mathbf{X}\text{-}0]$ ∪ $[2^23^1\mathbf{X}\text{-}1]$ ∪ ⋯ ∪ $[2^23^1\mathbf{X}\text{-}10]$ ∪ $[2^23^1\mathbf{X}\text{-}11]$ =

$[2^13^2\mathbf{X}\text{-}0]$ ∪ $[2^13^2\mathbf{X}\text{-}1]$ ∪ ⋯ ∪ $[2^13^2\mathbf{X}\text{-}16]$ ∪ $[2^13^2\mathbf{X}\text{-}17]$ =

⋯

$[2^r3^s\mathbf{X}\text{-}0]$ ∪ $[2^r3^s\mathbf{X}\text{-}1]$ ∪ ⋯ ∪ $[2^r3^s\mathbf{X}\text{-}(2^r3^s\text{-}\mathbf{P})]$ ∪ ⋯ ∪ $[2^r3^s\mathbf{X}\text{-}(2^r3^s\text{-}\mathbf{2})]$ ∪ $[2^r3^s\mathbf{X}\text{-}(2^r3^s\text{-}\mathbf{1})]$ = ℕ

⋯ Notice that for constant A = $2^r3^s$ there exist $2^r3^s$ different *pairs* (**P**+B) = $2^r3^s$, **P** ∊ ℕ, B ∊ $ℕ_0$

Even though a ***Covering System*** is usually required to have a **finite** amount of residue-classes, it is possible to **construct** *Covering Systems*, which have an *infinite* amount of residue-classes, and there **also** exist *Covering Systems*, for which the *countably infinite* $[2^r\mathbf{X}\text{-}B]$ have *different* **moduli** $2^r$, **r** ∊ ℕ. This does **not** contradict the Mirsky-Newman theorem, which states that for a **finite** amount of residue-classes that are **both** (*pairwise*) *disjoint* **and** *distinct* (have *different* **moduli**) it is **impossible** for the **Union** of the residue-classes to *Cover* **all** N ∊ ℤ (**all** integers).

**Comment**

**Appendix Alpha** is a **List** containing the first **u** ≤ 142 **Converging* **IV**-classes $^*_{\mathbf{u}}[2^{*r_s}\mathbf{X}\text{-}{}^*B_{Con}]_{\mathbf{q}}$ where the $_{\mathbf{u}}$index **u** ∊ ℕ is used to identify a **u**nique **reducing* **IV**-class and the ***pair*** $_{\mathbf{u}}(\mathbf{P_q}+{}^*B_{Con})$, and the index$_{\mathbf{q}}$ **q** ∊ ℕ is used to identify the $\mathbf{P_q}$-values; $\mathbf{P_1} < \mathbf{P_2} < ... < \mathbf{P_{q_{Max}}}$ in a **Group** for A = $2^{*r_s}$. The **List** also includes the corresponding **TV**-classes $[3^s\mathbf{X}\text{-}B_{j=s+{}^*r_s}]_{TV}$ from the *Class-series. The $^*_{u\le142}[2^{*r_s}\mathbf{X}\text{-}{}^*B_{Con}]_{\mathbf{q}}$ were found by constructing the Class-series $[2^r\mathbf{X}\text{-}B_0]_{IV}$ → ... → $[3^s\mathbf{X}\text{-}B_{j=s+r}]_{TV}$ for all **IV**-classes $[2^r\mathbf{X}\text{-}B_0]_{IV}$ for all $B_0$ ≡ **1** (mod 4), $B_0 < 2^r$, where **modulo** $A_{IV} = 2^r$ for all **r** ≤ 13, thereby confirming the **Converging* **IV**-classes and the $^{\Phi}$**Diverging IV**-classes. All subsets of **Converging* **IV**-classes are $^{\Delta}$**Redundant** and $^{\Delta}B_0 \equiv {}^*B_{Con}$ (mod $2^{*r_s}$) for all **r** > *$\mathbf{r_s}$ in $^{\Delta}[2^r\mathbf{X}\text{-}{}^{\Delta}B_0]_{IV}$.

In the present work it is shown, that there exist a *countable infinity* of **Converging* **IV**-classes $^*_{u\to\infty}[2^{*r_s}\mathbf{X}\text{-}{}^*B_{Con}]_{q_{Max\to\infty}}$, **u** ∊ ℕ, and that the **Union** of **all** $^*_{\mathbf{u}}[2^{*r_s}\mathbf{X}\text{-}{}^*B_{Con}]_{\mathbf{q}}$ for **all u** ∊ ℕ is in itself a *Covering System*, as all the **Converging* **IV**-classes $^*_{\mathbf{u}}[2^{*r_s}\mathbf{X}\text{-}{}^*B_{Con}]_{\mathbf{q}}$ are *pairwise disjoint*.

It is also shown that there exist a *countable infinity* of $^{\Phi}$**Diverging IV**-classes $^{\Phi}_{u'}[2^r\mathbf{X}\text{-}B_0]_{q'}$ **u'** ∊ ℕ, but very important; the $^{\Phi}$**Diverging IV**-classes are **not** *pairwise disjoint* because all the $^{\Phi}[2^r\mathbf{X}\text{-}B_0]_{IV}$ are proper subsets of a $^{\Phi}$**Diverging** $^{\Phi}[2^k\mathbf{X}\text{-}B_0]_{IV}$, **k** < **r**, with **higher density**; $^1/_{2^k} > {}^1/_{2^r}$.

*Also all* $^*_{\mathbf{u}}[2^{*r_s}\mathbf{X}\text{-}{}^*B_{Con}]_{\mathbf{q}}$ *are proper subsets of a* $^{\Phi}$**Diverging** $^{\Phi}[2^k\mathbf{X}\text{-}B_0]_{IV}$, **k** < *$\mathbf{r_s}$ (*the proper subsets of the *Converging* **IV**-*classes are considered* $^{\Delta}$**Redundant IV**-*classes* $^{\Delta}[2^r\mathbf{X}\text{-}{}^{\Delta}B_0]_{IV}$, **r** > *$\mathbf{r_s}$).

**NB**! *Also shown is the existence of an* **Uncountable Infinity** *of proper subsets of* $^*_{\mathbf{u}}[2^{*r_s}\mathbf{X}\text{-}{}^*B_{Con}]_{\mathbf{q}}$, *actually an* **Uncountable Infinity** *of proper subsets of* any **IV**-*class* $[2^r\mathbf{X}\text{-}\mathbf{B_0}]_{IV}$ *in the used* **model**.

## 2. The Flowcharts (*flowgraphs*) for Collatz-iterations

**Collatz-operations on values N in Branches**

Suppose $N \in \mathbb{N}$ is described by the **odd** part $(2C-1)$, $C \in \mathbb{N}$ & **even** part $2^r$, $r \in \mathbb{N}_0$ $\boxed{N = (2C-1)2^r}$

For **r** > 0 in $(2C-1)2^r$ a Type Even operation is observed: $(2C-1)2^r \rightleftarrows (2C-1)2^{r-1}$

For **r** = 0 in $(2C-1)2^r$ a Type Odd operation is observed: $(2C-1) \rightleftarrows \mathbf{3}(2C-1)+\mathbf{1} = (6C-2)$

**OBS**! Notice that the *Dynamics* in The Problem is the result of the two above shown *bijections*.

A **Branch** $(2C-1)2^r$ is uniquely defined by the value $(2C-1)$ and contains only Even operations.

| ⋯ | 4 | | 3 | | 2 | | 1 | | r = 0 | | | |
|---|---|---|---|---|---|---|---|---|---|---|---|---|
| ⋯ | $(2C-1)2^4$ | ⇄ | $(2C-1)2^3$ | ⇄ | $(2C-1)2^2$ | ⇄ | $(2C-1)2^1$ | ⇄ | (2C-1) | ⇄ | (6C-2) | C |
| ⋯ | 16 | ⇄ | 8 | ⇄ | 4 | ⇄ | 2 | ⇄ | 1 | ⇄ | 4 | 1 |
| ⋯ | 48 | ⇄ | 24 | ⇄ | 12 | ⇄ | 6 | ⇄ | 3 | ⇄ | 10 | 2 |
| ⋯ | 80 | ⇄ | 40 | ⇄ | 20 | ⇄ | 10 | ⇄ | 5 | ⇄ | 16 | 3 |
| ⋯ | 112 | ⇄ | 56 | ⇄ | 28 | ⇄ | 14 | ⇄ | 7 | ⇄ | 22 | 4 |
| ⋯ | 144 | ⇄ | 72 | ⇄ | 36 | ⇄ | 18 | ⇄ | 9 | ⇄ | 28 | 5 |
| ⋯ | 176 | ⇄ | 88 | ⇄ | 44 | ⇄ | 22 | ⇄ | 11 | ⇄ | 34 | 6 |
| ⋯ | ⋯ | ⋯ | ⋯ | ⋯ | ⋯ | ⋯ | ⋯ | ⋯ | ⋯ | ⋯ | ⋯ | ⋯ |

**NB**! *The above* (*horizontal*) Branches *are the vertical* Branches *depicted in* **Figure 1a** *on page* 1.
It is observed for the **Result** of a Type Odd operation that: $\mathbf{3}(2C_a-1)+\mathbf{1} = (6C_a-2) = (2C_b-1)2^{r_b}$.
The value $N_b = (6C_a-2) = (2C_b-1)2^{r_b}$ is called the *connection-point* for the odd value $N_a = (2C_a-1)$.

**Lemma 1**

For the **Result** of a Type Odd operation $(2C_b-1)2^{r_b}$ one of the following two possibilities is true:

$\boxed{(2C_b-1) \equiv \mathbf{1} \ (\text{mod } 6) \ \wedge \ r_b = \mathbf{2n},\ n \in \mathbb{N}}$ **or** $\boxed{(2C_b-1) \equiv \mathbf{5} \ (\text{mod } 6) \ \wedge \ r_b = \mathbf{2n-1},\ n \in \mathbb{N}}$

**Proof for Lemma 1**

| For **X** ∈ ℕ, **n** ∈ ℕ | $(6X-5)2^{2n-1} \equiv \mathbf{2} \ (\text{mod } 6)$ | $(6X-1)2^{2n-1} \equiv \mathbf{4} \ (\text{mod } 6)$ |
|---|---|---|
| $(6X-3)2^n \equiv \mathbf{0} \ (\text{mod } 6)$ | $(6X-5)2^{2n} \equiv \mathbf{4} \ (\text{mod } 6)$ | $(6X-1)2^{2n} \equiv \mathbf{2} \ (\text{mod } 6)$ ■ |

**Theorem 1**

There exist **one and only one** *simple loop* where it is true that: $\boxed{\mathbf{3}(2C-1)+\mathbf{1} = (2C-1)2^r}$

**Proof for Theorem 1**

$\mathbf{3}(2C-1)+\mathbf{1} = (2C-1)2^r \Leftrightarrow (2C-1)2^r - \mathbf{3}(2C-1) = \mathbf{1} \Leftrightarrow (2C-1)(2^r-\mathbf{3}) = \mathbf{1} \Rightarrow (2C-1) = \mathbf{1} \wedge (2^r-\mathbf{3}) = \mathbf{1}$.
$(2C-1) = \mathbf{1} \Rightarrow \underline{C = 1}$ & $(2^r-\mathbf{3}) = \mathbf{1} \Rightarrow \underline{r = 2}$ so in *the loop* it is true that: $\mathbf{3(}2(1)-1\mathbf{)+1} = \mathbf{(}2(1)-1\mathbf{)}2^2$ ■

**OBS**! The above mentioned *simple loop* is ***The Loop*** $\boxed{\cdots 1_j \nearrow 4_{j+1} \searrow 2_{j+2} \searrow 1_{j+3} \cdots}$

**Corollary 2**

For all $C_a > 1$ in $(2C_a-1)$ it is true for $\boxed{\mathbf{3}(2C_a-1)+\mathbf{1} = (2C_b-1)2^{r_b}}$ that: $C_a \neq C_b \Leftrightarrow (2C_a-1) \neq (2C_b-1)$.

Notice that $C_a > 1$ does *not* imply that $C_b > 1$ e.g. for $C_a = 3$: $\mathbf{3}(2(3_a)-1)+\mathbf{1} = (2C_b-1)2^{4_b} \Rightarrow C_b = 1$.

## Collatz-operations on classes $[2^r3^s X\text{-}B]$ in Branches

Remember from the previous page that operations on **values** N in Branches is described by:

| … | 4 | | 3 | | 2 | | 1 | | r = 0 | | |
|---|---|---|---|---|---|---|---|---|---|---|---|
| … | $(2C-1)2^4$ | ⇄ | $(2C-1)2^3$ | ⇄ | $(2C-1)2^2$ | ⇄ | $(2C-1)2^1$ | ⇄ | $(2C-1)$ | ⇄ | $(6C-2)$ |

For odd $B_j$ in $[2^r3^s X\text{-}B_j]$ there exist a parallel description for operations on **classes** in Branches:

| … | r+2 | | r+1 | | r > 0 | | s+1 |
|---|---|---|---|---|---|---|---|
| … | $[2^{r+2}3^s\mathbf{X}\text{-}(2^2B_j)_{j-2}]$ | ⇄ | $[2^{r+1}3^s\mathbf{X}\text{-}(2^1B_j)_{j-1}]$ | ⇄ | $[2^r3^s\mathbf{X}\text{-}B_j]$ | ⇄ | $[2^r3^{s+1}\mathbf{X}\text{-}(3B_j\text{-}1)_{j+1}]$ |

While a Value-Branch $(2C-1)2^r$ is valid for a particular odd $(2C-1)$, $C \in \mathbb{N}$, $\mathbf{r} \in \mathbb{N}_0$, a Class-Branch is valid for a particular odd class $[2^r3^s\mathbf{X}\text{-}B_j]$, $B_j$ odd, $\mathbf{r} > 0$, i.e. for all $N_X \in [2^r3^s\mathbf{X}\text{-}B_j]$, $\mathbf{X} \in \mathbb{N}$.

There exist three different types of **Main**-branches:

| | | | | | | | | | |
|---|---|---|---|---|---|---|---|---|---|
| … | $[2^33^1\mathbf{X}\text{-}8_{j-3}]$ | ⇄ | $[2^33^1\mathbf{X}\text{-}4_{j-2}]$ | ⇄ | $[2^23^1\mathbf{X}\text{-}2_{j-1}]$ | ⇄ | $[2^13^1\mathbf{X}\text{-}1_j]$ | ⇄ | $[2^13^2\mathbf{X}\text{-}2_{j+1}]$ |
| … | $[2^33^1\mathbf{X}\text{-}24_{j-3}]$ | ⇄ | $[2^33^1\mathbf{X}\text{-}12_{j-2}]$ | ⇄ | $[2^23^1\mathbf{X}\text{-}6_{j-1}]$ | ⇄ | $[2^13^1\mathbf{X}\text{-}3_j]$ | ⇄ | $[2^13^2\mathbf{X}\text{-}8_{j+1}]$ |
| … | $[2^33^1\mathbf{X}\text{-}40_{j-3}]$ | ⇄ | $[2^33^1\mathbf{X}\text{-}20_{j-2}]$ | ⇄ | $[2^23^1\mathbf{X}\text{-}10_{j-1}]$ | ⇄ | $[2^13^1\mathbf{X}\text{-}5_j]$ | ⇄ | $[2^13^2\mathbf{X}\text{-}14_{j+1}]$ |

Notice the similarities with the illustration on the previous page and the proof for **Lemma 1**.

- Notice that [**6X-5**] ∪ [**6X-3**] ∪ [**6X-1**] = $^\Phi[2\mathbf{X}\text{-}1]_{IV}$ i.e. all odd natural numbers
- and that [**18X-14**] ∪ [**18X-8**] ∪ [**18X-2**] = [**6X-2**] i.e. all possible *connection-points*

  For the Main-branches it is observed, that the **modulo** in the odd classes is $A = 2^13^1$. A branch which have a higher **modulo** in the odd class $[2^r3^s\mathbf{X}\text{-}B_j]$, $B_j$ odd, for $\mathbf{r} > 1$ and/or $\mathbf{s} > 1$, contains a **subset** of one of the three different types of Main-branches.

Branches can be used to construct **Flowcharts** (*flowgraphs*) showing the *connections* between different classes. The below illustration shows the types of **nodes** used in the Flowcharts.

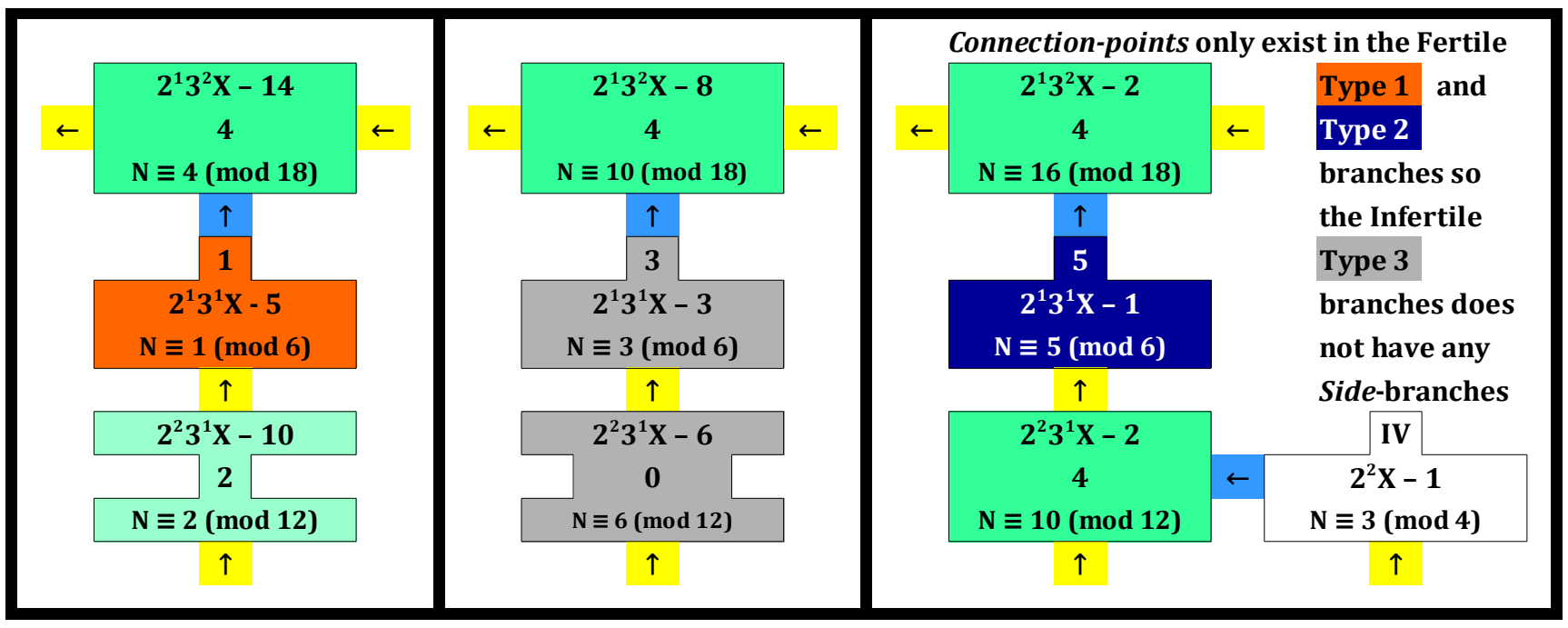

Overview showing the three **Types** of Main-Branches {**1**,**2**,**3**}

→ is an Even operation.
→ is an Odd operation.
The arrows indicates the *Flow-direction*. When the arrows are followed the *direction* is "**DOWN**" and vice versa the *direction* "**UP**" is **against** the direction of arrows.

The Class-branches only contains type Even operations while type Odd operations are **links** (*connections*) between *different* branches in the Flowcharts/*flowgraphs*. In the Flowcharts shape and colour of the **nodes** reflects the *congruence* **modulo** 6. Three type **even** nodes are used and one type **odd** node. Type of branch depends on the *congruence* **modulo** 3 for the **odd** class in the branch. It is possible to colour an **odd** node **when** the *congruence* **modulo** 3 is **known**; (**6X-5**) ≡ **1** (mod 3), (**6X-1**) ≡ **2** (mod 3), (**6X-3**) ≡ **3** (mod 3). In the Initial Value classes $[2^r\mathbf{X}\text{-}B_{j=0}]_{IV}$ the **modulo** does ***not*** have 3 as a factor, so it is *not possible* to determine the *congruence* **modulo** 6 i.e. the **IV**-classes contains **subsets** from all the three different **types**.

All the **nodes** contain the **Result** of a Type Even operation, so it is possible to ***extrapolate*** **against** type **E**-arrows and construct "*the **even** previous node*"; $[2^r3^s\mathbf{X}\text{-}B_j]$ $\leftarrow_E$ $[2^{r+1}3^s\mathbf{X}\text{-}(2B_j)_{j-1}]$ Notice that it is possible to extrapolate "*to the right*" from the *connections-points* on top and extrapolate "*down on the page*" without any limitations when extrapolating **against** the $\leftarrow_E$.

For a *connection-point* it is possible to extrapolate **against** type **O**-arrows and construct "*the **odd** previous node*" $[6\mathbf{X}\text{-}2_j]$ $\leftarrow_O$ $[2\mathbf{X}\text{-}1_{j-1}]$; $[2^r3^s\mathbf{X}\text{-}B_j]$ $\leftarrow_O$ $[2^r3^{s-1}\mathbf{X}\text{-}(^{B_j+1}/_3)_{j-1}]$ for **r** > 0. *Eventually* an **IV**-class $[2^r\mathbf{X}\text{-}B_j]_{IV}$ is reached - where the *congruence* **modulo** 3 depends on **X**.

The nodes contain well-defined ***sets*** of values $N_X$, and because two different descriptions of the ***set*** are shown; $N_X \in [2^r3^s\mathbf{X}\text{-}B] \Leftrightarrow N_X \equiv \mathbf{P} \pmod{2^r3^s}$ where $(\mathbf{P}+B) = 2^r3^s$, all the important information about the ***set*** is ***known***. The smallest value $N_1$ in each ***set*** is $\mathbf{P} = (2^r3^s\text{-}B)$ and the second-smallest value is $N_2 = (\mathbf{P}+2^r3^s) = (2^{r+1}3^s\text{-}B)$. The **modulo** and the constant B are ***known***. When the arrows are followed from node to node the **Class-series** describes the iteration-series for all $N_X$-values in the classes and in **particular** for the **P**-value (for **X** = 1) in each class.

## Comment

It is observed that the $^{\Phi}$**Diverging IV**-class $^{\Phi}[2^2\mathbf{X}\text{-}1]_{IV}$ is found as the **first** *Side*-branch *connecting* to the first *connection-point* in Type 2 branches, and the **IV**-class $^{\Phi}[2^2\mathbf{X}\text{-}1]_{IV}$ is named **The Domain-class** (which is a proper subset of **The Universe-class**) as the ***set*** $\{\mathbf{3}_1,7_2,11_3,15_4, \ldots ,N_X, \ldots\}$ contains all $N_X$ **not** in the *Union* $*_1[2^{*1_0}\mathbf{X}\text{-}*0_{Con}]_1 \cup *_2[2^{*2_1}\mathbf{X}\text{-}*3_{Con}]_1$.

The below "Minimal Flowcharts" demonstrates the minimum amount of ***extrapolation*** which is required to have all odd N *Covered* i.e. included in *at least one* class $[2^r3^s\mathbf{X}\text{-}B]$ for odd B, **r** > 0.

The Flowcharts illustrates the earlier described **Class-series** $[2^r\mathbf{X}\text{-}B_{j=0}]_{IV} \rightarrow \cdots \rightarrow [3^s\mathbf{X}\text{-}B_{j=s+r}]_{TV}$

## Illustrations of the two Minimal Flowcharts

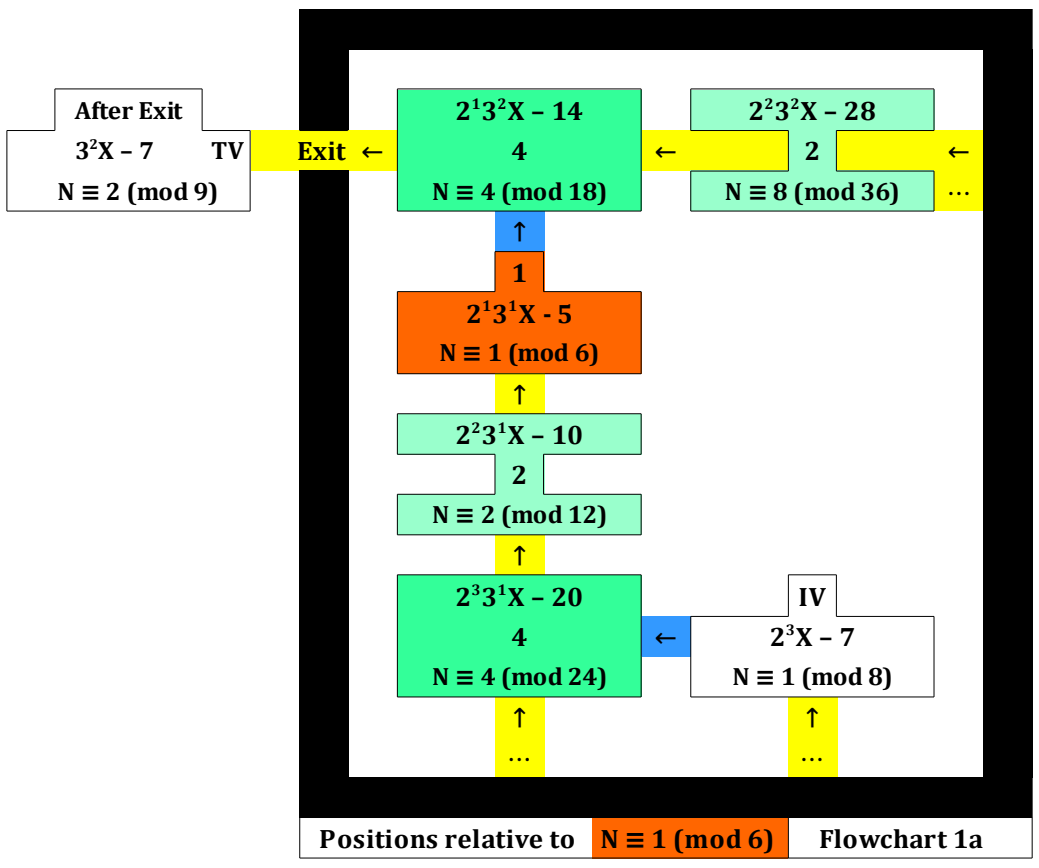

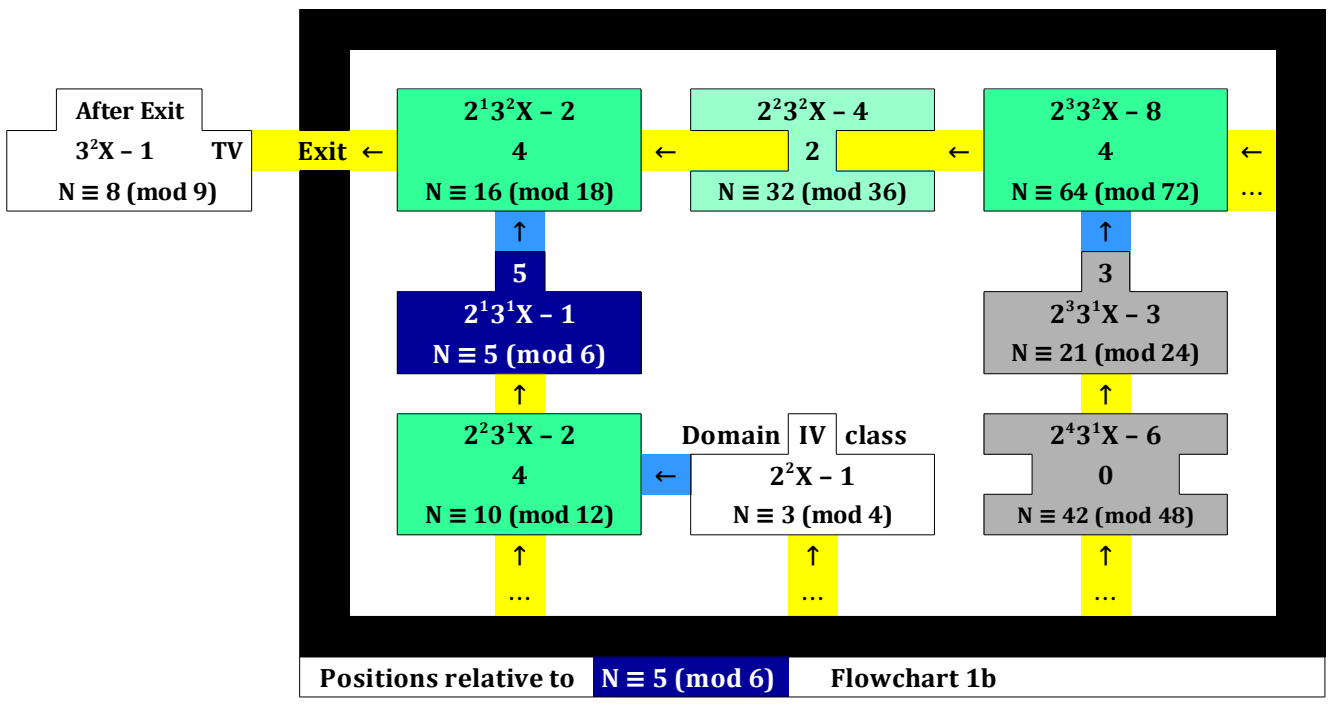

All elements in $N \equiv 10 \pmod{18}$ are included in the Union of the congruence-classes $N \equiv 10 \pmod{12} \cup N \equiv 4 \pmod{24} \cup N \equiv 64 \pmod{72}$ i.e all $N \equiv 4 \pmod 6$ are included.

**Flowchart 1a & 1b**. The Minimal Type 1 and Type 2 Flowcharts.

**Nodes** *are identified with a* **class** $[2^r3^s\mathbf{X}\text{-}B]$ *and parity-sequence i.e. order of* **O** & **E** *operations.*
In both the **Flowcharts** the horizontal **Exit**-branch is for **X** = 1 for the **odd** value 1; $(2(\mathbf{1})\text{-}1)2^r$.

The Flowcharts contains all $N_X \equiv \mathbf{1} \pmod 6 \Leftrightarrow N_X \in \{1_1,7_2,13_3,19_4,25_5, \ldots\}$ in Type **1** Flowcharts and all $N_X \equiv \mathbf{5} \pmod 6 \Leftrightarrow N_X \in \{5_1,11_2,17_3,23_4,29_5, \ldots\}$ in Type **2** Flowcharts. Because the Union of $[2^2\mathbf{X}\text{-}1] \cup [2^3\mathbf{X}\text{-}7] \cup [2^33^1\mathbf{X}\text{-}3]$ contains all $N_X \equiv \mathbf{3} \pmod 6$ it follows that the **Union** of the **Minimal Flowcharts** for all $\mathbf{X} \in \mathbb{N}$ contains all the *possible* (2C-1) implicitly all $(2C\text{-}1)2^r$, $\mathbf{r} \in \mathbb{N}_0$.
$[2^2\mathbf{X}\text{-}1] = \{\mathbf{3}_1,7_2,11_3,\mathbf{15}_4,19_5, \ldots\}$ & $[2^3\mathbf{X}\text{-}7] = \{1_1,9_2,17_3,25_4,\mathbf{33}_5, \ldots\}$ & $[2^33^1\mathbf{X}\text{-}3] = \{\mathbf{21}_1,\mathbf{45}_2,\mathbf{69}_3, \ldots\}$

In the Flowcharts/*flowgraphs* the classes **after** the **Exit** are actually the previously described **T**erminal **V**alue classes $[3^s\mathbf{X}\text{-B}]_{TV}$ having odd **modulo** $A_{TV} = 3^s$ i.e. the parity depends on **X**. Please notice that for a **X** ∊ ℕ the nodes in a Flowchart contains specific *values* for a *constant* **X**.

The two **Minimal Flowcharts** contains two different **IV**-classes i.e. two different Class-series:

**1b**; $[2^2\mathbf{X}\text{-}1_0]_{IV}$ ⇄O $[2^2 3^1\mathbf{X}\text{-}2_1]$ ⇄E $[2^1 3^1\mathbf{X}\text{-}1_2]$ ⇄O $[2^1 3^2\mathbf{X}\text{-}2_3]$ ⇄E $[3^2\mathbf{X}\text{-}1_4]_{TV}$

$^\Phi[2^2\mathbf{X}\text{-}1_0]_{IV}$ is **The Domain-class** having the $^\Phi$**Diverging** $^\Phi$parity-sequence $^\Phi$ **R O E - O E - T**.

**1a**; $[2^3\mathbf{X}\text{-}7_0]_{IV}$ ⇄O $[2^3 3^1\mathbf{X}\text{-}20_1]$ ⇄E $[2^2 3^1\mathbf{X}\text{-}10_2]$ ⇄E $^*[2^1 3^1\mathbf{X}\text{-}5_3]$ ⇄O $[2^1 3^2\mathbf{X}\text{-}14_4]$ ⇄E $[3^2\mathbf{X}\text{-}7_5]_{TV}$

$^\Delta[2^3\mathbf{X}\text{-}7_0]_{IV}$ is a **subset** of the **Converging* **IV**-class $^*{}_2[2^{*2_1}\mathbf{X}\text{-}^*3_{Con}]_{IV}$ where **X** = 1 ⇒ ***The Loop***.

For **X** = 1 the vertical Type **1** Main-branch in **Flowchart 1a**, including the odd class [**6X-5**], is a representation of "***The Tree before***" (as seen in **Figure 1b,** page 1), because the **P**-values in the branch includes values from ***The Loop*** i.e. from the Value-Branch **(**2(**1**)-1**)**$2^r$. Similarly the first *Side*-branch, including the odd **IV**-class [8**X**-7], represents "***The Tree before that***".

**Conclusions** can be made from comparing the **modulo** $A_{IV} = 2^r$ and $A_{TV} = 3^s$ in Class-series, *or* by observing the **P**-values (or the B-values) in the Class-series; $^*\mathbf{P}_{TV} < \mathbf{P}_{IV}$ (i.e. $B_{TV} < B_{IV}$) ⇒ **Converging* **IV**-class $^*[2^{*r_s}\mathbf{X}\text{-}^*B_{Con}]_{IV}$ & v.v. $\mathbf{P}_{TV} > \mathbf{P}_{IV}$ ($B_{TV} \geq B_{IV}$) ⇒ $^\Phi$**Diverging IV**-class $^\Phi[2^r\mathbf{X}\text{-}B_0]_{IV}$

## Note

In the below **Extended** Flowcharts it is demonstrated that the *Union* of *Side*-branches on the Main-branches *Covers* **all** odd N, and that **IV**-classes in *Side*-branches are *pairwise disjoint*.

## Illustrations of the Extended Flowcharts for modulo A = $2^r 3^s$

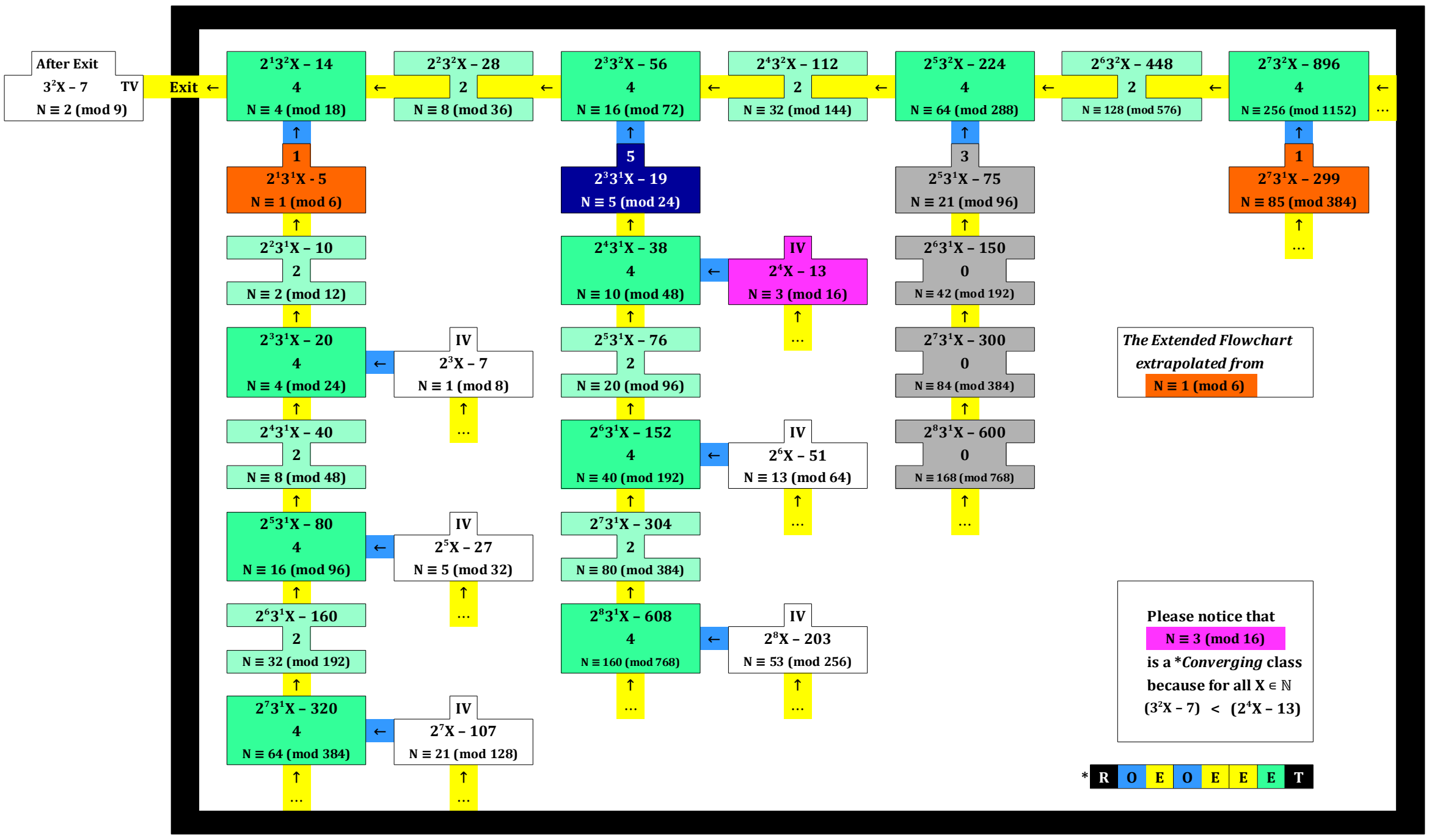

**Flowchart 2a**. The Extended **Type 1** Flowchart. For **X** = 1 the horizontal branch is **(**2(**1**)-1**)**$2^r$.
**OBS**! The vertical Type **1** Main-branch is only relevant for **X** > 1 due to ***The Loop***.

Notice that **Flowchart 2a** have a **node** containing the **Converging* **IV**-class $^*{}_3[2^{*4_2}\mathbf{X}\text{-}^*13_{Con}]_1$:

$[2^{*4_2}\mathbf{X}\text{-}13_0]_{IV}$ ⇄O $[2^4 3^1\mathbf{X}\text{-}38_1]$ ⇄E $[2^3 3^1\mathbf{X}\text{-}19_2]$ ⇄O $[2^3 3^2\mathbf{X}\text{-}56_3]$ ⇄E $[2^2 3^2\mathbf{X}\text{-}28_4]$ ⇄E $[2^1 3^2\mathbf{X}\text{-}14_5]$ ⇄E $[3^2\mathbf{X}\text{-}7_6]_{TV}$

Please also notice that *after* the first type Odd operation it is **not** possible for the **P**-values and the B-values to have **3** as a factor i.e. the $B_j$ is **either** $B_j \equiv$ **1** (mod 3) **or** $B_j \equiv$ **2** (mod 3).

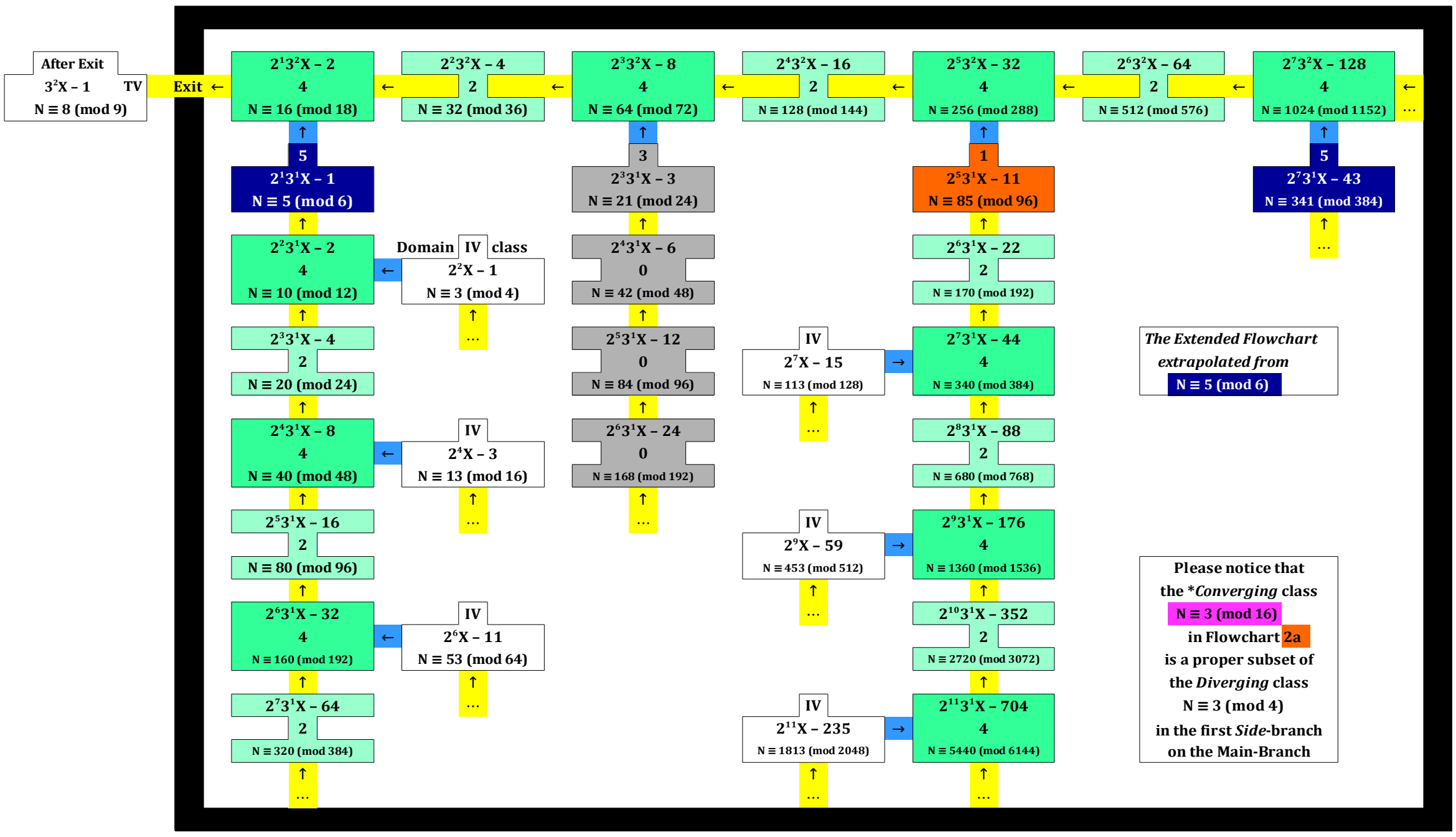

**Flowchart 2b**. The Extended Type 2 Flowchart. For **X** = 1 the horizontal branch is **(**2(**1**)-1**)**$2^r$**.**
**OBS**! The vertical Type **2** Main-branch is relevant for all **X** ∊ ℕ (no ***Loop***-values).

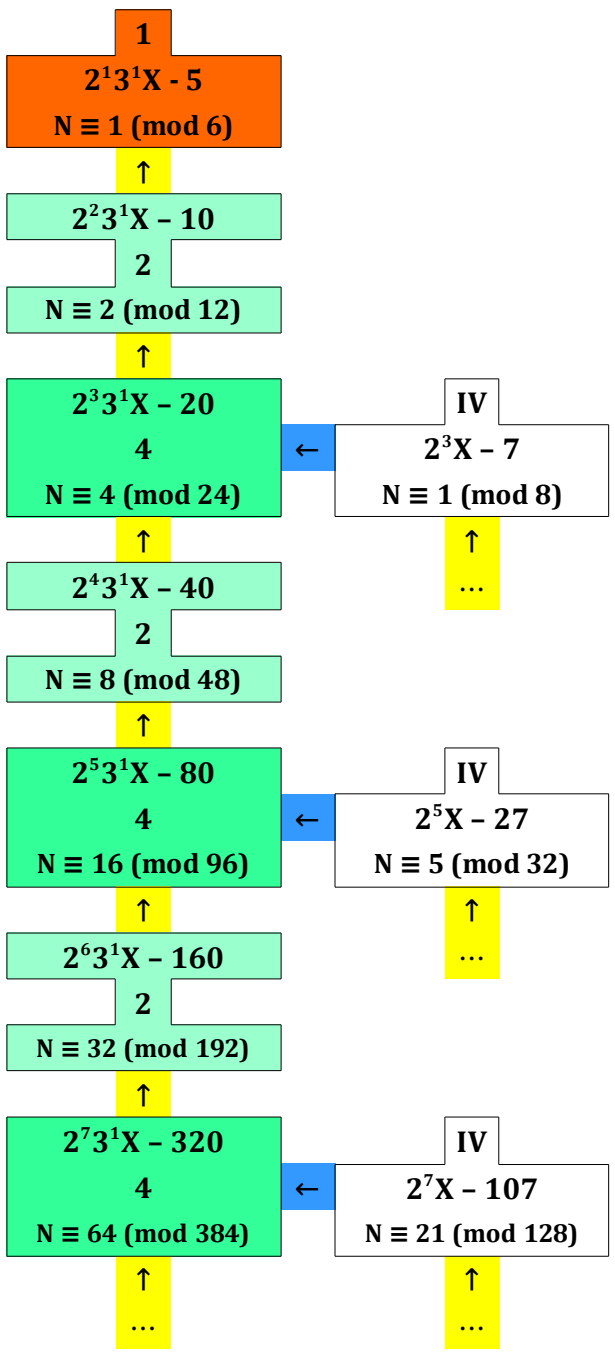

Type **1** Main-branch

The two Fertile Main-branches contains all the possible *connection-points* [**6X-2**] and the Infertile Type **3** branches can not contain any *connection-point* at all, because 3 is a factor in all classes in a Type **3** branch, and a *connection-point* always is congruent to **1** (mod 3), which is a direct consequence of the rule for Type Odd operations; 3N+1 for odd N. **NB**! See **Lemma 1** on page 6 for details.

Because all odd values N *connects* to one of the *connection-points*, then the *Union* of all the **IV**-classes in the *Side*-branches *connecting* to one of the Fertile Main-branches **must** include **all** (2C-1), C ∊ ℕ.

**Patterns** can be observed for neighbour **IV**-classes in the shown *Side*-branches.

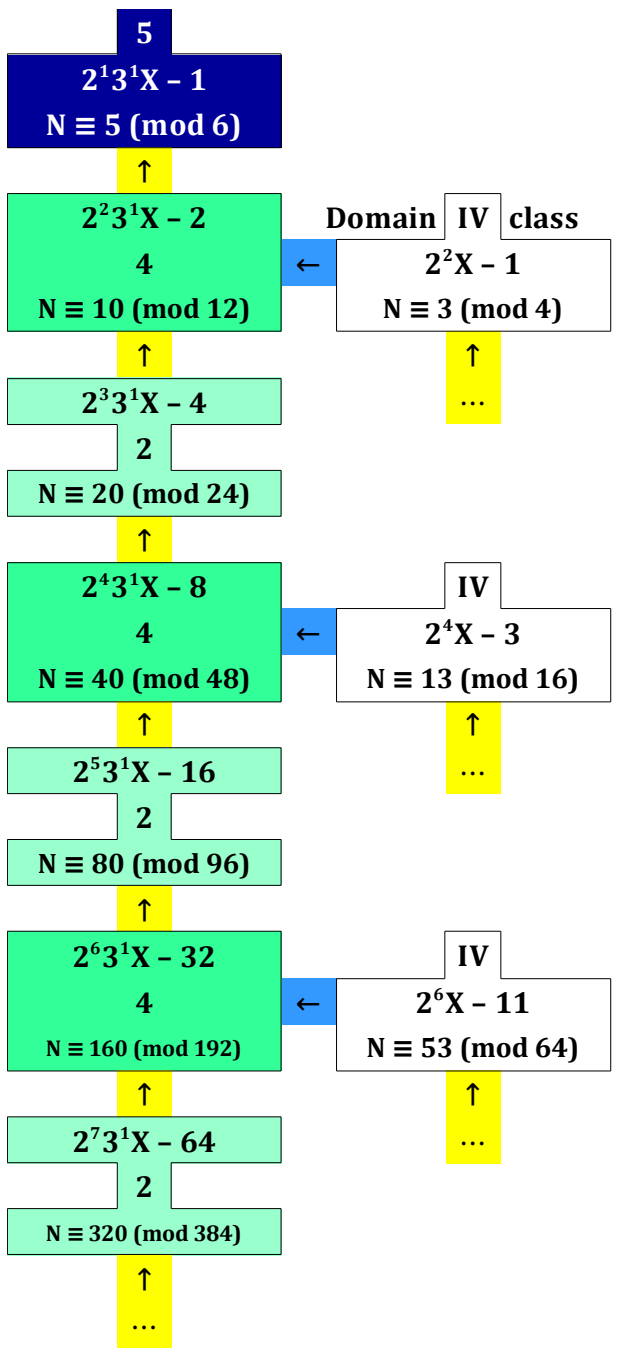

Type **2** Main-branch

The *Union* of **IV**-classes in the *Side*-branches on the two Main-branches contains **all** odd N:

[$2^2$**X**-1] U [$2^3$**X**-7] U [$2^4$**X**-3] U [$2^5$**X**-27] U [$2^6$**X**-11] U [$2^7$**X**-107] U ⋯ → [$2^1$**X**-1]

In the above *sequence* it is true that the **next** class (A = $2^{r+1}$) contains *half of the rest* of odd N.

| ***Side*-branches on Type 1 Main-branch** | **n∈ℕ** | ***Side*-branches on Type 2 Main-branch** |
|---|---|---|
| $[2^3\mathbf{X}\text{-}7] \Leftrightarrow N \equiv \mathbf{1} \pmod{2^3}$ | **1** | $[2^2\mathbf{X}\text{-}1] \Leftrightarrow N \equiv \mathbf{3} \pmod{2^2}$ |
| $[2^5\mathbf{X}\text{-}27] \Leftrightarrow N \equiv \mathbf{5} \pmod{2^5}$ | 2 | $[2^4\mathbf{X}\text{-}3] \Leftrightarrow N \equiv \mathbf{13} \pmod{2^4}$ |
| $[2^7\mathbf{X}\text{-}107] \Leftrightarrow N \equiv \mathbf{21} \pmod{2^7}$ | 3 | $[2^6\mathbf{X}\text{-}11] \Leftrightarrow N \equiv \mathbf{53} \pmod{2^6}$ |
| $[2^9\mathbf{X}\text{-}427] \Leftrightarrow N \equiv \mathbf{85} \pmod{2^9}$ | 4 | $[2^8\mathbf{X}\text{-}43] \Leftrightarrow N \equiv \mathbf{213} \pmod{2^8}$ |
| … ⇔ … | … | … ⇔ … |
| $[2^{2n+1}\mathbf{X}\text{-}B] \Leftrightarrow N \equiv \mathbf{P} \pmod{2^{2n+1}}$ | **n** | $[2^{2n}\mathbf{X}\text{-}B] \Leftrightarrow N \equiv \mathbf{P} \pmod{2^{2n}}$ |
| $[2^{2n+3}\mathbf{X}\text{-}(4B\text{-}1)] \Leftrightarrow N \equiv (4\mathbf{P}+1) \pmod{2^{2n+3}}$ | **n**+1 | $[2^{2n+2}\mathbf{X}\text{-}(4B\text{-}1)] \Leftrightarrow N \equiv (4\mathbf{P}+1) \pmod{2^{2n+2}}$ |
| … ⇔ … | … | … ⇔ … |

| Exact formula for $[2^r\mathbf{X}\text{-}B]$ on Type **1** | | Exact formula for $[2^r\mathbf{X}\text{-}B]$ on Type **2** |
|---|---|---|
| $[2^{2n+1}\mathbf{X} - (\sum_{n=2}^{n} 2^{2(n-1)}5^1 + 7_{n=1})]$ | **n** | $[2^{2n}\mathbf{X} - (\sum_{n=2}^{n} 2^{2(n-1)-1} + 1_{n=1})]$ |

The *Union* of **IV**-classes in the *Side*-branches *Covers* all odd N = (2C-1), C ∈ ℕ **exactly** once, and in **Union** with the **Converging* **IV**-class $*_1[2^{*1_0}\mathbf{X}\text{-}*0_{Con}]_{IV}$ it is a *Covering System* where all the *countably infinite* **IV**-classes are *pairwise disjoint* **and** also *distinct* as all **moduli** $A_{IV} = 2^r$, **r** ∈ ℕ, are *different* and observe the ***geometric series*** for the **Sum** of *densities*; $\lim \sum_{r=1}^{\infty} \mathbf{1/2^r = 1}$

**NB**! The *density* $1/2^r$ for an **IV**-class $[2^r\mathbf{X}\text{-}B_0]_{IV}$ is **defined** to be the *inverse* of the **modulo** $2^r$.

Notice that the **Domain-class** $^{\Phi}[2^2\mathbf{X}\text{-}1]_{IV} \Leftrightarrow N \equiv \mathbf{3} \pmod 4$ is found in the first *Side*-branch (**n** = 1) on the Type **2** Main-branch i.e. **all** the **other** *Side*-branches contains **subsets** of the **Converging* **IV**-class $*N_X \in *_2[2^{*2_1}\mathbf{X}\text{-}*3_{Con}]_{IV} \Leftrightarrow *N_X \equiv \mathbf{1} \pmod 4 \Leftrightarrow *N_X \in \{1_1, 5_2, 9_3, 13_4, 17_5, \ldots\}$ and notice that parity-sequences for these **subsets** *starts* with the **reducing* **R O E - E** -... while the parity-sequence for all $N_X \in {}^{\Phi}[2^2\mathbf{X}\text{-}1]_{IV}$ *starts* with the $^{\Phi}$**Diverging R O E - O E** -...

**Comment**

It is **significant** that all **IV**-classes in *Side*-branches are *pairwise disjoint* as for all (2C-1), C ∈ ℕ: N = (2C-1) have a **specific** *position* **X** ∈ ℕ in a unique **IV**-class $N_X \in [2^r\mathbf{X}\text{-}B_0]_{IV}$ in a **specific** *Side*-branch **n** ∈ ℕ from the above *Lists* of *Side*-branches **n** on the Type **1** branch & Type **2** branch. Notice that all odd values in the classes [**6X-5**] and [**6X-1**] are *also* elements in **one** of the *Side*-branches from the above *Lists*, and from **Corollary 2** (p. 6) it is known, that no (2C-1) > 1 can possibly be *connected* to a *connection-point* in "its own Branch" $(2C\text{-}1)2^r$. This is a very strong indication, that an alternative loop (other than the known ***Loop*** ...1↗4↘2↘1...) can **not** exist. The later **Discussion** contains more details about *why* an alternative loop can not exist (**NB**! It is possible to prove this **after** the analysis of parity-sequences in the following main-sections).

**Note**

In the below Flowcharts for ***triple* modulo** it is demonstrated, that it is possible to extrapolate **against** an **additional** Type **Odd** operation $\rightleftarrows_o$ every time the **modulo** A is ***tripled*** to 3A, so it is possible to determine the *congruence* **modulo** 3 in more of the (now coloured) odd classes. **NB**! The Flowcharts for *triple* **modulo** are *illustrations* of the *proof* for the following **Lemma 2.**

**Appendix Beta** contains larger versions of the six below shown Flowcharts. The focus-point at the present time is the fact, that all charts are **different** (and the **Union** contains **all** odd N).

### Illustrations of the six Flowcharts 3a-f for *triple* modulo 3A = $2^r3^{s+1}$

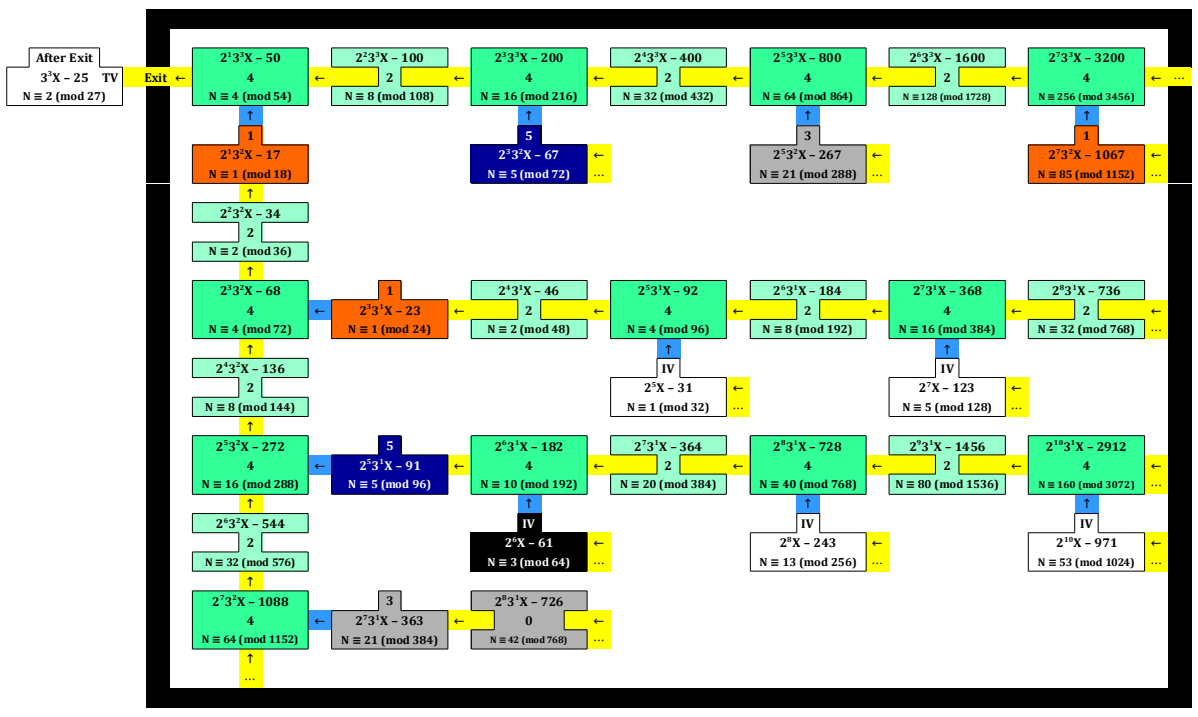

**3a**. Extrapolated from N ≡ **1** (mod $2^13^2$)
**TV**-class is N ∊ [$3^3$**X**-25] ⇔ N ≡ **2** (mod 27)

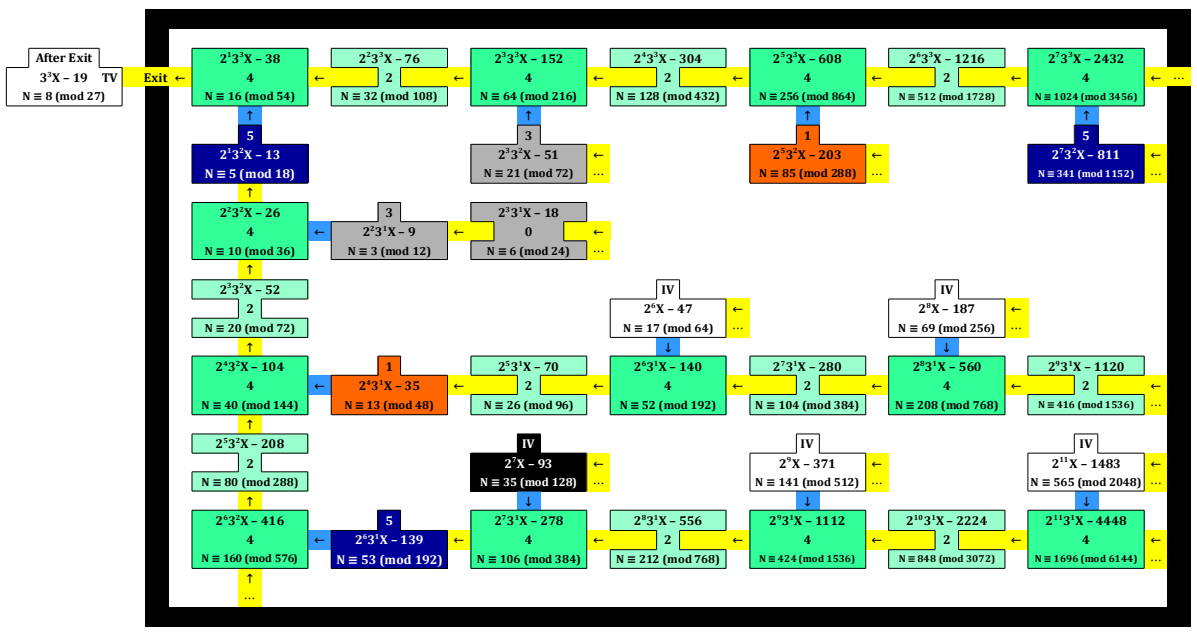

**3b**. Extrapolated from N ≡ **5** (mod $2^13^2$)
**TV**-class is N ∊ [$3^3$**X**-19] ⇔ N ≡ **8** (mod 27)

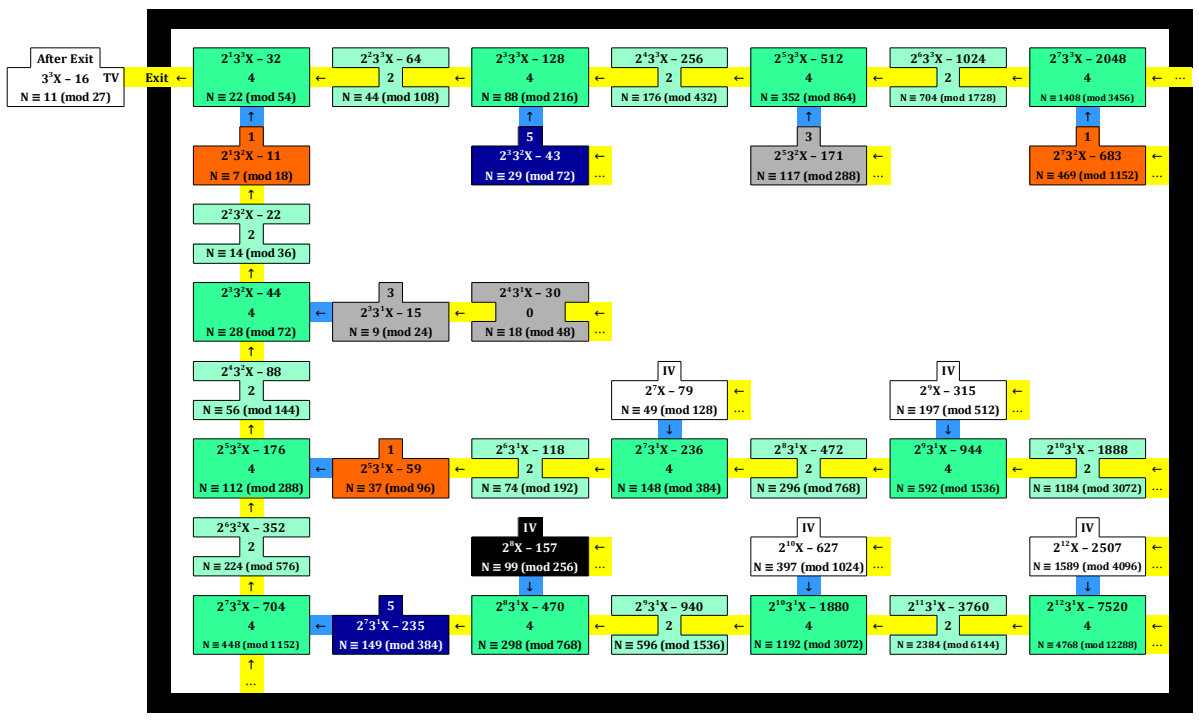

**3c**. Extrapolated from N ≡ **7** (mod $2^13^2$)
**TV**-class is N ∊ [$3^3$**X**-16] ⇔ N ≡ **11** (mod 27)

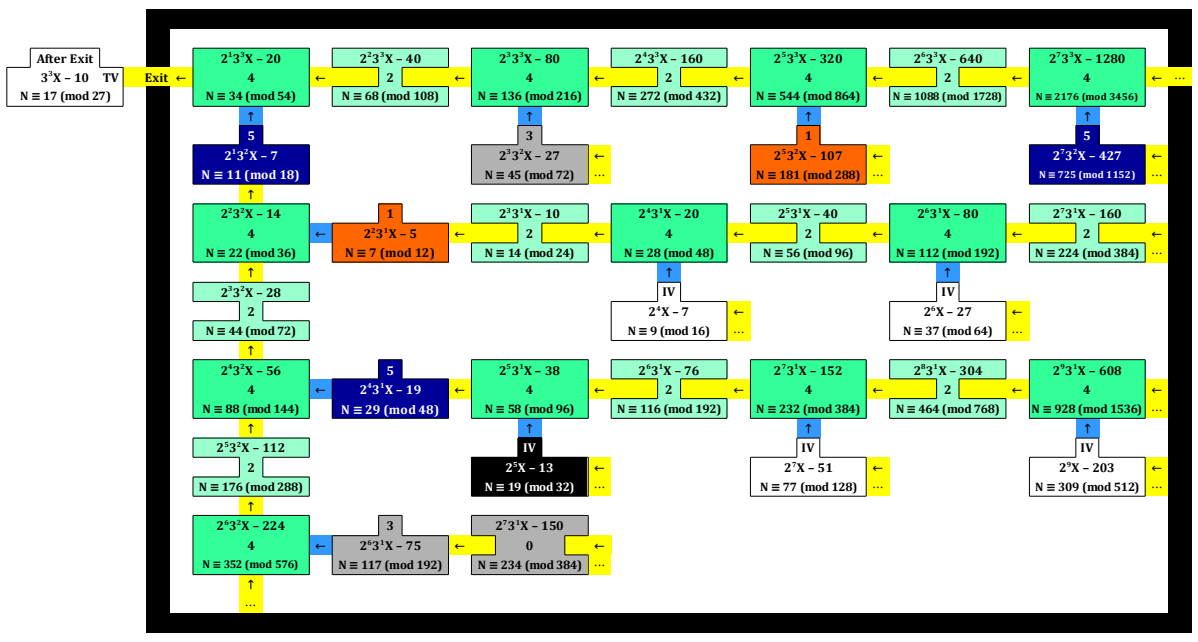

**3d**. Extrapolated from N ≡ **11** (mod $2^13^2$)
**TV**-class is N ∊ [$3^3$**X**-10] ⇔ N ≡ **17** (mod 27)

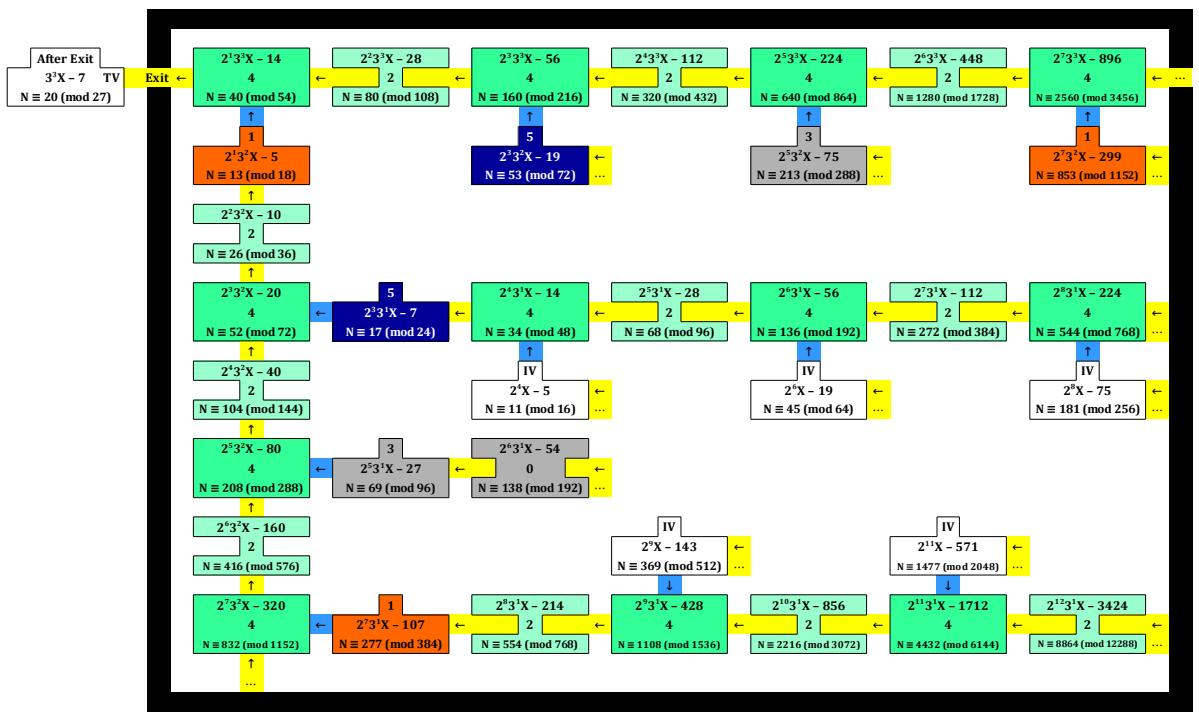

**3e**. Extrapolated from N ≡ **13** (mod $2^13^2$)
**TV**-class is N ∊ [$3^3$**X**-7] ⇔ N ≡ **20** (mod 27)

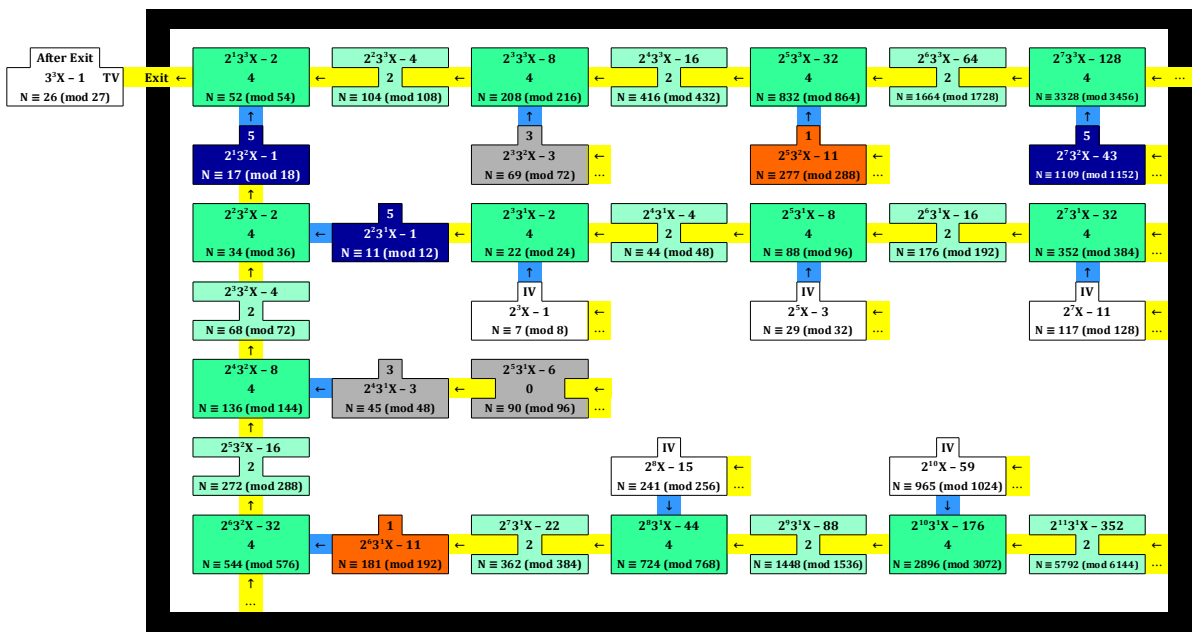

**3f**. Extrapolated from N ≡ **17** (mod $2^13^2$)
**TV**-class is N ∊ [$3^3$**X**-1] ⇔ N ≡ **26** (mod 27)

Be aware that the iterationseries for 7 as the **I**nitial **V**alue $N_0$ is:

$7_0$ ↗ $22_1$ ↘ $11_2$ ↗ $34_3$ ↘ $17_4$ ↗ $52_5$ ↘ $26_6$ ↘ $13_7$ ↗ $40_8$ ↘ $20_9$ ↘ $10_{10}$ ↘ ***5**$_{11}$ ↗ $16_{12}$ ↘ $8_{13}$ ↘ $4_{14}$ ↘ $2_{15}$ ↘ †$1_{16}$

Notice that all **P** ∊ {$7_{0c}$,$11_{2d}$,$17_{4f}$,$13_{7e}$,$5_{11b}$,$1_{16a}$} from the **Fertile** N ≡ **P** (mod 18) are in the series.

Notice that all **P** ∊ {$2_a$,$8_b$,$11_c$,$17_d$,$20_e$,$26_f$} from the **TV**-classes N ≡ **P**$_{TV}$ (mod $3^3$) are in the series.

This implies that all **P**-values in above shown classes N ≡ **P** (mod A) are in fact *connected* to †$1_j$.

It is observed in the **Flowcharts 3a-f** that all the **TV**-classes now have *higher* **modulo** $A_{TV} = 3^3$, and that it is now possible to determine the *congruence* **modulo** 3 in **more** of the odd **classes** *connecting* to **Fertile** branches, as all the **IV**-classes from **Flowchart 2a** & **2b** have now been separated into the three **types** {**1**,**2**,**3**}. This would also be possible for other *Side*-branches on the *horizontal* **Exit**-branches containing **subsets** of the **classes** from **Flowchart 2a** & **2b**.

Notice that from the six **Flowcharts 3a-3f**, it can be deduced exactly which **subsets** of the *Side*-branches on the two Main-branches in **Flowchart 2a** & **2b**, that are **Fertile** (have *Side*-branches) and the *one of three* **Infertile** *Side*-branches (which do **not** have any *Side*-branches).

While it is definitely true that **more** information is gained by **tripling** the **modulo**, it *seems* to come at some kind of cost, as a "**loss of generality**" is observed concerning the **Converging* **IV**-classes. For example for the **reducing* **IV**-class $^*N \equiv \mathbf{3} \pmod{\mathbf{16}}$ **subsets** are now observed in several charts; $N \equiv \mathbf{19} \pmod{2^5}$ & $N \equiv \mathbf{3} \pmod{2^6}$ & $N \equiv \mathbf{35} \pmod{2^7}$ & $N \equiv \mathbf{99} \pmod{2^8}$

Observations in the shown **Flowcharts 2a & 2b** and **Flowcharts 3a-f** leads to the following:

**Lemma 2**

When *Side*-branches on **Fertile branches** are observed, any three **consecutive** branches are a **triple** of one *each* of the three **types** of branches {**1**,**2**,**3**}.

**Proof for Lemma 2**

It is observed that an **odd** value (2C-1) in a *Side*-branch have the *connection-point* (6C-2) in the **Parent**-branch. For the **next** *Side*-branch on the **Parent**-branch it is **true** that **modulo** have been **doubled twice** *ergo* the **next** *connection-point* is **4**(6C-2) = (24C-8) and the **odd** value for **this** *point* is (**2C'-1**) = **(**(24C-8)**-1)**/**3** = (24C-9)/**3** = (8C-3) = (8C-4)**+1** = **4**(2C-1)**+1**

Suppose N = (2C-1) = (3**X**-a), a < 3, implying **4**N+**1** = **4**(3**X**-a)+**1** = (**12X+1**)-**4**a

a = 0 ⇒ **4**a ≡ **0** (mod 3) ∧ a = 1 ⇒ **4**a ≡ **1** (mod 3) ∧ a = 2 ⇒ **4**a ≡ **2** (mod 3) so the (**12X+1**) ≡ **1** (mod 3) - but not the **4**a - **changes** the *congruence* i.e. it is a **cycle** ■

**Comment**

Remember that it is possible to extrapolate against type $\leftarrow_E$ operations (in the direction "**UP**") without any limitation in the Flowcharts and all Class-Branches contains an infinity of nodes. Suppose a **Flowchart 2a** and a **Flowchart 2b** are printed for all **X** ∊ ℕ. In this case each of the printed Flowcharts shows the *connections* and relative positions between the Value-Branches included, as each **node** now contain **one** value $N_X = (2^r3^s\mathbf{X}\text{-}B)$ for a particular **X**. The infinity of printed Flowcharts is a ***collection*** which includes all odd $N_X$ ∊ [2**X**-1] and all even $N_X$ ∊ [**6X-2**]. If all the six different **Flowchart 3a** to **3f** are printed for all **X** ∊ ℕ, then the result is the **same** infinite ***collection*** of Flowcharts for Value-Branches where iterationseries for N are observed. Imagine treating the ***collection*** of printed Flowcharts as a jigsaw-puzzle, where the odd value, N = (2C-1), in a specific Value-Branch in a particular Flowchart can be observed to **overlap** with an odd value in *another* Flowchart. Remember that all the odd $N_X$ ∊ [**6X-5**] ∪ $N_X$ ∊ [**6X-1**] are **also** an element in **one** of the *Side*-branches found in the *Lists* of *Side*-branches (page 11).

*When* it is proved that **all** the printed Flowcharts in the ***collection*** are *connected* to **all** the *other* Flowcharts in the ***collection***, then ***The Collatz Conjecture*** will be proven **true**. If **all** the Flowcharts are *connected*, they are actually part of **one** giant Flowchart, which includes all odd N = (2C-1) implicitly all $N = (2C-1)2^r$, C ∊ ℕ, **r** ∊ $\mathbb{N}_0$, i.e. **all** natural numbers N ∊ ℕ.

## Note

In the **Flowcharts** on the following pages it is demonstrated, how it is possible to ***extend*** the parity-sequences, when the **modulo** is **doubled** and all classes are ***split*** into two *Child*-classes, and how this can bring knowledge about **more** **Converging* **subsets** of **The Domain-class**, when the **modulo** is **doubled** *ad infinitum*. (**NB**! *Which is the subject of the* ***rest*** *of the analysis.*)

**Appendix Delta** contains larger versions of the four below shown **Flowcharts 4a-d**.

### Illustrations of the four Flowcharts for *double* modulo 2A = $2^{r+1}3^s$

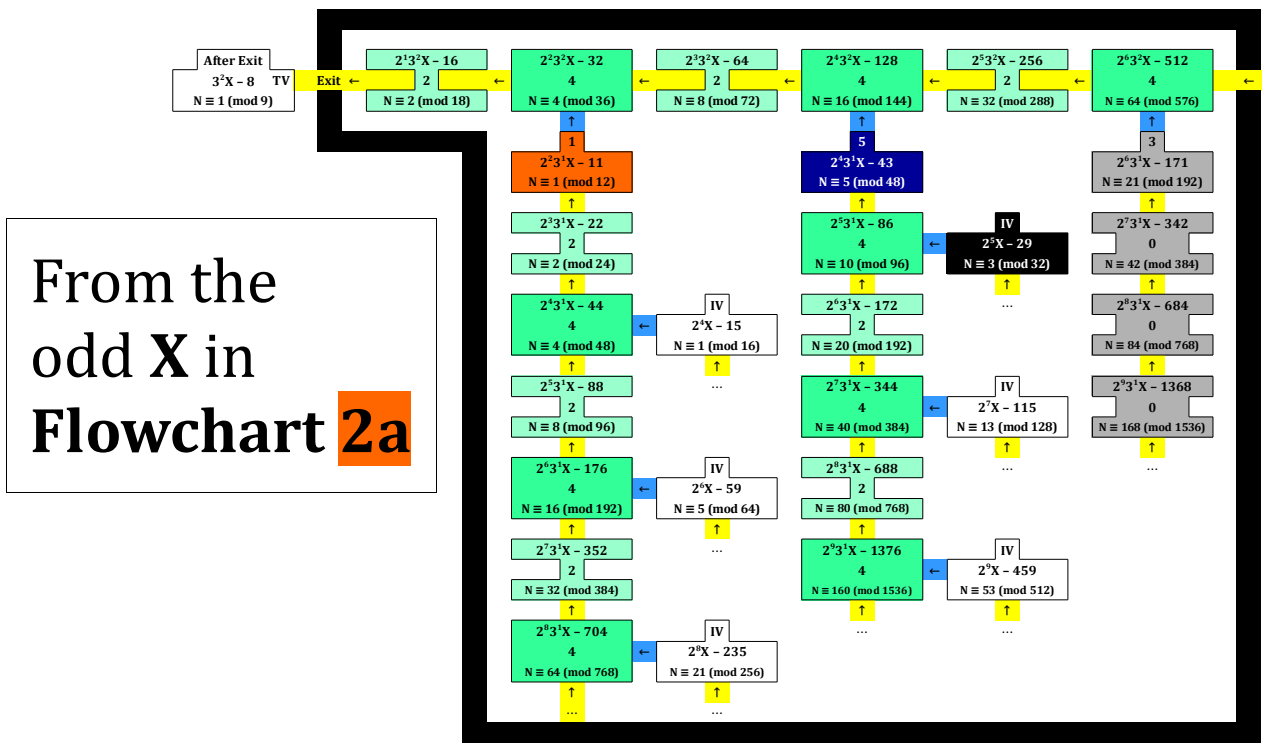

**4a**. Extrapolated from N ≡ 1 (mod $2^2 3^1$)
**TV**-class [$3^2$**X**-8] ⇔ N ≡ **1** (mod $3^2$)

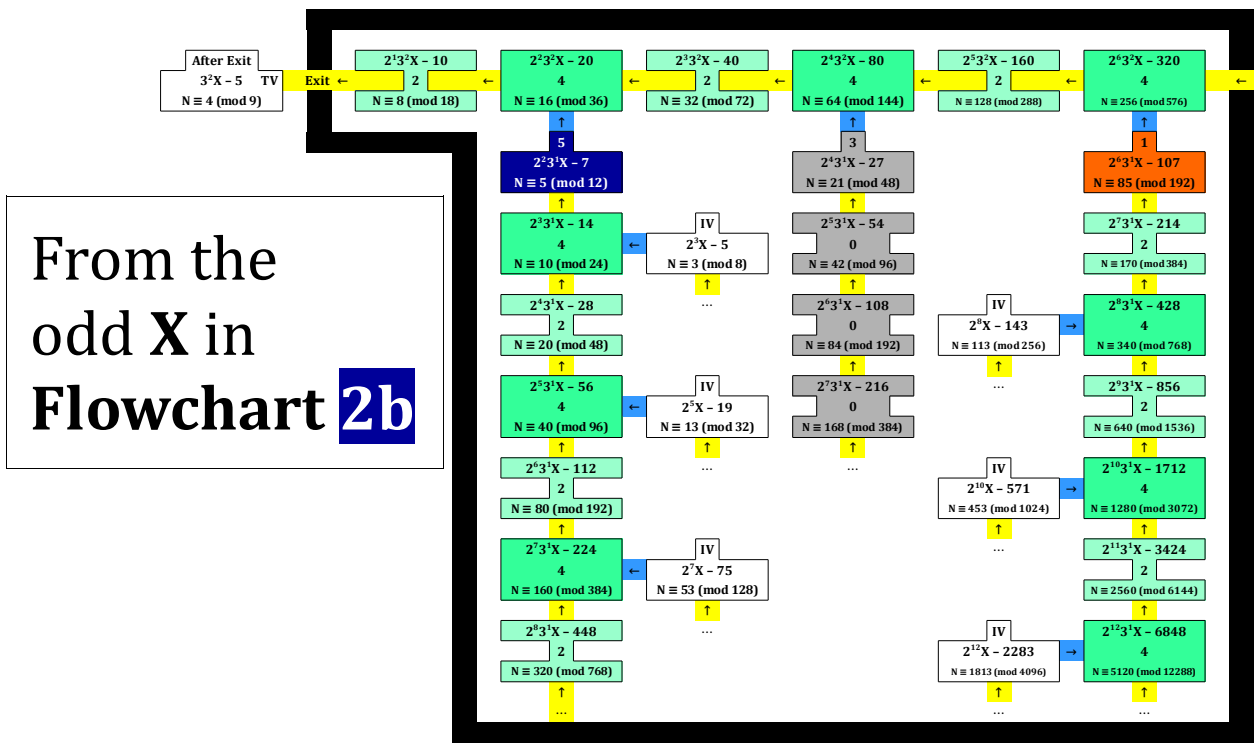

**4b**. Extrapolated from N ≡ 5 (mod $2^2 3^1$)
**TV**-class [$3^2$**X**-5] ⇔ N ≡ **4** (mod $3^2$)

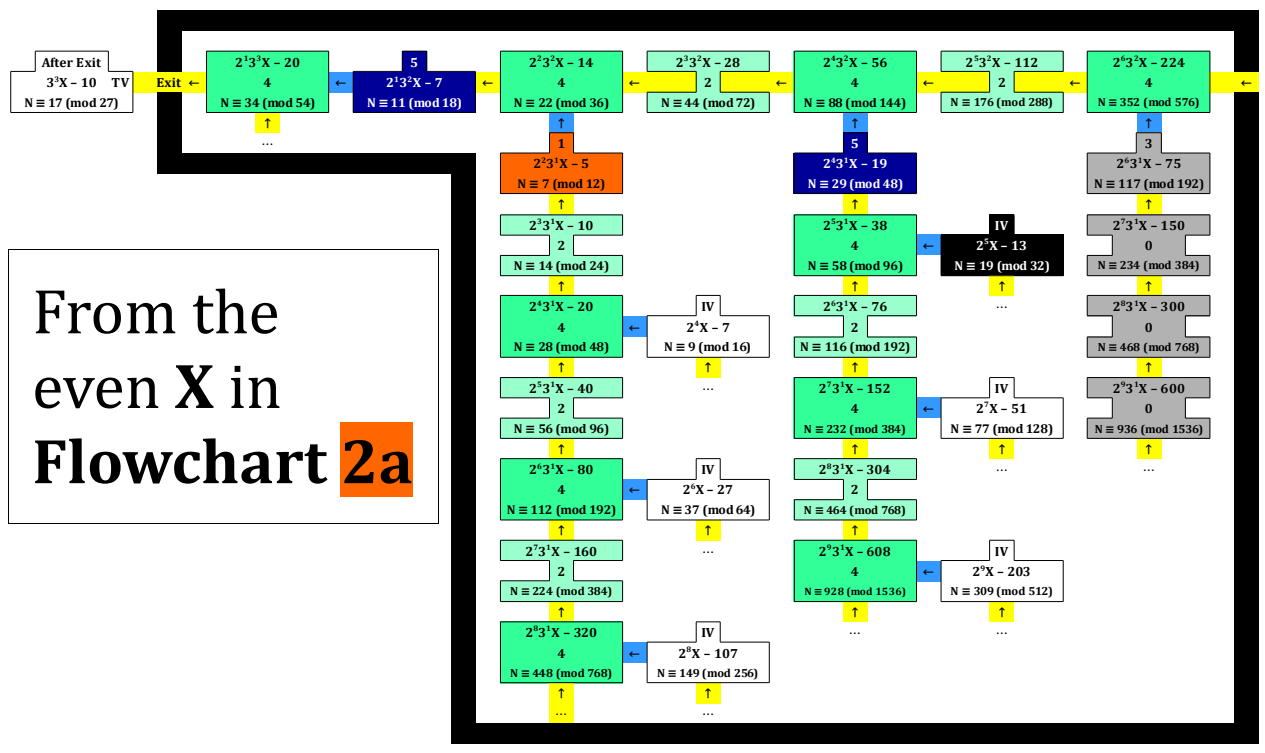

**4c**. Extrapolated from N ≡ 7 (mod $2^2 3^1$)
**TV**-class [$3^3$**X**-10] ⇔ N ≡ **17** (mod $3^3$)

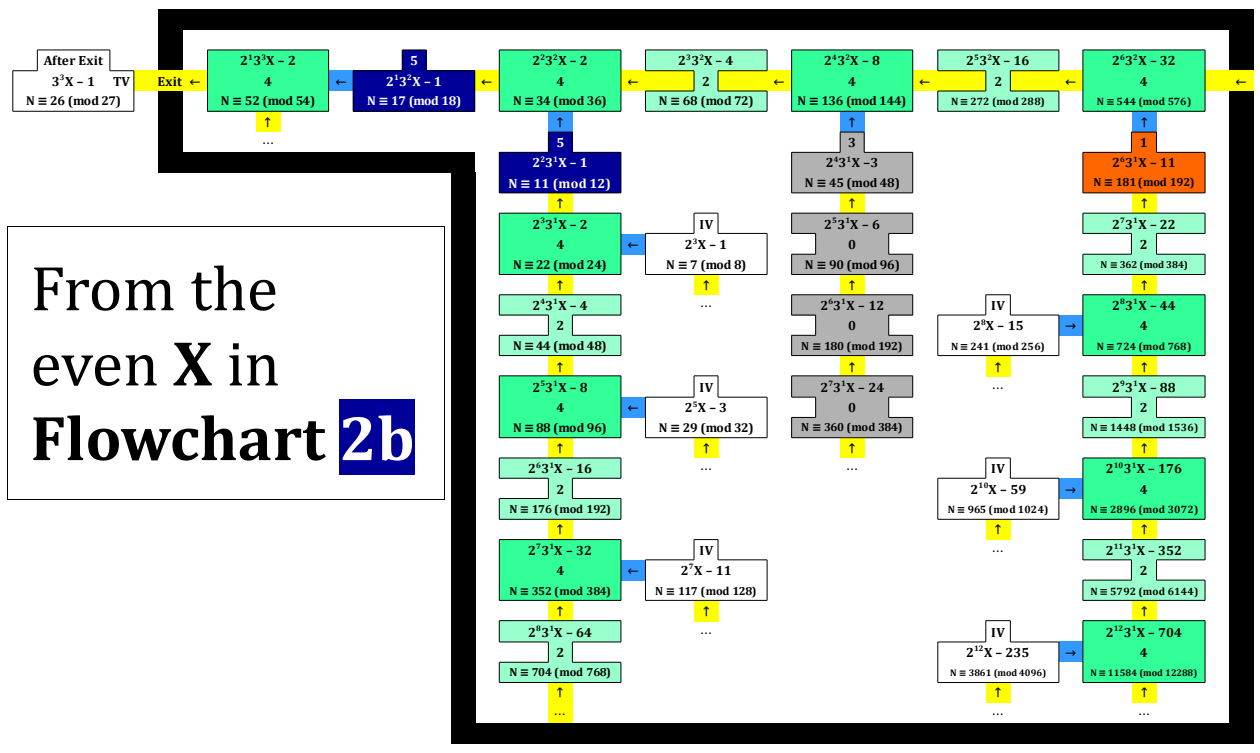

**4d**. Extrapolated from N ≡ 11 (mod $2^2 3^1$)
**TV**-class [$3^3$**X**-1] ⇔ N ≡ **26** (mod $3^3$)

Notice in **4a** that the Class-series for [$2^5$**X**-$29_{\mathbf{End}}$]$_3$ terminates at [$3^2$**X**-$8_{\mathbf{j=2+5}}$]$_{\mathbf{TV}}$ & ($3^2$-8) = $^\dagger\mathbf{1}_7$
$^\dagger$[$2^5$**X**-($2^5$-$3_3$)$_{\mathbf{End}}$]$_3$ → ... → [$3^2$**X**-($3^2$-$^\dagger\mathbf{1}$)$_7$]$_{\mathbf{TV}}$ is the $^\dagger$**End** Class-series for the Initial **Value** $N_N$ = $3_3$

For all $N_N$ which reach $^\dagger\mathbf{1}_j$ the $^\dagger$**End** Class-series exist: $^\dagger$[$2^d$**X**-($2^d$-$N_N$)$_{\mathbf{End}}$]$_N$ → ... → [$3^{s_d}$**X**-($3^{s_d}$-$^\dagger\mathbf{1}$)$_j$]$_{\mathbf{TV}}$

A *unique* $^\dagger$**End** Class-series for $N_N \in \mathbb{N}$ have a *unique* **End** $^\dagger$parity-sequence: $^\dagger$ **R** ... - **E** - ... - **E** - **D**
containing $\mathbf{s}_d$ pcs of **O** $_{\text{Operation}}$ i.e. $\mathbf{s}_d$ pcs of **O E** **-tiles** and **d** pcs of **E** $_{\text{Operation}}$ i.e. (**d**-$\mathbf{s}_d$) **E** **-tiles.**

## Comment

*Notice that all possible* **End TV**-classes [$3^{s_d}$**X**-($3^{s_d}$-1)$_{\mathbf{j=s_d+d}}$]$_{\mathbf{TV}}$ ⇒ $\mathbf{P}_{\mathbf{TV}}$ = $^\dagger\mathbf{1}_j$ *certainly exist for* $\mathbf{s}_d \in \mathbb{N}$.
*It can be deduced from the behaviour of* ***The Loop*** *- an infinity of* $^\Delta$**Redundant** *Class-series exist*:

| | | |
|---|---|---|
| $1_1 \in {}^{*\dagger}$[$2^2$**X**-*$3_{Loop}$]$_1$ ≥ [$3^{1_d}$**X**-$21_d$+2]$_{\mathbf{TV}}$ | 1 ∈ $^\Delta$[$2^4$**X**-$15_{Loop}$]$_{\mathbf{IV}}$ ≥ [$3^{2_d}$**X**-$82_d$+4]$_{\mathbf{TV}}$ | 1 ∈ $^\Delta$[$2^6$**X**-$63_{Loop}$]$_{\mathbf{IV}}$ ≥ [$3^{3_d}$**X**-$263_d$+6]$_{\mathbf{TV}}$ |
| 1 ∈ $^\Delta$[$2^8$**X**-$255_{Loop}$]$_{\mathbf{IV}}$ ≥ [$3^{4_d}$**X**-$804_d$+8]$_{\mathbf{TV}}$ | … | 1 ∈ $^\Delta$[$2^{2s_d}$**X**-($2^{2s_d}$-1)$_{Loop}$]$_{\mathbf{IV}}$ ≥ [$3^{s_d}$**X**-($3^{s_d}$-1)$_{s_d+2s_d}$]$_{\mathbf{TV}}$ |

**Appendix Gamma** contains larger versions of the eight below shown **Flowcharts 5a-h**.

**Illustrations of the eight Flowcharts for *quadruple* modulo 4A = $2^{r+2}3^s$**

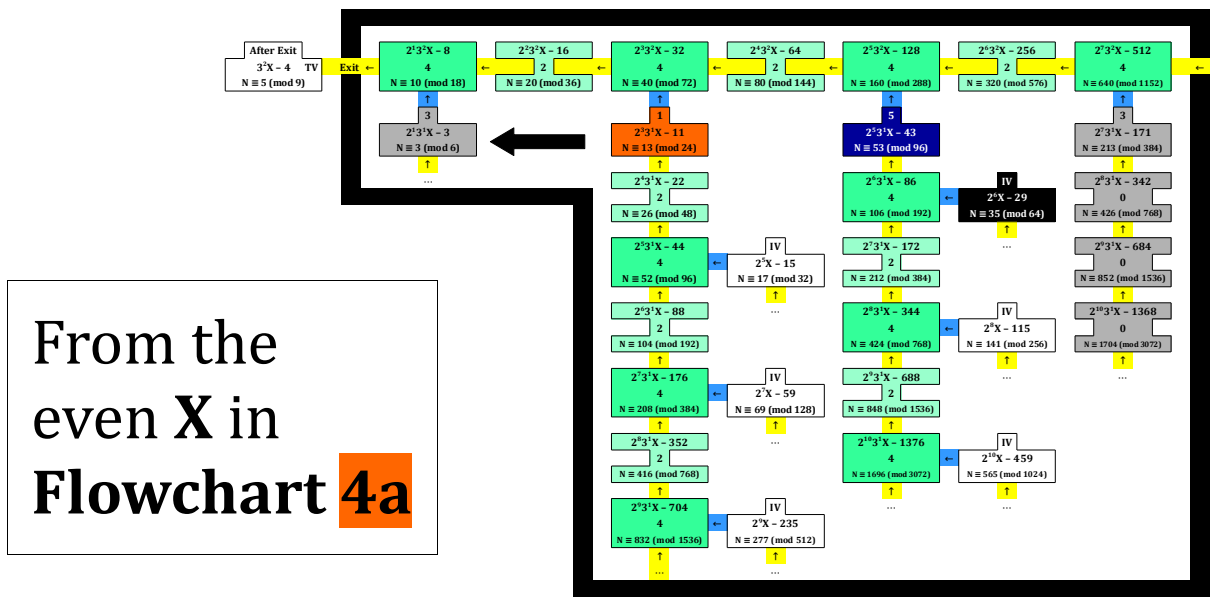

**5a**. Extrapolated from $N \equiv 13 \pmod{2^3 3^1}$
**TV**-class $[3^2$**X**$-4] \Leftrightarrow N \equiv$ **5** $\pmod{3^2}$

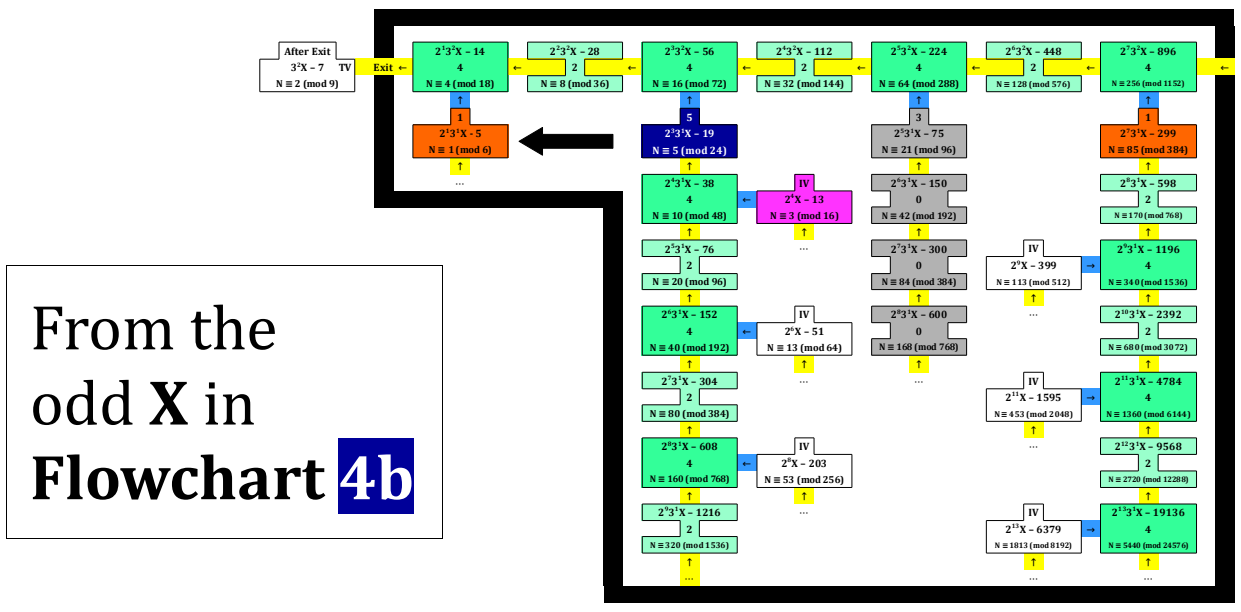

**5b**. Extrapolated from $N \equiv 5 \pmod{2^3 3^1}$
**TV**-class $[3^2$**X**$-7] \Leftrightarrow N \equiv$ **2** $\pmod{3^2}$

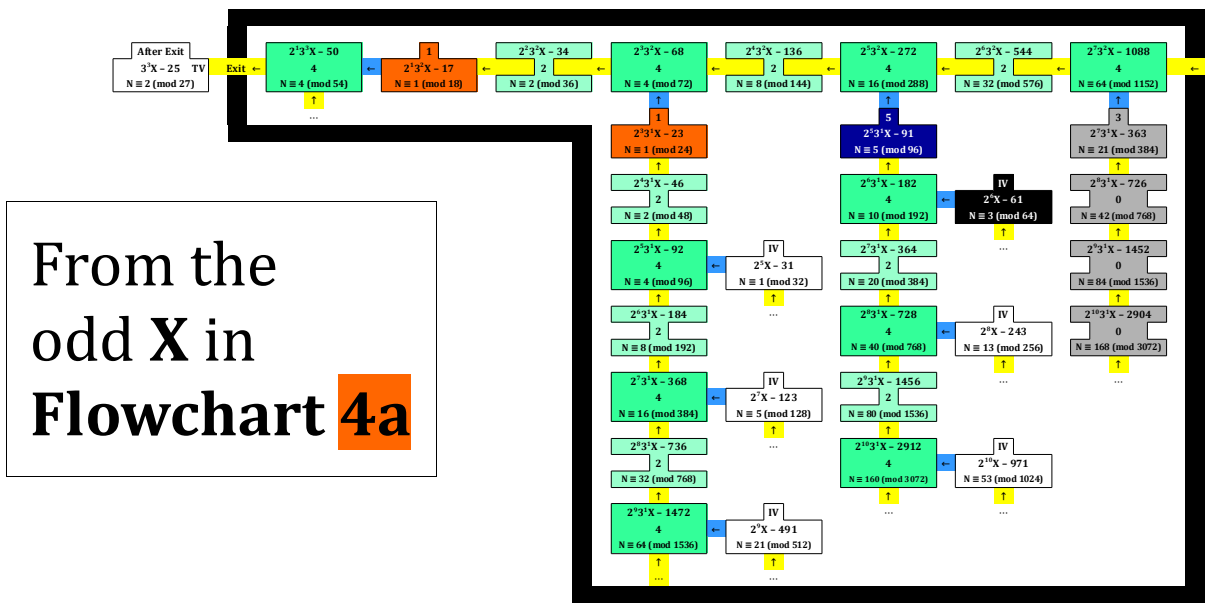

**5c**. Extrapolated from $N \equiv 1 \pmod{2^3 3^1}$
**TV**-class $[3^3$**X**$-25] \Leftrightarrow N \equiv$ **2** $\pmod{3^3}$

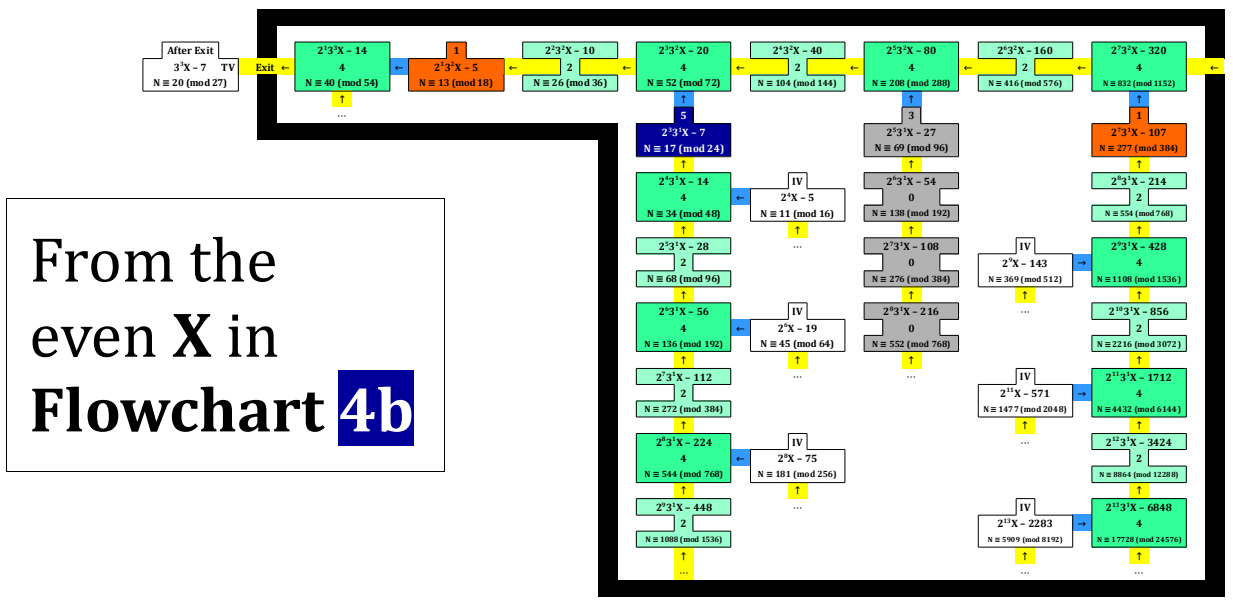

**5d**. Extrapolated from $N \equiv 17 \pmod{2^3 3^1}$
**TV**-class $[3^3$**X**$-7] \Leftrightarrow N \equiv$ **20** $\pmod{3^3}$

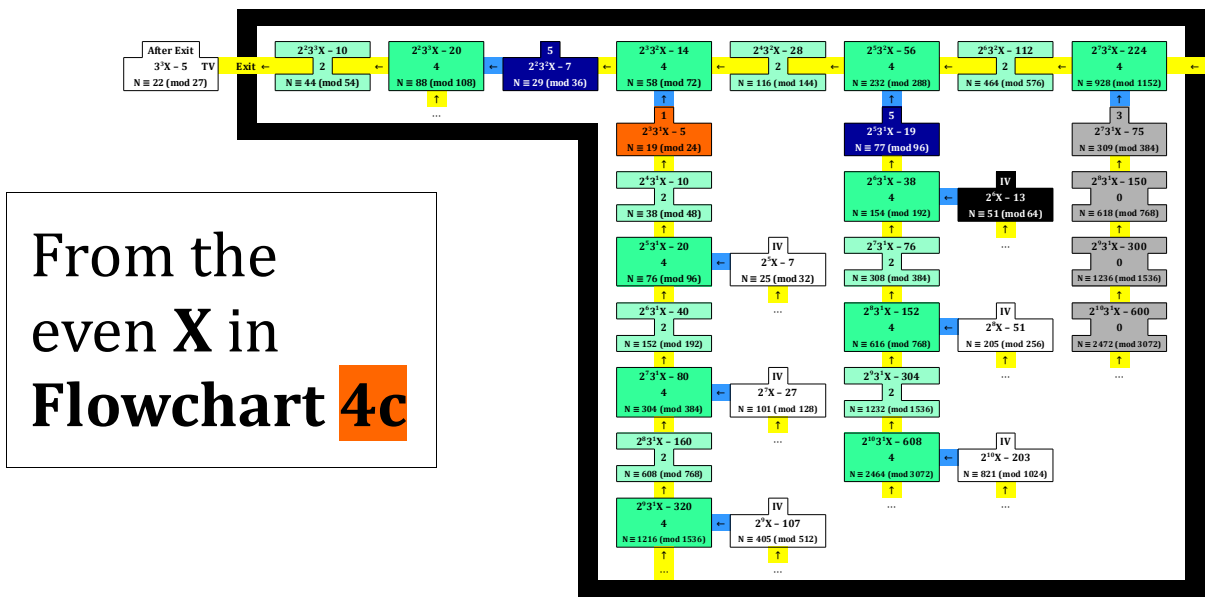

**5e**. Extrapolated from $N \equiv 19 \pmod{2^3 3^1}$
**TV**-class $[3^3$**X**$-5] \Leftrightarrow N \equiv$ **22** $\pmod{3^3}$

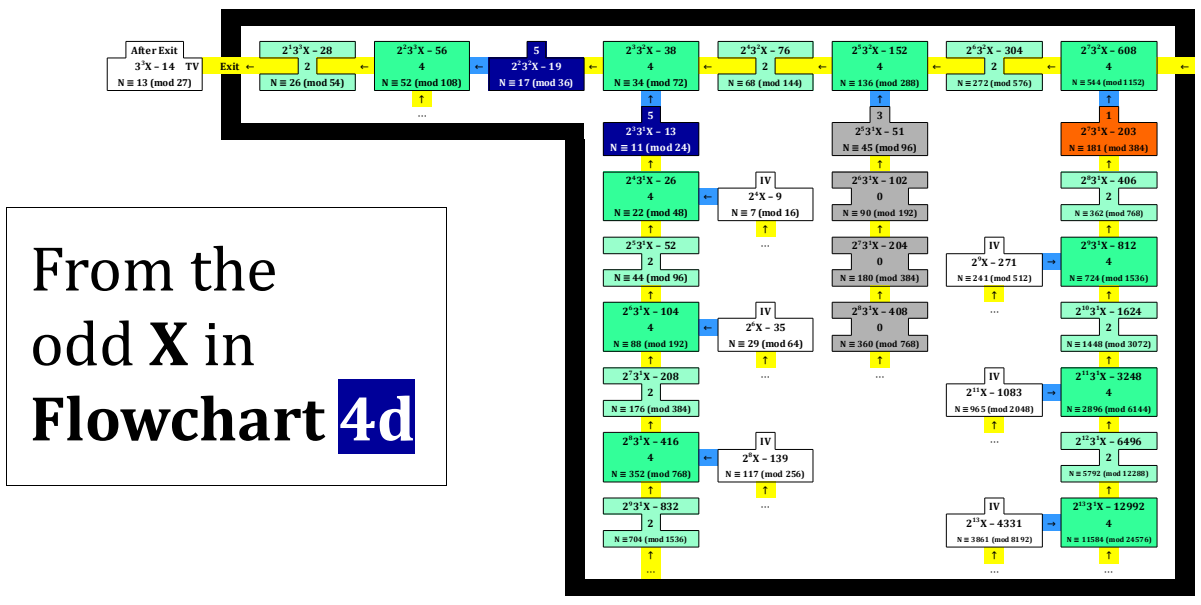

**5f**. Extrapolated from $N \equiv 11 \pmod{2^3 3^1}$
**TV**-class $[3^3$**X**$-14] \Leftrightarrow N \equiv$ **13** $\pmod{3^3}$

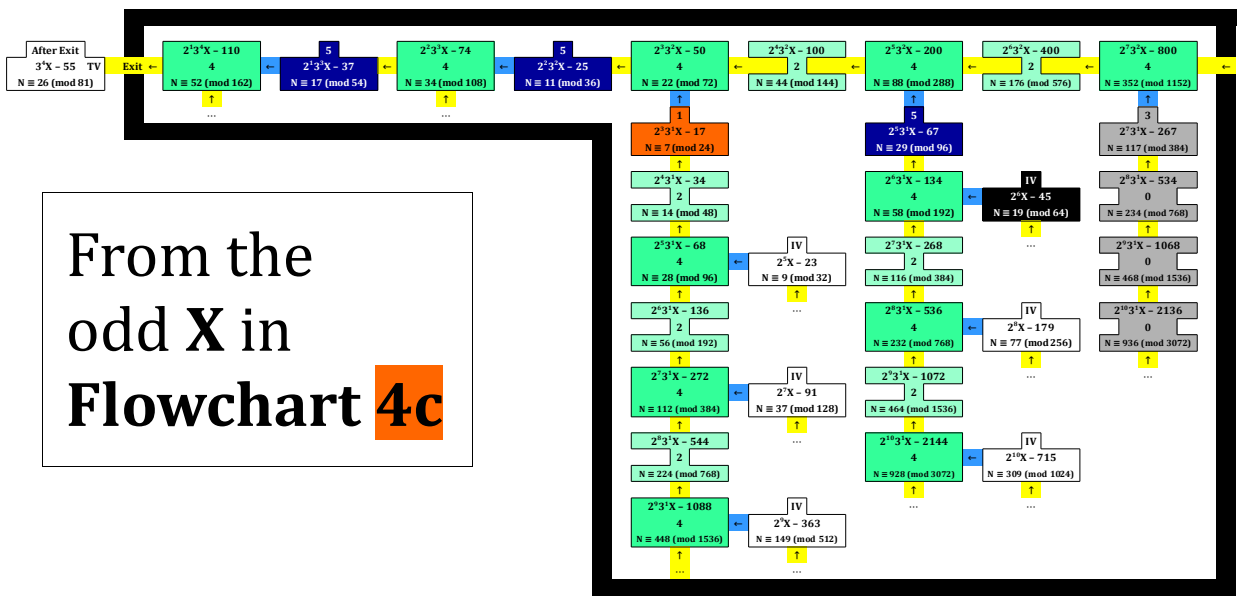

**5g**. Extrapolated from $N \equiv 7 \pmod{2^3 3^1}$
**TV**-class $[3^4$**X**$-55] \Leftrightarrow N \equiv$ **26** $\pmod{3^4}$

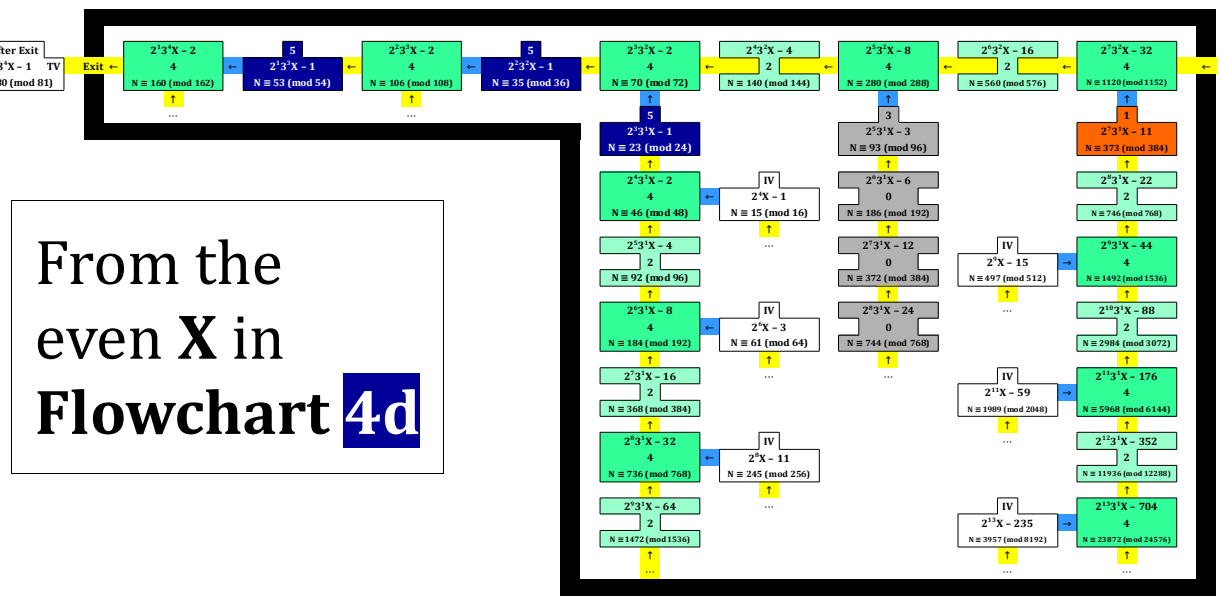

**5h**. Extrapolated from $N \equiv 23 \pmod{2^3 3^1}$
**TV**-class $[3^4$**X**$-1] \Leftrightarrow N \equiv$ **80** $\pmod{3^4}$

Notice that the **modulo** $A_{TV} = 3^s$ depends on the different **extensions** of the parity-sequences. It is possible to gain information about the **extensions** of the parity-sequences when **modulo** A is doubled to 2A hence information about **Converging* Class-series for **subsets** of $^{\Phi}[2^2$**X**$-1]_{IV}$.

**Comment**

From observations in the shown Flowcharts it is possible to make some interesting (and some very important) deductions/conclusions about the development of Class-series for *increasing* **modulo** A = $2^r3^s$. Notice that in **Flowchart** **5b** the nodes actually contains the same (*identical*) classes [$2^r3^s$**X**-B] as found in the nodes from **Flowchart** **2a** i.e. all N ≡ **1** (mod 6) *connects* to one of the four Type **2** Flowcharts for *quadruple* **modulo** 4A. Also notice that in **Flowchart** **5a** it is observed that all N ≡ **3** (mod 6) *connects* to one of the four Type **1** Flowcharts, **4**A=$2^{r+2}3^s$. The two above observations/facts implies, that it is possible to have **Flowchart** **2b** as the *only* **Start-chart**, as the sixteen (*split*) Type **2** charts for **modulo 16**A=$2^{r+4}3^s$ ***certainly*** includes all N ≡ **5** (mod 6) & for ***certain*** contains (are *connected* to) all N ≡ **1** (mod 6) & all N ≡ **3** (mod 6).

The two Type **2** **Flowcharts** **4b,d** (from **Flowchart** **2b**) contains the Class-series for the two subsets of the $^{\Phi}$**Diverging Domain-class**; $^{\Phi}[2^3\mathbf{X}\text{-}1_0]_{IV} \cup {}^{\Phi}[2^3\mathbf{X}\text{-}5_0]_{IV} = {}^{\Phi}[2^2\mathbf{X}\text{-}1_0]_{IV}$. Please notice the two different **extensions** {**OE** , **E**} of the $^{\Phi}$**Diverging** $^{\Phi}$parity-sequence $^{\Phi}$**R O E - O E - T**.

$^{\Phi}[2^3\mathbf{X}\text{-}1_0]_{IV}$ →→ $[3^3\mathbf{X}\text{-}1_6]_{TV}$  $\mathbf{P}_{IV} = 7$ & $\mathbf{P}_{TV} = 26$  $^{\Phi}$**R O E - O E - O E - T**

$^{\Phi}[2^3\mathbf{X}\text{-}5_0]_{IV}$ →→ $[3^2\mathbf{X}\text{-}5_5]_{TV}$  $\mathbf{P}_{IV} = 3$ & $\mathbf{P}_{TV} = 4$  $^{\Phi}$**R O E - O E - E - T**

The four Type **2** **Flowcharts** **5b,d,f,h** contains the Class-series for the four different subsets of **The Domain-class**; $^{\Phi}[2^4\mathbf{X}\text{-}1_0]_{IV} \cup {}^{\Phi}[2^4\mathbf{X}\text{-}5_0]_{IV} \cup {}^{\Phi}[2^4\mathbf{X}\text{-}9_0]_{IV} \cup {}^{*}[2^{*4_2}\mathbf{X}\text{-}{*13}_{Con}]_{IV} = {}^{\Phi}[2^2\mathbf{X}\text{-}1_0]_{IV}$.
Please notice the four different **extensions** of the parity-sequences {**OE-OE, E-OE, OE-E, E-E**}

$^{\Phi}[2^4\mathbf{X}\text{-}1_0]_{IV}$ →→ $[3^4\mathbf{X}\text{-}1_8]_{TV}$  $\mathbf{P}_{IV} = 15$ & $\mathbf{P}_{TV} = 80$  $^{\Phi}$**R O E - O E - O E - O E - T**

$^{\Phi}[2^4\mathbf{X}\text{-}5_0]_{IV}$ →→ $[3^3\mathbf{X}\text{-}7_7]_{TV}$  $\mathbf{P}_{IV} = 11$ & $\mathbf{P}_{TV} = 20$  $^{\Phi}$**R O E - O E - E - O E - T**

$^{\Phi}[2^4\mathbf{X}\text{-}9_0]_{IV}$ →→ $[3^3\mathbf{X}\text{-}14_7]_{TV}$  $\mathbf{P}_{IV} = 7$ & $\mathbf{P}_{TV} = 13$  $^{\Phi}$**R O E - O E - O E - E - T**

$*_3[2^{*4_2}\mathbf{X}\text{-}{*13}_{Con}]_1$ →→ $[3^2\mathbf{X}\text{-}7_6]_{TV}$  $\mathbf{P}_{IV} = 3$ & $\mathbf{P}_{TV} = *\mathbf{2}$  $^{*}$**R O E - O E - E - E - T**

Notice that $*_3[2^{*4_2}\mathbf{X}\text{-}{*13}_{Con}]_1 = \{\mathbf{3}_1,19_2,35_3, ...\}$ is a subset of $^{\Phi}[2^3\mathbf{X}\text{-}5_0]_{IV} = \{\mathbf{3}_1,11_2,19_3,27_4,35_5, ...\}$

Here is shown the first five **Converging* **IV**-classes and the **reducing* **parity-sequences*.

$*_{u=1}[2^{*1_0}\mathbf{X}\text{-}{*0}_{Con}]_{q=1}$ →→ $[3^0\mathbf{X}\text{-}0_1]_{TV}$  $\mathbf{P}_{IV} = 2$ & $\mathbf{P}_{TV} = *\mathbf{1}$  $^{*}$**R E - T**

$*_2[2^{*2_1}\mathbf{X}\text{-} \ {*3}_{Con}]_1$ →→ $[3^1\mathbf{X}\text{-}2_3]_{TV}$  $\mathbf{P}_{IV} = 1$ & $\mathbf{P}_{TV} = *\mathbf{1}$  $^{*}$**R O E - E - T** (**X** = 1 ⇒ ***The Loop***)

$*_3[2^{*4_2}\mathbf{X}\text{-}{*13}_{Con}]_1$ →→ $[3^2\mathbf{X}\text{-}7_6]_{TV}$  $\mathbf{P}_{IV} = 3$ & $\mathbf{P}_{TV} = *\mathbf{2}$  $^{*}$**R O E - O E - E - E - T**

$*_4[2^{*5_3}\mathbf{X}\text{-}{*21}_{Con}]_1$ →→ $[3^3\mathbf{X}\text{-}17_8]_{TV}$  $\mathbf{P}_{IV} = 11$ & $\mathbf{P}_{TV} = *\mathbf{10}$  $^{*}$**R O E - O E - E - O E - E - T**

$*_5[2^{*5_3}\mathbf{X}\text{-} \ {*9}_{Con}]_2$ →→ $[3^3\mathbf{X}\text{-}7_8]_{TV}$  $\mathbf{P}_{IV} = 23$ & $\mathbf{P}_{TV} = *\mathbf{20}$  $^{*}$**R O E - O E - O E - E - E - T**

Notice that $*_4[2^{*5_3}\mathbf{X}\text{-}{*21}_{Con}]_1$ is a subset of $^{\Phi}[2^4\mathbf{X}\text{-}5_0]_{IV}$ & $*_5[2^{*5_3}\mathbf{X}\text{-}{*9}_{Con}]_2$ is a subset of $^{\Phi}[2^4\mathbf{X}\text{-}9_0]_{IV}$.

$^{\Phi}$**Diverging IV**-classes for **modulo** $A_{IV} = 2^5$ contains all $N_X$ **not** included in a $*_{u\leq5}[2^{*r_s}\mathbf{X}\text{-}{*B}_{Con}]_q$:

$^{\Phi}[2^5\mathbf{X}\text{-}1_0]_{IV}$ →→ $[3^5\mathbf{X}\text{-}1_{10}]_{TV}$  $\mathbf{P}_{IV} = 31$ & $\mathbf{P}_{TV} = 242$  $^{\Phi}$**R O E - O E - O E - O E - O E - T**

$^{\Phi}[2^5\mathbf{X}\text{-}5_0]_{IV}$ →→ $[3^4\mathbf{X}\text{-}10_9]_{TV}$  $\mathbf{P}_{IV} = 27$ & $\mathbf{P}_{TV} = 71$  $^{\Phi}$**R O E - O E - E - O E - O E - T**

$^{\Phi}[2^5\mathbf{X}\text{-}17_0]_{IV}$ →→ $[3^4\mathbf{X}\text{-}41_9]_{TV}$  $\mathbf{P}_{IV} = 15$ & $\mathbf{P}_{TV} = 40$  $^{\Phi}$**R O E - O E - O E - O E - E - T**

$^{\Phi}[2^5\mathbf{X}\text{-}25_0]_{IV}$ →→ $[3^4\mathbf{X}\text{-}61_9]_{TV}$  $\mathbf{P}_{IV} = 7$ & $\mathbf{P}_{TV} = 20$  $^{\Phi}$**R O E - O E - O E - E - O E - T**

These four $^{\Phi}[2^5\mathbf{X}\text{-}B_0]_{IV}$ contains **subsets** of the three $^{\Phi}$**Diverging IV**-classes $^{\Phi}[2^4\mathbf{X}\text{-}B_0]_{IV}$, $A_{IV} = 2^4$.

## The Collatz-constant Zeta

In $^*$*Converging* $^*$Class-series $\;^*{}_u[2^{*r_s}\mathbf{X}\text{-}^*B_{Con}]_q \to\to [3^s\mathbf{X}\text{-}B_j]_{TV}$ it **is** true that; $\tfrac{1}{2} \le 3^s/2^{*r_s} < \mathbf{1}$
For example the $^*$*reducing* $^*$Class-series $\;^*{}_3[2^{*4_2}\mathbf{X}\text{-}^*13_{Con}]_1 \to\to [3^2\mathbf{X}\text{-}7_6]_{TV} \Rightarrow \mathbf{P}_{IV} = 3_{IV} \;\wedge\; ^*2_6 = {}^*\mathbf{P}_{TV}$
where the ratio $A_{TV}/A_{IV} = 3^2/2^{*4_2}$ is 0.5625 i.e. it is true for the $^*$Class-series that $\tfrac{1}{2} < 3^2/2^{*4_2} < \mathbf{1}$

In $^\dagger$**End**-Class-series $^\dagger[2^d\mathbf{X}\text{-}(2^d\text{-}N_N)_{End}]_N \to\to [3^{s_d}\mathbf{X}\text{-}(3^{s_d}\text{-}\mathbf{1})_j]_{TV}$, $\mathbf{d} \in \mathbb{N}_0$, $\mathbf{s}_d \in \mathbb{N}_0$; $\tfrac{1}{2} < N_N 3^{s_d}/2^d \le \mathbf{1}$
For example in the $^\dagger$**End**-Class-series $^\dagger[2^5\mathbf{X}\text{-}29_{End}]_3 \to\to [3^{2_d}\mathbf{X}\text{-}8_7]_{TV} \;\Rightarrow\; \mathbf{P}_{IV} = 3_{IV} \;\wedge\; ^\dagger1_7 = {}^\dagger\mathbf{P}_{TV}$
where the ratio $A_{TV}/A_{IV} = 3^{2_d}/2^5$ is 0.28125 $\Rightarrow 3_{IV}3^{2_d}/2^5 = 0.84375$ i.e. $\tfrac{1}{2} < 3_3 3^{2_d}/2^5 < \mathbf{1}$ is true.

Notice that $^\dagger$**End-IV**-classes are *proper subsets* of $^*$*reducing* **IV**-classes & $B_{End} \equiv {}^*B_{Con} \pmod{2^{*r_s}}$.

**OBS**! From the iterationseries for $N_0 = 7$ (p. 12) it can be *deduced* that 7 belongs to these Class-series: $^*{}_6[2^{*7_4}\mathbf{X}\text{-}^*(2^7\text{-}7_0)_{Con}]_1 \to\to [3^4\mathbf{X}\text{-}(3^4\text{-}^*5_{11})]_{TV}$ **and** $^\dagger[2^{11}\mathbf{X}\text{-}(2^{11}\text{-}7_7)_{End}]_7 \to\to [3^{5_d}\mathbf{X}\text{-}(3^{5_d}\text{-}\mathbf{1})_{16}]_{TV}$
$^\dagger$**R O E - O E - O E - E - O E - E - E - O E - E - E - E - D**

## Theorem 2

For **s** any number of Type **O** $_{\text{Operation}}$ there **exist** a $^*$*reducing* number $^*\mathbf{r}_s$ of Type **E** $_{\text{Operation}}$

## Proof for Theorem 2

The proof comes from observations of $^*$Class-series for $^*$*Converging* **I**nitial **V**alue classes:

$\boxed{^*{}_u[2^{*r_s}\mathbf{X}\text{-}^*B_{Con}]_q > [3^s\mathbf{X}\text{-}B_{j=s+^*r}]}$ where it is true that $\boxed{2^{*r_s} > 3^s}$ i.e. that $\boxed{3^s/2^{*r_s} < 1}$

*Remember that* $A = (\mathbf{P}+\mathbf{B})$ *which is true for all the* **IV**-*classes* $2^r = (\mathbf{P}_{IV}+B_0)$ ***and*** *for all* **TV**-*classes* $3^s = (\mathbf{P}_{TV}+B_{j=s+r})$ *so* ***when*** $A_{TV} = A_{IV}3^s/2^{*r_s}$ *it* ***is*** *true that* $\boxed{(\mathbf{P}_{TV}+B_{j=s+^*r}) = (\mathbf{P}_{IV}+{}^*B_{Con})3^s/2^{*r_s}}$ *which* ***implies*** $\boxed{\mathbf{P}_{TV} \approx \mathbf{P}_{IV}3^s/2^{*r_s} \;\wedge\; B_{j=s+^*r} \approx {}^*B_{Con}3^s/2^{*r_s}}$ *and* ***this*** *delivers the desired* ***criterion*** *for existence of* $^*$*Converging* **IV**-classes $^*{}_u[2^{*r_s}\mathbf{X}\text{-}^*B_{Con}]_q$ *which* ***is*** *that* $\boxed{2^{*r_s-1} < 3^s < 2^{*r_s}}$ *is true for* $\mathbf{s} \in \mathbb{N}$

Suppose that $\mathbf{c} \in \mathbb{R}$ and that $\mathbf{s} \in \mathbb{N}$ is number of type **Odd**-operations $\mathbf{s} = \#$ **O**:

$2^c = 3^s \quad\Leftrightarrow\quad log(2^c) = log(3^s) \quad\Leftrightarrow\quad \mathbf{c}\,log\mathbf{2} = \mathbf{s}\,log\mathbf{3} \quad\Leftrightarrow\quad \mathbf{c} = \mathbf{s}\,log\mathbf{3}/log\mathbf{2}$

In this case $2^{\lceil c\rceil-1} < 3^s < 2^{\lceil c\rceil}$ is true when $^*\mathbf{r}_s = \lceil\mathbf{c}\rceil$ is true, i.e. $\boxed{2^{*r_s-1} < 3^s < 2^{*r_s}}$ is true.

**NB**! Here $\lceil\mathbf{c}\rceil$ is the *ceiling*-function (round **up** to nearest integer) e.g. $\lceil 15.84\rceil = 16$.

*Notice the* $^*$*Reducing Combinations* $(\mathbf{s},{}^*\mathbf{r}_s)$ in $3^0 < 2^{*1} < 3^1 < 2^{*2} < 2^3 < 3^2 < 2^{*4} < 3^3 < 2^{*5} < 2^6 < \ldots$
*which is "A Never Ending Sequence"* ***so*** $^*{}_u[2^{*r_s}\mathbf{X}\text{-}^*B_{Con}]_q$ ***exist*** *when* $^*\mathbf{r}_s = \lceil\mathbf{s}\,{}^{log3}/_{log2}\rceil$ *is true* ■

**OBS**! ***If*** $2^{*r-1} < 3^s < 2^{*r}$ *is true* ***then*** $2^{*r} < 3^{s+1} < 2^{*r+2}$ *is true and* ***sometimes*** $2^{*r} < 3^{s+1} < 2^{*r+1}$ *is true.*

## Corollary 3

The estimate $\mathbf{P}_{TV} \approx \mathbf{P}_{IV}3^s/2^{*r_s}$ is a little **to low** and *v.v.* the estimate $B_{j=s+^*r} \approx B_{Con}3^s/2^{*r_s}$ is a little **to high** $\Rightarrow$ $\boxed{B_{j=s+r} < B_0 \Leftrightarrow \textit{reducing}\text{ Combination }(\mathbf{s},\mathbf{r}) \;\&\; B_{j=s+r} \ge B_0 \Leftrightarrow \textbf{Diverging}\text{ Combination}}$

## Definition 2

The ratio $log\mathbf{3}/log\mathbf{2}$ is named ***The Collatz-constant* Zeta** $\boxed{\mathbf{Z} = log\mathbf{3}/log\mathbf{2}} \Rightarrow \boxed{^*\mathbf{r}_s = \lceil\mathbf{s}\mathbf{Z}\rceil,\ \mathbf{s} \in \mathbb{N}}$

The equation $\boxed{^*\mathbf{r}_s = \lceil\mathbf{s}\mathbf{Z}\rceil,\ \mathbf{s} \in \mathbb{N},\ \mathbf{Z} = log\mathbf{3}/log\mathbf{2}}$ is named **Formula (1)** for future reference**.**

Please notice that; **Z** is *irrational*, which ensures that **sZ** is *never* an integer for any $\mathbf{s} \in \mathbb{N}$.

Please notice that; $\boxed{\mathbf{Z} = log_{10}\mathbf{3}/log_{10}\mathbf{2} = ln_e\mathbf{3}/ln_e\mathbf{2} = log_2\mathbf{3}/log_2\mathbf{2} = log_2\mathbf{3} = \mathbf{1.5849625007...}}$

**NB**! *The ratio* $log\mathbf{3}/log\mathbf{2}$ *is well-known in relation to parity-sequences for congruence-classes e.g.* R. Terras[III] *uses the ratio to calculate densities for* **Diverging** *parity-sequences for* # **O** $\to \infty$.
**Formulas** *shown in the present work are based on the same method i.e. delivers identical results.*

**Corollary 4**

From **Definition 2** follows, that it is possible to **split** all $\mathbf{r} \in \mathbb{N}_0$ into the two **proper subsets**:
$^*\mathbf{r}_s \in \{1_0,2_1,4_2,5_3, ...\}$ where $^*\mathbf{r}_s$ = [**sZ**] is true & $^\Phi\mathbf{r}_v \in \{0_0,3_1,6_2,9_3, ...\}$ where $^\Phi\mathbf{r}_v \neq$ [**sZ**] is true.

**Comment**

The **Total** number of **IV**-classes for constant **modulo** $2^r$ in $[2^r\mathbf{X}\text{-}B_0]_{IV}$ is the **Sum** of the **three** *different* **types** of **IV**-classes; $|^\Phi\mathbf{Div}(\mathbf{r})| + |^*Con(\mathbf{r})| + |^\Delta\mathbf{Red}(\mathbf{r})| = 2^r$

Notice that $|^*Con(^\Phi\mathbf{r}_v)| = 0$ & $|^*Con(^*\mathbf{r}_s)| > 0$. The **Group** is named $\mathbf{G(s,{}^*r)} = |^*Con(^*\mathbf{r}_s)|$

**Consequences of splitting IV-classes**

All classes in the used **model** have the form $[2^r3^s\mathbf{X}\text{-}B_j]$ and for **TV**-classes *parity* depends on **X**.
All (*Parent*-)classes $[2^r3^s\mathbf{X}\text{-}B_j]$ are ***split*** into two *Child*-classes $[2^r3^s(\mathbf{2X\text{-}1})\text{-}B_j] \cup [2^r3^s(\mathbf{2X})\text{-}B_j]$

| | | | |
|---|---|---|---|
| $[2^r3^s(\mathbf{2X\text{-}1})\text{-}B_j] \Leftrightarrow [2^{r+1}3^s\mathbf{X}\text{-}(2^r3^s\mathbf{+}B_j)_j]$ | **and** | $[2^r3^s(\mathbf{2X})\text{-}B_j] \Leftrightarrow [2^{r+1}3^s\mathbf{X}\text{-}B_j]$ | **All** |
| $[2^r(\mathbf{2X\text{-}1})\text{-}B_0] \Leftrightarrow [2^{r+1}\mathbf{X}\text{-}(2^r\text{+}B_0)_0]$ | **and** | $[2^r(\mathbf{2X})\text{-}B_0] \Leftrightarrow [2^{r+1}\mathbf{X}\text{-}B_0]$ | **IV** |
| $[3^s(\mathbf{2X\text{-}1})\text{-}B_j] \Leftrightarrow [2^13^s\mathbf{X}\text{-}(3^s\text{+}B_j)_j]$ | **and** | $[3^s(\mathbf{2X})\text{-}B_j] \Leftrightarrow [2^13^s\mathbf{X}\text{-}B_j]$ | **TV** |

Notice that after *splitting* a **TV**-class it is **always possible** to perform ***operations*** on the two new *Child*-classes, because the **modulo** is now even, **but** the *Child*-classes are extended with **different** operations **O E** or **E** because $[2^13^s\mathbf{X}\text{-}(3^s\text{+}B_j)_j]$ & $[2^13^s\mathbf{X}\text{-}B_j]$ have **opposite** *parity*.
Notice that *Children* of $[2^r\mathbf{X}\text{-}B_0]_{IV}$ **are IV**-classes, but *Children* of $[3^s\mathbf{X}\text{-}B_j]_{TV}$ are **not TV**-classes, and that it is the **P**-value (and all **odd X**) in the *Parent*-class that gets the ***new*** higher B-value in **one** *Child*-class, while it is $N_2$ (and all **even X**) that keeps the ***same*** $B_j$ in the **other** *Child*-class, so it is the *parity* of $B_j$ in **TV**-classes that decides the extensions **OE**/**E** of the parity-sequences.
Even $\mathbf{P}_{TV}$ (i.e. Odd $B_{TV}$) means the **odd X** in the *Parent*-class gets the **E**-**tile** extension (and v.v.).

It is a very important point, that ***either*** the **P** ***or*** the B is "lucky" and gets the **E**-**tile** extension.

When the **modulo** is doubled all the (*Parent*-)Class-series $[2^r\mathbf{X}_a\text{-}B_0]_{IV} \rightarrow\rightarrow [3^s\mathbf{X}_a\text{-}B_j]_{TV}$ are ***split***:

$[2^{r+1}\mathbf{X}_b\text{-}B_0]_{IV} \rightarrow\rightarrow [2^13^s\mathbf{X}_b\text{-}B_j]$ (even $\mathbf{X}_a$) **and** $[2^{r+1}\mathbf{X}_c\text{-}(2^r\text{+}B_0)_0] \rightarrow\rightarrow [2^13^s\mathbf{X}_c\text{-}(3^s\mathbf{+}B_j)_j]$ (odd $\mathbf{X}_a$)

For **Odd** $B_j$ in the *Parent*-**TV**-class:

$[2^{r+1}\mathbf{X}_b\text{-}B_0]_{IV} \rightarrow\rightarrow [2^13^s\mathbf{X}_b\text{-}B_j]$ ↗ $[2^13^{s+1}\mathbf{X}_b\text{-}(3B_j\text{-}1)_{j+1}]$ ↘ $[3^{s+1}\mathbf{X}_b\text{-}\mathbf{(}(3B_j\text{-}1)/2\mathbf{)}_{j+2}]_{TV}$

$[2^{r+1}\mathbf{X}_c\text{-}(2^r\text{+}B_0)_0]_{IV} \rightarrow\rightarrow [2^13^s\mathbf{X}_c\text{-}(3^s\text{+}B_j)_j]$ ↘ $[3^s\mathbf{X}_c\text{-}\mathbf{(}(3^s\text{+}B_j)/2\mathbf{)}_{j+1}]_{TV}$ (**E**-**tile** for odd $\mathbf{X}_a$ & "lucky" **P**)

For **Even** $B_j$ in the *Parent*-**TV**-class:

$[2^{r+1}\mathbf{X}_b\text{-}B_0]_{IV} \rightarrow\rightarrow [2^13^s\mathbf{X}_b\text{-}B_j]$ ↘ $[3^s\mathbf{X}_b\text{-}(B_j/2)_{j+1}]_{TV}$ (**E**-**tile** for even $\mathbf{X}_a$ & "lucky" B)

$[2^{r+1}\mathbf{X}_c\text{-}(2^r\text{+}B_0)_0]_{IV} \rightarrow\rightarrow [2^13^s\mathbf{X}_c\text{-}(3^s\text{+}B_j)_j]$ ↗ $[2^13^{s+1}\mathbf{X}_c\text{-}\mathbf{(}3(3^s\text{+}B_j)\text{-}1\mathbf{)}_{j+1}]$ ↘ $[3^{s+1}\mathbf{X}_c\text{-}\mathbf{(}(3^{s+1}\text{-}1\text{+}3B_j)/2\mathbf{)}_{j+2}]_{TV}$

Notice that when $\mathbf{P}_{TV}$ is "very even" i.e. $\mathbf{P}_{TV} = (2C\text{-}1)2^r$ for a "high **r**", then **P** is "lucky" **r** times in a row, while B is "unlucky" **r** times in a row, before the odd $\mathbf{P}_{TV} = (2C\text{-}1)$ is reached. Similarly when $B_j$ is "very even" i.e. $B_j = (2C\text{-}1)2^r$ for a "high **r**", then B is "lucky" **r** times in a row, while **P** is "unlucky" **r** times in a row, before the odd $B_{j+r} = (2C\text{-}1)$ is reached. Remember that every time $A_{IV}$ is doubled from $2^r$ to $2^{r+1}$, all the $\mathbf{P} \leq 2^r$ are *paired* with a *new* B-value; $B = 2^{r+1}\text{-}\mathbf{P}$, where $2^r \leq B < 2^{r+1}$, and all the $B < 2^r$ are *paired* with a *new* **P**-value; $\mathbf{P} = 2^{r+1}\text{-}B$, where $2^r < \mathbf{P} \leq 2^{r+1}$.

Be aware that the ***pair*** {**P**,B} for $\mathbf{P}+B = 2^r$ **share** a parity-sequence with **r** pcs of **tiles** {**O E**, **E**}.

While the N-values (the **P**-values) follows the rules in the **3N plus 1 Problem**, the B-values are observed to follow the rules in the **3N minus 1 Problem**:

- For **r** > 0 and Odd $B_j$; $B_{j+1} = (3B_j-1)$ **and** $B_{j+2} = (3B_j-1)/2$
- For **r** > 0 and Even $B_j$; $B_{j+1} = (B_j/2)$

Iterationseries $N_0 \rightarrow N_1 \rightarrow N_2 \rightarrow \ldots \rightarrow N_{j=s+r} \rightarrow \ldots$ for different $N_0 \in \mathbb{N}$ are able to reveal which **Converging* **IV**-class $N_0$ belongs to, when it is observed that $*N_{*j=s+*r_s} < N_0$. When it is observed that $N_j = {}^{\dagger}1_{j=s_d+d}$ it is confirmed that $N_0$ belongs to ${}^{\dagger}[2^d\mathbf{X}\text{-}(2^d\text{-}N_N)_{End}]_N \rightarrow\rightarrow [3^{s_d}\mathbf{X}\text{-}(3^{s_d}\text{-}\mathbf{1})_j]_{TV}$

Class-series $[2^r\mathbf{X}\text{-}B_0]_{IV} \rightarrow\rightarrow [3^s\mathbf{X}\text{-}B_j]_{TV}$ are able to reveal which **Converging* **IV**-class $B_0$ belongs to, when it is observed that $\mathbf{B}_{*j=s+*r_s} < B_0 \Rightarrow B_0 \equiv *B_{Con} \pmod{2^{*r_s}}$. For example for $B_0 = 9$:

| | | | |
|---|---|---|---|
| $[2^{r\text{-}0}3^{s+0}\mathbf{X} - 9_0]_{IV}$ | ⇄ O | $[2^{r\text{-}0}3^{s+1}\mathbf{X} - 26_1]$ | ⇄ E |
| $[2^{r\text{-}1}3^{s+1}\mathbf{X} - 13_2]$ | ⇄ O | $[2^{r\text{-}1}3^{s+2}\mathbf{X} - 38_3]$ | ⇄ E |
| $[2^{r\text{-}2}3^{s+2}\mathbf{X} - 19_4]$ | ⇄ O | $[2^{r\text{-}2}3^{s+3}\mathbf{X} - 56_5]$ | ⇄ E |
| | | $[2^{r\text{-}3}3^{s+3}\mathbf{X} - 28_6]$ | ⇄ E |
| | | $[2^{r\text{-}4}3^{s+3}\mathbf{X} - 14_7]$ | ⇄ E |
| $[2^{r\text{-}5}3^{s+3}\mathbf{X} - 7_8]_{TV}$ | **Reducing Combination* $(\mathbf{s},*\mathbf{r}_s) = (3,*5_3)$ | | |

It is observed that $\mathbf{7}_8 < 9_0$ **so** $B_0 = 9$ belongs to the **reducing* *Class-series;

$*_5[2^{*5_3}\mathbf{X}\text{-}*9_{Con}]_2 \rightarrow\rightarrow [3^3\mathbf{X}\text{-}\mathbf{7}_8]_{TV}$

Whenever it is observed that $B_{j=s+r} < B_0$, it can be concluded that $3^s < 2^r$ *and* $\mathbf{P}_{TV} < \mathbf{P}_{IV}$.

**Comment**

Notice that it has been shown, that there exist an infinity of **Converging* **IV**-classes, as all the **Reducing Combinations* $(\mathbf{s},*\mathbf{r}_s)$ results in a non-empty ***set*** of $*_u[2^{*r_s}\mathbf{X}\text{-}*B_{Con}]_q$, $\mathbf{q} \le \mathbf{q}_{Max}$, which all share the same **modulo** $2^{*r_s}$. What remains to be shown is, that **all** N ∈ ℕ belongs to one of the $*_{u\rightarrow\infty}[2^{*r_s}\mathbf{X}\text{-}*B_{Con}]_{q\rightarrow\infty}$ i.e. that **no** $N_{ce} \in \mathbb{N}$ *can* exist, which do **not** belong to a **Converging* **IV**-class.

In the below **Figure** the top left square is the **Universe-class** $N_X \equiv \mathbf{1} \pmod 1 \Leftrightarrow {}^{\Phi}[2^0\mathbf{X}\text{-}0_0]_{IV/TV}$
The bottom left square is the **Domain-class** $N_X \equiv \mathbf{3} \pmod{2^2} \Leftrightarrow {}^{\Phi}[2^2\mathbf{X}\text{-}1_0]_{IV}$
The large square contains the four ${}^{\Phi}$**Diverging IV**-classes found in the bottom left square.
*The* ***Figure*** *is illustrating that* $*_u[2^{*r_s}\mathbf{X}\text{-}*B_{Con}]_q$ *eventually Covers half of the remaining* ${}^{\Phi}[2^r\mathbf{X}\text{-}B_0]_{IV}$

The **Universe-class** ${}^{\Phi}[2^0\mathbf{X}\text{-}0_0]_{IV}$

| | |
|---|---|
| $*_2[2^2X\text{-}*3]_1$ | $[2^2X\text{-}1]$ |
| *N ≡ 1 (mod 4) | N ≡ 3 (mod 4) |
| B ≡ *3 (mod 4) | B ≡ 1 (mod 4) |
| N ∈ {1,5,9,13, ...} | N ∈ {3,7,11,15, ...} |
| $*_1[2^1X\text{-}*0]_1$ | |
| *N ≡ 2 (mod 2) | |
| B ≡ *0 (mod 2) | |
| N ∈ {2,4,6,8,10, ...} | |

Remember that for **IV**-classes $(P_{IV} + B_{IV}) = 2^r$

The **Domain-class** ${}^{\Phi}[2^2\mathbf{X}\text{-}1_0]_{IV}$

| | |
|---|---|
| N ≡ 7 (mod 32) | N ≡ 31 (mod 32) |
| B ≡ 25 (mod 32) | B ≡ 1 (mod 32) |
| N ≡ 15 (mod 32) | N ≡ 27 (mod 32) |
| B ≡ 17 (mod 32) | B ≡ 5 (mod 32) |
| *N ≡ 11 (mod 32) | $*_3[2^4X\text{-}*13]_1$ |
| B ≡*21 (mod 32) | *N ≡ 3 (mod 16) |
| *N ≡ 23 (mod 32) | B ≡ *13 (mod 16) |
| B ≡ *9 (mod 32) | *N ∈ {3,19,35, ...} |

| | | | |
|---|---|---|---|
| *N ≡ 7 (mod $2^7$) | N ≡ 71 (mod $2^8$) | N ≡ 31 (mod $2^8$) | N ≡ 127 (mod $2^8$) |
| B ≡*121 (mod $2^7$) | P = 583 | P = 287 | |
| | N ≡ 167 (mod $2^8$) | N ≡ 63 (mod $2^8$) | N ≡ 159 (mod $2^8$) |
| | P = 423 | P = 575 | |
| *N ≡ 39 (mod $2^8$) | N ≡ 231 (mod $2^8$) | N ≡ 223 (mod $2^8$) | *N ≡ 95 (mod $2^8$) |
| B ≡ *217 (mod $2^8$) | P = 999 | P = 735 | B ≡ *161 (mod $2^8$) |
| *N ≡ 199 (mod $2^8$) | N ≡ 103 (mod $2^8$) | N ≡ 191 (mod $2^8$) | N ≡ 255 (mod $2^8$) |
| B ≡ *57 (mod $2^8$) | | | B ≡ 1 (mod $2^8$) |
| *N ≡ 15 (mod $2^7$) | N ≡ 239 (mod $2^8$) | N ≡ 27 (mod $2^8$) | *N ≡ 59 (mod $2^7$) |
| B ≡*113 (mod $2^7$) | | | B ≡ *69 (mod $2^7$) |
| | N ≡ 47 (mod $2^8$) | N ≡ 91 (mod $2^8$) | |
| | P = 815 | P = 347 | |
| *N ≡ 79 (mod $2^8$) | N ≡ 111 (mod $2^8$) | N ≡ 155 (mod $2^8$) | *N ≡ 123 (mod $2^8$) |
| B ≡ *177 (mod $2^8$) | P = 367 | P = 923 | B ≡ *133 (mod $2^8$) |
| *N ≡ 175 (mod $2^8$) | N ≡ 207 (mod $2^8$) | N ≡ 251 (mod $2^8$) | *N ≡ 219 (mod $2^8$) |
| B ≡ *81 (mod $2^8$) | P = 975 | P = 507 | B ≡ *37 (mod $2^8$) |

"Half-way-zooms" on split subsets of the Universe-class N ≡ 1 (mod 1)

| $\lvert *Con(r) \rvert$ |
|---|
| $\lvert *Con(*1_0) \rvert = 1$ |
| $\lvert *Con(*2_1) \rvert = 1$ |
| $\lvert *Con(3) \rvert = 0$ |
| $\lvert *Con(*4_2) \rvert = 1$ |
| $\lvert *Con(*5_3) \rvert = 2$ |
| $\lvert *Con(6) \rvert = 0$ |
| $\lvert *Con(*7_4) \rvert = 3$ |
| $\lvert *Con(*8_5) \rvert = 7$ |
| $\lvert *Con(9) \rvert = 0$ |
| $\lvert *Con(*10_6) \rvert = 12$ |
| $\lvert *Con(11) \rvert = 0$ |
| $\lvert *Con(*12_7) \rvert = 30$ |
| $\lvert *Con(*13_8) \rvert = 85$ |
| $\lvert *Con(14) \rvert = 0$ |
| ... |

*Emperically confirmed*

**Figure 2.** *Illustrating *Converging* **IV**-*classes (literally) Covering more of* ℕ *for* $*\mathbf{r}_s \rightarrow \infty$.

*Notice that for* A = $2^5$ *only four* ${}^{\Phi}$**Div**(5) *exist*; ${}^{\Phi}[2^5\mathbf{X}\text{-}1_0]_{IV} \cup {}^{\Phi}[2^5\mathbf{X}\text{-}5_0]_{IV} \cup {}^{\Phi}[2^5\mathbf{X}\text{-}17_0]_{IV} \cup {}^{\Phi}[2^5\mathbf{X}\text{-}25_0]_{IV}$

**NB**! The (**1**+**1**+**1**+**2**+**3**+**7**+**12**+30+85) = **142** $*Con_u(*\mathbf{r}_s)_q$ can be found in **Appendix Alpha**.

## *The Collatz Tree* as Flowchart

It is well known, that it is **not** possible to prove *The Collatz Conjecture* by extrapolating **values** $N \in \mathbb{N}$ from the **Origo** $N = \mathbf{1!}$ and generate the ***Collatz-Tree*** illustrated in **Figure 1a**. The ***Tree*** obviously contains a *countable infinity* of **Branches** $(2C-1)2^r$, $C \in \mathbb{N}$, $\mathbf{r} \in \mathbb{N}_0$ but it is **not** obvious that the ***Tree*** contains *all possible* **Branches** for **all** odd $(2C-1)$ for **all** $C \in \mathbb{N}$. The below shown **Flowchart 6** is equivalent to the ***Tree*** in **Figure 1a**. *The nodes in* **Flowchart 6** *contains the unique* $^{\dagger}$**End-IV**-*classes for the* **P**-*values in the* **IV**-*classes* $N \equiv \mathbf{P} \pmod{2^d} \Rightarrow \mathbf{P} \in {}^{\dagger}[2^d\mathbf{X}\text{-}(2^d\text{-}\mathbf{P}_P)_{\mathbf{End}}]_{\mathbf{P}}$

**$^{\dagger}$End-IV-classes for the $^{\dagger}$End-Class-series** $^{\dagger}[2^{\mathbf{d}}\mathbf{X}\text{-}(2^{\mathbf{d}}\text{-}N_N)_{\mathbf{End}}]_{\mathbf{N}} \rightarrow\rightarrow [3^{s_d}\mathbf{X}\text{-}(3^{s_d}\text{-}\mathbf{1})_j]_{\mathbf{TV}}, \mathbf{j} = (\mathbf{s}_d + \mathbf{d})$

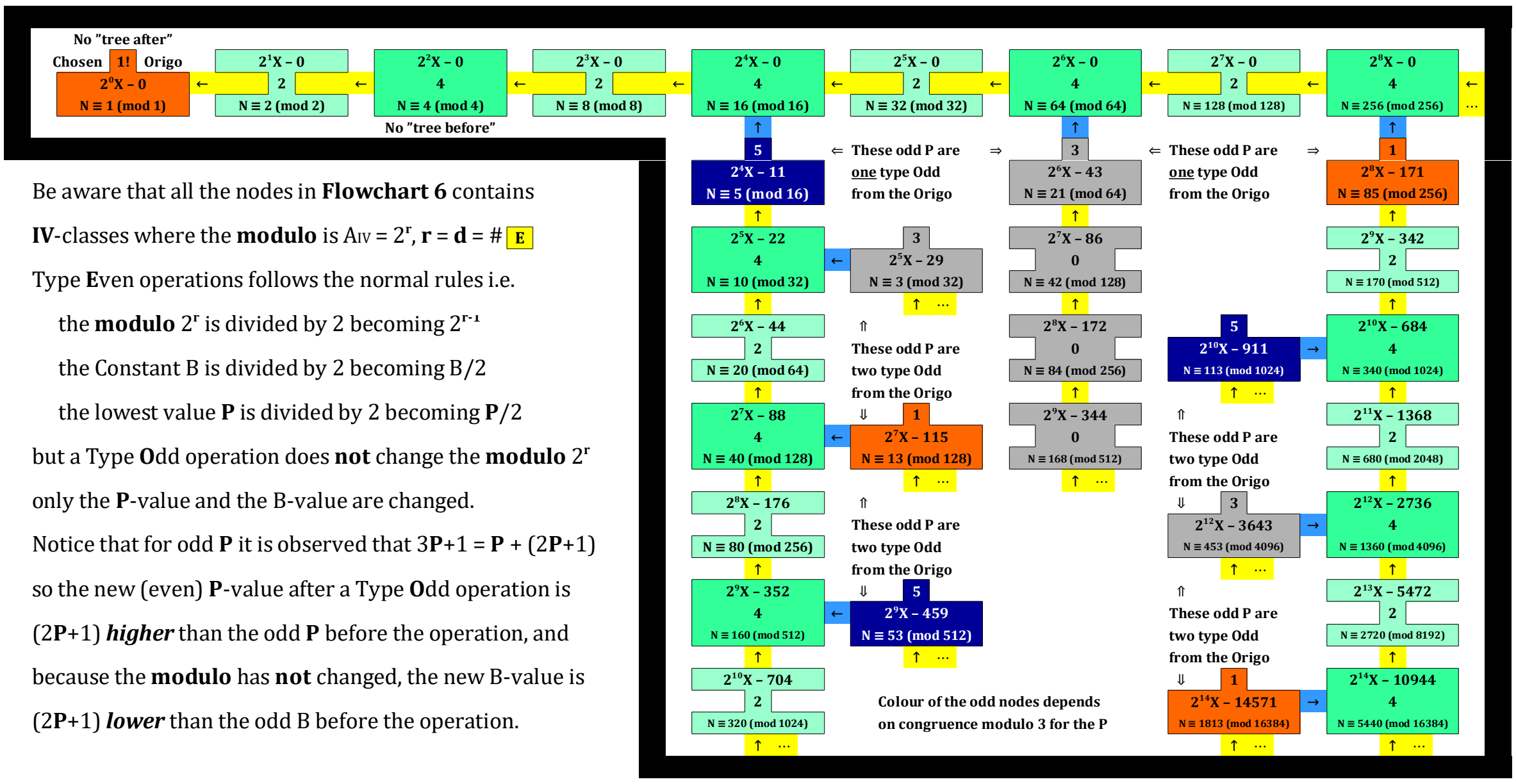

Be aware that all the nodes in **Flowchart 6** contains

**IV**-classes where the **modulo** is $A_{IV} = 2^r$, $\mathbf{r} = \mathbf{d} = \#$ **E**

Type **E**ven operations follows the normal rules i.e.

the **modulo** $2^r$ is divided by 2 becoming $2^{r-1}$

the Constant B is divided by 2 becoming B/2

the lowest value **P** is divided by 2 becoming **P**/2

but a Type **O**dd operation does **not** change the **modulo** $2^r$

only the **P**-value and the B-value are changed.

Notice that for odd **P** it is observed that 3**P**+1 = **P** + (2**P**+1)

so the new (even) **P**-value after a Type **O**dd operation is

(2**P**+1) ***higher*** than the odd **P** before the operation, and

because the **modulo** has **not** changed, the new B-value is

(2**P**+1) ***lower*** than the odd B before the operation.

**Flowchart 6.** The $^{\dagger}$**End-IV**-classes $^{\dagger}[2^{\mathbf{d}}\mathbf{X}\text{-}(2^{\mathbf{d}}\text{-}N_N)_{\mathbf{End}}]_{\mathbf{N}}$ incl. $^{\dagger}[2^{\mathbf{0}}\mathbf{X}\text{-}0_{\mathbf{j=0}}]_{\mathbf{N=1}}$ for **The Origo 1!**

## Summary

An iterationseries exist for all Initial Values $N_{j=0} \in \mathbb{N}$; $N_0 \rightarrow N_1 \rightarrow \ldots \rightarrow N_{j=s+r} \rightarrow \ldots$

If ***Collatz Conjecture*** is true all iterationseries contains a **Reduced Value* and the **End**-value $1_j$.

$$N_0 \rightarrow N_1 \rightarrow \ldots \xrightarrow{}_{E} {}^{*}\mathbf{N}_{*j=s+*r} \rightarrow N_{j=s+*r+1} \rightarrow \ldots \xrightarrow{}_{E} 2_{j=s+d-1} \xrightarrow{}_{E} 1_{j=s+d}$$

A *Class-series* exist for all **I**nitial **V**alue-classes $[2^r\mathbf{X}\text{-}B_{j=0}]_{IV}$, $B_0 < 2^r$; $[2^r\mathbf{X}\text{-}B_0]_{IV} \rightarrow \ldots \rightarrow [3^s\mathbf{X}\text{-}B_{j=s+r}]_{TV}$

If ***Collatz Conjecture*** is true all $N \in \mathbb{N}$ belongs to a **Converging* Class-series.

$${}^{*}_{u}[2^{*r_s}\mathbf{X}\text{-}{}^{*}B_{Con}]_q \rightarrow\rightarrow [3^s\mathbf{X}\text{-}B_{j=s+*r}]_{TV} \text{ for a unique } \mathbf{u} \in \mathbb{N},\ {}^{*}\mathbf{r}_s = {}^{*}[\mathbf{sZ}],\ \mathbf{s} \in \mathbb{N},\ \mathbf{X} \in \mathbb{N},\ {}^{*}B_{Con} \in \mathbb{N}_0.$$

Notice that a **Converging* **IV**-class ${}^{*}_{u}[2^{*r_s}\mathbf{X}\text{-}{}^{*}B_{Con}]_q$ contains the ***pair*** $_u(\mathbf{P}_q + {}^{*}B_{Con}) = 2^{*r_s}$. The ***pair*** implies two **modular *sets*** ${}^{*}N \equiv \mathbf{P}_q \pmod{2^{*r_s}}$ and ${}^{\Delta}B_{\mathbf{Redundant}} \equiv {}^{*}\mathbf{B}_{Con} \pmod{2^{*r_s}}$ for all $2^r > 2^{*r_s}$. Please notice that any **subset** ${}^{\Delta}[2^r\mathbf{X}\text{-}{}^{\Delta}B_{\mathbf{Redundant}}]_{IV}$ of ${}^{*}_{u}[2^{*r_s}\mathbf{X}\text{-}{}^{*}B_{Con}]_q$ is $^{\Delta}$**Redundant** for all $\mathbf{r} > {}^{*}\mathbf{r}_s$.

If ***Collatz Conjecture*** is true all $N_N \in \mathbb{N}$ belongs to a unique $^{\dagger}$**End** Class-series.

$${}^{\dagger}[2^{\mathbf{d}}\mathbf{X}\text{-}(2^{\mathbf{d}}\text{-}N_N)_{\mathbf{End}}]_{\mathbf{N}} \rightarrow\rightarrow [3^{s_d}\mathbf{X}\text{-}(3^{s_d}\text{-}{}^{\dagger}1)_{j=s_d+d}]_{TV},\ \mathbf{s} \in \mathbb{N},\ \mathbf{X} = \mathbf{1},\ B_{\mathbf{End}} = (2^{\mathbf{d}}\text{-}N_N),\ \mathbf{P}_{IV} = N_N,\ \mathbf{P}_{TV} = {}^{\dagger}1.$$

The **End** $^{\dagger}$parity-sequence for an $^{\dagger}$**End IV**-Class is unique - it is only valid for $\mathbf{P}_{IV} = N_N$. Notice a **reducing* **parity-sequence* for a **Converging* **IV**-class is valid for all $N_X \in {}^{*}_{u}[2^{*r_s}\mathbf{X}\text{-}{}^{*}B_{Con}]_q$, $\mathbf{X} \in \mathbb{N}$.

**NB**! The following main-sections contains a deeper analysis of (**reducing*) parity-sequences.

# 3. The Infinite Complete Binary Trees

The scope of this main-section is to show that all Class-series for **IV**-classes can be identified with a unique parity-sequence **R**···**E**-**T** where the **Total** number of **tiles** is **r** from the **modulo** in **IV**-classes $2^r$ and the number of **O E**-**tiles** is **s** from the **modulo** in **TV**-classes $3^s$:

| | | | |
|---|---|---|---|
| $[2^r\mathbf{X}\text{-}B_{j=0}]_{IV} \rightarrow\rightarrow [3^s\mathbf{X}\text{-}B_{j=s+r}]_{TV}$ | **r** = #**Total**$_{tile}$ = # **E** $_{Operation}$ | **s** = # **O E** $_{tile}$ = # **O** $_{Operation}$ | **r-s** = # **E** $_{tile}$ |

## The Domino-Tree

An Infinite Binary Domino-Tree {**O E** $_{Left}$, **E** $_{Right}$, **R** $_{Root}$} is defined using **O E**-**tiles** & **E**-**tiles**:

The **R**oot-node **R** is in layer **r** = 0 and any node/**tile** is linked to **two** nodes/**tiles** in layer **r+1**

- The **Left**-child is an **O E**-**tile** & The **Right**-child is an **E**-**tile**

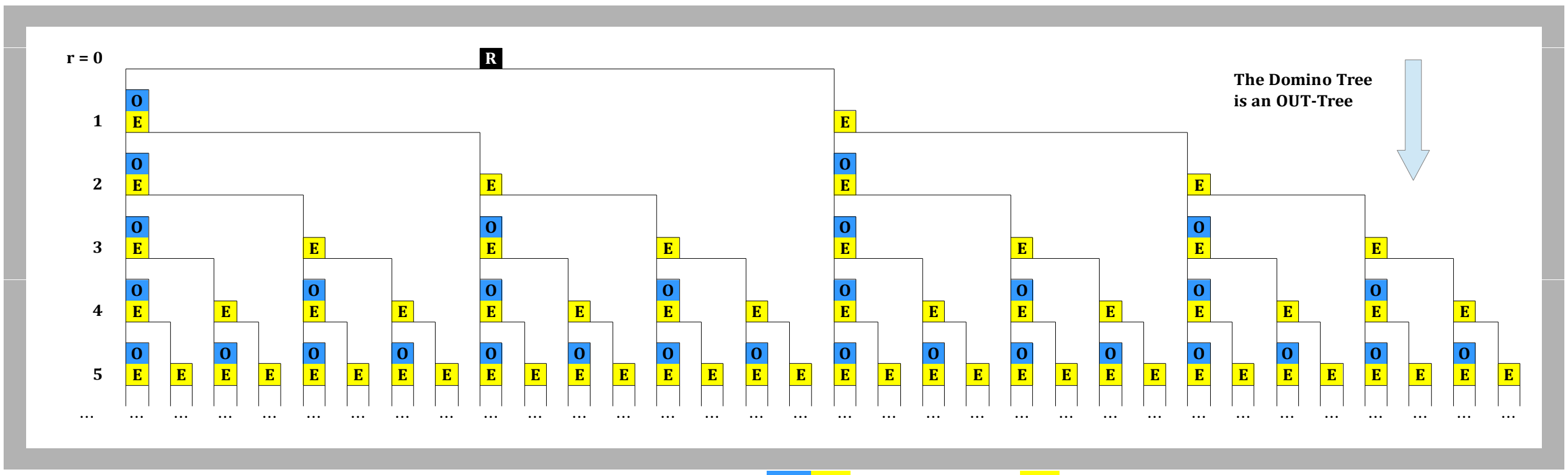

**Graph 1**. The Domino-Tree constructed using **O E**-**tiles** and **E**-**tiles** as vertices/nodes.

*When a **route** is followed **from** the Root-**tile*** **R** ***any** unique **tile** in the **Tree** can be reached and **when** it is, a parity-sequence (minus the **Terminating tile*** **T**) *for the unique **tile** was followed.*

***Routes/tiles** are associated with a Class-series for an* **IV**-*class* **&** *the corresponding* **TV**-*class*:

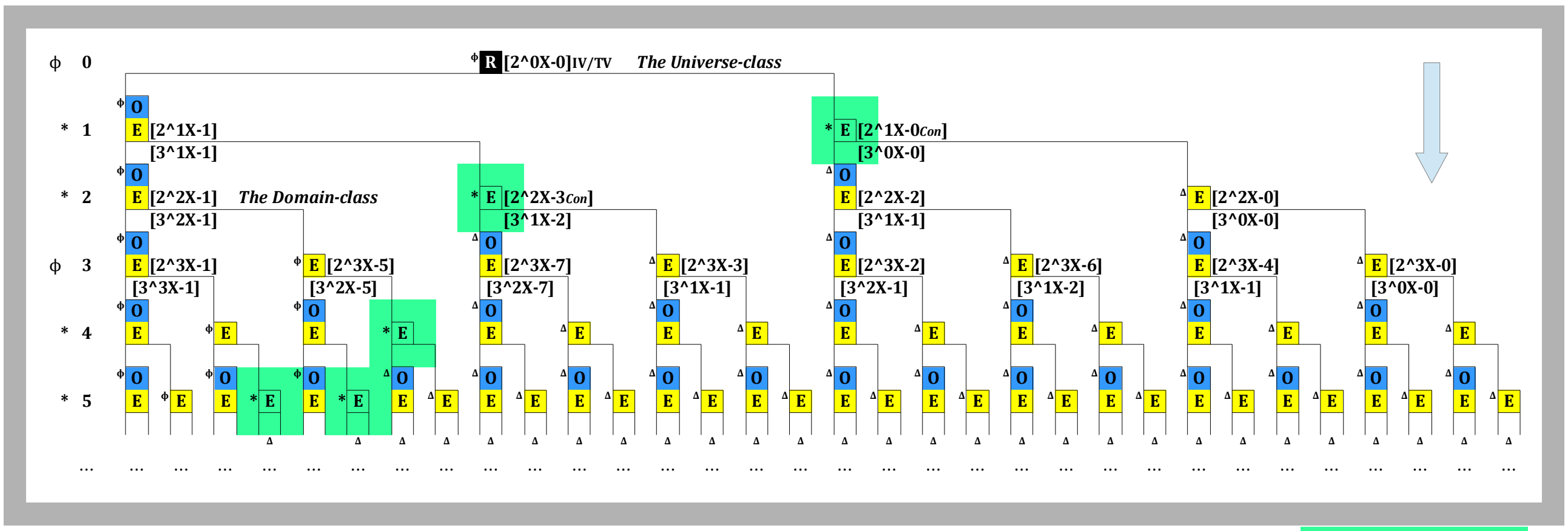

**Graph 2**. Domino-Tree with $[2^r\mathbf{X}\text{-}B_0]_{IV}$ & $[3^s\mathbf{X}\text{-}B_{j=s+r}]_{TV}$ complete to layer **r** = 3 $*[2^{*r}\mathbf{X}\text{-}{*B}_{Con}]_{IV}$

*Notice that the lowest element* **P** *and the second-lowest element* (**P**+$2^r$) *in a* $^{\phi}$**Diverging IV**-*class can **not** belong to the **same** *Converging* **IV**-*class, as all* **IV**-*classes are **split** in the next layer* **r+1**.

**NB**! *Technically it is the* **E**-*operation in* **O E**-**tiles** *that can be associated with the Class-series.*

**OBS**! *All *reducing *parity-sequences **ends** with an* * **E**-**tile** (*as an* **O E**-**tile** ⇒ ***larger*** $N_{j+2}$).

As all nodes/**tiles** in **Graph 2** have a sub-tree **identical** to the (sub-) Tree for the Root R in **Graph 1**, we can ***Zoom in*** on any node/**tile** and picture it as the "New Root" in the sub-tree.

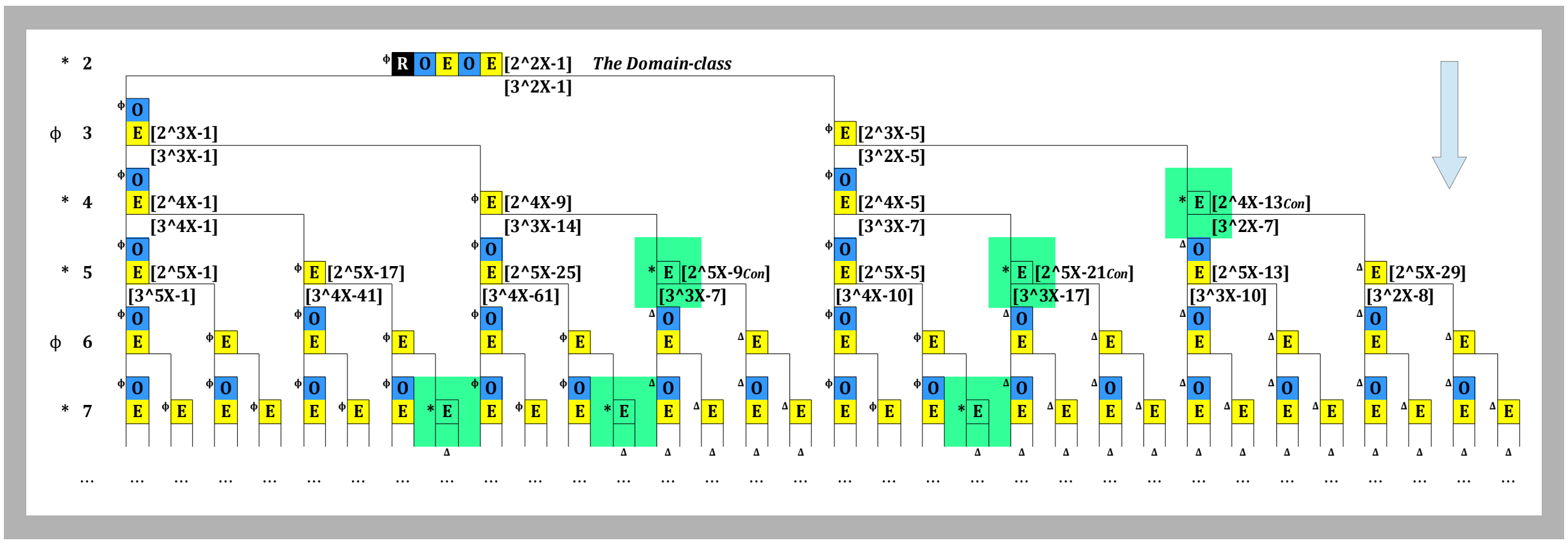

**Graph 3**. Domino-Tree ***Zoom*** on the sub-tree for **The Domain-class** $^{\phi}[2^2\mathbf{X}\text{-}1_0]_{\mathbf{IV}}$

In layer **r** = 5 is found four $^{\phi}$**Div(5)** i.e. $^{\phi}[2^5\mathbf{X}\text{-}1_0]_{\mathbf{IV}}$ & $^{\phi}[2^5\mathbf{X}\text{-}5_0]_{\mathbf{IV}}$ & $^{\phi}[2^5\mathbf{X}\text{-}17_0]_{\mathbf{IV}}$ & $^{\phi}[2^5\mathbf{X}\text{-}25_0]_{\mathbf{IV}}$ and notice, that *these* are the four **IV**-classes *also* contained in the **large square** in **Figure 2** (p. 19).

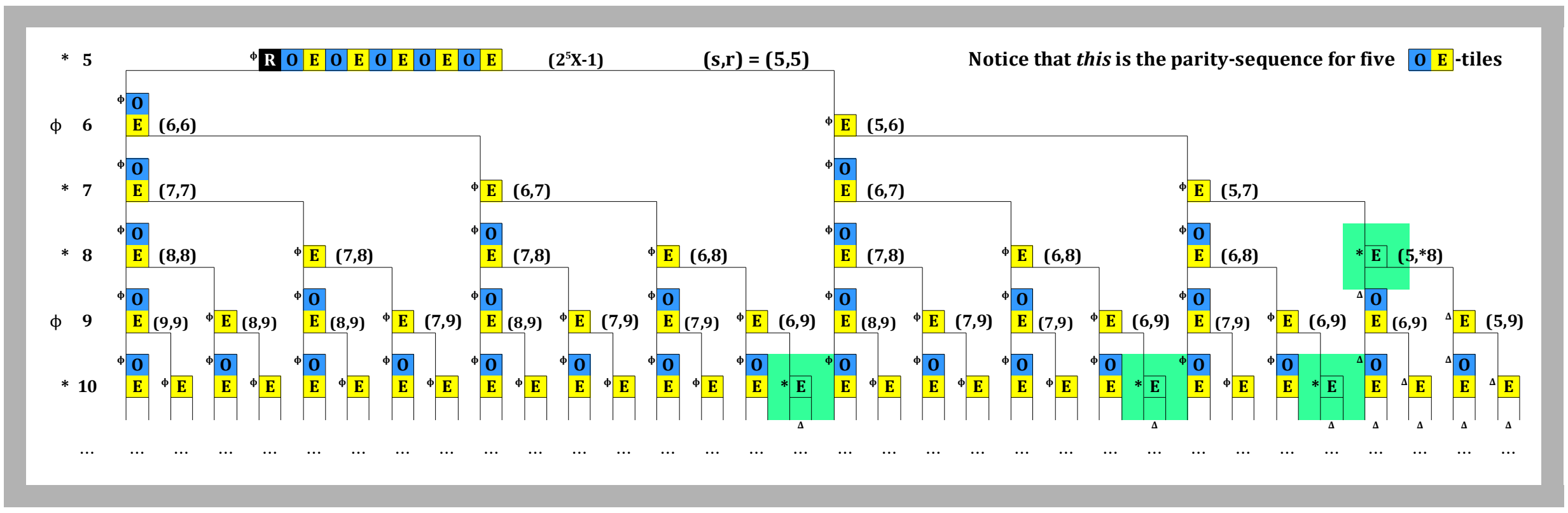

**Graph 4**. Domino-Tree ***Zoom*** on the sub-tree for the **IV**-class $^{\phi}[2^5\mathbf{X}\text{-}1_0]_{\mathbf{IV}}$ in layer **r** = 5.

*Notice that in the sub-trees for* $B_0 \in \{5,17,25\}$ *the **positions** for* *Con(*$\mathbf{r_s}$) *have identical **routes***:

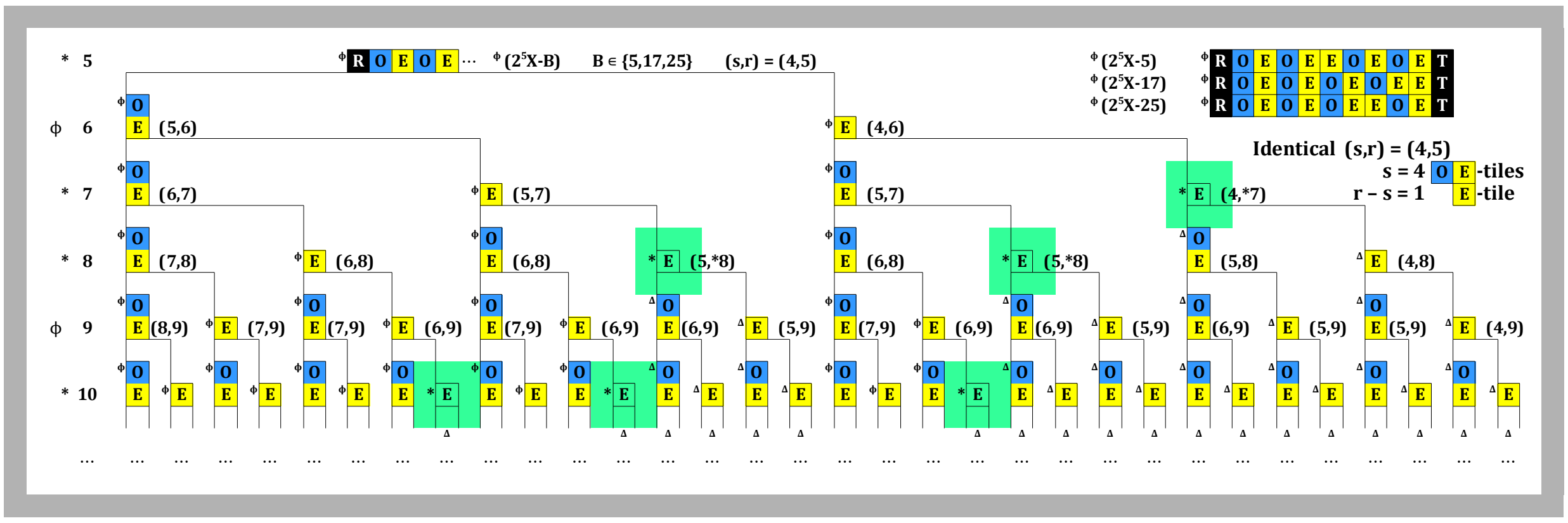

**Graph 5**. The "identical" Domino sub-trees for the **IV**-classes $^{\phi}[2^5\mathbf{X}\text{-}B_0]_{\mathbf{IV}}$ $B_0 \in \{5,17,25\}$, **r** = 5.
Notice in **Figure 2**: $[2^5\mathbf{X}\text{-}5_0]_{\mathbf{IV}} \cup [2^5\mathbf{X}\text{-}17_0]_{\mathbf{IV}} \cup [2^5\mathbf{X}\text{-}25_0]_{\mathbf{IV}}$ have same $|{}^*[2^{*r_s}\mathbf{X}\text{-}{}^*B_{Con}]_{\mathbf{IV}}|$ for $2^{*r_s} > 2^5$.

## Comment

*Please notice the Combinations of Operations* (**s**,**r**) = (# **O** Operation, # **E** Operation) *as these are of high importance in following argumentation* (**NB**! *The distribution of amount of the* **O**-*Operations in* ***routes*** *in* ***layers*** **r** ∈ ℕ *follows the Binomial-distribution, which makes later* **Formulas** *possible*).

## The Base-2 Tree

An Infinite Binary **Base-2 Tree** {1 *Left*, 0 *Right*, **P** *Root*} containing *binary strings* is defined:

The Root-node **P** (the "*head*") is in layer **r** = 0 & all nodes have two Children in layer **r+1**

- The **Left**-child have a 1 added between **P** and the current *binary string* (the "*tail*")
- The **Right**-child have a 0 added between **P** and the current *binary string*

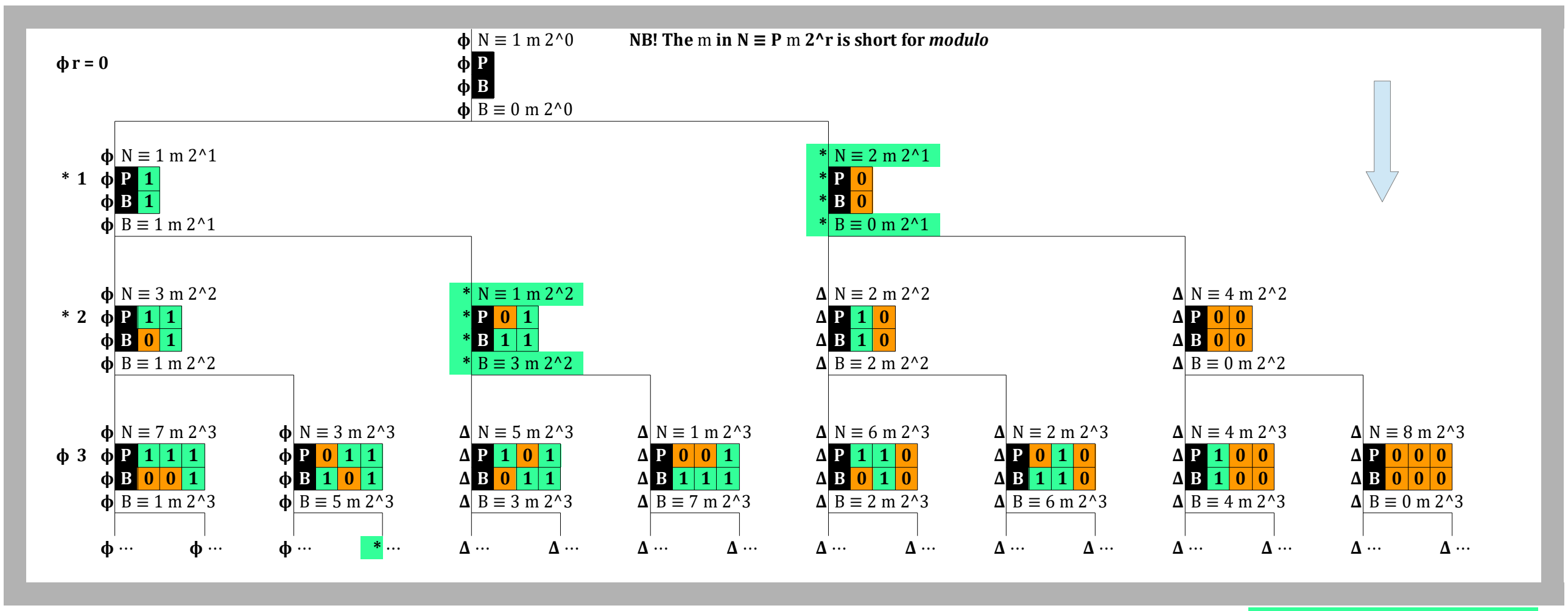

**Graph 6**. **Base-2 Tree** containing the *binary strings* for **P** and $B_0 = (2^r\text{-}\mathbf{P})$ $*N_X \equiv \mathbf{P} \pmod{2^{*r}}$

*The layer* **r** *have the* ***same*** $N_X \equiv (\mathbf{2^r\text{-}B_0}) \pmod{2^r} \Leftrightarrow N_X \in [2^r\mathbf{X}\text{-}B_0]_{IV}$ *as found in The Domino-Tree*.

*Observe in layer* **r** = 2; *one* **Con*(***2**$_1$), *two* $^\Delta$**Red**(**2**) & *one* $^\phi$**Div**(**2**), *the* **Domain-class** $^\phi[2^2\mathbf{X}\text{-}1_0]_{IV}$

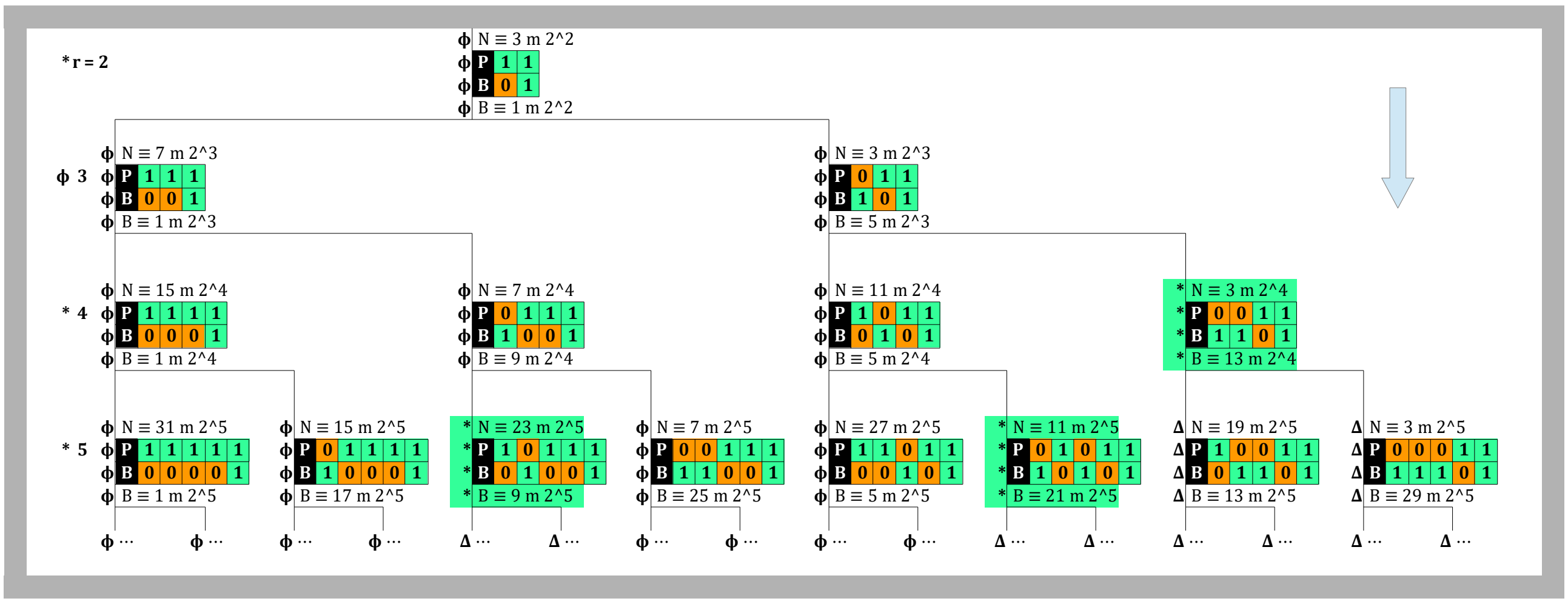

**Graph 7**. **Base-2** sub-tree for the $^\phi$**Diverging Domain-class** $^\phi[2^2\mathbf{X}\text{-}1_0]_{IV}$ in layer **r** = 2.

## Lemma 3

Any node in The Domino-Tree can be identified with a *unique* node in **The Base-2 Tree**.
Any node in **The Base-2 Tree** can be identified with a *unique* node in The Domino-Tree.
***Ergo*** there **exist** a *one-to-one-correspondence* between **IV**-classes and parity-sequences.

## Proof for Lemma 3

The proof follows from ***The Pigeon-hole Principle*** and the fact that in both **Trees** *each* of the layers **r = k** contains $2^k$ *different* nodes. In The Domino-Tree is found $2^k$ *different* ***routes****/ parity-sequences* and in **The Base-2 Tree** is found $2^k$ *different* $[2^k\mathbf{X}\text{-}B_0]$ with *different* $B_0 < 2^k$ ■

## The Standard Binary Tree and The Fibonacci Trees

A Standard Infinite Binary **Tree** {**O** *Left* , **E** *Right* , **R** *Root*} is defined using **O-nodes** & **E-nodes**:

- The **Left**-child is an **O-node/**operation & The **Right**-child is an **E-node/**operation

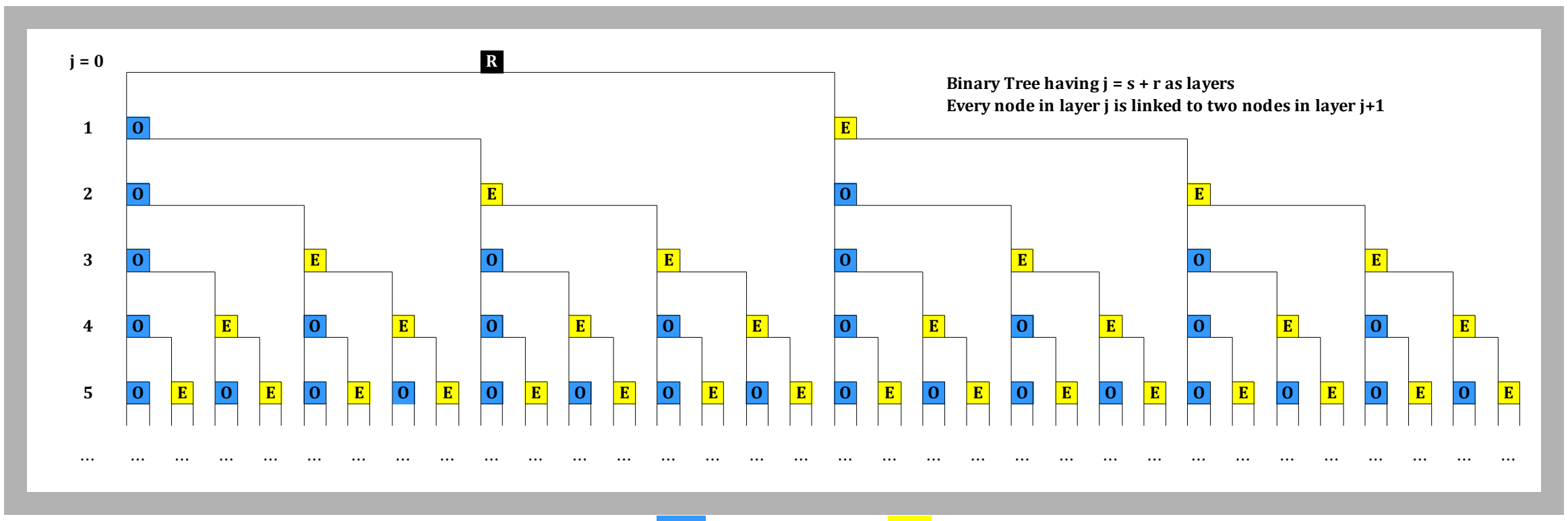

**Graph 8**. Standard Infinite Binary **Tree**. **O-nodes** & **E-nodes** as *vertices* and **j = s+r** as *layers*.

Remember that **all** $(3N_j+1) \equiv$ **0** (mod 2) for **all** $N_j \equiv$ **1** (mod 2) so this **restriction** exist:

| |
|---|
| **No O↔O-link** is allowed as **all O** are followed by an **E** in the rules for Collatz-Iterations |

This is achieved by ***pruning*** the **Tree**. When a *Parent*-node in layer **j** is ***pruned*** the **link** to a *Child*-node in layer **j+1** is **deleted** i.e. the *Child*-node no longer have a ***route*** from the Root **R**.

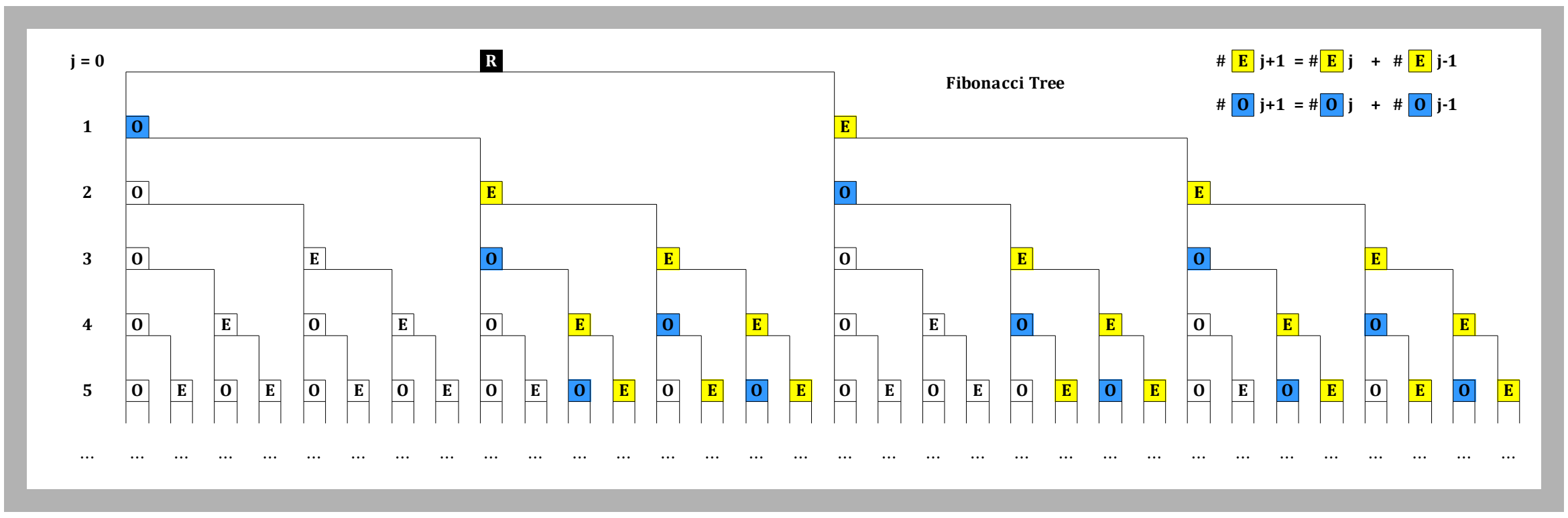

**Graph 9. Fibonacci IN/OUT Tree** with colourless "ghosts" of disconnected sub-trees.

**Notice** that in layer **j = s+r** = # **O** + # **E**: | # **O** (j) + # **E** (j) = #Total(j) = # **E** (j+1) = # **O** (j+2) |

*The above figure is included to emphasize that all nodes in* **Graph 10** *have a position in* **Graph 8**.

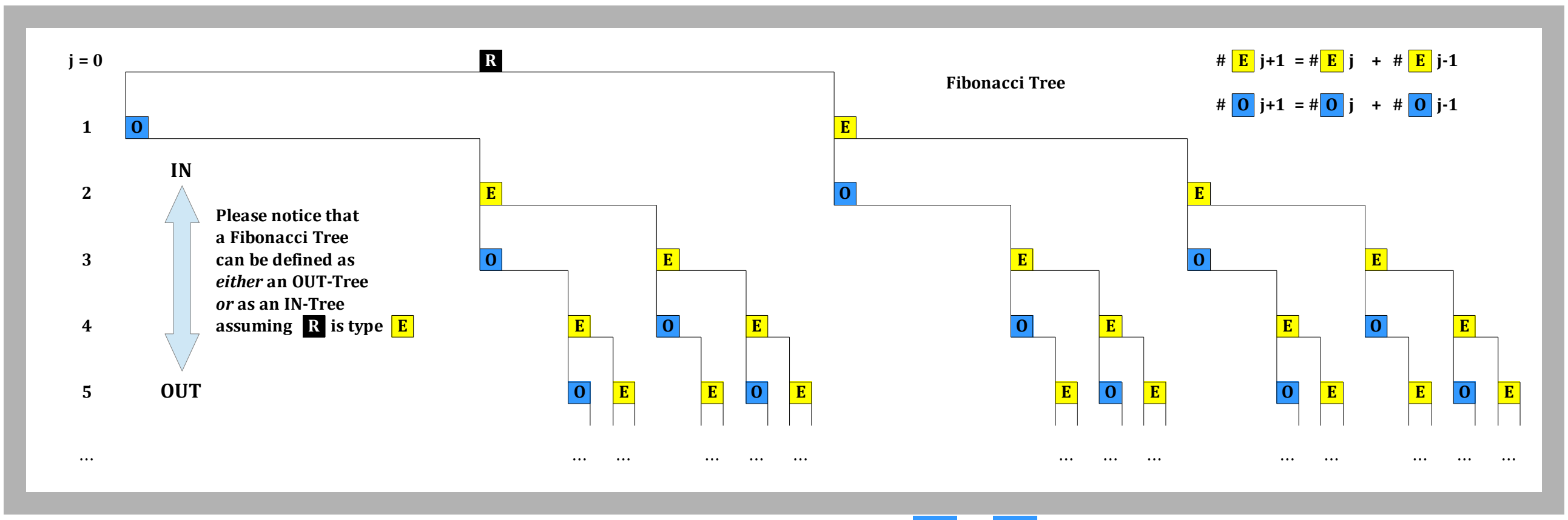

**Graph 10. Fibonacci IN/OUT Tree** after removing all **O↔O-links** from **Graph 8**.

## The Fibonacci Tree

Suppose that **Graph 10** is the starting-point for a **Fibonacci Tree** defined as an **OUT-Tree**. In **this** case (all) the ***routes*** (*from* the Root) found in **Graph 11** are identical to ***routes*** found in **Graph 2**, implying that the **Total** amount of *unique* ***routes***/nodes is an **Uncountable Infinity**. **NB**! *As in The Domino-Tree it is* ***only*** **E***-operations that can be associated with the* ***Class-series***.

When **Converging *parity-sequences* in **The Fibonacci Tree** are ***pruned***, it is found that all the **reducing* **IV**-classes **Con*(*$\mathbf{r_s}$) are grouped together in the layer ***j = s +** *$\mathbf{r_s}$

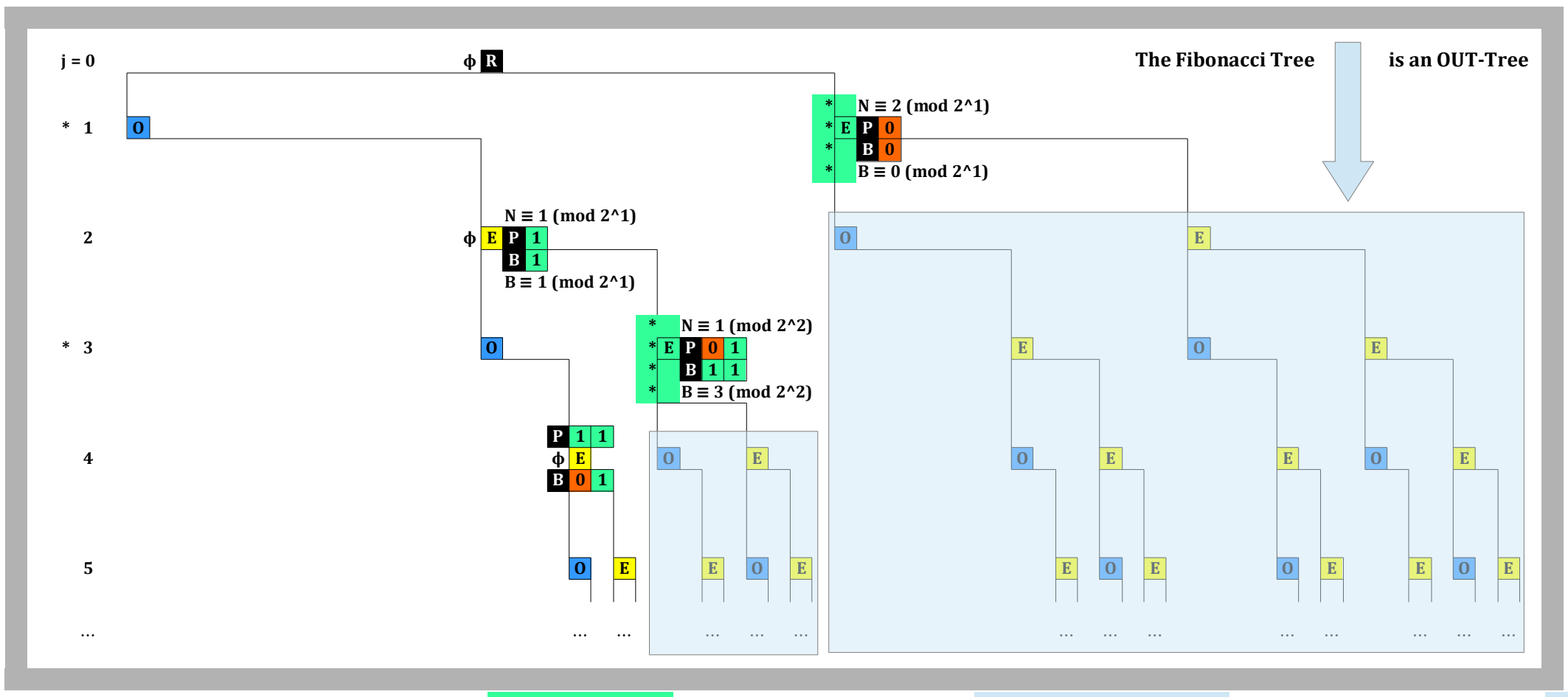

**Graph 11. Fibonacci Tree** with **Con*(*$\mathbf{r_s}$) & $^{\Phi}$**Div(r)** but $^{\Delta}$**Red(r >***$\mathbf{r_s}$**)** ***pruned away*** (shaded).

**OBS**! Notice that this ***pruning*** of *N ≡ $P_{IV}$ (mod $2^{*r}$) ***could*** have been done in the Domino-Tree.

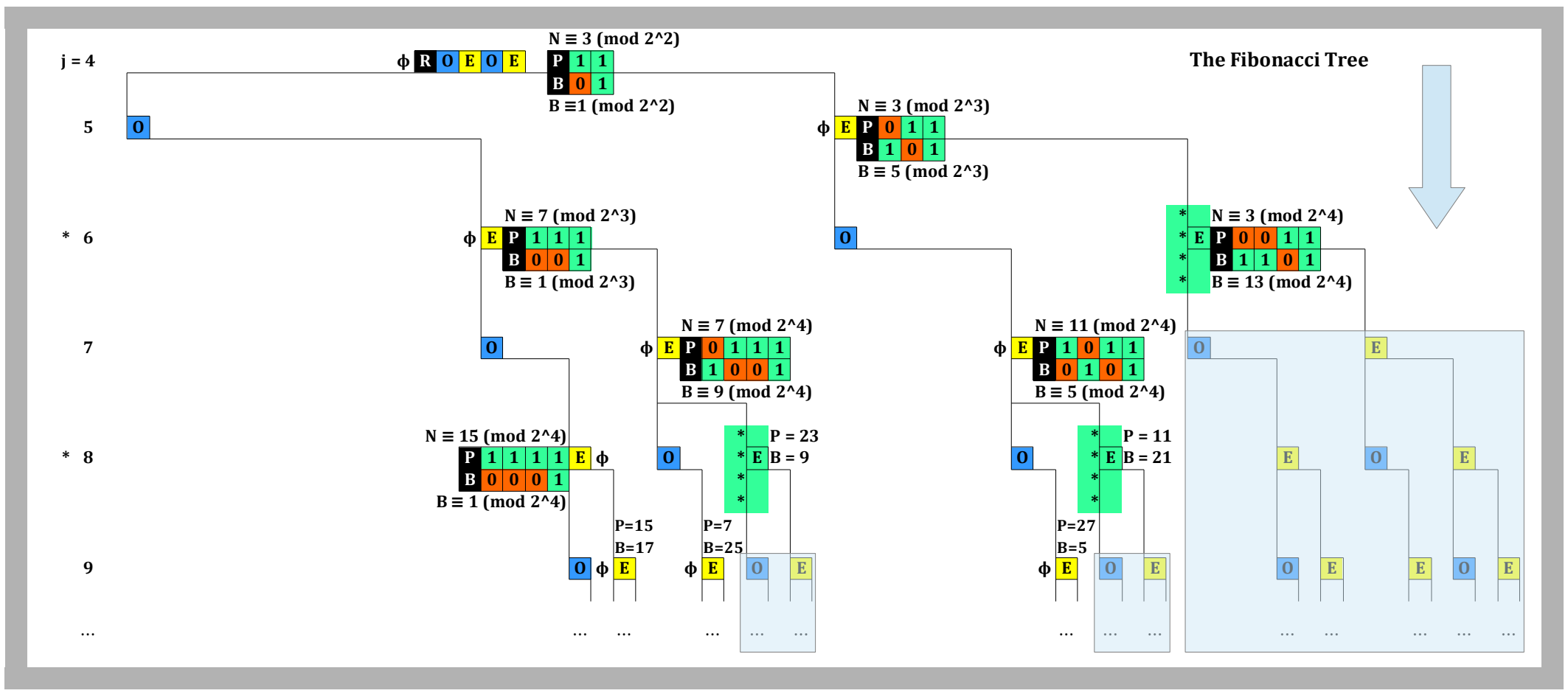

**Graph 12.** ***Zoom*** on the **Fibonacci sub-tree** for the ***Domain-class*** $^{\Phi}[2^2\mathbf{X}\text{-}1_0]_{IV}$ in layer **j** = 4.

### Comment

*After pruning all* $*_{\mathbf{u}}[2^{*r_s}\mathbf{X}\text{-}*B_{Con}]_{\mathbf{q}}$, **u** ∊ ℕ *the* ***Tree*** *contains a countable infinity of *reducing routes.*

*Any* ***sub-tree*** (*which have* ***not*** *been pruned*) *in the* ***Tree*** *contains an* **Uncountable Infinity** *of different* ***routes****/nodes. This implies, that only after pruning* ***all*** $*_{\mathbf{u}\to\infty}[2^{*r_s}\mathbf{X}\text{-}*B_{Con}]_{\mathbf{q}\to\infty}$, **u** ∊ ℕ, *is it* ***possible*** *for the* ***Tree*** *to contain* ***less*** *than an* **Uncountable Infinity** *of different* ***routes****. Please notice that the *Converging *parity-sequences are* ***extensions*** *of* $^{\Phi}$**Diverging** $^{\Phi}$parity-sequences.

*In the following sections in the present work it is shown, that there* ***exist*** *a countable infinity of the *reducing *parity-sequences i.e. a countable infinity of *Converging* **IV***-classes, and it is* ***also*** *shown, that there* ***exist*** *a countable infinity of* $^{\Phi}$**Diverging IV**-*classes*/$^{\Phi}$parity-sequences/***routes***.

## The *Reverse* Fibonacci Tree

Suppose that **Graph 10** is the starting-point for a **Fibonacci Tree** defined as an **IN-Tree**. In **this** case all {**O** *Left*, **E** *Right*, **D** *Root*} can be associated with a **Value** N ∊ ℕ *and* an †**End-IV-class.** The new **Root-tile D** have the **Value** N = **1** (***The Origo* 1!**) and as **1** is Odd it is necessary to ***prune*** the odd-odd-link **O**↔**D** *disconnecting* the entire left half of **Graph 10** to get **Graph 13.**

- Any **node**/sub-tree that **can not** be identified with **Values** N ∊ ℕ are ***pruned away*** i.e.
- All **links** to the **Left**-child **O** are deleted for all N ≡ **0** (mod 6) **and** all N ≡ **2** (mod 6)
- Also ***pruned*** is the **Left**-child for the **E** -**node** with the **Value** N = **4** (to avoid ***The Loop***)

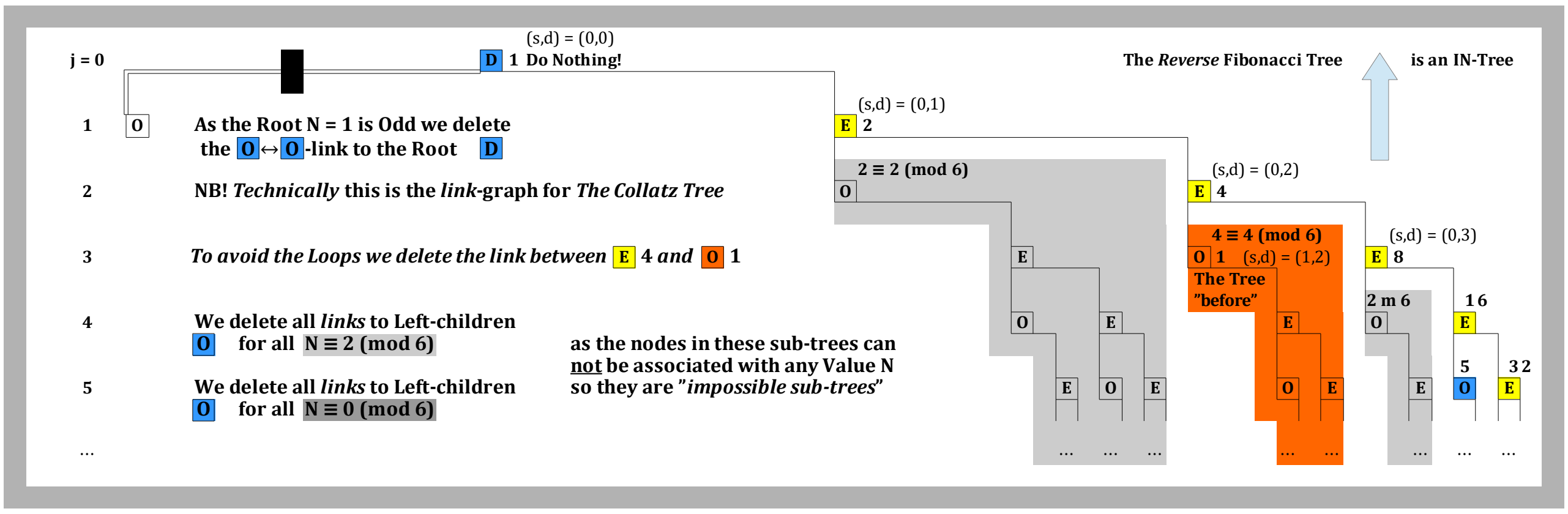

**Graph 13. The *Reverse* Fibonacci Tree a.k.a. *The Collatz Tree*.**

Please notice that all ***routes*** (towards ***The Origo* 1!**) in **The *Reverse* Fibonacci Tree** *definitely* exist in **The (OUT) Fibonacci Tree** as **all possible** parity-sequences do exist in **Graph 11**, but **not** all ***routes*** in **Graph 11** exist in **Graph 13**, as only **some** nodes in **Graph 11** represents †**End-Class-series**, but **all** nodes in **Graph 13** represents a: $^\dagger[2^d X\text{-}(2^d\text{-}N_N)_{End}]_N >> [3^{s_d} X\text{-}(3^{s_d}\text{-}1)]$

*Hopefully it is perfectly clear, that* **Graph 13** *is* ***equivalent*** *to* **Figure 1a** *hence to* **Flowchart 6.**

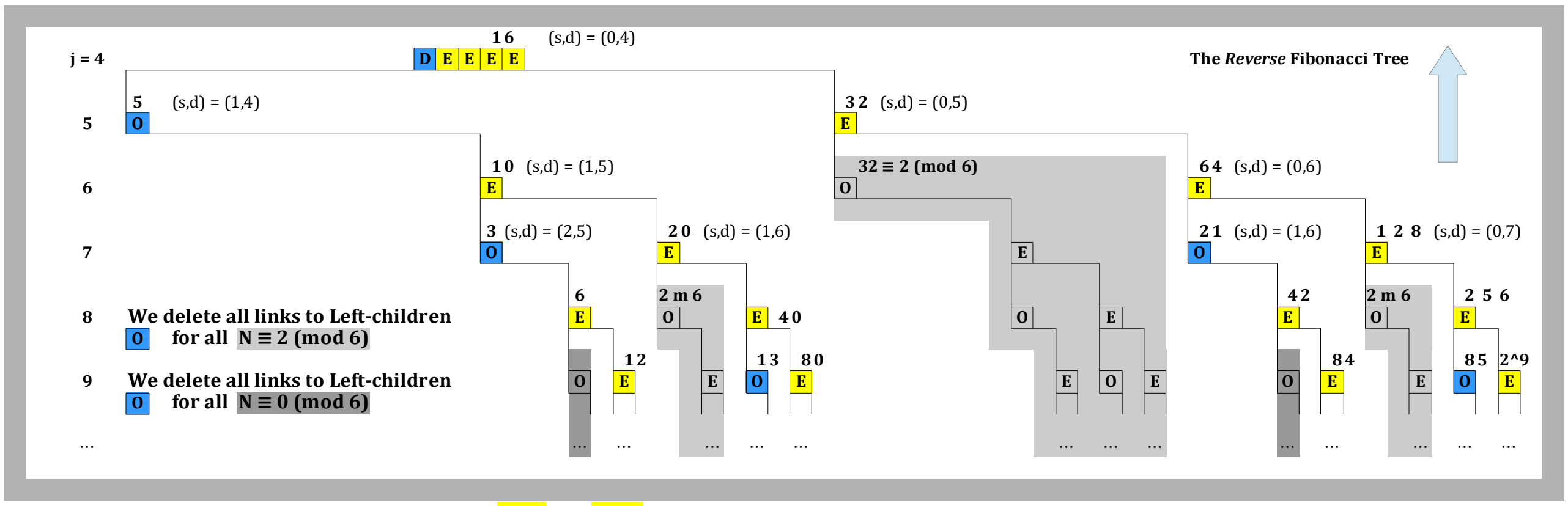

**Graph 14. Zoom** on the node **E** ~ 16 in layer **j** = 4 in **The *Reverse* Fibonacci Tree.**

### Comment

*After pruning* ***all*** [6**X**-0] ∪ [6**X**-4] *the* ***Reverse Tree*** *have a countable infinity of* ***nodes*** N ∊ ℕ.

*Any* ***sub-tree*** (*which have* ***not*** *been pruned*) *in the* ***Tree*** *contains an* **Uncountable Infinity** *of different* ***routes***/*nodes. This implies, that only after pruning* ***all*** [6**X**-0] ∪ [6**X**-4] *is it* ***possible*** *for the* ***Tree*** *to contain* ***less*** *than an* **Uncountable Infinity** *of different* ***routes***, *which* ***implies*** *that all* [6**X**-3] *are included in the* ***Tree***, *because all* [6**X**-0] *have been pruned, and that all* [6**X**-2] *are included in the* ***Tree,*** *because all* [6**X**-4] *are pruned & all* [6**X**-2] ***implies*** *all* [6**X**-5] ∪ [6**X**-1].

# 4. Possible *_Converging_ *$Con_u$(s,*r)$_q$ and possible $^\Phi$Diverging $^\Phi$Div$_{u'}$(s,r)$_{q'}$

*The scope of **this** main-section is to **show** how analysis of **parity-sequences** leads to the **valid*** **Formulas** *for the Densities of the three **types*** [2$^r$**X**-B$_0$]$_{IV}$ {$^\Phi$**Diverging**, **Converging*, $^\Delta$**Redundant**}

## Arranging *_Converging_ *$Con_u$(*r$_s$)$_q$ by *_Reducing Combinations_ (s,*r$_s$) and ascending $P_{IV}$

*All *Con*$_u$(*r$_s$)$_q$ *are arranged in **Groups** for constant* **s** *in* (s,*[sZ]). *The* $_u$*index* **u** ∊ ℕ *is used as a **counting** variable for the* **u**nique **Converging* **IV**-classes *i.e.* ***each*** **u**nique ***pair*** $_u$($P_{IV}$+*$B_{IV}$) = 2$^{*r_s}$.
*Inside **Groups** for constant* (**s**,*r$_s$) *the index*$_q$ **q** ∊ ℕ *is used to identify* $P_q < P_{q+1} \Rightarrow q_{Max}$ = |*Con(*r$_s$)|
*The **first** *Converging* **IV**-class *$_1$[2$^{*1_0}$**X**-*0$_{Con}$]$_1$ *is associated with the index* **u** = 1 & **q** = 1.
*The **second** *Converging* **IV**-class *$_2$[2$^{*2_1}$**X**-*3$_{Con}$]$_1$ *is associated with the index* **u** = 2 & **q** = 1, ***etc.***

**$P_q$** + *$B_{Con}$ = **2$^{*r_s}$** & **Con*$_u$(*r$_s$)$_q$ & **parity-sequence for* *R... E - E - T *implied*:

| | | | | | |
|---|---|---|---|---|---|
| s: 6 5 4 3 2 1 0 ← s | | | #O = s → 0 1 2 3 4 5 6 | | |
| r: … 11 10 9 8 7 6 5 4 3 2 1 0 ← r | | N ≡ P (mod 2^*r) | #E = *r → 1 2 4 5 7 8 10 | | |
| | | N ∈ [2^*rX-BCon] | j = s + r → 1 2 3 4 5 6 7 8 9 10 11 12 13 14 15 16 17 … | | |
| P 0 | u =1 | N ≡ 2 (mod 2^1) | q = 1 | *R E | |
| B 0 | | N ∈ [2^1X-0] | | | \|*Con(*1$_0$)\| = 1 |
| P 0 1 | 2 | N ≡ 1 (mod 2^2) | 1 | *R O E E | |
| B 1 1 | | N ∈ [2^2X-3] | | | \|*Con(*2$_1$)\| = 1 |
| P 0 0 1 1 | 3 | N ≡ 3 (mod 2^4) | 1 | *R O E O E E E | |
| B 1 1 0 1 | | N ∈ [2^4X-13] | | | \|*Con(*4$_2$)\| = 1 |
| P 0 1 0 1 1 | 4 | N ≡ 11 (mod 2^5) | 1 | *R O E O E E O E E | |
| B 1 0 1 0 1 | | N ∈ [2^5X-21] | | | \|*Con(*5$_3$)\| = 2 |
| P 1 0 1 1 1 | 5 | N ≡ 23 (mod 2^5) | 2 | *R O E O E O E E E | |
| B 0 1 0 0 1 | | N ∈ [2^5X-9] | | | |
| P 0 0 0 0 1 1 1 | 6 | N ≡ 7 (mod 2^7) | 1 | *R O E O E O E E O E E E | |
| B 1 1 1 1 0 0 1 | | N ∈ [2^7X-121] | | | \|*Con(*7$_4$)\| = 3 |
| P 0 0 0 1 1 1 1 | 7 | N ≡ 15 (mod 2^7) | 2 | *R O E O E O E O E E E E | |
| B 1 1 1 0 0 0 1 | | N ∈ [2^7X-113] | | | |
| P 0 1 1 1 0 1 1 | 8 | N ≡ 59 (mod 2^7) | 3 | *R O E O E E O E O E E E | |
| B 1 0 0 0 1 0 1 | | N ∈ [2^7X-69] | | | |
| P 0 0 1 0 0 1 1 1 | 9 | N ≡ 39 (mod 2^8) | 1 | *R O E O E O E E O E O E E E | |
| B 1 1 0 1 1 0 0 1 | | N ∈ [2^8X-217] | | | \|*Con(*8$_5$)\| = 7 |
| P 0 1 0 0 1 1 1 1 | 10 | N ≡ 79 (mod 2^8) | 2 | *R O E O E O E O E E E O E E | |
| B 1 0 1 1 0 0 0 1 | | N ∈ [2^8X-177] | | | |
| P 0 1 0 1 1 1 1 1 | 11 | N ≡ 95 (mod 2^8) | 3 | *R O E O E O E O E O E E E E | |
| B 1 0 1 0 0 0 0 1 | | N ∈ [2^8X-161] | | | |
| P 0 1 1 1 1 0 1 1 | 12 | N ≡ 123 (mod 2^8) | 4 | *R O E O E E O E O E E O E E | |
| B 1 0 0 0 0 1 0 1 | | N ∈ [2^8X-133] | | | |
| P 1 0 1 0 1 1 1 1 | 13 | N ≡ 175 (mod 2^8) | 5 | *R O E O E O E O E E O E E E | |
| B 0 1 0 1 0 0 0 1 | | N ∈ [2^8X-81] | | | |
| P 1 1 0 0 0 1 1 1 | 14 | N ≡ 199 (mod 2^8) | 6 | *R O E O E O E E O E E O E E | |
| B 0 0 1 1 1 0 0 1 | | N ∈ [2^8X-57] | | | |
| P 1 1 0 1 1 0 1 1 | 15 | N ≡ 219 (mod 2^8) | 7 | *R O E O E E O E O E O E E E | |
| B 0 0 1 0 0 1 0 1 | | N ∈ [2^8X-37] | | | |
| … | … | … | … | … | |

**List 1. The first 15** *$_{u\le15}$[2$^{*r_s}$**X**-*B$_{Con}$]$_q$ **in ascending order for u → ∞.**

*This* is **The Infinite List** of **reducing* **IV**-classes *$_u$[2$^{*r_s}$**X**-*B$_{Con}$]$_q$ for **u** ∊ ℕ, **u** → ∞.

**For example**:

For the **Reducing Combination* (**s**,*r$_s$) = (**4**,*7$_4$) ***three*** pcs of *$_u$[2$^{*7_4}$**X**-*B$_{Con}$]$_q$ are observed for *B$_{Con}$ ∊ {121,113,69} and $u_q$ ∊ {6$_1$,7$_2$,8$_3$} as the **IV**-classes are here arranged according to *size* of **$P_q$**; 7$_1$ < 15$_2$ < 59$_3$

**OBS**! *The **three** *parity-sequences are found in layer* **r** = *7 *in* **Graph 5**.

**NB**! *Remember that*

**$P_q$** + *$B_{Con}$ = **2$^{*r_s}$**

**Notice** that for **u ≥ 3** it is observed, that **all** the **reducing* **IV**-classes are **proper subsets** of the **Domain-class** $^\Phi$[2$^2$**X**-1$_0$]$_{IV}$ which have the $^\Phi$**Diverging** $^\Phi$**parity-sequence** $^\Phi$**R O E - O E - T** & also **notice** that a **Converging* **IV**-class with a **parity-sequence* ***R ... E'- E - T** contains **half** of the N$_X$ for **either** odd **X or** even **X** in the $^\Phi$**Diverging IV**-class with the $^\Phi$**sequence** $^\Phi$**R ... E'- T**.

## Comment

The purpose of **List 1** is to **demonstrate**, that **each** **u**nique **reducing* *$_u$[2$^{*r_s}$**X**-*B$_{Con}$]$_q$, **u** ∊ ℕ, **can** be identified with a **u**nique **reducing *parity-sequence* and with the *pair* $_u$($P_q$+*$B_{Con}$) = 2$^{*r_s}$.
In **List 2** on the next page, the **G**roups for constant **Reducing Combination* (**s**,*r$_s$) have been **rearranged** according to **patterns** in the **u**nique **reducing *parity-sequences* ***R ... E - E - T**.
The **G**roups in **List 2** contains the same **amount** of **Converging* **IV**-classes *$_u$[2$^{*r_s}$**X**-*B$_{Con}$]$_q$ as in **List 1** i.e. |*Con*$_u$(*r$_s$)$_q$| = $q_{Max}$ for a specific $A_{IV}$ = 2$^{*r_s}$, but the **IV**-classes are **not** in the same *order*.
**NB**! **Graph 4** and **Graph 5** are found on page 22. These graphs are important when **subsets** of the four $^\Phi$**Diverging IV**-classes $^\Phi$[2$^5$**X**-B$_0$]$_{IV}$ for B$_0$ ∊ {1,5,17,25}, $A_{IV}$ = 2$^5$, are observed in **List 2**.

## Arranging *$Converging$* *$Con_u(*r_s)_q$* by $(s,*r_s)$ and *patterns* in the **parity-sequences*

For all the **Reducing Combination* $(s,*r_s)$ a **finite number** $|*Con_u(*r_s)_q|$ is observed for $*_u[2^{*r_s}X\text{-}*B_{Con}]_{q \le q_{Max}}$ so it is possible to **rearrange** $*_u[2^{*r_s}X\text{-}*B_{Con}]_q$ in any alternative order. The following **List 2** uses **u'** & **q'** to arrange $*Con(*r_s)$ in the **Groups** for constant $(s,*r_s)$:

| | | | | | |
|---|---|---|---|---|---|
| s: | 6, 5, 4, 3, 2, 1, 0 ← s | | #O = s → | 0, 1, 2, 3, 4, 5, 6 | |
| r: | … 11 10 9 8 7 6 5 4 3 2 1 0 ← r | N ≡ **P** (mod 2^*r) / N ∈ [2^*rX-BCon] | #E = *r → | 1, 2, 4, 5, 7, 8, 10 | |
| | | | j = s + r → | 1 2 3 4 5 6 7 8 9 10 11 12 13 14 15 16 17 … | |

| New index u' | q' | P / B | u | N | q | parity-sequence | |
|---|---|---|---|---|---|---|---|
| u' = 1 | q' = 1 | P: 0 / B: 0 | u = 1 | N ≡ **2** (mod 2^1); N ∈ [2^1**X**-0] | q = 1 | * R E | \|*Con(*1₀)\| = 1 |
| 2 | 1 | P: 0 1 / B: 1 1 | 2 | N ≡ **1** (mod 2^2); N ∈ [2^2**X**-3] | 1 | * R O E E (1) | \|*Con(*2₁)\| = 1 |
| 3 | 1 | P: 0 0 1 1 / B: 1 1 0 1 | 3 | N ≡ **3** (mod 2^4); N ∈ [2^4**X**-13] | 1 | * R O E O E E E (1 1) | \|*Con(*4₂)\| = 1 |
| 4 | 1 | P: 1 0 1 1 1 / B: 0 1 0 0 1 | 5 | N ≡ **23** (mod 2^5); N ∈ [2^5**X**-9] | 2 | * R O E O E O E E E (1 2) | \|*Con(*5₃)\| = 2 |
| 5 | 2 | P: 0 1 0 1 1 / B: 1 0 1 0 1 | 4 | N ≡ **11** (mod 2^5); N ∈ [2^5**X**-21] | 1 | * R O E O E E O E E | # E after last O = 3 |
| 6 | 1 | P: 0 0 0 1 1 1 1 / B: 1 1 1 0 0 0 1 | 7 | N ≡ **15** (mod 2^7); N ∈ [2^7**X**-113] | 2 | * R O E O E O E O E E E E (1 3 3) | \|*Con(*7₄)\| = 3 |
| 7 | 2 | P: 0 0 0 0 1 1 1 / B: 1 1 1 1 0 0 1 | 6 | N ≡ **7** (mod 2^7); N ∈ [2^7**X**-121] | 1 | * R O E O E O E E O E E E | |
| 8 | 3 | P: 0 1 1 1 0 1 1 / B: 1 0 0 0 1 0 1 | 8 | N ≡ **59** (mod 2^7); N ∈ [2^7**X**-69] | 3 | * R O E O E E O E O E E E (1 3 3) | # E = 7 |
| 9 | 1 | P: 0 1 0 1 1 1 1 1 / B: 1 0 1 0 0 0 0 1 | 11 | N ≡ **95** (mod 2^8); N ∈ [2^8**X**-161] | 3 | * R O E O E O E O E O E E E E (1 4 7) | \|*Con(*8₅)\| = 7 |
| 10 | 2 | P: 1 0 1 0 1 1 1 1 / B: 0 1 0 1 0 0 0 1 | 13 | N ≡ **175** (mod 2^8); N ∈ [2^8**X**-81] | 5 | * R O E O E O E O E E O E E E | |
| 11 | 3 | P: 0 0 1 0 0 1 1 1 / B: 1 1 0 1 1 0 0 1 | 9 | N ≡ **39** (mod 2^8); N ∈ [2^8**X**-217] | 1 | * R O E O E O E E O E O E E E | |
| 12 | 4 | P: 1 1 0 1 1 0 1 1 / B: 0 0 1 0 0 1 0 1 | 15 | N ≡ **219** (mod 2^8); N ∈ [2^8**X**-37] | 7 | * R O E O E E O E O E O E E E | |
| 13 | 5 | P: 0 1 0 0 1 1 1 1 / B: 1 0 1 1 0 0 0 1 | 10 | N ≡ **79** (mod 2^8); N ∈ [2^8**X**-177] | 2 | * R O E O E O E O E E E O E E | |
| 14 | 6 | P: 1 1 0 0 0 1 1 1 / B: 0 0 1 1 1 0 0 1 | 14 | N ≡ **199** (mod 2^8); N ∈ [2^8**X**-57] | 6 | * R O E O E O E E O E E O E E | |
| 15 | 7 | P: 0 1 1 1 1 0 1 1 / B: 1 0 0 0 0 1 0 1 | 12 | N ≡ **123** (mod 2^8); N ∈ [2^8**X**-133] | 4 | * R O E O E E O E O E E O E E (1 4 7) | # E = 12 |
| … | … | … | … | … | … | … | |

**List 2. The first 15 $*_{u \le 15}[2^{*r_s}X\text{-}*B_{Con}]_q$ in ascending order for u' → ∞**

In *this* **List** the $*Con(*r_s)$ have been *rearranged* according to *patterns* in **parity-sequences*.

The *first* entry **q'** = **1** in *each* of the **Groups** for constant $(s,*r_s)$ is the **IV**-class having a **parity-sequence* ***starting*** with only the ***increasing*** **O E-tiles** and *ending* with only the *decreasing* **E-tiles**.

The *next* **sub**-group are the *sequences* with one **E-tile** *before* the *final* **O E-tile**, & the *next* **sub**-group those *sequences* with two **E-tiles** *before* the *final* **O E-tile** *etc.*

The **TRICK** is, that just by ***counting*** all ***theoretically possible*** $^{\phi}$**Diverging** $^{\phi}$parity-sequences (for same **s** = #**O**) in a **Group**, it is ***possible*** to *deduce how many* $*Con(*r_{s+1})$ **exist** in the *next* **Group** for $G(s+1,*r_{s+1})$ by ***counting*** the type **E** Operation after the last type **O** Operation before * **E**-**T**

**Example:** Constructing the $|Con(*8_5)| = 7$ **reducing* **parity-sequences* from the **3** = $|Con(*7_4)|$

| #O = s → | 0 | 1 | | 2 | 3 | | 4 | 5 | | 6 | | | | | | | | |
|---|---|---|---|---|---|---|---|---|---|---|---|---|---|---|---|---|---|---|
| #E = *r → | 1 | 2 | | 4 | 5 | | 7 | 8 | | 10 | | | | | | | | |
| j = s + r → | 1 | 2 | 3 | 4 | 5 | 6 | 7 | 8 | 9 | 10 | 11 | 12 | 13 | 14 | 15 | 16 | 17 | … |

| u' |  | q' |
|---|---|---|
| 6 | R O E O E O E O E | 1 |

The $^{\phi}$parity-sequence with four **O E-tiles.**

| u' |  | q' |
|---|---|---|
| 9 | R O E O E O E O E | 1 |

For **u'** = 9 the copy of the $^{\phi}$sequence in **u'** = 6.

| u' |  | from q' |
|---|---|---|
| 6 | R O E O E O E O E | 1 |
| 6 | R O E O E O E O E E | 1 |
| 7 | R O E O E O E E O E | 2 |
| 8 | R O E O E E O E O E | 3 |

Including the three $^{\phi}$parity-sequences containing four **O E-tiles** and one **E-tile** i.e. for the combination **(s,r)** = (#**O**,#**E**) = (4,5)

| u' |  | to q' |
|---|---|---|
| 9 | R O E O E O E O E | 1 |
| 10 | R O E O E O E O E E | 2 |
| 11 | R O E O E O E E O E | 3 |
| 12 | R O E O E E O E O E | 4 |

For **u'** ∈ {10,11,12} copy of $^{\phi}$parity-sequences from **u'** ∈ {6,7,8}. The following *List* includes **u'** ∈ {13,14,15} from **u'** ∈ {6,7,8}, **(s,r)** = (4,6)

| u' |  | from q' | to q' |
|---|---|---|---|
| 9 | R O E O E O E O E | 1 | 1 |
| 10 | R O E O E O E O E E | 1 | 2 |
| 11 | R O E O E O E E O E | 2 | 3 |
| 12 | R O E O E E O E O E | 3 | 4 |
| 13 | R O E O E O E O E E E | 1 | 5 |
| 14 | R O E O E O E E O E E | 2 | 6 |
| 15 | R O E O E E O E O E E | 3 | 7 |

The seven $^{\phi}$**Diverging** $^{\phi}$parity-sequences are the complete collection of $^{\phi}$parity-sequences which contains exactly **s** = 4 pcs of **O E-tiles.**

| u' | sequence (1 4 7) | q' |
|---|---|---|
| 9 | * R O E O E O E O E O E E E E | 1 |
| 10 | * R O E O E O E O E E O E E E | 2 |
| 11 | * R O E O E O E E O E O E E E | 3 |
| 12 | * R O E O E E O E O E O E E E | 4 |
| 13 | * R O E O E O E O E E E O E E | 5 |
| 14 | * R O E O E O E E O E E O E E | 6 |
| 15 | * R O E O E E O E O E E O E E | 7 |

The seven $^{\phi}$**Diverging** $^{\phi}$parity-sequences are ***extended*** with **1 O E-tile** & an ***appropriate amount*** of **E-tiles** to reach $(s,*r_s) = (5,*8_5)$

## How to read List 2

Notice that the **Infinite Lists** contains The Complete **Collection** of all existing $*_u[2^{*r_s}\mathbf{X}\text{-}*B_{Con}]_q$ and that **all** $*Con_u(*r_s)_q$ have a **u**nique **parity-sequence* ***R...- E - T Terminating** at the ***r**$_s$'th **Type E-operation**, as the unique **parity-sequence* consists of **s O E -tiles** and (***r**$_s$**-s**) **E -tiles.**

**s = # O E -tiles** | ***r**$_s$ **- s = *⌈sZ⌉ - s = ⌈(Z-1)s⌉ = # E -tiles** (incl. the * **E -tile**)

Each **Group** (**s**,***r**$_s$) in the **List** contains the *complete set* of **reducing* **parity-sequences* which have **s** = # **O** type **Odd** Operations and notice that **Groups** (**s**,***r**$_s$), **s** ≥ 2, in the **List** contains the **complete set** of $^\Phi$**Diverging** $^\Phi$parity-sequences which have (**s-1**) type **Odd** Operations. **These** $^\Phi$parity-sequences in **List 2** have (**s-1**) **O E -tiles** & (**r-s**) *colourless* **E -tiles**; (**s-1**) ≤ **r** < ⌈(**s-1**)**Z**⌉

Notice in the **Group** for (**s**,***r**$_s$) = (4,*7$_4$) where **G**(4,*7$_4$) = 3 for $*B_{Con}$ ∊ {*121,*113,*69} that:

*69 ≡ **5** (mod 32) *113 ≡ **17** (mod 32) *121 ≡ **25** (mod 32)

and also notice that the **parity-sequences* are indeed those found in layer **r** = *7 in **Graph 5** i.e. the *Identical extensions* to the **three** $^\Phi$parity-sequences having 4 **O E -tiles** & 1 **E -tile** in **r** = 5.

In the *Next* **Group** (5,*8$_5$) exist **two** *triplets* of **parity-sequences* which are found in **Graph 5** & one **sequence* is found in **Graph 4** i.e. **the** one in **layer r** = *8 for $*B_{Con}$; *161 ≡ **1** (mod 32). In all **Groups** (**s**,***r**$_s$) for all higher **s** *triplets* of **parity-sequences* having *Identical extensions* exist:

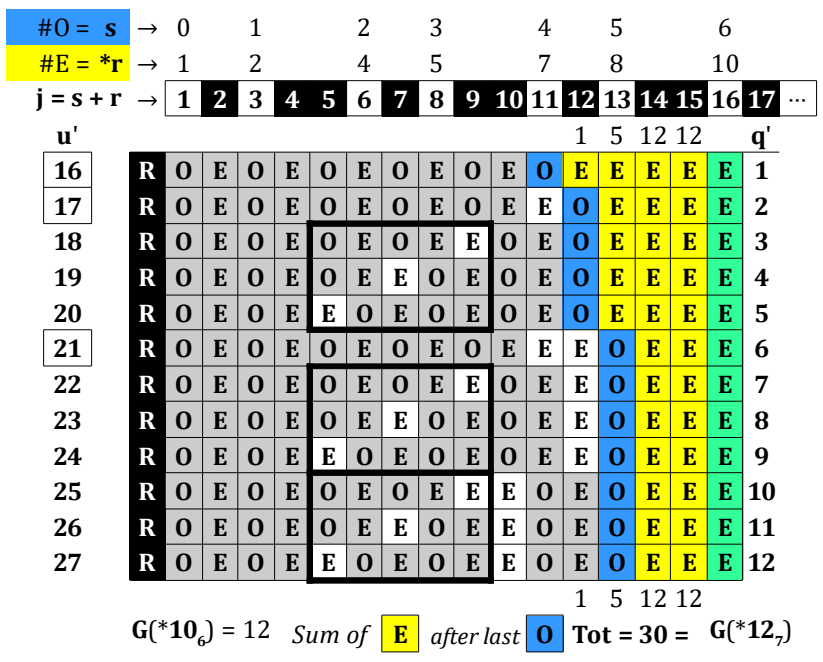

**List 2***10$_6$. (**s**,***r**$_s$) = (6,*10$_6$)

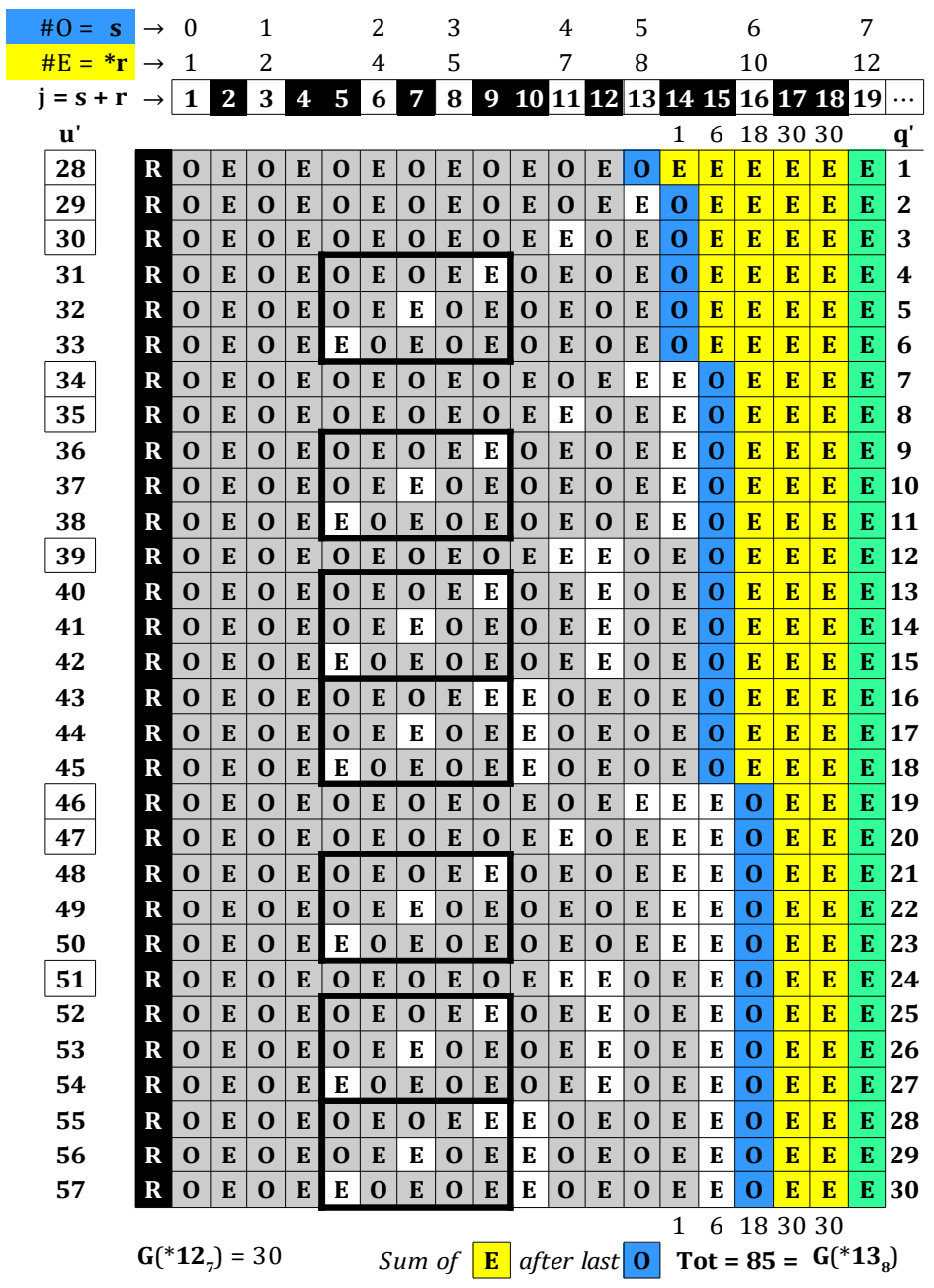

**List 2***12$_7$. (**s**,***r**$_s$) = (7,*12$_7$)

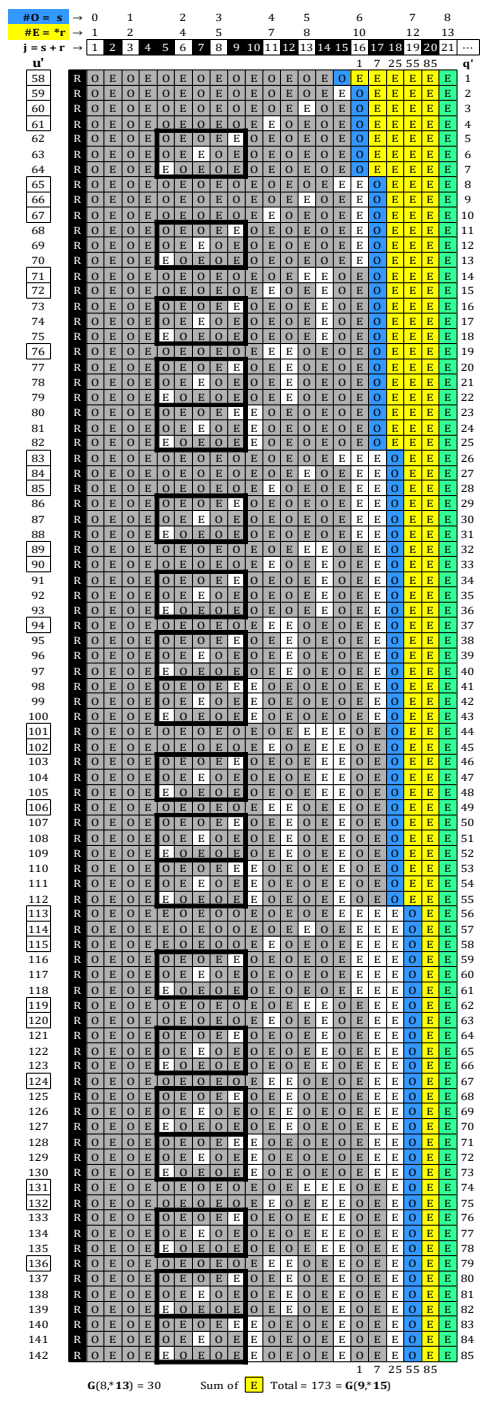

**List 2***13$_8$. *Next* **List**

### Definition 3

The **parity-sequences* in the *Next* **List** are ***generated*** by ***copying*** the $^\Phi$**Diverging** $^\Phi$parity-sequences incl. the final **O E -tile** in a **List 2***r$_s$ and in the *Next* **List 2***r$_{s+1}$ extend **each** $^\Phi$parity-sequence with **1 O E -tile** & an ***appropriate amount*** of **E -tiles.**

**OBS**! The ***copying*** is done *column by column* to ***ensure*** **uniqueness.**

If the *final* **E -tile** is *removed* from all **reducing* **parity-sequences* in a *collection* for a **constant** **Reducing Combination* (**s**,***r**$_s$) what is ***left*** is a *collection* of *longest possible* $^\Phi$**Diverging** $^\Phi$**parity-sequences**, as the $^\Phi$**Diverging** Combination $^\Phi$(**s**,***r**$_s$**-1**) is observed i.e. there exist a ***one-to-one-****correspondence* between the $^\Phi$**Diverging IV**-classes $^\Phi{}_{u'}$(**s**,***r**$_s$**-1**)$_{q'}$ and the **Converging* **IV**-classes $*_{u'}$(**s**,***r**$_s$)$_{q'}$. Please notice that the mentioned *longest possible* $^\Phi$**parity-sequences** $^\Phi$**R... E'-T** for the $^\Phi$**Diverging IV**-classes represents **subsets** of **IV**-classes identified with shorter $^\Phi$**Diverging** $^\Phi$sequences and that the *longest possible* $^\Phi$**parity-sequences** are re-found in the *last* sub-group (*max* **E** *before last* **O E**) in the *Next* **List** for (**s+1**,***r**$_{s+1}$). As observed above there exist a **1-*to*-1-***correspondence* for *sequences*/**IV**-classes $*Con_{u'\geq 3}(*\mathbf{r}_s)_{q'}$ ↔ $^\Phi\mathbf{Div}_{u'\geq 3}(\mathbf{r}_{s\text{-}1})_{q'}$ for (**s-1**) ≤ **r** < ⌈(**s-1**)**Z**⌉.

## ***The Sequence-formula* and *The Series-formula***

**Comment**

For the analysis it is very fortunate, that when *splitting* **T**erminal **V**alue classes $[3^s\mathbf{X}\text{-}B_{j=s+r}]_{TV}$ *one* of the parity-sequences for *one* of the **Child**-classes is extended with one **E-tile** and the *other* parity-sequence for the *other* **Child**-class is extended with one **O E-tile**. It is observed, that the **Combination** (**s**,**r**) *splits* into (**s**,**r+1**) & (**s+1**,**r+1**) and **this** is in actual fact **identical** to the definition of ***Pascal's Triangle*** i.e. there exist a ***Binomial-distribution*** of (**s**,**r**). The following **Figure** is a (left-justified) ***Pascal's Triangle*** with **Reducing Combinations* marked.

| *r | 1 | 2 | 4 | 5 | 7 | 8 | 10 | 12 | 13 | 15 | 16 | 18 | 20 | ··· |
|---|---|---|---|---|---|---|---|---|---|---|---|---|---|---|
| **s** | **0** | **1** | **2** | **3** | **4** | **5** | **6** | **7** | **8** | **9** | **10** | **11** | **12** | ··· |
| **r** | | | | | | | | | | | | | | |
| **0** | 1 | | | | | | | | | | | | | |
| **1** | 1 | 1 | | | | | | | | | | | | |
| **2** | 1 | 2 | 1 | | | | | | | | | | | |
| **3** | 1 | 3 | 3 | 1 | | | | | | | | | | |
| **4** | 1 | 4 | 6 | 4 | 1 | | | | | | | | | |
| **5** | 1 | 5 | 10 | 10 | 5 | 1 | | | | | | | | |
| **6** | 1 | 6 | 15 | 20 | 15 | 6 | 1 | | | | | | | |
| **7** | 1 | 7 | 21 | 35 | 35 | 21 | 7 | 1 | | | | | | |
| **8** | 1 | 8 | 28 | 56 | 70 | 56 | 28 | 8 | 1 | | | | | |
| **9** | 1 | 9 | 36 | 84 | 126 | 126 | 84 | 36 | 9 | 1 | | | | |
| **10** | 1 | 10 | 45 | 120 | 210 | 252 | 210 | 120 | 45 | 10 | 1 | | | |
| **11** | 1 | 11 | 55 | 165 | 330 | 462 | 462 | 330 | 165 | 55 | 11 | 1 | | |
| **12** | 1 | 12 | 66 | 220 | 495 | 792 | 924 | 792 | 495 | 220 | 66 | 12 | 1 | |
| **13** | 1 | 13 | 78 | 286 | 715 | 1287 | 1716 | 1716 | 1287 | 715 | 286 | 78 | 13 | ··· |
| **14** | 1 | 14 | 91 | 364 | 1001 | 2002 | 3003 | 3432 | 3003 | 2002 | 1001 | 364 | 91 | ··· |
| **15** | 1 | 15 | 105 | 455 | 1365 | 3003 | 5005 | 6435 | 6435 | 5005 | 3003 | 1365 | 455 | ··· |
| **16** | 1 | 16 | 120 | 560 | 1820 | 4368 | 8008 | 11440 | 12870 | 11440 | 8008 | 4368 | 1820 | ··· |
| **17** | 1 | 17 | 136 | 680 | 2380 | 6188 | 12376 | 19448 | 24310 | 24310 | 19448 | 12376 | 6188 | ··· |
| **18** | 1 | 18 | 153 | 816 | 3060 | 8568 | 18564 | 31824 | 43758 | 48620 | 43758 | 31824 | 18564 | ··· |
| **19** | 1 | 19 | 171 | 969 | 3876 | 11628 | 27132 | 50388 | 75582 | 92378 | 92378 | 75582 | 50388 | ··· |
| **20** | 1 | 20 | 190 | 1140 | 4845 | 15504 | 38760 | 77520 | 125970 | 167960 | 184756 | 167960 | 125970 | ··· |
| ··· | ··· | ··· | ··· | ··· | ··· | ··· | ··· | ··· | ··· | ··· | ··· | ··· | ··· | ··· |

**Figure 3a.** ***Pascal's Triangle*** for **r** *choose* **s** with **Reducing Combinations* marked.

## ***The Sequence-formula***

It is possible to calculate the number of **Converging* **IV**-classes/**parity-sequences* for all **modulo** $A_{IV} = 2^{*r_s}$ by the help of the ***Binomial-distribution***. The above marked *binomial-coefficients* for the **Reducing Combinations* (**s**,***r**$_s$) are the **Totals** in the below **Formula (2a)**.

$$|{*}Con({*}r_{*s})| = G({*}s,{*}r_{*s}) = \binom{\text{※}r_s}{\text{※}s}_{\text{Total}} - \sum_{s=0}^{s=\text{※}s-1} G(s,\text{※}r)\binom{\text{※}r_s-\text{※}r}{\text{※}s-s}_{\text{Redundant}} \quad , s \ge 0 \qquad \textbf{(2a)}$$

The *second* term in (**2a**) is ***The Weighted Sum*** of **Redundant Combinations** from the earlier **Converging* ***Combinations*** i.e. for lower (**s**,***r**$_s$) < (***s**,***r**$_{*s}$) which have to be ***subtracted*** from the **Total** to find the **correct** amount of **Converging* **IV**-classes for the current **modulo** $2^{*r_{*s}}$.

**Comment**

It is already *known* that $*_1[2^{*1_0}\mathbf{X}\text{-}{*}0_{Con}]_1$ with the **reducing* **parity-sequence* ***R E - T** exist, and that $^{\Phi}[2^1\mathbf{X}\text{-}1_0]_{IV}$ with the $^{\Phi}$parity-sequence $^{\Phi}$**R O E - T** is a $^{\Phi}$**Diverging IV**-class, so it is possible to ignore the **IV**-classes having $*N_X \equiv \mathbf{2} \pmod{2^{*1_0}} \Leftrightarrow {}^{\Delta}B \equiv *\mathbf{0}_{Con} \pmod{2^1}$ in **Formula (2b)**:

$$G({*}s,{*}\lceil sZ\rceil) = \binom{\text{※}\lceil sZ\rceil-1}{\text{※}s-1}_{\text{Total}} - \sum_{s=1}^{s=\text{※}s-1} G(s,\lceil sZ\rceil)\binom{\text{※}\lceil sZ\rceil-\lceil sZ\rceil}{\text{※}s-s}_{\text{Redundant}} \quad , s \ge 1 \qquad \textbf{(2b)}$$

Remember **Formula (1)**: $*r_s = *\lceil sZ\rceil,\ s \in \mathbb{N},\ Z = \log 3/\log 2 = 1.5849625007\ldots \qquad \textbf{(1)}$

**OBS**! Please notice that **Formula (2b)** is a *function* of number of type Odd operations **s** = # **O**.

It is also *known* that $*_2[2^{*2_1}\mathbf{X}\text{-}{*3}_{Con}]_1$ with the **reducing* *parity-sequence* *R O E - E - T exist and that $^{\phi}[2^2\mathbf{X}\text{-}1_0]_{IV}$ with $^{\phi}$parity-sequence $^{\phi}$R O E - O E - T is a $^{\phi}$**Diverging IV**-class, so it is **also** possible to ignore **IV**-classes, $*N_X \equiv \mathbf{1} \pmod{2^{*2_1}} \Leftrightarrow {}^{\Delta}B \equiv *\mathbf{3}_{Con} \pmod{2^2}$ in **Formula (2c)**:

$$\mathbf{G}(*\mathbf{s},*\lceil \mathbf{sZ}\rceil) = \binom{\text{※}\lceil \mathbf{sZ}\rceil - 2}{\text{※}\mathbf{s}-2}_{\text{Total}} - \sum_{\mathbf{s}=2}^{\mathbf{s}=\text{※}\mathbf{s}-1} \mathbf{G}(\mathbf{s},\lceil \mathbf{sZ}\rceil)\binom{\text{※}\lceil \mathbf{sZ}\rceil - \lceil \mathbf{sZ}\rceil}{\text{※}\mathbf{s}-\mathbf{s}}_{\text{Redundant}} \quad ,\mathbf{s} \geq 2 \qquad \mathbf{(2c)}$$

In **Formula (2c)** the **Totals** are found in **Figure 3b** where *subsets* of $^{\phi}[2^2\mathbf{X}\text{-}1_0]_{IV}$ are observed:

| | *r | 4 | 5 | 7 | 8 | 10 | 12 | 13 | 15 | 16 | 18 | 20 | ··· |
|---|---|---|---|---|---|---|---|---|---|---|---|---|---|
| | **s** | **2** | **3** | **4** | **5** | **6** | **7** | **8** | **9** | **10** | **11** | **12** | ··· |
| **r** | | | | | | | | | | | | | |
| **2** | | 1 | | | | | | | | | | | |
| **3** | | 1 | 1 | | | | | | | | | | |
| **4** | | 1 | 2 | 1 | | | | | | | | | |
| **5** | | 1 | 3 | 3 | 1 | | | | | | | | |
| **6** | | 1 | 4 | 6 | 4 | 1 | | | | | | | |
| **7** | | 1 | 5 | 10 | 10 | 5 | 1 | | | | | | |
| **8** | | 1 | 6 | 15 | 20 | 15 | 6 | 1 | | | | | |
| **9** | | 1 | 7 | 21 | 35 | 35 | 21 | 7 | 1 | | | | |
| **10** | | 1 | 8 | 28 | 56 | 70 | 56 | 28 | 8 | 1 | | | |
| **11** | | 1 | 9 | 36 | 84 | 126 | 126 | 84 | 36 | 9 | 1 | | |
| **12** | | 1 | 10 | 45 | 120 | 210 | 252 | 210 | 120 | 45 | 10 | 1 | |
| **13** | | 1 | 11 | 55 | 165 | 330 | 462 | 462 | 330 | 165 | 55 | 11 | ··· |
| **14** | | 1 | 12 | 66 | 220 | 495 | 792 | 924 | 792 | 495 | 220 | 66 | ··· |
| **15** | | 1 | 13 | 78 | 286 | 715 | 1287 | 1716 | 1716 | 1287 | 715 | 286 | ··· |
| **16** | | 1 | 14 | 91 | 364 | 1001 | 2002 | 3003 | 3432 | 3003 | 2002 | 1001 | ··· |
| **17** | | 1 | 15 | 105 | 455 | 1365 | 3003 | 5005 | 6435 | 6435 | 5005 | 3003 | ··· |
| **18** | | 1 | 16 | 120 | 560 | 1820 | 4368 | 8008 | 11440 | 12870 | 11440 | 8008 | ··· |
| **19** | | 1 | 17 | 136 | 680 | 2380 | 6188 | 12376 | 19448 | 24310 | 24310 | 19448 | ··· |
| **20** | | 1 | 18 | 153 | 816 | 3060 | 8568 | 18564 | 31824 | 43758 | 48620 | 43758 | ··· |
| ··· | | ··· | ··· | ··· | ··· | ··· | ··· | ··· | ··· | ··· | ··· | ··· | ··· |

**Figure 3b.** ***Pascal's Triangle*** for $^{\phi}[2^2\mathbf{X}\text{-}1_0]_{IV}$. The **Totals** for **Formula (2c)** are marked.

**Formula (2c)** is named ***The Sequence-formula*** as the result is a ***sequence*** $\mathbf{G}(\mathbf{s},*\mathbf{r_s})$ for $\mathbf{s} \to \infty$:

$$|*Con(*\mathbf{r_s})| = \mathbf{G}(*\mathbf{s},*\mathbf{r_s}) \in \{(\mathbf{1}{*1_0}), (\mathbf{1}{*2_1}), \mathbf{1}{*4_2}, \mathbf{2}{*5_3}, \mathbf{3}{*7_4}, \mathbf{7}{*8_5}, \mathbf{12}{*10_6}, \mathbf{30}{*12_7}, \ldots\}$$

## Example of using *The Sequence-formula*

$$|*Con(*\mathbf{7}_{*4})| = \mathbf{G}(*\mathbf{4},*\mathbf{7}_{*4}) = \binom{\text{※}\mathbf{7_4}-\mathbf{2}}{\text{※}\mathbf{4}-\mathbf{2}}_{\text{Total}} - \sum_{\mathbf{s}=2}^{\mathbf{s}=\text{※}\mathbf{s}-1} \mathbf{G}(\mathbf{s},\text{※}\mathbf{r})\binom{\text{※}\mathbf{7_4}-\text{※}\mathbf{r}}{\text{※}\mathbf{4}-\mathbf{s}}_{\text{Redundant}}$$

$$= \binom{\mathbf{5}}{\mathbf{2}}_{\text{Total}} - \mathbf{G}(\mathbf{2},\text{※}\,\mathbf{4})\binom{\mathbf{3}}{\mathbf{2}}_{\text{Redundant}} - \mathbf{G}(\mathbf{3},\text{※}\,\mathbf{5})\binom{\mathbf{2}}{\mathbf{1}}_{\text{Redundant}}$$

$$\Rightarrow \quad \mathbf{G}(*\mathbf{4},*\mathbf{7}_{*4}) = (\mathbf{10})_{\text{Total}} - \big(1(\mathbf{3})_{\text{Redundant}} + 2(\mathbf{2})_{\text{Redundant}}\big) = \mathbf{10} - (\mathbf{3}+\mathbf{4}) = \underline{\underline{3}}$$

## Comment

In **Resultlist 1** on the following page the first results $\mathbf{G}(\mathbf{s},*\mathbf{r_s}) = |*Con_{\mathbf{u}}(*\mathbf{r_s})_{\mathbf{q}}|$ are observed.

Remember that $|*Con_{\mathbf{u}}(*\mathbf{r_s})_{\mathbf{q}}| > 0$ and that $|*Con(^{\phi}\mathbf{r_v})| = 0$ for $^{\phi}\mathbf{r_v} \neq \lceil \mathbf{sZ}\rceil$, $\mathbf{s} \in \mathbb{N}$, $\mathbf{Z} = log\mathbf{3}/log\mathbf{2}$.

From $|*Con_{\mathbf{u}}(*\mathbf{r_s})_{\mathbf{q}}|$ it is possible to calculate all remaining data for $|^{\phi}\mathbf{Div}(\mathbf{r})|$ and $|^{\Delta}\mathbf{Red}(\mathbf{r})|$ as;

$$|^{\phi}\mathbf{Div}(\mathbf{r+1})| = 2|^{\phi}\mathbf{Div}(\mathbf{r})| - |*Con_{\mathbf{u}}(\mathbf{r+1})_{\mathbf{q}}|$$

$$|^{\Delta}\mathbf{Red}(\mathbf{r+1})| = 2|^{\Delta}\mathbf{Red}(\mathbf{r})| + 2|*Con_{\mathbf{u}}(\mathbf{r})_{\mathbf{q}}|$$

**OBS**! In The Domino-Tree **Converging* **IV**-classes $*_{\mathbf{u}}[2^{*\mathbf{r_s}}\mathbf{X}\text{-}*B_{Con}]_{\mathbf{q}}$ **exist** in all ***layers*** $*\mathbf{r_s} = \lceil \mathbf{sZ}\rceil$

## The number of different *types* of IV-classes for increasing **Reducing Combinations* (s,*r$_s$)

| s = #O | *r$_s$ = #E | *r$_s$-s | *r$_s$+s | Total = $2^{*r_s}$ | $\vert{}^{\Delta}$Red(*r$_s$)$\vert$ | $\vert{}^{\Phi}$Div(*r$_s$)$\vert$ | $\vert$*Con(*r$_s$)$\vert$ | $\sum \vert$*Con(*r$_s$)$\vert$ |
|---|---|---|---|---|---|---|---|---|
| **0** | 1 | 1 | 1 | 2 | **0** | 1 | **1** | 1 |
| **1** | 2 | 1 | 3 | 4 | 2 | 1 | **1** | 2 |
| **2** | 4 | 2 | 6 | 16 | 12 | 3 | **1** | 3 |
| **3** | 5 | 2 | 8 | 32 | 26 | 4 | **2** | 5 |
| **4** | 7 | 3 | 11 | 128 | 112 | 13 | **3** | 8 |
| **5** | 8 | 3 | 13 | 256 | 230 | 19 | **7** | 15 |
| **6** | 10 | 4 | 16 | 1024 | 948 | 64 | **12** | 27 |
| **7** | 12 | 5 | 19 | 4096 | 3840 | 226 | **30** | 57 |
| **8** | **13** | 5 | 21 | **8192** | **7740** | 367 | **85** | **142** |
| **9** | 15 | 6 | 24 | 32768 | 31300 | 1295 | **173** | 315 |
| **10** | 16 | 6 | 26 | 65536 | 62946 | 2114 | **476** | 791 |
| **11** | 18 | 7 | 29 | 262144 | 253688 | 7495 | **961** | 1752 |
| **12** | 20 | 8 | 32 | 1048576 | 1018596 | 27328 | **2652** | 4404 |
| **13** | 21 | 8 | 34 | 2097152 | 2042496 | 46611 | **8045** | 12449 |
| **14** | 23 | 9 | 37 | 8388608 | 8202164 | 168807 | **17637** | 30086 |
| **15** | 24 | 9 | 39 | 16777216 | 16439602 | 290496 | **47118** | 77204 |
| **16** | 26 | 10 | 42 | 67108864 | 65946880 | 1074149 | **87835** | 165039 |
| **17** | 27 | 10 | 44 | 134217728 | 132069430 | 1852478 | **295820** | 460859 |
| **18** | 29 | 11 | 47 | 536870912 | 529461000 | 6840772 | **569140** | 1029999 |
| **19** | 31 | 12 | 50 | 2147483648 | 2120120560 | 25841433 | **1521655** | 2551654 |
| **20** | 32 | 12 | 52 | 4294967296 | 4243284430 | 46010008 | **5672858** | 8224512 |
| **21** | 34 | 13 | 55 | 17179869184 | 16995829152 | 172315631 | **11724401** | 19948913 |
| **22** | 35 | 13 | 57 | 34359738368 | 34015107106 | 306244032 | **38387230** | 58336143 |
| … | … | … | … | … | … | … | … | … |
| **s → ∞** | **⌈sZ⌉** | **⌈(Z-1)s⌉** | **⌈(Z+1)s⌉** | $2^{\lceil sZ \rceil}$ | $\vert{}^{\Delta}$**Red(*r$_s$)**$\vert \to 2^{\lceil sZ \rceil}$ | **q' → ∞** | **q → ∞** | **u → ∞** |
| … | … | … | … | … | … | … | … | … |

**Resultlist 1.** The number of different **IV**-classes for increasing **Reducing Combinations* (**s**,***r$_s$**)

For *sufficiently high* **modulo** $A_{IV} = 2^r$ it is true for the ***Union*** of **IV**-classes $[2^r\mathbf{X}\text{-}B_0]_{IV}$, $B_0 < 2^r$, that:

$2^{r-*1_0}/2^r = {}^1/_2$ *of the possible residue-classes for* $B_0 < 2^r$ *can be associated with* $*_1[2^{*1_0}\mathbf{X}\text{-}*0_{Con}]_1$

$2^{r-*2_1}/2^r = {}^1/_4$ *of the possible residue-classes for* $B_0 < 2^r$ *can be associated with* $*_2[2^{*2_1}\mathbf{X}\text{-}*3_{Con}]_1$

$2^{r-*4_2}/2^r = {}^1/_{16}$ *of the possible residue-classes for* $B_0 < 2^r$ *can be associated with* $*_3[2^{*4_2}\mathbf{X}\text{-}*13_{Con}]_1$

$2^{r-*5_3}/2^r = {}^1/_{32}$ *of the possible residue-classes for* $B_0 < 2^r$ *can be associated with* $*_4[2^{*5_3}\mathbf{X}\text{-}*21_{Con}]_1$

$2^{r-*5_3}/2^r = {}^1/_{32}$ *of the possible residue-classes for* $B_0 < 2^r$ *can be associated with* $*_5[2^{*5_3}\mathbf{X}\text{-}*9_{Con}]_2$

$2^{r-*7_4}/2^r = {}^1/_{128}$ *of the possible residue-classes for* $B_0 < 2^r$ *can be associated with* $*_6[2^{*7_4}\mathbf{X}\text{-}*121_{Con}]_1$

$2^{r-*7_4}/2^r = {}^1/_{128}$ *of the possible residue-classes for* $B_0 < 2^r$ *can be associated with* $*_7[2^{*7_4}\mathbf{X}\text{-}*113_{Con}]_2$

$2^{r-*7_4}/2^r = {}^1/_{128}$ *of the possible residue-classes for* $B_0 < 2^r$ *can be associated with* $*_8[2^{*7_4}\mathbf{X}\text{-}*69_{Con}]_3$

…

$2^{r-*r_s}/2^r = {}^1/2^{*r_s}$ *of the possible residue-classes for* $B_0 < 2^r$ *can be associated with* $*_u[2^{*r_s}\mathbf{X}\text{-}*B_{Con}]_q$

…

A ***density*** ${}^1/2^{*r_s}$ can be found for all $*_u[2^{*r_s}\mathbf{X}\text{-}*B_{Con}]_q$, **u** ∈ ℕ, for all $2^{*r_s} \leq 2^r$.

## ***The Series-Formula***

### **Definition 4**

The result of **Formula (3)** is the **F**raction of **IV**-classes $[2^r\mathbf{X}\text{- }B_0]_{IV}$, $B_0 < 2^r$, for a **Reducing Combination* (**s**,*r$_s$) which are **reduced* in **exactly *j = s+*r$_s$** Operations **Odd** & **Even**.

| | | |
|---|---|---|
| **F(s,*r$_s$) = G(s,*r$_s$)/$2^{*r_s}$** | this **F**raction (***ratio***) is a ***density*** | **(3)** |

The result of **Formula (4)** is the **S**um of the above **F**ractions i.e. of the *subsets* of **IV**-classes $*_u[2^{*r_s}\mathbf{X}\text{-}*B_{Con}]_q$ for *constant* (**s**,*r$_s$) which are **reduced* in **at most *j = s+*r$_s$** Operations **O** & **E**

| | | |
|---|---|---|
| **S(s,*r$_s$)** = $\sum_{s=0}^{s} \mathbf{F}(\mathbf{s}, * \mathbf{r}_s)$ | this **S**um (of **F**ractions) is a ***Density*** | **(4)** |

**Formula (4)** is named ***The Series-formula*** as the result is a ***series*** $\mathbf{S(s,{}^*r_s) \rightarrow 1 \text{ for } s \rightarrow \infty}$

The first results from **Formula (3)** and **Formula (4)** are found in the following **Resultlist 2**:

| s = #O | $^*r_s$ = #E | $^*r_s$–s | $^*r_s$+s | Total = $2^{*r_s}$ | G(s,$^*r_s$) | F(s,$^*r_s$) | S(s,$^*r_s$) | 1 - S(s,$^*r_s$) | F(s,$^*r_s$) / (1 - S(s,$^*r_s$)) |
|---|---|---|---|---|---|---|---|---|---|
| **0** | 1 | 1 | 1 | 2 | **1** | 0.500000 | 0.500000 | 0.500000 | 1.000000 |
| **1** | 2 | 1 | 3 | 4 | **1** | 0.250000 | **0.750000** | 0.250000 | 1.000000 |
| **2** | 4 | 2 | 6 | 16 | **1** | 0.062500 | 0.812500 | 0.187500 | 0.333333 |
| **3** | 5 | 2 | 8 | 32 | **2** | 0.062500 | 0.875000 | 0.125000 | 0.500000 |
| **4** | 7 | 3 | 11 | 128 | **3** | 0.023438 | 0.898438 | 0.101563 | 0.230769 |
| **5** | 8 | 3 | 13 | 256 | **7** | 0.027344 | **0.925781** | 0.074219 | 0.368421 |
| **6** | 10 | 4 | 16 | 1024 | **12** | 0.011719 | 0.937500 | 0.062500 | 0.187500 |
| **7** | 12 | 5 | 19 | 4096 | **30** | 0.007324 | 0.944824 | 0.055176 | 0.132743 |
| **8** | **13** | 5 | 21 | **8192** | **85** | **0.010376** | **0.955200** | **0.044800** | 0.231608 |
| **9** | 15 | 6 | 24 | 32768 | **173** | 0.005280 | 0.960480 | 0.039520 | 0.133591 |
| **10** | 16 | 6 | 26 | 65536 | **476** | 0.007263 | 0.967743 | 0.032257 | 0.225166 |
| **11** | 18 | 7 | 29 | 262144 | **961** | 0.003666 | 0.971409 | 0.028591 | 0.128219 |
| **12** | 20 | 8 | 32 | 1048576 | **2652** | 0.002529 | 0.973938 | 0.026062 | 0.097043 |
| **13** | 21 | 8 | 34 | 2097152 | **8045** | 0.003836 | 0.977774 | 0.022226 | 0.172599 |
| **14** | 23 | 9 | 37 | 8388608 | **17637** | 0.002102 | 0.979877 | 0.020123 | 0.104480 |
| **15** | 24 | 9 | 39 | 16777216 | **47118** | 0.002808 | 0.982685 | 0.017315 | 0.162198 |
| **16** | 26 | 10 | 42 | 67108864 | **87835** | 0.001309 | 0.983994 | 0.016006 | 0.081772 |
| **17** | 27 | 10 | 44 | 134217728 | **295820** | 0.002204 | 0.986198 | 0.013802 | 0.159689 |
| **18** | 29 | 11 | 47 | 536870912 | **569140** | 0.001060 | 0.987258 | 0.012742 | 0.083198 |
| **19** | 31 | 12 | 50 | 2147483648 | **1521655** | 0.000709 | 0.987967 | 0.012033 | 0.058884 |
| **20** | 32 | 12 | 52 | 4294967296 | **5672858** | 0.001321 | 0.989287 | 0.010713 | 0.123296 |
| **21** | 34 | 13 | 55 | 17179869184 | **11724401** | 0.000682 | 0.989970 | 0.010030 | 0.068040 |
| **22** | 35 | 13 | 57 | 34359738368 | **38387230** | 0.001117 | **0.991087** | 0.008913 | 0.125348 |
| … | … | … | … | … | … | … | … | … | … |
| **#O=#OE** | **#E=#OE+#E$_t$** | **#E$_t$=#E-#OE** | **#E+#O** | $2^{\lceil sZ \rceil}$ | $q_s \rightarrow \infty$ | 0.00⋯5 → **0** | 0.99⋯5 → **1** | 0.00⋯5 → **0** | **NB**! New **min** for **s** → ∞ |
| … | … | … | … | … | … | … | … | … | |

The index t in #E$_t$ indicates number of E-**tiles** in parity-sequences.

**Resultlist 2**. The **F**ractions and **S**ums for the **Converging* ***IV**-classes/**parity-sequences*.

**OBS**! Please notice that $\mathbf{S(s+1,{}^*r_{s+1}) > S(s,{}^*r_s)}$ & $\mathbf{F(s,{}^*r_s) > 0}$ for $\mathbf{s} \geq 0$

**OBS**! Please notice that for **s** ≥ 2 it is observed that $\mathbf{G(s+1,{}^*r_{s+1}) > G(s,{}^*r_s)}$

**OBS**! $|{}^\Phi\mathbf{Div}(\mathbf{r})| + |{}^*Con(\mathbf{r})| + |{}^\Delta\mathbf{Red}(\mathbf{r})| = 2^r \Rightarrow \left(|{}^\Phi\mathbf{Div}({}^*\mathbf{r}_s)| + |{}^*Con({}^*\mathbf{r}_s)| + |{}^\Delta\mathbf{Red}({}^*\mathbf{r}_s)|\right)/2^{*r_s} = \mathbf{1}$

which implies:

$$\mathbf{S(s,{}^*r_s)} = \left(|{}^*Con({}^*\mathbf{r}_s)| + |{}^\Delta\mathbf{Red}({}^*\mathbf{r}_s)|\right)/2^{*r_s} \qquad \& \qquad \mathbf{1-S(s,{}^*r_s)} = |{}^\Phi\mathbf{Div}({}^*\mathbf{r}_s)|/2^{*r_s}$$

$$\Rightarrow \quad \mathbf{1-S(k,{}^*r_k)} = \mathbf{1} - \sum_{s=0}^{s=k} \mathbf{F(s,{}^*r_s)} = |{}^\Phi\mathbf{Div}({}^*\mathbf{r}_k)|/2^{*r_k} = \sum_{s=k+1}^{\infty} \mathbf{F(s,{}^*r_s)} \qquad \textbf{Formula (5)}$$

## Comment

In **Resultlist 2** the values in the column **1-S(*r$_s$)** are **identical** to values for the **Union** of the $^\Phi$**Diverging** $^\Phi$parity-sequences/$^\Phi$**IV**-classes from R. Terras[III]. Notice the ***Asymptotic Densities***;

$\mathbf{S(s,{}^*r_s) \rightarrow 1_-}$ for $^*\mathbf{r_s} \rightarrow \infty$ $\qquad$ $\mathbf{1\text{-}S(s,{}^*r_s) \rightarrow 0_+}$ for $^*\mathbf{r_s} \rightarrow \infty$

**NB**! The **Discussion** contains more details about the ***Asymptotic Densities*** and **Formula (5)**.

Notice that the **Sum**/***Density*** **S(s,*r$_s$)** is *simultaneously* an expression for the *fraction* of N ≤ $2^{*r_s}$ which are **reduced* in **at most** ***j** = (**s**+***r$_s$**) = (# **O** +*# **E** ) operations, **and** an expression for the *fraction* of $N_X \in {}^\Phi[2^0\mathbf{X}\text{-}0_0]_{IV/TV}$ (all N ∈ ℕ) which are **reduced* in **at most** ***j** = (**s**+***r$_s$**) operations.

The **Union** of the $2^{*r_s}$ pcs of $[2^{*r_s}\mathbf{X}\text{-}B_0]_{IV}$, $0 \leq B_0 < 2^{*r_s}$, is a ***Covering System*** and all the **IV**-classes have ***different*** parity-sequences - but with the ***same*** number of **tiles**; ***r$_s$** = (# **O** **E**$_{tile}$ + # **E**$_{tile}$)

## The *Inverse* Densities

When the *numbers* of the three different **types** of **IV**-classes from **Resultlist 1** are divided by the **Total** $2^{*r_s}$, the ***ratios/densities*** found in **Resultlist 2** are generated. If the *inverse* operation is performed i.e. $2^{*r_s}$ is divided by the three *numbers*, the values in **Resultlist 3** are generated:

| s = #O | $*r_s$ = #E | $*r_s$–s | $*r_s$+s | Total = $2^{*r_s}$ | $2^{*r}$/ $\|^{\Delta}$Red($*r_s$)\| | $2^{*r}$/ $\|^{\phi}$Div($*r_s$)\| | $2^{*r}$/ \|*Con($*r_s$)\| |
|---|---|---|---|---|---|---|---|
| **0** | **1** | 1 | 1 | 2 | | **2.00** | 2.00 |
| **1** | 2 | 1 | 3 | 4 | 2.000000 | **4.00** | 4.00 |
| **2** | 4 | 2 | 6 | 16 | 1.333333 | **5.33** | 16.00 |
| **3** | 5 | 2 | 8 | 32 | 1.230769 | **8.00** | 16.00 |
| **4** | 7 | 3 | 11 | 128 | 1.142857 | **9.85** | 42.67 |
| **5** | 8 | 3 | 13 | 256 | 1.113043 | **13.47** | 36.57 |
| **6** | 10 | 4 | 16 | 1024 | 1.080169 | **16.00** | 85.33 |
| **7** | 12 | 5 | 19 | 4096 | 1.066667 | **18.12** | 136.53 |
| **8** | **13** | 5 | **21** | **8192** | **1.058398** | **22.32** | **96.38** |
| **9** | 15 | 6 | 24 | 32768 | 1.046901 | **25.30** | 189.41 |
| **10** | 16 | 6 | 26 | 65536 | 1.041146 | **31.00** | 137.68 |
| **11** | 18 | 7 | 29 | 262144 | 1.033332 | **34.98** | 272.78 |
| **12** | 20 | 8 | 32 | 1048576 | 1.029433 | **38.37** | 395.39 |
| **13** | 21 | 8 | 34 | 2097152 | 1.026759 | **44.99** | 260.68 |
| **14** | 23 | 9 | 37 | 8388608 | 1.022731 | **49.69** | 475.63 |
| **15** | 24 | 9 | 39 | 16777216 | 1.020537 | **57.75** | 356.07 |
| **16** | 26 | 10 | 42 | 67108864 | 1.017620 | **62.48** | 764.03 |
| **17** | 27 | 10 | 44 | 134217728 | 1.016266 | **72.45** | 453.71 |
| **18** | 29 | 11 | 47 | 536870912 | 1.013995 | **78.48** | 943.30 |
| **19** | 31 | 12 | 50 | 2147483648 | 1.012906 | **83.10** | 1411.28 |
| **20** | 32 | 12 | 52 | 4294967296 | 1.012180 | **93.35** | 757.11 |
| **21** | 34 | 13 | 55 | 17179869184 | 1.010829 | **99.70** | 1465.31 |
| **22** | 35 | 13 | 57 | 34359738368 | 1.010132 | **112.20** | 895.08 |
| … | … | … | … | … | … | … | … |
| s → ∞ | ⌈sZ⌉ → ∞ | ⌈(Z-1)s⌉ → ∞ | ⌈(Z+1)s⌉ → ∞ | $2^{\lceil sZ \rceil} \to 2^{\infty}$ | 1/0.99···5 → **1** | 1/0.00···5 → ∞ | 1/0.00···5 → ∞ |
| … | … | … | … | … | … | … | … |

**Resultlist 3**. ***The Average Distance*** between the **IV**-classes $\{^{\Delta},^{\phi},^{*}\}$

In **Resultlist 3** it is observed for $*\mathbf{r_s} \to \infty$ that:

$2^{*r_s}/|^{\Delta}$**Red**($*\mathbf{r_s}$)| is *strictly decreasing*

$2^{*r_s}/|^{\phi}$**Div**($*\mathbf{r_s}$)| is **strictly increasing**

$2^{*r_s}/|$**Con*($*\mathbf{r_s}$)| is increasing (not strictly)

"***The Average Distance***" is higher for **Converging* than for $^{\phi}$**Diverging IV**-classes for all $*\mathbf{r_s}$ **> 2.**

## Lemma 4

***The Average Distance*** in the Union of ***elements*** $N_X \in {}^{\phi}{}_{u'}[2^{*r_s}X\text{-}B_0]_{q'}$ grows *much* faster than $*\mathbf{r_s}$

## Proof for Lemma 4

In the individual $^{\phi}$**Diverging IV**-classes the *distance* between ***neighbour-elements*** $N_X$ & $N_{X+1}$ is the **modulo** $A_{IV} = 2^r$ and obviously $2^r$ tends to infinity *much* faster than **r** ∈ ℕ for **r** → ∞. In **Resultlist 3** it can be observed that ***The Average Distance*** between the $^{\phi}$**Diverging IV**-classes $2^{*r_s}/|^{\phi}$**Div**($*\mathbf{r_s}$)| *also* tends to infinity *much* faster than $*\mathbf{r_s}$ ∈ ℕ for $*\mathbf{r_s} \to \infty$. That both the *distance* between the *possible* **IV**-classes and the *distance* between individual *elements* in the **finite set** of **IV**-classes for constant $*\mathbf{r_s}$ tends to infinity *much* faster than $*\mathbf{r_s}$ *implies* that the *distance* between *elements* in the Union of ${}^{\phi}{}_{u'}[2^{*r_s}X\text{-}B_0]_{q'}$ tends to infinity *much* faster than $*\mathbf{r_s}$ ■

## Comment

In the **extreme** for "*very high*" **r** ∈ ℕ it is observed that the *distance* between ***elements*** in the **Union** of ${}^{\phi}{}_{u'}[2^rX\text{-}B_0]_{q'}$ is "*much higher*". For the hypothetical "**r** *close to* ∞" it is observed that the *distance* between ***elements*** in the **Union** of ${}^{\phi}{}_{u'}[2^rX\text{-}B_0]_{q'}$ tends to "*equal to infinity*" and for "**r** *equal to* ∞" the *distance* apparently becomes "*higher than infinity*" i.e. non-existent. This is an indication that a hypothetical **c**ounter**e**xample $N_{ce}$ can not exist, as $N_{ce}$ would need to be higher than **any** N ∈ ℕ, because **all** N ∈ ℕ belongs to an existing $*[2^{*r_s}X\text{-}*B_{Con}]_{IV}$ for ***some*** $*\mathbf{r_s}$ ∈ ℕ.

Please notice that the foundation for two important **conclusions** have been established:

- *There* ***exist*** *a countable infinity of *Converging* **IV**-*classes* $*_u[2^{*r_s}X\text{-}*B_{Con}]_q$ *for* **u** ∈ ℕ
- *There* ***exist*** *a countable infinity of* $^{\phi}$**Diverging IV**-*classes* ${}^{\phi}{}_{u'}[2^rX\text{-}B_0]_{q'}$ *for* **u'** ∈ ℕ, **u' ≥ 3**
  *These infinities of* **IV**-*classes are subsets of the* $^{\phi}$**Diverging *Universe-class*** ${}^{\phi}[2^0X\text{-}0_0]_{IV/TV}$
  *Notice that *Converging* **IV**-*classes (but not* $^{\phi}$**Diverging IV**-*classes) are pairwise* ***disjoint.***

***The Algorithm***

**Appendix Omega** contains the *code* (*Thonny-Python*) for ***The Algorithm*** incl. some examples. **NB**! It is strongly recommended to *have a look* at the two examples for N = **15** & N = ($2^{14}$**+15**) as they illustrate some important properties about the ***pair*** ($P_{IV}$+*$B_{Con}$) = $2^{*r_s}$ in *[$2^{*r_s}$**X**-*$B_{Con}$]$_{IV}$

***The Algorithm*** generates for ***any*** Test Value $N_N \in \mathbb{N}$:

- The **Converging* **IV**-class *[$2^{*r_s}$**X**-*$B_{Con}$]$_{IV}$ **NB**! ***The Algorithm*** is **not** able to find the **u** & **q**
- The corresponding **TV**-class in the ***Class-series*** *[$2^{*r_s}$**X**-*$B_{Con}$]$_{IV}$ →→ [$3^s$**X**-$B_{j=s+*r_s}$]$_{TV}$
- The **position X** i the **reducing* **IV**-class $N_X \in$ *[$2^{*r_s}$**X**-*$B_{Con}$]$_{IV}$ $N_X = N_N$
- The **Reducing Combination* (**s**,*$r_s$) & The **Ending Combination** ($s_d$,**d**)
- The **End-IV**-class †[$2^d$**X**-$B_{End}$]$_N$ →→ [$3^{s_d}$**X**-$B_{j=s_d+d}$] $B_{End} = (2^d - N_N)$ & $B_{j=s_d+d} = (3^{s_d} - {}^{\dagger}1_j)$

**Optional (if** print-lines activated):

- *The *reducing iterationseries* $N_0 \rightarrow ... \rightarrow {}^*N_{*j=s+*r_s} < N_0$ *for* $N_j$ *in* Base-**10** *and in* Base-**2**
- *The **shared** *reducing *parity-sequence* *R O E - O E -...- E - T *for* *[$2^{*r_s}$**X**-*$B_{Con}$]$_{IV}$, **X** ∈ ℕ
- *The iterationseries from the **Reduced Value*** $N_{j=s+*r_s+1} \rightarrow ... \rightarrow {}^{\dagger}1_{j=s_d+d}$ *to **The Loop***
- *The* (*Unique **Total***) *Parity-sequence from the **Reduced Value*** *$N_{j=s+*r_s}$ *to **The Loop*** †$1_j$

Also included is an **Algorithm** to Test B-values. Here rules from the **3N-1 Problem** are used. This **Algorithm** is able to find the **reducing* **IV**-class associated with the Test-B and also print the **reducing *parity-sequence* for *[$2^{*r_s}$**X**-*$B_{Con}$]$_{IV}$. All odd B *except* {1,5,17} are *reduced* at **some** point, so **all** B-values ***eventually*** gets caught in a B-loop (see the **Discussion** for more details).

**Note**

The keen and/or curious reader is highly encouraged to copy the *code* on the final pages of **Appendix Omega** and test different odd $N_0$ and $B_0$. **NB**! *On a normal laptop the code/* ***The Algorithm*** *is accurate for* N < $10^{13}$ (*but **would** have* "PERFECT *accuracy" on a* "PERFECT *computer*").

**Examples (from *The Algorithm*) of high P versus high B in the *Converging IV-classes:**

| | |
|---|---|
| The Test-value 999 N in binary: 1111100111<br>is an Element in the *reducing IV-class *[*AX-*B] where:<br>*A = 2^r = 1024 *A in binary: 10000000000<br>*B = *A-*P = 25 *B in binary: 00000110 01<br>*P = 2^r-*B = 999 *P in binary: 1111100111<br>The InitialValue-class [2^ 10 X - 25 ] | The Test-series for *[*AX-*B] Terminates<br>at the TV-class [AX-B] where:<br>A = 3^s = 729 A in Binary: 1011011001<br>B = A- P = 17 B in Binary: 10001<br>P = 3^s- B = 712 P in Binary: 1011001000<br>The TerminalValue-class [3^ 6 X - 17 ] |
| The Test-value 255 N in binary: 11111111<br>is an Element in the *reducing IV-class *[*AX-*B] where:<br>*A = 2^r = 8192 *A in binary: 10000000000000<br>*B = *A-*P = 7937 *B in binary: 1111100000001<br>*P = 2^r-*B = 255 *P in binary: 0000011111111<br>The InitialValue-class [2^ 13 X - 7937 ] | The Test-series for *[*AX-*B] Terminates<br>at the TV-class [AX-B] where:<br>A = 3^s = 6561 A in Binary: 1100110100001<br>B = A- P = 6356 B in Binary: 1100011010100<br>P = 3^s- B = 205 P in Binary: 11001101<br>The TerminalValue-class [3^ 8 X - 6356 ] |

*Remember that in* **IV**-*classes the **pair*** ($P_{IV}$+$B_{IV}$) = $A_{IV}$ = $\mathbf{2}^{r_{IV}}$ *is observed so when the **focus** is on the **binary strings** for* ***P** & *B *in the *reducing* **IV**-*classes, it is observed that the* "***head***" *for **one** of the values* {***P**,*B} *have a sequence of* **one's** *and the* "***head***" *for the **other** one of the values* {*B,***P**} *have a sequence of* **zero's** *in the **string** and **notice** that the **strings*** "*fits like a zipper*" *because the **last*** **1** *in the* "***tail***" (for $2^0$) *have to **carry** all the way to the* "***head***" $2^{*r_s}$ = *$P_{IV}$ + *$B_{Con}$

# 5. Discussion

*The scope of this last main-section is* (from the previous shown **Figures** and **argumentation**) *to extract all the clues for the* ***explanation*** *of The Governing Dynamics controlling* The Problem.

## ***The Governing Dynamics*** **in The Problem is the** ***result*** **of two** ***bijections***

For $\mathbf{r_a} > 0$ in $(2C_a-1)2^{r_a}$ a Type Even Operation exist: $(2C_a-1)2^{r_a}$ $\rightleftarrows_E$ $(2C_a-1)2^{r_a-1}$
For $\mathbf{r_a} = 0$ in $(2C_a-1)2^{r_a}$ a Type Odd Operation exist: $(2C_a-1)$ $\rightleftarrows_O$ $(2C_b-1)2^{r_b}$, $(2C_a-1) \neq {}^{\dagger}\mathbf{1}_{Loop}$

The **Result** of an Operation is *absolutely determined*, but it **is** possible for two different values N (one odd and one even) to reach the same **Result** e.g. $5_j$ $\rightleftarrows_O$ $16_{j+1}$ and $32_j$ $\rightleftarrows_E$ $16_{j+1}$, so the *rules* in The Problem - taken as a **whole** - is **not** a *bijection* while the *rules* individually **are.**

*As stated in the beginning the scope of the present work is to explain* ***why*** *all* N *are in* **Figure 1a**:

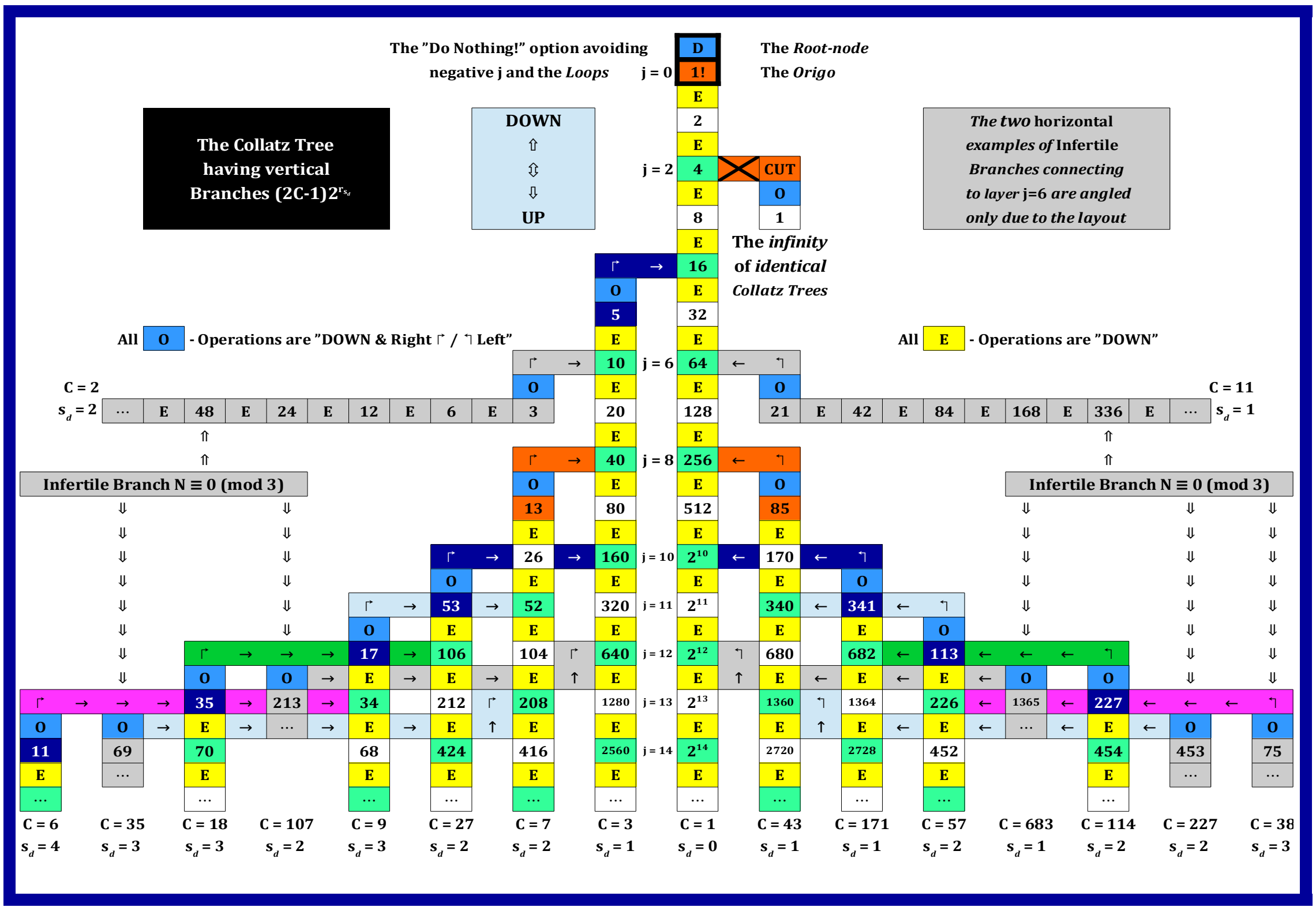

**Figure 1a** which is **equivalent** to **Flowchart 6** and to **The *Reverse* Fibonacci Tree.**
**OBS**! Each Branch $(2Cs_d-1)2^{r_{s_d}}$ is *uniquely defined* by the odd value $(2Cs_d-1)$ or by $Cs_d \in \mathbb{N}$

There exist *one and only one* Branch for which it is true that the †**End Class-series** for all N is:

| | | | | |
|---|---|---|---|---|
| $^{\dagger}[2^d\mathbf{X}\text{-}0_{\mathbf{End}}]_{\mathbf{N}}$ >> $[3^{0_d}\mathbf{X}\text{-}0_{\mathbf{j=0_d+d}}]_{\mathbf{N}}$ | namely The Branch | $N_X \in (2Cs_d-1)2^{r_{s_d}}$ | for | $Cs_d = \mathbf{1}_{0_d}$ & $\mathbf{d} = \mathbf{r}_{s_d}$ |

The Branch which includes ***The Origo* 1!** is the unique **Root Trunk** in ***The Collatz Tree***. All Branches in the **Right side** of the ***Tree*** *connects* to the unique **Root Trunk**. All Branches in the **Left side** of the ***Tree*** *connects* to the unique **Twin Trunk**; $(2Cs_d-1)2^{r_{s_d}}$ for $(23_{1_d}-1)$ = **5** $\rightleftarrows_O$ **16**

The †**End IV-classes** $^{\dagger}[2^d\mathbf{X}\text{-}0_{\mathbf{End}}]_{\mathbf{N}}$ for the N-values in the unique **Root Trunk** $(21_{0_d}-1)2^{r_{0_d}}$:

| | | |
|---|---|---|
| $1_1 \in {}^{\dagger}[2^0\mathbf{X}\text{-}0_{\mathbf{End}}]_1 = [3^{0_d}\mathbf{X}\text{-}0_{0_d+0}]$ | $2_2 \in {}^{*\dagger}[2^1\mathbf{X}\text{-}{}^*0_{\mathbf{End}}]_2 > [3^{0_d}\mathbf{X}\text{-}0_{0_d+1}]$ | $4_4 \in {}^{\dagger}[2^2\mathbf{X}\text{-}0_{\mathbf{End}}]_4 >> [3^{0_d}\mathbf{X}\text{-}0_{0_d+2}]$ |
| $8_8 \in {}^{\dagger}[2^3\mathbf{X}\text{-}0_{\mathbf{End}}]_8 >> [3^{0_d}\mathbf{X}\text{-}0_{0_d+3}]$ | $16_{16} \in {}^{\dagger}[2^4\mathbf{X}\text{-}0_{\mathbf{End}}]_{16} >> [3^{0_d}\mathbf{X}\text{-}0_{0_d+4}]$ | … |

**OBS**! The corresponding †**parity-sequences** for $\mathbf{s_d} = 0$ do *not* contain any type Odd operations.

There exist **one** infinity of Branches for which it is true that the †**End Class-series** for all N is:

| | | | | |
|---|---|---|---|---|
| $^{\dagger}[2^{d}\mathbf{X}\text{-}B_{End}]_{N}$ >> $[3^{1_d}\mathbf{X}\text{-}B_{j=1_d+d}]_{TV}$ | i.e. The Branches | $N_X \in (2C_{s_d}\text{-}1)2^{r_{s_d}}$ | for | $s_d = 1_d$ & **d** = $r_{1_d}$ |

These Branches having $A_{TV} = 3^{1_d}$ are the *Side*-branches on the **Root Trunk** $(21\mathbf{0}_d\text{-}1)2^{r_{0_d}}$ which have *connection-points* $(21\mathbf{0}_d\text{-}1)2^{r_{0_d}} \equiv$ **4** (mod 6) for all even **d** = $r_{0_d}$. †**End IV**-classes $^{\dagger}[2^{d}\mathbf{X}\text{-}B_{End}]_{N}$:

| | | |
|---|---|---|
| $1_1 \in {}^{*\dagger}[2^{2}\mathbf{X}\text{-}{}^{*}3_{Loop}]_1$ ≥ $[3^{1_d}\mathbf{X}\text{-}21_d\text{+}2]$ | $5_5 \in {}^{\dagger}[2^{4}\mathbf{X}\text{-}11_{End}]_5$ >> $[3^{1_d}\mathbf{X}\text{-}21_d\text{+}4]$ | $21_{21} \in {}^{\dagger}[2^{6}\mathbf{X}\text{-}43_{End}]_{21}$ >> $[3^{1_d}\mathbf{X}\text{-}21_d\text{+}6]$ |
| $85_{85} \in {}^{\dagger}[2^{8}\mathbf{X}\text{-}171_{End}]_{85}$ >> $[3^{1_d}\mathbf{X}\text{-}21_d\text{+}8]$ | $341_{341} \in {}^{\dagger}[2^{10}\mathbf{X}\text{-}683_{End}]_{341}$ >> $[3^{1_d}\mathbf{X}\text{-}21_d\text{+}10]$ | ··· *connecting* to **1** i.e. to $(21\mathbf{0}_d\text{-}1)2^{r_{0_d}}$ |

The infinity of odd values in *Side*-branches having $s_d$ = **1** on the *Parent*-branch $(21\mathbf{0}_d\text{-}1)2^{r_{0_d}}$

As any three *consecutive Side*-branches on **Fertile** Branches consists of a **triple** of ***one each*** of the three ***types*** of Branches { **1** , **2** , 3 } **two** out of **three** *Side*-branches on $(21\mathbf{0}_d\text{-}1)2^{r_{0_d}}$ are **Fertile** and in the *Side*-branches on *these* **Fertile** Branches all **End**-series with $s_d$ = **2** are found.

| | | |
|---|---|---|
| $3_3 \in {}^{\dagger}[2^{5}\mathbf{X}\text{-}29_{End}]_3$ >> $[3^{2_d}\mathbf{X}\text{-}82_d\text{+}5]$ | $13_{13} \in {}^{\dagger}[2^{7}\mathbf{X}\text{-}115_{End}]_{13}$ >> $[3^{2_d}\mathbf{X}\text{-}82_d\text{+}7]$ | $53_{53} \in {}^{\dagger}[2^{9}\mathbf{X}\text{-}459_{End}]_{53}$ >> $[3^{2_d}\mathbf{X}\text{-}82_d\text{+}9]$ |
| $213_{213} \in {}^{\dagger}[2^{11}\mathbf{X}\text{-}1835_{End}]_{213}$ >> $[3^{2_d}\mathbf{X}\text{-}82_d\text{+}11]$ | $853_{853} \in {}^{\dagger}[2^{13}\mathbf{X}\text{-}7339_{End}]_{853}$ >> $[3^{2_d}\mathbf{X}\text{-}82_d\text{+}13]$ | ··· *connecting* to **5** i.e. to $(23\mathbf{1}_d\text{-}1)2^{r_{1_d}}$ |

The infinity of *Side*-branches on the *Parent*-branch $(23\mathbf{1}_d\text{-}1)2^{r_{1_d}}$ for $(2C_{s_d}\text{-}1)$ = **5** ⇄$_O$ **16**

| | | |
|---|---|---|
| $113 \in {}^{\dagger}[2^{10}\mathbf{X}\text{-}911_{End}]_{113}$ >> $[3^{2_d}\mathbf{X}\text{-}82_d\text{+}10]$ | $453 \in {}^{\dagger}[2^{12}\mathbf{X}\text{-}3643_{End}]_{453}$ >> $[3^{2_d}\mathbf{X}\text{-}82_d\text{+}12]$ | $1813 \in {}^{\dagger}[2^{14}\mathbf{X}\text{-}14571_{End}]_{1813}$ >> $[3^{2_d}\mathbf{X}\text{-}82_d\text{+}14]$ |
| $7253 \in {}^{\dagger}[2^{16}\mathbf{X}\text{-}58283_{End}]_{7253}$ >> $[3^{2_d}\mathbf{X}\text{-}82_d\text{+}16]$ | $29013 \in {}^{\dagger}[2^{18}\mathbf{X}\text{-}233131_{End}]_{29013}$ >> $[3^{2_d}\mathbf{X}\text{-}82_d\text{+}18]$ | ··· *connecting* to **85** i.e. to $(243\mathbf{1}_d\text{-}1)2^{r_{1_d}}$ |

The infinity of *Side*-branches on the *Parent*-branch $(243\mathbf{1}_d\text{-}1)2^{r_{1_d}}$ for $(2C_{s_d}\text{-}1)$ = **85** ⇄$_O$ **256**

Remember that for an **IV**-class N ≡ **P** (mod $2^{r}$) ⇔ $[2^{r}\mathbf{X}\text{-}B_0]_{IV}$ which is a *Side*-branch on a **Fertile** Branch, the **IV**-class for the **next** *Side*-branch is N ≡ **4P+1** (mod $2^{r+2}$) ⇔ N ∈ $[2^{r+2}\mathbf{X}\text{-}(4B_0\text{-}1)_0]_{IV}$

**Comment**

The ***infinity*** of *Side-branches*$_1$ connecting to the **Root Trunk** contains all values $N_N \in \mathbb{N}$ for which it is true that the †**End-Class-series terminates** at $[3^{1_d}\mathbf{X}\text{-}2_{j=1_d+d}]_{TV}$ i.e. for $s_d$ = $\mathbf{1}_d$. The **infinity** of ***infinities*** of *Side-branches*$_2$ connecting to the **Fertile** *Side-branches*$_1$ contains all values $N_N \in \mathbb{N}$ for which it is true that †**End Class**-series ends at $[3^{2_d}\mathbf{X}\text{-}8_{j=2_d+d}]_{TV}$ i.e. for $s_d$ = $\mathbf{2}_d$. The ***infinity*** of **infinities** of ***infinities*** of *Side-branches*$_3$ connecting to **Fertile** *Side-branches*$_2$ contains all values $N_N \in \mathbb{N}$ for which the †**End Class**-series ends at $[3^{3_d}\mathbf{X}\text{-}26_{j=3_d+d}]_{TV}$, $s_d$ = $\mathbf{3}_d$ **etc.**

*An implication of this is, that in the below **general*** **Class-series** *it is possible to go* DOWN *from any* N = $(2C_{s_d}\text{-}1)2^{r_{s_d}=0}$ & *reach* †$\mathbf{1}\mathbf{0}_d$. *Vice versa; it is possible to go* UP *from* †**1** & *reach any* N > †**1**.

| | | | | | | | | |
|---|---|---|---|---|---|---|---|---|
| $(2C_{s_d}\text{-}1)2^{r_{s_d}=0}$ | ⇄$_O$ | $(2C_{s_d}\text{-}1\text{-}1)2^{r_{s_d}\text{-}1}$ | ⇄$_E$ | $(2C_{s_d}\text{-}1\text{-}1)2^{r_{s_d}\text{-}1\text{-}1}$ | ⇄$_E$ ··· ⇄$_E$ | $(2C_{s_d}\text{-}1\text{-}1)2^{r_{s_d}\text{-}1=1}$ | ⇄$_E$ | |
| $(2C_{s_d}\text{-}1\text{-}1)2^{r_{s_d}\text{-}1=0}$ | ⇄$_O$ | $(2C_{s_d}\text{-}2\text{-}1)2^{r_{s_d}\text{-}2}$ | ⇄$_E$ | $(2C_{s_d}\text{-}2\text{-}1)2^{r_{s_d}\text{-}2\text{-}1}$ | ⇄$_E$ ··· ⇄$_E$ | $(2C_{s_d}\text{-}2\text{-}1)2^{r_{s_d}\text{-}2=1}$ | ⇄$_E$ | |
| $(2C_{s_d}\text{-}2\text{-}1)2^{r_{s_d}\text{-}2=0}$ | ⇄$_O$ | $(2C_{s_d}\text{-}3\text{-}1)2^{r_{s_d}\text{-}3}$ | ⇄$_E$ | $(2C_{s_d}\text{-}3\text{-}1)2^{r_{s_d}\text{-}3\text{-}1}$ | ⇄$_E$ ··· ⇄$_E$ | $(2C_{s_d}\text{-}3\text{-}1)2^{r_{s_d}\text{-}3=1}$ | ⇄$_E$ | |
| ··· ↓DOWN↓ | ⇄$_O$ | ··· ↑UP↑ | ⇄$_E$ | ··· | ⇄$_E$ ··· ⇄$_E$ | ··· | ⇄$_E$ | |
| $(2C_{s_d}\text{-}t\text{-}1)2^{r_{s_d}\text{-}t=0}$ | = | $(2C_{s_d}\text{-}t\text{-}1)$ = †$\mathbf{1}\mathbf{0}_d$ | | **t** = # **O** Operation | | and $\sum r_{s_d}$ = # **E** Operation | | |

**The *general* Class-series for Branches** (**NB**! Please notice the *similarity* with **Flowchart 6**)

## Some details about End Initial Value classes

It has now been established, that there exist an *infinity* of $^{*}{}_{u}[2^{*r_s}\mathbf{X}\text{-}^{*}B_{Con}]_q$ for $\mathbf{u} \in \mathbb{N}$. For $N_N > \mathbf{2}$ it is true for all $^{\dagger}$**End IV**-classes $^{\dagger}[2^{d}\mathbf{X}\text{-}(2^{d}\text{-}N_N)_{End}]_N$ >> $[3^{s_d}\mathbf{X}\text{-}(3^{s_d}\text{-}^{\dagger}\mathbf{1})_{j=s_d+d}]_{TV}$ that they are ***proper subsets*** of a **Converging* **IV**-class $N_X \in {}^{*}{}_{u}[2^{*r_s}\mathbf{X}\text{-}^{*}B_{Con}]_q$ i.e. $N_N = N_X$ & $(2^{d}\text{-}N_N)_{End} \equiv {}^{*}\mathbf{B}_{Con}$ (mod $2^{*r_s}$)

E.g. the **Converging* **IV**-class $^{*}{}_{3}[2^{*4_2}\mathbf{X}\text{-}^{*}13_{Con}]_1 \Leftrightarrow {}^{*}N_X \equiv \mathbf{3}$ (mod $2^{*4_2}$) $\Leftrightarrow {}^{*}N_X \in \{3_1, 19_2, 35_3, 51_4, \ldots\}$

| | | | |
|---|---|---|---|
| $3 \in {}^{\dagger}[2^{5}\mathbf{X}\text{-}(2^{5}\text{-}3)]_3$ | $19 \in {}^{\dagger}[2^{14}\mathbf{X}\text{-}(2^{14}\text{-}19)]_{19}$ | $35 \in {}^{\dagger}[2^{10}\mathbf{X}\text{-}(2^{10}\text{-}35)]_{35}$ | $51 \in {}^{\dagger}[2^{17}\mathbf{X}\text{-}(2^{17}\text{-}51)]_{51}$ |
| $67 \in {}^{\dagger}[2^{19}\mathbf{X}\text{-}(2^{19}\text{-}67)]_{67}$ | $83 \in {}^{\dagger}[2^{70}\mathbf{X}\text{-}(2^{70}\text{-}83)]_{83}$ | $99 \in {}^{\dagger}[2^{18}\mathbf{X}\text{-}(2^{18}\text{-}99)]_{99}$ | $115 \in {}^{\dagger}[2^{23}\mathbf{X}\text{-}(2^{23}\text{-}115)]_{115}$ |
| $131 \in {}^{\dagger}[2^{20}\mathbf{X}\text{-}(2^{20}\text{-}131)]_{131}$ | $147 \in {}^{\dagger}[2^{74}\mathbf{X}\text{-}(2^{74}\text{-}147)]_{147}$ | $163 \in {}^{\dagger}[2^{17}\mathbf{X}\text{-}(2^{17}\text{-}163)]_{163}$ | … |

| j | 0 | 1 | 2 | 3 | 4 | 5 | 6 | $N_X$ | 7 | 8 | 9 | 10 | 11 | 12 | 13 | 14 | 15 | 16 | 17 | 18 | 19 | 20 | 21 | 22 | 23 | 24 | 25 | 26 | 27 | 28 | … |
|---|---|---|---|---|---|---|---|---|---|---|---|---|---|---|---|---|---|---|---|---|---|---|---|---|---|---|---|---|---|---|---|
| * | R | O | E | O | E | E | E | 3 | E | D | $(s_d,d) = (2,5), j_{End} = 7$ | | | | | | | | | | | | | | | | | | | | |
| *N ≡ 3 (mod 16) | | | | | | | | 19 | O | E | O | E | E | O | E | E | E | O | E | E | E | E | D | $(s_d,d) = (6,14), j_{End} = 20$ | | | | | | | |
| $(s_d,d) =$ (3,10), $j_{End} =$ 13 | | | | | | | | 35 | E | E | O | E | E | E | E | D | | | | | | | | | | | | | | | |
| $(s_d,d) =$ (7,17), $j_{End} =$ 24 | | | | | | | | 51 | O | E | E | E | O | E | O | E | E | O | E | E | E | O | E | E | E | E | D | | | | |
| $(s_d,d) =$ (8,19), $j_{End} =$ 27 | | | | | | | | 67 | E | O | E | O | E | E | E | O | E | O | E | E | O | E | E | E | O | E | E | E | E | D | |
| $(s_d,d) =$ (40,70), $j_{End} =$ 110 | | | | | | | | 83 | O | E | … | | | | | | | | | | | | | | | | | | | | |
| $(s_d,d) =$ (7,18), $j_{End} =$ 25 | | | | | | | | 99 | E | E | E | O | E | O | E | O | E | E | O | E | E | E | O | E | E | E | E | D | | | |
| $(s_d,d) =$ (10,23), $j_{End} =$ 33 | | | | | | | | 115 | O | E | … | | | | | | | | | | | | | | | | | | | | |
| $(s_d,d) =$ (8,20), $j_{End} =$ 28 | | | | | | | | 131 | E | … | | | | | | | | | | | | | | | | | | | | | |
| $(s_d,d)$ , $j_{End}$ | | | | | | | | … | … | | | | | | | | | | | | | | | | | | | | | | |

| j | 0 | 1 | 2 | 3 | 4 | 5 | 6 |
|---|---|---|---|---|---|---|---|
| * | R | O | E | O | E | E | E |

| | 6 | 7 | 8 | 9 | 10 | 11 | 12 | 13 | $N_X$ | 14 | 15 | 16 | 17 | 18 | 19 | 20 | 21 | 22 | 23 | 24 | 25 | 26 | 27 | 28 | 29 | 30 | 31 | 32 | 33 | 34 | 35 | … |
|---|---|---|---|---|---|---|---|---|---|---|---|---|---|---|---|---|---|---|---|---|---|---|---|---|---|---|---|---|---|---|---|---|
| | … | E | E | O | E | E | E | E | 35 | D | $(s_d,d) = (3,10), j_{End} = 13$ | | | | | N = 1059 ⇒ | | | $(s_d,d) = (8,23), j_{End} = 31$ | | | | | | | | | | | | | |
| N ≡ 35 (mod 2^10) | | | | | | | | | 1059 | E | E | O | E | O | E | O | E | E | O | E | E | E | O | E | E | E | E | D | | | | |
| $(s_d,d) =$ (44,81), $j_{End} =$ 125 | | | | | | | | | 2083 | O | E | O | E | O | E | E | E | O | E | O | E | O | E | O | E | E | O | E | … | D | | |
| $(s_d,d) =$ (43,80), $j_{End} =$ 123 | | | | | | | | | 3107 | E | O | E | E | O | E | O | E | O | E | O | E | O | E | E | O | E | … | D | | | | |
| $(s_d,d) =$ (44,82), $j_{End} =$ 126 | | | | | | | | | 4131 | O | E | E | E | O | E | E | O | E | O | E | O | E | O | E | O | E | E | O | E | … | D | |
| $(s_d,d) =$ (6,22), $j_{End} =$ 28 | | | | | | | | | 5155 | E | E | E | O | E | E | O | E | E | E | O | E | E | E | E | D | | | | | | | |
| $(s_d,d) =$ (9,27), $j_{End} =$ 36 | | | | | | | | | 6179 | O | E | … | | | | | | | | | | | | | | | | | | | | |
| $(s_d,d) =$ (41,78), $j_{End} =$ 119 | | | | | | | | | 7203 | E | … | | | | | | | | | | | | | | | | | | | | | |
| $(s_d,d)$ , $j_{End}$ | | | | | | | | | … | … | | | | | | | | | | | | | | | | | | | | | | |

Theses examples are an illustration of the fact, that a *countable infinity* of **Convergent End** $^{\dagger}$**parity-sequences** exist, and that for N > 2 all **End** $^{\dagger}$**parity-sequences** are extensions of **Converging *parity-sequences* i.e. that all $^{\dagger}$**End IV**-classes are ***proper subsets*** of **Converging* **IV**-classes $^{*}{}_{u}[2^{*r_s}\mathbf{X}\text{-}^{*}B_{Con}]_q$.

Notice that it is possible to focus on any **End** $^{\dagger}$**parity-sequence** and find an *infinity* of **extensions** for all $\mathbf{X} > 1$ for the $N_X \in {}^{\dagger}[2^{d}\mathbf{X}\text{-}(2^{d}\text{-}N_N)_{End}]_N$

**End IV**-classes $^{\dagger}[2^{d}\mathbf{X}\text{-}(2^{d}\text{-}N_N)_{End}]_N$ are valid for the **lowest** value $\mathbf{P}_{IV} = N_N$ i.e. all the $N_X$ for $\mathbf{X} > 1$ in an **End IV**-class **also** belongs to **its own** class where **they** are the lowest $\mathbf{P}_{IV} = N_N$ **&** index$_N$ **N.**

**The *three* B-loops for B ∈ {1,5,17}** **NB**! See Wikipedia[I] *Collatz Conjecture* for more info.

| The B-loop for $B = 1_{j+2}$ | | | |
|---|---|---|---|
| $[2^{r-0}3^{s+0}\mathbf{X} - 1_0]$ | $\rightleftarrows_O$ | $[2^{r-0}3^{s+1}\mathbf{X} - 2_1]$ | $\rightleftarrows_E$ |
| $[2^{r-1}3^{s+1}\mathbf{X} - 1_2]$ | $\rightleftarrows_O$ | $[2^{r-1}3^{s+2}\mathbf{X} - 2_3]$ | $\rightleftarrows_E$ |
| … | | … | |

| The B-loop for $B = 5_{j+5}$ | | | |
|---|---|---|---|
| $[2^{r-0}3^{s+0}\mathbf{X} - 5_0]$ | $\rightleftarrows_O$ | $[2^{r-0}3^{s+1}\mathbf{X} - 14_1]$ | $\rightleftarrows_E$ |
| $[2^{r-1}3^{s+1}\mathbf{X} - 7_2]$ | $\rightleftarrows_O$ | $[2^{r-1}3^{s+2}\mathbf{X} - 20_3]$ | $\rightleftarrows_E$ |
| | | $[2^{r-2}3^{s+2}\mathbf{X} - 10_4]$ | $\rightleftarrows_E$ |
| $[2^{r-3}3^{s+2}\mathbf{X} - 5_5]$ | $\rightleftarrows_O$ | $[2^{r-3}3^{s+3}\mathbf{X} - 14_6]$ | $\rightleftarrows_E$ |
| … | | … | |

B = **1** ⇒ **(s,r) is (1,1) = (1,*2₁-1)**
B = **5** ⇒ **(s,r) is (2,3) = (2,*4₂-1)**
B = **17** ⇒ **(s,r) is (7,11) = (7,*12₇-1)**

| The B-loop for $B = 17_{j+18}$ | | | |
|---|---|---|---|
| $[2^{r-0}3^{s+0}\mathbf{X} - 17_0]$ | $\rightleftarrows_O$ | $[2^{r-0}3^{s+1}\mathbf{X} - 50_1]$ | $\rightleftarrows_E$ |
| $[2^{r-1}3^{s+1}\mathbf{X} - 25_2]$ | $\rightleftarrows_O$ | $[2^{r-1}3^{s+2}\mathbf{X} - 74_3]$ | $\rightleftarrows_E$ |
| $[2^{r-2}3^{s+2}\mathbf{X} - 37_4]$ | $\rightleftarrows_O$ | $[2^{r-2}3^{s+3}\mathbf{X} - 110_5]$ | $\rightleftarrows_E$ |
| $[2^{r-3}3^{s+3}\mathbf{X} - 55_6]$ | $\rightleftarrows_O$ | $[2^{r-3}3^{s+4}\mathbf{X} - 164_7]$ | $\rightleftarrows_E$ |
| | | $[2^{r-4}3^{s+4}\mathbf{X} - 82_8]$ | $\rightleftarrows_E$ |
| $[2^{r-5}3^{s+4}\mathbf{X} - 41_9]$ | $\rightleftarrows_O$ | $[2^{r-5}3^{s+5}\mathbf{X} - 122_{10}]$ | $\rightleftarrows_E$ |
| $[2^{r-6}3^{s+5}\mathbf{X} - 61_{11}]$ | $\rightleftarrows_O$ | $[2^{r-6}3^{s+6}\mathbf{X} - 182_{12}]$ | $\rightleftarrows_E$ |
| $[2^{r-7}3^{s+6}\mathbf{X} - 91_{13}]$ | $\rightleftarrows_O$ | $[2^{r-7}3^{s+7}\mathbf{X} - 272_{14}]$ | $\rightleftarrows_E$ |
| | | $[2^{r-8}3^{s+7}\mathbf{X} - 136_{15}]$ | $\rightleftarrows_E$ |
| | | $[2^{r-9}3^{s+7}\mathbf{X} - 68_{16}]$ | $\rightleftarrows_E$ |
| | | $[2^{r-10}3^{s+7}\mathbf{X} - 34_{17}]$ | $\rightleftarrows_E$ |
| $[2^{r-11}3^{s+7}\mathbf{X} - 17_{18}]$ | $\rightleftarrows_O$ | $[2^{r-11}3^{s+8}\mathbf{X} - 50_{19}]$ | $\rightleftarrows_E$ |
| … | | … | |

Please notice that it is possible to **draw conclusions** from $B_{*j=s+*r} < B_j$ in all the **Class-series** shown here and in all the **Flowcharts.**

$17_{18} < 25_2 \Rightarrow$
$^{*}{}_{27}[2^{*10_6}\mathbf{X}\text{-}^{*}25_{Con}]_{12}$

$34_{17} < 37_4 \Rightarrow$
$^{*}{}_{15}[2^{*8_5}\mathbf{X}\text{-}^{*}37_{Con}]_7$

$41_9 < 55_6 \Rightarrow$
$^{*}{}_{u=2}[2^{*2_1}\mathbf{X}\text{-}^{*}3_{Con}]_{q=1}$

***etc.***

**Figure 4.** The ***three*** different B-loops for $B_{r_s} \in \{1_{1_1}, 5_{3_2}, 17_{11_7}\}$

## Conclusions from the Flowcharts

If a **Flowchart** 2a and a **Flowchart** 2b is printed for all $\mathbf{X} \in \mathbb{N}$ the individual Flowcharts now contains Value-Branches as each node now have a distinct value $N_X = (2C-1)2^r$, $C \in \mathbb{N}$, $\mathbf{r} \in \mathbb{N}_0$, where $N_X \in [2^r3^s\mathbf{X}\text{-B}]$, $B < 2^r3^s$, $\mathbf{X} \in \mathbb{N}$. Notice that the ***Domain-class*** $^{\Phi}[2^2\mathbf{X}\text{-}1_0]_{IV}$ is found in the first *Side*-branch in **Flowchart** 2b, so it is **guaranteed** that **each** printed type 2 Flowchart for a constant $\mathbf{X} \in \mathbb{N}$ contains a **specific** $N_X$ from $\{3_1,7_2,11_3,15_4,19_5,23_6,27_7, \dots ,N_X, \dots\} = {}^{\Phi}[2^2\mathbf{X}\text{-}1_0]_{IV}$. Remember that it is the ***same*** *infinity* of Flowcharts (with Value-Branches) that is printed, when Flowchart 2a&2b is printed for all $\mathbf{X} \in \mathbb{N}$, Flowchart 3a&3b&3c&3d&3e&3f is printed for all **X**, when 4a&4b&4c&4d is printed or when 5a&5b&5c&5d&5e&5f&5g&5h is printed. When the **modulo** is *doubled* it is possible to print Flowcharts with *extended* parity-sequences.

Remember that; $\quad {}^*{}_1[2^{*1_0}\mathbf{X}\text{-}{}^*0_{Con}]_1 \cup {}^*{}_2[2^{*2_1}\mathbf{X}\text{-}{}^*3_{Con}]_1 \cup {}^{\Phi}[2^2\mathbf{X}\text{-}1_0]_{IV} = {}^{\Phi}[2^0\mathbf{X}\text{-}0_0]_{IV/TV} = \mathbb{N}$

All **subsets** $^{\Delta}[2^r\mathbf{X}\text{-}{}^{\Delta}B_0]_{IV}$, $^{\Delta}B_0 \equiv {}^*\mathbf{0}_{Con} \pmod{2^1}$, $\mathbf{r} > {}^*\mathbf{1}_0$, of ${}^*{}_1[2^{*1_0}\mathbf{X}\text{-}{}^*0_{Con}]_1$ are $^{\Delta}$**Redundant**.

All **subsets** $^{\Delta}[2^r\mathbf{X}\text{-}{}^{\Delta}B_0]_{IV}$, $^{\Delta}B_0 \equiv {}^*\mathbf{3}_{Con} \pmod{2^2}$, $\mathbf{r} > {}^*\mathbf{2}_1$, of ${}^*{}_2[2^{*2_1}\mathbf{X}\text{-}{}^*3_{Con}]_1$ are $^{\Delta}$**Redundant**.

For $\mathbf{r} > 2$ then all ${}^*{}_u[2^{*r_s}\mathbf{X}\text{-}{}^*B_{Con}]_q$ and all $^{\Phi}[2^r\mathbf{X}\text{-}B_0]_{IV}$ are ***subsets*** of the ***Domain-class*** $^{\Phi}[2^2\mathbf{X}\text{-}1_0]_{IV}$.

If the focus is on **Converging* subsets of the ***Domain-class*** it is possible to start with only one **Flowchart** 2b, as the first *Side*-branch on the branch with the odd class [6X-1] is $^{\Phi}[2^2\mathbf{X}\text{-}1_0]_{IV}$. Every time the **modulo** is *doubled* an **IV**-class $[2^r\mathbf{X}\text{-}B_0]_{IV}$ for $B_0 \equiv \mathbf{1} \pmod{2^2}$ is split into two **IV**-classes; $[2^{r+1}\mathbf{X}\text{-}(2^r+B_0)_0]_{IV} \cup [2^{r+1}\mathbf{X}\text{-}B_0]_{IV} = [2^r\mathbf{X}\text{-}B_0]_{IV}$ and the parity-sequences are ***extended***. It is possible to identify the **Converging* subsets ($\mathbf{u} \geq 3$) of the **Domain-class** when it is observed that the **TV**-class $[3^s\mathbf{X}\text{-}B_j]_{TV}$ have smaller **modulo** $3^s$ than the **IV**-class ${}^*{}_u[2^{*r_s}\mathbf{X}\text{-}{}^*B_{Con}]_q$ in which case $2^{*r_s} > 3^s$ & ${}^*B_{Con} > B_{*j=s+*r_s}$ & $(2^{*r_s}\text{-}{}^*B_{Con}) > (3^s\text{-}B_{*j=s+*r_s})$ i.e. $\mathbf{P}_{IV} > \mathbf{P}_{TV}$.

In the **extreme** for $\mathbf{r} \to \infty$ the pile of Flowcharts for **Converging* **IV**-classes tends to infinity, $\mathbf{u} \to \infty$, and the pile of charts containing $^{\Phi}$**Diverging IV**-classes *Covers* less & less of $N \in \mathbb{N}$.

An alternative loop is **not** possible, because there exist only **one** *flow-direction* for any node in any of the *infinity* of printed Flowcharts; the *direction* "**DOWN**" which is **towards** the **Origo 1**!

## Conclusions from the Binary Trees and the analysis of parity-sequences

The most important property of The Domino-Tree is, that for all layers **r** the **Union** of the $2^r$ different **IV**-classes $[2^r\mathbf{X}\text{-}B_0]_{IV}$, $B_0 < 2^r$, is a ***Covering System*** i.e. **all** $N \in \mathbb{N}$ are included in a class. There exist an **Uncountable Infinity** of different ***routes/****parity-sequences* in The Domino-Tree and of these a (comparatively small) *countable infinity* are **reducing* and can be identified with a **Converging* **Class-series** ${}^*{}_u[2^{*r_s}\mathbf{X}\text{-}{}^*B_{Con}]_{IV} > [3^s\mathbf{X}\text{-}B_{j=s+*r_s}]_{TV}$, $\mathbf{u} \in \mathbb{N}$. In the Domino-Tree there also exist a *countable infinity* of $^{\dagger}$**Convergent** ***routes***, which can be identified with **End Class-series** $^{\dagger}[2^d\mathbf{X}\text{-}(2^d\text{-}N_N)_{End}]_N >> [3^{s_d}\mathbf{X}\text{-}(3^{s_d}\text{-}\mathbf{1})_{j=s_d+d}]_{TV}$, $\mathbf{N} \in \mathbb{N}$. From **List 2** it is known, that for any $^{\Phi}$**Diverging** parity-sequence there exist a **reducing* extension having the **same** amount of the increasing **O E -tiles** and also a **reducing* extension having **one** more **O E -tile**. In all $^{\Phi}$**Diverging IV**-classes $N_X \in {}^{\Phi}[2^k\mathbf{X}\text{-}{}^*B_0]_{IV}$, $\mathbf{k} < {}^*\mathbf{r}_s$, there exist a *countable infinity* of **Converging* **IV**-classes and for *all* $N_X$ it is *possible* to find the $^{\dagger}$**Convergent** $^{\dagger}$**End IV**-class $^{\dagger}[2^d\mathbf{X}\text{-}(2^d\text{-}N_N)_{End}]_{N=N_X}$

Please notice in **Resultlist 2** that for all $(\mathbf{s},{}^*\mathbf{r}_s)$, $\mathbf{s} \to \infty$; $|{}^*Con({}^*\mathbf{r}_s)| \,/\, |{}^{\Phi}\mathbf{Div}({}^*\mathbf{r}_s)| > |{}^{\Phi}\mathbf{Div}({}^*\mathbf{r}_s)|/2^{*r_s}$
The **ratio** for $|{}^*Con({}^*\mathbf{r}_s)| \,/\, |{}^{\Phi}\mathbf{Div}({}^*\mathbf{r}_s)|$ is ***more significant*** than the **density** for $^{\Phi}\mathbf{Div}({}^*\mathbf{r}_s)$ i.e. a larger *fraction* of $^{\Phi}\mathbf{Div}({}^*\mathbf{r}_s)$ are **Converging* than the *fraction* of the **Total** which are **Diverging**.

## <u>Observations about *collections* of Class-series</u>

In the analysis leading to the present work all Class-series for all **IV**-classes $[2^r\mathbf{X}\text{-}B_0]_{IV}$, $\mathbf{r} \leq 13$ and $B_0 \equiv \mathbf{1} \pmod 4$ were constructed. For constant **r** the Class-series including all $B_0 < 2^r$ is a *complete collection* as it contains **all** the *possible* **subsets** of **The Domain-class** $^{\Phi}[2^2\mathbf{X}\text{-}1_0]_{IV}$. Besides generating all the **reducing* **IV**-classes listed in **Appendix Alpha** the *collections* for constant **r** reveals, that (some) Class-series *merge* i.e. have *identical* **T**erminal **V**alue-classes.

The *first* example is found in the *collection* for **r** = 9, where it is observed that the $^{\Delta}$Class-series for $^{\Delta}[2^9\mathbf{X}\text{-}121_0]_{IV}$ & $^{\Delta}[2^9\mathbf{X}\text{-}125_0]_{IV}$ both terminates at the **TV**-class $[3^4\mathbf{X}\text{-}19_{13}]_{TV}$. What this means is, that the iterationseries for 387 & 391 reach the **same** N after 13 Operations, and it implies that for all **r** ≥ 9 the $^{\Delta}$Class-series for $^{\Delta}[2^r\mathbf{X}\text{-}121_0]_{IV}$ & $^{\Delta}[2^r\mathbf{X}\text{-}125_0]_{IV}$ have identical $[3^s\mathbf{X}\text{-}B_{j=s+r}]_{TV}$:

$387_0$ ↗ $1162_1$ ↘ $581_2$ ↗ $1744_3$ ↘ $872_4$ ↘ $436_5$ ↘ $*\mathbf{218}_6$ ↘ $109_7$ ↗ $328_8$ ↘ $164_9$ ↘ $82_{10}$ ↘ $41_{11}$ ↗ $124_{12}$ ↘ $\mathbf{62}_{13}$

$391_0$ ↗ $1174_1$ ↘ $587_2$ ↗ $1762_3$ ↘ $881_4$ ↗ $2644_5$ ↘ $1322_6$ ↘ $661_7$ ↗ $1984_8$ ↘ $992_9$ ↘ $496_{10}$ ↘ $*\mathbf{248}_{11}$ ↘ $124_{12}$ ↘ $\mathbf{62}_{13}$

More examples are found in the *collections* for **r** > 9:

$^{\Delta}[2^{10}\mathbf{X}\text{-}469_0]_{IV} \rightarrow\rightarrow [3^4\mathbf{X}\text{-}37_{14}]_{TV}$ & $^{\Delta}[2^{10}\mathbf{X}\text{-}477_0]_{IV} \rightarrow\rightarrow [3^4\mathbf{X}\text{-}37_{14}]_{TV}$

$^{\Delta}[2^{10}\mathbf{X}\text{-}1017_0]_{IV} \rightarrow\rightarrow [3^5\mathbf{X}\text{-}241_{15}]_{TV}$ & $^{\Delta}[2^{10}\mathbf{X}\text{-}1021_0]_{IV} \rightarrow\rightarrow [3^5\mathbf{X}\text{-}241_{15}]_{TV}$

$^{\Delta}[2^{11}\mathbf{X}\text{-}137_0]_{IV} \rightarrow\rightarrow [3^5\mathbf{X}\text{-}16_{16}]_{TV}$ & $^{\Delta}[2^{11}\mathbf{X}\text{-}141_0]_{IV} \rightarrow\rightarrow [3^5\mathbf{X}\text{-}16_{16}]_{TV}$

$^{\Delta}[2^{11}\mathbf{X}\text{-}213_0]_{IV} \rightarrow\rightarrow [3^5\mathbf{X}\text{-}25_{16}]_{TV}$ & $^{\Delta}[2^{11}\mathbf{X}\text{-}221_0]_{IV} \rightarrow\rightarrow [3^5\mathbf{X}\text{-}25_{16}]_{TV}$

$^{\Delta}[2^{11}\mathbf{X}\text{-}249_0]_{IV} \rightarrow\rightarrow [3^6\mathbf{X}\text{-}88_{17}]_{TV}$ & $^{\Delta}[2^{11}\mathbf{X}\text{-}253_0]_{IV} \rightarrow\rightarrow [3^6\mathbf{X}\text{-}88_{17}]_{TV}$

… & …

In the *collections* for higher **r** it is observed that *more and more* $^{\Delta}$Class-series *merge*.

The explanation for the *merging* $^{\Delta}$Class-series is, that the **TV**-classes "run out" of $B_{TV} < A_{TV}$.

| *r | 4 | 5 | 7 | 8 | 10 | 12 | 13 | 15 | 16 | 18 | 20 | … |
|---|---|---|---|---|---|---|---|---|---|---|---|---|
| s | **2** | **3** | **4** | **5** | **6** | **7** | **8** | **9** | **10** | **11** | **12** | … |
| **r** | | | | | | | | | | | | |
| **2** | 1 | | | | | | | | | | | |
| **3** | 1 | 1 | | | | | | | | | | |
| **4** | 1 | 2 | 1 | | | | | | | | | |
| **5** | 1 | 3 | 3 | 1 | | | | | | | | |
| **6** | 1 | 4 | 6 | 4 | 1 | | | | | | | |
| **7** | 1 | 5 | 10 | 10 | 5 | 1 | | | | | | |
| **8** | 1 | 6 | 15 | 20 | 15 | 6 | 1 | | | | | |
| **9** | 1 | 7 | 21 | 35 | 35 | 21 | 7 | 1 | | | | |
| **10** | 1 | 8 | 28 | 56 | 70 | 56 | 28 | 8 | 1 | | | |
| **11** | 1 | 9 | 36 | 84 | 126 | 126 | 84 | 36 | 9 | 1 | | |
| **12** | 1 | 10 | 45 | 120 | 210 | 252 | 210 | 120 | 45 | 10 | 1 | |
| **13** | 1 | 11 | 55 | 165 | 330 | 462 | 462 | 330 | 165 | 55 | 11 | … |
| **14** | 1 | 12 | 66 | 220 | 495 | 792 | 924 | 792 | 495 | 220 | 66 | … |
| **15** | 1 | 13 | 78 | 286 | 715 | 1287 | 1716 | 1716 | 1287 | 715 | 286 | … |
| **16** | 1 | 14 | 91 | 364 | 1001 | 2002 | 3003 | 3432 | 3003 | 2002 | 1001 | … |
| **17** | 1 | 15 | 105 | 455 | 1365 | 3003 | 5005 | 6435 | 6435 | 5005 | 3003 | … |
| **18** | 1 | 16 | 120 | 560 | 1820 | 4368 | 8008 | 11440 | 12870 | 11440 | 8008 | … |
| **19** | 1 | 17 | 136 | 680 | 2380 | 6188 | 12376 | 19448 | 24310 | 24310 | 19448 | … |
| **20** | 1 | 18 | 153 | 816 | 3060 | 8568 | 18564 | 31824 | 43758 | 48620 | 43758 | … |
| … | … | … | … | … | … | … | … | … | … | … | … | … |

**Figure 3b** for the **Domain-class**

After a type Odd operation the new value is $N_{j+1}=3N_j+1$ and it is observed, that $N_{j+1} \equiv \mathbf{1} \pmod 3$ and also that $N_{j+1} \equiv \mathbf{1} \pmod{N_j}$. In a Class-series it is observed that $[2^r3^s\mathbf{X}\text{-}B_j]$ $\rightleftarrows_{\mathbf{o}}$ $[2^r3^{s+1}\mathbf{X}\text{-}(3B_j\text{-}1)_{j+1}]$, $(3B_j\text{-}1)_{j+1} \equiv \mathbf{2} \pmod 3$ so after the first type Odd operation all the B-values are $B_j \equiv \mathbf{1} \pmod 3$ or $B_j \equiv \mathbf{2} \pmod 3$ i.e. only two of three B-values lower than $A_{TV} = 3^{s_{TV}}$ are *possible* in the $^{\Delta}$Class-series for $^{\Delta}$**Redundant Combinations** (**s**,**r**), **r** > *$\mathbf{r}_s$.

Notice in **Figure 3b** that the coefficients to the left of the marked Totals are **Redundant**, as the combinations (**s**,**r**) have *fewer* type Odd operations than *necessary* for a **reducing combination* (**s**,*$\mathbf{r}_s$) and the coefficients below the marked Totals are also **Redundant**, as the combinations have *more* type Even operations than *necessary* for a **reducing combination* (**s**,*$\mathbf{r}_s$) e.g. for the combination (**s**,**r**) = (5,15) there exist 286 Class-series $^{\Delta}[2^{15}\mathbf{X}\text{-}B_0]_{IV} \rightarrow\rightarrow [3^5\mathbf{X}\text{-}B_{20}]_{TV}$ where all the $B_{20}$ are *lower* than $A_{TV} = 3^5 = 243$ so ***some*** Class-series for *different* $B_0$ **do** have *identical* $B_{20}$.

Also observe that for growing **r** an increasing fraction of the $^{\Phi}$**Diverging Combinations** (to the right of the marked Totals) can be identified with $^{\Delta}$**Redundant** *Class-series*. For **r** = 5 all four $^{\Phi}$**Diverging Combinations** represents $^{\Phi}$**Diverging** Class-series, but for **r** = 7 three of the 16 are $^{\Delta}$**Redundant** as only 13 $^{\Phi}$**Diverging** Class-series exist. For **r** = 8 three of the 15+6+1 = 22 are $^{\Delta}$**Redundant** as only 19 $^{\Phi}$**Diverging** Class-series exist. For **r** = 10 twenty-nine of the 93 are $^{\Delta}$**Redundant** as only 64 $^{\Phi}$**Diverging** Class-series exist etc. For increasing **r** the fraction of the $^{\Phi}$**Diverging Combinations** which can be identified with $^{\Delta}$**Redundant** Class-series *grows*. **OBS**! $^{\Phi}$**Diverging** *Class-series* ***do*** *exist for all* **modulo** $A_{IV} = 2^r$, $\mathbf{r} \in \mathbb{N}$, *but* $|^{\Phi}\mathbf{Div}(\mathbf{r})|/2^r \rightarrow 0$, $\mathbf{r} \rightarrow \infty$.

## Observations about well-ordered sets

All **IV**-classes are **well-ordered** in the sense that they contain a *lowest* element **P** = $(2^r\text{-}B_0)$, so the *Union* of **Converging* **IV**-classes $*_u[2^{*r_s}\mathbf{X}\text{-}*B_{Con}]_{q\le q_{Max}}$ for a constant **reducing* combination (**s**,***r**$_s$) also contain a *lowest* element $\mathbf{P}*_{r_s}$ which is found in the **IV**-class for **q** = 1 as the classes are arranged according to size of $\mathbf{P_1} < \mathbf{P_2} < \ldots < \mathbf{P}_{q_{Max}}$ so $\mathbf{P_1}$ is the *representative* $\mathbf{P}*_{r_s}$ for the **Group**. There exist a *countable infinity* of $\mathbf{P}*_{r_s} \in \{2*1_0, 1*2_1, 3*4_2, 11*5_3, 7*7_4, 39*8_5, 287*10_6, \ldots\}$ for **s** → ∞.

For a constant **reducing* combination (**s**,***r**$_s$) the **Union** of $^\Phi$**Diverging IV**-classes $^\Phi[2^{*r_s}\mathbf{X}\text{-}^\Phi B_0]_{IV}$ also contains a *lowest* element i.e. the lowest N which is **not yet** included in a **Converging* **IV**-class having a **modulo** lower than or equal to the current **modulo** $2^{*r_s}$. As the **modulo** grows more and more of these **Threshold-values** $N_T \in \{1*2_1, 3*4_2, 7*7_4, 27*59_{37}, 703*81_{51}, \ldots\}$ are found to belong to the shown **set** of **Group**-*representatives* $\mathbf{P}*_{r_s}$ and in actual fact the **set** of $N_T$ is **also** a *countable infinity*, as for any **reducing* combination a *lowest* N exist in the **Union** of $^\Phi$**Div**(***r**$_s$).
**OBS**! **Appendix Omega** contains more data about the Class-series for **Threshold-values** $N_T$.

## The *lowest* value $P_{IV}$ in a $^\Phi$Diverging IV-class can *not* be unlucky forever

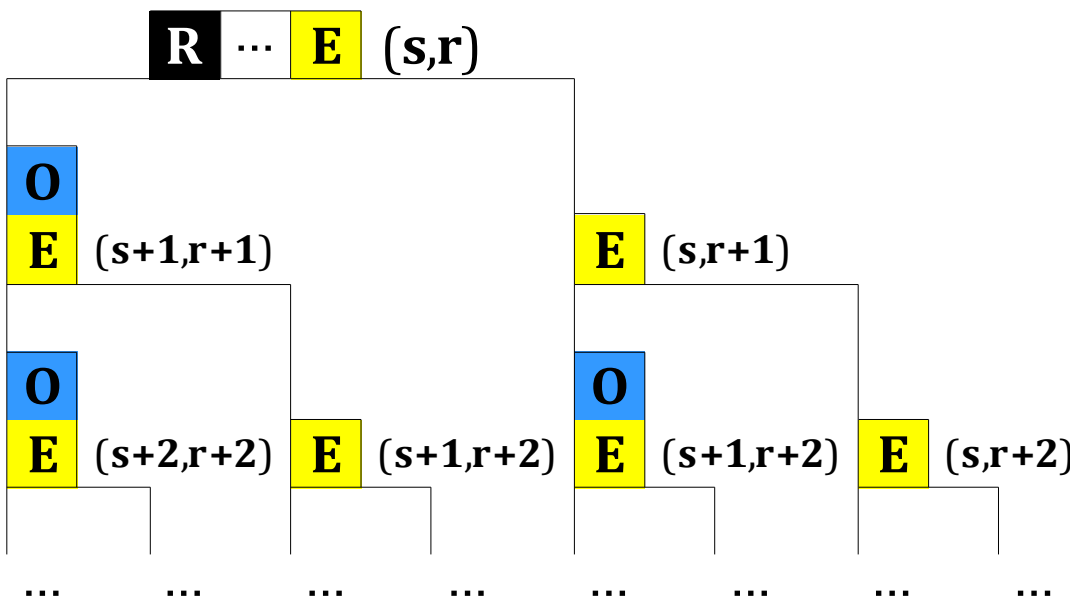

Every time we go **two layers UP** in The Domino-Tree the four lowest $N_X$ in any **IV**-class $[2^r\mathbf{X}\text{-}B_0]_{IV}$ will be the lowest value **P** in one of the *Child*-**IV**-classes $[2^r\mathbf{X}\text{-}B_0]_{IV} = [2^{r+2}\mathbf{X}\text{-}B_0]_{IV} \cup [2^{r+2}\mathbf{X}\text{-}(2^r+B_0)_0]_{IV} \cup [2^{r+2}\mathbf{X}\text{-}(2^{r+1}+B_0)_0]_{IV} \cup [2^{r+2}\mathbf{X}\text{-}(2^{r+1}+2^r+B_0)_0]_{IV}$ and three of the *parity-sequences* are **more** *reducing* with the **extensions** {**E**-**E**, **E**-**O E**, **O E**-**E**} and only one of the four is **less** *reducing* {**O E**-**O E**}.

**Figure 5.** Only one of the four $N_{X\le 4}$ in $[2^r\mathbf{X}\text{-}B_0]_{IV}$, $A_{IV} = 2^r$ can be unlucky when $A_{IV}$ becomes $2^{r+2}$.

Any $^\Phi$**Diverging** Class-series $^\Phi[2^k\mathbf{X}\text{-}^\Phi B_0]_{IV} \to\to [3^s\mathbf{X}\text{-}B_{j=s+k}]_{TV}$ have a $^\Phi$**Diverging** $^\Phi$parity-sequence consisting of **s** **O E**-**tiles** & (**k**-**s**) **E**-**tiles** and the *lowest* value in the **IV**-class is $P_{IV} = (2^k\text{-}^\Phi B_0)$. The first layer **r** in The Domino-Tree containing a **reducing* **subset** of $^\Phi[2^k\mathbf{X}\text{-}^\Phi B_0]_{IV}$ is ***r**$_s$ = ⌈**sZ**⌉ where the **IV**-class with a **parity-sequence* extended only with the decreasing **E**-**tiles** is found and as $2^{(*r_s-k)}$ **subsets** of $^\Phi[2^k\mathbf{X}\text{-}^\Phi B_0]_{IV}$ exist in layer ***r**$_s$ one of the lowest $N_X$ for $\mathbf{X} \le 2^{(*r_s-k)}$ is found as the **P** = $(2^{*r_s}\text{-}*B_{Con})$ in the **Converging* **IV**-class $*[2^{*r_s}\mathbf{X}\text{-}*B_{Con}]_{IV}$ where $*B_{Con} \equiv {}^\Phi B_0 \pmod{2^k}$. In the next **reducing* layer ***r**$_{s+1}$ there exist (***r**$_s$-**k**) **reducing* **subsets** of $^\Phi[2^k\mathbf{X}\text{-}^\Phi B_0]_{IV}$ and in all the **reducing* layers ***r** for higher **modulo** a fast-growing number of **reducing* **subsets** exist. *It is inevitable that the lowest value in the* $^\Phi$**Diverging IV**-class **P** = $(2^k\text{-}^\Phi B_0)$ ***eventually*** *reach a layer for a* **modulo** *high* enough *for the parity-sequence to be *reducing. No restrictions on* **r** ∈ ℕ *exist. But there* ***do*** *exist restrictions on how much higher than* **P** *it is* ***possible*** *for* ***r**$_s$ *to get, before the* **P** *is shown to be an element in a *Converging* **IV**-*class e.g. for all* **P** > 1 *it is observed that* ***r**$_s < 2^P$*. No* N *can visit the* ***power-set*** *of* N (*i.e.* $2^N$) *different layers/classes before they are *reduced. Not enough* $^\Phi$**Diverging** *routes/classes are found in the* Tree (*not enough* **O E**-**tiles** *are available*).

## Observations about sifted sets

If **Converging *parity-sequences* are *pruned* in The Domino-Tree (with $\mathbf{r}_{IV}$ as layers) or in **The** (OUT) **Fibonacci-Tree** (with $\mathbf{j}_{TV} = \mathbf{s}_{TV}+\mathbf{r}_{IV}$ as layers) the **Converging* **IV**-classes are in essence **sifting** the **set** of natural numbers, and each **reducing* layer "*catches*" **all** $N_X \in *Con_u(*r_s)_{q\le q_{Max}}$. *When all* $*_u[2^{*r_s}\mathbf{X}\text{-}*B_{Con}]_q$ *for all* **u** ∈ ℕ *from* **List 1** *have been pruned* no $^\Phi$**Diverging** *subset is left.*
**OBS**! An **Uncountable Infinity** of $^\Delta$**Redundant** routes/sub-classes is removed in each *pruning*.

**Final Comment**

Notice that in **TV**-classes $[3^s\mathbf{X}\text{-}B_{j=s+r}]_{TV}$ it is the parity of $B_{j_{TV}=s_{TV}+r_{IV}}$ that decides if it is the odd or the even **X** that are *lucky* and have the parity-sequence extended with an **E**-tile when the classes are ***split*** as the **modulo** is **doubled** and which subset is **unlucky** and get the **O E**-tile:

- Even $B_{j=s+r}$ ⇒ odd **X** (incl. **X** = 1) are **unlucky** because $[2^1 3^s\mathbf{X}\text{-}(3^s+B_{j=s+r})j]$ is **Odd**
- Odd $B_{j=s+r}$ ⇒ odd **X** (incl. **X** = 1) are *lucky* because $[2^1 3^s\mathbf{X}\text{-}(3^s+B_{j=s+r})j]$ is **Even**
- The $B_{TV}$ is *sometimes* odd and sometimes even, so *sometimes* $N_X$ for odd **X** are *lucky*

**NB**! Remember that $B_j$-values in Class-series follows the rules in the **3N minus 1** problem.

*As a growing number of *Converging* **IV**-*classes are confirmed, an equivalent number of elements disappears from the* **Union** *of* ᶲ**Diverging IV**-*classes and "The Average distance" between the individual* ᶲ**Diverging IV**-*classes i.e.* $2^r/|{}^{\Phi}\mathbf{Div}(\mathbf{r})|$ *tends to infinity for* $\mathbf{r} \to \infty$. *Eventually "The Average distance" will be higher than* **any** *previous* **Total** $2^r$ *for* $*\mathbf{r}_s >> \mathbf{r}$.

**End-IV**-*classes* $^{\dagger}[2^d\mathbf{X}\text{-}(2^d\text{-}N_N)_{End}]_N$ ***only*** *depends on* $\mathbf{d}$ = # **E** $_{Operation}$ *so eventually* **d** *is high* enough.

## Final Theorem

***The Collatz Conjecture*** *can not be false, because it is* ***impossible*** *for a counterexample to exist.*

## Proof for the Final Theorem

- **Formula (1)** (**reducing* number of **E**-operations $*\mathbf{r}_s = *\lceil \mathbf{sZ} \rceil$ for a given **s** = # **O** ) is exact
- ***The Sequence-formula*** (number of **Converging* **IV**-classes $_u(*\mathbf{r}_s)_q$ for a given **s**) is exact
- **Formula (3)** (*density* of **Converging* **IV**-classes $_u(*\mathbf{r}_s)_q$ for a given **s**) is exact
- ***The Series-formula*** (**S**um of *densities* for all $*Con_u(*\mathbf{r}_s)_q$ lower or equal to **s**) is exact
- That all the **Formulas** are exact *implies* that **Formula (5)** is also exact

*Density* of ᶲ**Diverging IV**-classes

$$|{}^{\Phi}\mathbf{Div}(*\mathbf{r}_k)|/2^{*r_k} = \sum_{s=k+1}^{\infty} \mathbf{F}(\mathbf{s}, *\mathbf{r}_s) \qquad \textbf{Formula (5)}$$

For all $\mathbf{k} \in \mathbb{N}$ the *density* of ᶲ**Diverging IV**-classes **equals** the **S**um for **Converging* **IV**-classes which have *higher* **modulo** $A_{IV} = 2^{*r_s}$ than $2^{\lceil kZ \rceil}$ (a *countable infinity* of **s** > **k** exist for all $\mathbf{k} \in \mathbb{N}$).

Suppose that in **Figure 2** all $*_u[2^{*r_s}\mathbf{X}\text{-}*B_{Con}]_q$ for **all** $\mathbf{u} \in \mathbb{N}$ are marked by colour. In **this** case there is *absolutely* NO ***uncovered*** area ***left*** in **Figure 2**. For $\mathbf{s} \to \infty$ the **Union** of **Converging* **IV**-classes $*_{u\to\infty}[2^{*r_s}\mathbf{X}\text{-}*B_{Con}]_{q\to\infty}$ for **all** $\mathbf{u} \in \mathbb{N}$ is a ***Covering System***. All the **Converging* **IV**-classes are *pairwise disjoint*, as they are *proper subsets* of $^{\Phi}[2^0\mathbf{X}\text{-}0_0]_{IV/TV}$

For all $N \in \mathbb{N}$ it is true for **modulo** $2^{*r_s}$ in the **Converging* **IV**-class $N \in *_u[2^{*r_s}\mathbf{X}\text{-}*B_{Con}]_q$ that the exponent $*\mathbf{r}_s$ can not be higher than $2^N$ i.e. not higher than the ***power-set*** of N. In reality it is observed that $*\mathbf{r}_s$ is much lower than $2^N$ and amongst "low N" it is observed that $*\mathbf{r}_s < 3N$ & the value $27 = 3^3$ represents a *local maximum* as $27 \in *_u[2^{*59_{37}}\mathbf{X}\text{-}*(2^{59}\text{-}27)_{Con}]_{q=1}$, $^{59}/_{27}$ = 2.185...

For any $\mathbf{k} \in \mathbb{N}$ the **Union** of ᶲ**Diverging IV**-classes $^{\Phi}[2^{*r_k}\mathbf{X}\text{-}B_0]_{IV}$ contain a ***lowest element***, as all classes are **well-ordered** and for a *combination* $(\mathbf{k}, *\mathbf{r}_k)$ a **finite** number of $^{\Phi}[2^{*r_k}\mathbf{X}\text{-}B_0]_{IV}$ exist. At ***some*** point for a higher **modulo** the ***lowest element*** is shown to belong to a $*_u[2^{*r_s}\mathbf{X}\text{-}*B_{Con}]_q$:
$*N_X \equiv \mathbf{P}_q \pmod{2^{*r_s}} \Leftrightarrow {}^{\Delta}B_0 \equiv *\mathbf{B}_{Con} \pmod{2^{*r_s}} \Rightarrow$ ***Union*** *of *Converging* $N_X \in *_{u\to\infty}[2^{*r_s}\mathbf{X}\text{-}*B_{Con}]_q \to \mathbb{N}$

A *hypothetical* lowest counter**e**xample $2^{r_{ce}-1} < N_{ce} < 2^{r_{ce}}$ can **not** exist which **implies** that ***The Collatz Conjecture*** can **not** be false! **QED**!

## Suggestions for future work

The ***method*** used in the present work can also be used on the **Q**N+1 Problem for other **Q** > 3:

***For example*** for **Q** = 5 it is simple to prove that more than one N-*loop* exist:

- $\mathbf{1}_0$ ↗ $6_1$ ↘ $3_2$ ↗ $16_3$ ↘ $8_4$ ↘ $4_5$ ↘ $2_6$ ↘ $\mathbf{1}_7$
- $\mathbf{13}_0$ ↗ $66_1$ ↘ $33_2$ ↗ $166_3$ ↘ $83_4$ ↗ $416_5$ ↘ $208_6$ ↘ $104_7$ ↘ $52_8$ ↘ $26_9$ ↘ $\mathbf{13}_{10}$

The ***higher*** the **Q** - the ***lower*** is the **S**um of **F**ractions **S**(**s**,*⌈**s***log***Q**/*log***2**⌉) for **s** → ∞.
It is **only** for **Q** = 3 (and the *degenerate case* for **Q** = 1) that it is **possible** for the ***Asymptotic Density*** for the **S**um of **F**ractions to have the **limit 1**; **S**(**s**,*⌈**s***log***3**/*log***2**⌉) → **1** for **s** → ∞.

Notice; $Z_5 = log5/log2 = 2.3219...$ $Z_7 = log7/log2 = 2.7475...$ $Z_9 = log3^2/log2 = 3.1699...$ ***etc***.

*In the end* ***The Collatz Conjecture*** *must be* ***true*** *because* $Z_3 < 2$, $Z_3 = log3/log2 = \mathbf{1.5849625}007...$

## References

The original inspiration for analysing the details in ***The Collatz Conjecture*** came from two YouTube videos available on Mr. Brady Haran's ***π*****Numberphile** channel (*thanks* Brady:)

## **Appendix Alpha**. *The first* **142** **Converging* **IV**-*classes* $*_{u}[2^{*r}\mathbf{X}\text{-}B_{Con}]_q$

$[2^{*r}\mathbf{X}\text{-}B_{Con}] - [3^s\mathbf{X}\text{-}B_{j=*r+s}] = [(2^{*r}-3^s)_D\mathbf{X}\text{-}(B_{Con}\text{-}B_{j=*r+s})_D]$ are called the **D**ifference-classes.

**Converging* **IV**-classes are arranged in **Groups** for $*\mathbf{r}_s$ and inside **Groups** by *size* of $\mathbf{P}_1 < \dots < \mathbf{P}_{q_{Max}}$

| **IV**-*class* | **TV**-*class* | **D**-*class* | **u** | **q** | $N_X \equiv \mathbf{2} \pmod{2^{*1}}$ |
|---|---|---|---|---|---|
| $[2^{*1}\mathbf{X}\text{-}0_{Con}]$ | $[3^0\mathbf{X}\text{-}0_{j=1}]$ | $[1_D\mathbf{X}\text{-}0_D]$ | 1 | **1** | *{**2**,4,6,8,10,12,14,16, ...} |

**F**(**0**,*1) = **1**/2^{*1} = 0.50 — **S**(**0**,*1) = **F**(**0**,*1) = 0.50

| **IV**-*class* | **TV**-*class* | **D**-*class* | **u** | **q** | $N_X \equiv \mathbf{1} \pmod{2^{*2}}$ |
|---|---|---|---|---|---|
| $[2^{*2}\mathbf{X}\text{-}3_{Con}]$ | $[3^1\mathbf{X}\text{-}2_3]$ | $[1_D\mathbf{X}\text{-}1_D]$ | 2 | **1** | *{**1**,5,9,13,17,21,25,29, ...} |

**F**(**1**,*2) = **1**/$2^{*2}$ = 0.25 — **S**(**1**,*2) = 0.75

*From here* **all** $*_{u}[2^{*r}\mathbf{X}\text{-}B_{Con}]_q$ *are* ***subsets*** *of the* ***Domain****-class* $^{\Phi}[2^2\mathbf{X}\text{-}1_0]_{IV}$ = $^{\Phi}${**3**,7,11,15,19, ...}

*Please notice that from here we observe that the Constant* $B_{Con} \equiv \mathbf{1} \pmod 4$

| **IV**-*class* | **TV**-*class* | **D**-*class* | **u** | **q** | $N_X \equiv \mathbf{3} \pmod{2^{*4}}$ |
|---|---|---|---|---|---|
| $[2^{*4}\mathbf{X}\text{-}13_{Con}]$ | $[3^2\mathbf{X}\text{-}7_6]$ | $[7_D\mathbf{X}\text{-}6_D]$ | 3 | **1** | *{**3**,19,35,51,67,83,99, ...} |

**F**(**2**,*4) = **1**/$2^{*4}$ = 0.0625 — **S**(**2**,*4) = 0.8125

| **IV**-*class* | **TV**-*class* | **D**-*class* | **u** | **q** | $N_X \equiv \mathbf{P} \pmod{2^{*5}}$ |
|---|---|---|---|---|---|
| $[2^{*5}\mathbf{X}\text{-}21_{Con}]$ | $[3^3\mathbf{X}\text{-}17_8]$ | $[5_D\mathbf{X}\text{-}4_D]$ | 4 | 1 | *{**11**,43,75,107,139,171, ...} |
| $[2^{*5}\mathbf{X}\text{- }9_{Con}]$ | $[3^3\mathbf{X}\text{-}7_8]$ | $[5_D\mathbf{X}\text{-}2_D]$ | 5 | **2** | *{**23**,55,87,119,151,183, ...} |

**F**(**3**,*5) = **2**/$2^{*5}$ = 0.0625 — **S**(**3**,*5) =0.875

| **IV**-*class* | **TV**-*class* | **D**-*class* | **u** | **q** | $N_X \equiv \mathbf{P} \pmod{2^{*7}}$ |
|---|---|---|---|---|---|
| $[2^{*7}\mathbf{X}\text{-}121_{Con}]$ | $[3^4\mathbf{X}\text{-}76_{11}]$ | $[47_D\mathbf{X}\text{-}45_D]$ | 6 | 1 | *{**7**,135,263,391,519,647, ...} |
| $[2^{*7}\mathbf{X}\text{-}113_{Con}]$ | $[3^4\mathbf{X}\text{-}71_{11}]$ | $[47_D\mathbf{X}\text{-}42_D]$ | 7 | 2 | *{**15**,143,271,399,527,655, ...} |
| $[2^{*7}\mathbf{X}\text{- }69_{Con}]$ | $[3^4\mathbf{X}\text{-}43_{11}]$ | $[47_D\mathbf{X}\text{-}26_D]$ | 8 | **3** | *{**59**,187,315,443,571,699, ...} |

**F**(**4**,*7) = **3**/$2^{*7}$ = 0.0234375 — **S**(**4**,*7) = 0.8984375

| **IV**-*class* | **TV**-*class* | **D**-*class* | **u** | **q** | $N_X \equiv \mathbf{P} \pmod{2^{*8}}$ |
|---|---|---|---|---|---|
| $[2^{*8}\mathbf{X}\text{-}217_{Con}]$ | $[3^5\mathbf{X}\text{-}205_{13}]$ | $[13_D\mathbf{X}\text{-}12_D]$ | 9 | 1 | *{**39**,295,551,807,1063, ...} |
| $[2^{*8}\mathbf{X}\text{-}177_{Con}]$ | $[3^5\mathbf{X}\text{-}167_{13}]$ | $[13_D\mathbf{X}\text{-}10_D]$ | 10 | 2 | *{**79**,335,591,847,1103, ...} |
| $[2^{*8}\mathbf{X}\text{-}161_{Con}]$ | $[3^5\mathbf{X}\text{-}152_{13}]$ | $[13_D\mathbf{X}\text{-}9_D]$ | 11 | 3 | *{**95**,351,607,863,1119, ...} |
| $[2^{*8}\mathbf{X}\text{-}133_{Con}]$ | $[3^5\mathbf{X}\text{-}125_{13}]$ | $[13_D\mathbf{X}\text{-}8_D]$ | 12 | 4 | *{**123**,379,635,891,1147, ...} |
| $[2^{*8}\mathbf{X}\text{-}81_{Con}]$ | $[3^5\mathbf{X}\text{-}76_{13}]$ | $[13_D\mathbf{X}\text{-}5_D]$ | 13 | 5 | *{**175**,431,687,943,1199, ...} |
| $[2^{*8}\mathbf{X}\text{-}57_{Con}]$ | $[3^5\mathbf{X}\text{-}53_{13}]$ | $[13_D\mathbf{X}\text{-}4_D]$ | 14 | 6 | *{**199**,455,711,967,1223, ...} |
| $[2^{*8}\mathbf{X}\text{-}37_{Con}]$ | $[3^5\mathbf{X}\text{-}34_{13}]$ | $[13_D\mathbf{X}\text{-}3_D]$ | 15 | **7** | *{**219**,475,731,987,1243, ...} |

**F**(**5**,*8) = **7**/$2^{*8}$ = 0.02734375 — **S**(**5**,*8) = 0.92578125

| **IV**-*class* | **TV**-*class* | **D**-*class* | **u** | **q** | $N_X \equiv$ **P** (mod $2^{*10}$) |
|---|---|---|---|---|---|
| $[2^{*10}\mathbf{X}-737_0]$ | $[3^6\mathbf{X}-524_{16}]$ | $[295_D\mathbf{X}-213_D]$ | 16 | 1 | *{**287**,1311,2335,3359,4383, ...} |
| $[2^{*10}\mathbf{X}-677_0]$ | $[3^6\mathbf{X}-481_{16}]$ | $[295_D\mathbf{X}-196_D]$ | 17 | 2 | *{**347**,1371,2395,3419,4443, ...} |
| $[2^{*10}\mathbf{X}-657_0]$ | $[3^6\mathbf{X}-467_{16}]$ | $[295_D\mathbf{X}-190_D]$ | 18 | 3 | *{**367**,1391,2415,3439,4463, ...} |
| $[2^{*10}\mathbf{X}-601_0]$ | $[3^6\mathbf{X}-427_{16}]$ | $[295_D\mathbf{X}-174_D]$ | 19 | 4 | *{**423**,1447,2471,3495,4519, ...} |
| $[2^{*10}\mathbf{X}-517_0]$ | $[3^6\mathbf{X}-367_{16}]$ | $[295_D\mathbf{X}-150_D]$ | 20 | 5 | *{**507**,1531,2555,3579,4603, ...} |
| $[2^{*10}\mathbf{X}-449_0]$ | $[3^6\mathbf{X}-319_{16}]$ | $[295_D\mathbf{X}-130_D]$ | 21 | 6 | *{**575**,1599,2623,3647,4671, ...} |
| $[2^{*10}\mathbf{X}-441_0]$ | $[3^6\mathbf{X}-313_{16}]$ | $[295_D\mathbf{X}-128_D]$ | 22 | 7 | *{**583**,1607,2631,3655,4679, ...} |
| $[2^{*10}\mathbf{X}-289_0]$ | $[3^6\mathbf{X}-205_{16}]$ | $[295_D\mathbf{X}-84_D]$ | 23 | 8 | *{**735**,1759,2783,3807,4831, ...} |
| $[2^{*10}\mathbf{X}-209_0]$ | $[3^6\mathbf{X}-148_{16}]$ | $[295_D\mathbf{X}-61_D]$ | 24 | 9 | *{**815**,1839,2863,3887,4911, ...} |
| $[2^{*10}\mathbf{X}-101_0]$ | $[3^6\mathbf{X}-71_{16}]$ | $[295_D\mathbf{X}-30_D]$ | 25 | 10 | *{**923**,1947,2971,3995,5019, ...} |
| $[2^{*10}\mathbf{X}-49_0]$ | $[3^6\mathbf{X}-34_{16}]$ | $[295_D\mathbf{X}-15_D]$ | 26 | 11 | *{**975**,1999,3023,4047,5071, ...} |
| $[2^{*10}\mathbf{X}-25_0]$ | $[3^6\mathbf{X}-17_{16}]$ | $[295_D\mathbf{X}-8_D]$ | 27 | **12** | *{**999**,2023,3047,4071,5095, ...} |

**F(6,*10) = 12**/$2^{*10}$ = 0.01171875    **S(6,*10)** = 0.9375

The (*infinite*) list continues on the next page with the **30** **reducing* **IV**-*classes* having $A_{IV} = 2^{*12}$ and the following three pages with the **85** **reducing* **IV**-*classes* having $A_{IV} = 2^{*13}$.

Remember that $*\mathbf{r}_s$ **= [sZ]**, **s** ∊ ℕ, **Z** = *log*3/*log*2 = **1.5849625**007...

*An overview of the* **F**ractions, **F(s,*r)**, *and* **S**ums, **S(s,*r)**, *for the *Reducing Combinations* **(s,*r)**:

| **s** | $*\mathbf{r}_s$ **= [sZ]** | **j = s+*r** | **Total =** $\mathbf{2^{*r}}$ | $\mathbf{q_{Max}}$ | **F(s,*r) =** $\mathbf{q_{Max}/2^{*r}}$ | **S(s,*r)** |
|---|---|---|---|---|---|---|
| **0** | $1_0$ | 1 | 2 | **1** | 0.5 | 0.5 |
| **1** | 2 | 3 | 4 | **1** | 0.25 | 0.75 |
| **2** | 4 | 6 | 16 | **1** | 0.0625 | 0.8125 |
| **3** | 5 | 8 | 32 | **2** | 0.0625 | 0.875 |
| **4** | 7 | 11 | 128 | **3** | 0.0234375 | 0.8984375 |
| **5** | 8 | 13 | 256 | **7** | 0.02734375 | 0.**9**2578125 |
| **6** | 10 | 16 | 1024 | **12** | 0.01171875 | 0.9375 |
| **7** | 12 | 19 | 4096 | **30** | 0.00732421875 | 0.**9**4482421875 |
| **8** | 13 | 21 | 8192 | **85** | 0.0103759765625 | 0.**9**552001953125 |

Notice that these are **exact** values for the **F**ractions and **S**ums.

For (**6**,***10**) we *observe* that <u>93.75%</u> of N ∊ ℕ are **reduced* in <u>at most</u> (6+10) = 16 operations.

**Exactly** <u>95.52001953125%</u> of N ∊ ℕ are **reduced* in <u>at most</u> (8+13) = 21 operations and as explained in the work the ***same part*** of 1 ≤ N ≤ $2^{*13}$ are **reduced* in <u>at most</u> 21 operations.

| **IV**-*class* | **TV**-*class* | **D**-*class* | **u** | **q** | $N_X \equiv$ **P** (mod $2^{*12}$) |
|---|---|---|---|---|---|
| $[2^{*12}\mathbf{X}\text{-}3865_{0}]$ | $[3^{7}\mathbf{X}\text{-}2063_{19}]$ | $[1909_{D}\mathbf{X}\text{-}1802_{D}]$ | 28 | 1 | *{**231**,4327,8423,12519,16615, ...} |
| $[2^{*12}\mathbf{X}\text{-}3713_{0}]$ | $[3^{7}\mathbf{X}\text{-}1982_{19}]$ | $[1909_{D}\mathbf{X}\text{-}1731_{D}]$ | 29 | 2 | *{**383**,4479,8575,12671,16767, ...} |
| $[2^{*12}\mathbf{X}\text{-}3633_{0}]$ | $[3^{7}\mathbf{X}\text{-}1939_{19}]$ | $[1909_{D}\mathbf{X}\text{-}1694_{D}]$ | 30 | 3 | *{**463**,4559,8655,12751,16847, ...} |
| $[2^{*12}\mathbf{X}\text{-}3481_{0}]$ | $[3^{7}\mathbf{X}\text{-}1858_{19}]$ | $[1909_{D}\mathbf{X}\text{-}1623_{D}]$ | 31 | 4 | *{**615**,4711,8807,12903,16999, ...} |
| $[2^{*12}\mathbf{X}\text{-}3217_{0}]$ | $[3^{7}\mathbf{X}\text{-}1717_{19}]$ | $[1909_{D}\mathbf{X}\text{-}1500_{D}]$ | 32 | 5 | *{**879**,4975,9071,13167,17263, ...} |
| $[2^{*12}\mathbf{X}\text{-}3161_{0}]$ | $[3^{7}\mathbf{X}\text{-}1687_{19}]$ | $[1909_{D}\mathbf{X}\text{-}1474_{D}]$ | 33 | 6 | *{**935**,5031,9127,13223,17319, ...} |
| $[2^{*12}\mathbf{X}\text{-}3077_{0}]$ | $[3^{7}\mathbf{X}\text{-}1642_{19}]$ | $[1909_{D}\mathbf{X}\text{-}1435_{D}]$ | 34 | 7 | *{**1019**,5115,9211,13307,17403, ...} |
| $[2^{*12}\mathbf{X}\text{-}3009_{0}]$ | $[3^{7}\mathbf{X}\text{-}1606_{19}]$ | $[1909_{D}\mathbf{X}\text{-}1403_{D}]$ | 35 | 8 | *{**1087**,5183,9279,13375,17471, ...} |
| $[2^{*12}\mathbf{X}\text{-}2865_{0}]$ | $[3^{7}\mathbf{X}\text{-}1529_{19}]$ | $[1909_{D}\mathbf{X}\text{-}1336_{D}]$ | 36 | 9 | *{**1231**,5327,9423,13519,17615, ...} |
| $[2^{*12}\mathbf{X}\text{-}2661_{0}]$ | $[3^{7}\mathbf{X}\text{-}1420_{19}]$ | $[1909_{D}\mathbf{X}\text{-}1241_{D}]$ | 37 | 10 | *{**1435**,5531,9627,13723,17819, ...} |
| $[2^{*12}\mathbf{X}\text{-}2449_{0}]$ | $[3^{7}\mathbf{X}\text{-}1391_{19}]$ | $[1909_{D}\mathbf{X}\text{-}1058_{D}]$ | 38 | 11 | *{**1647**,5743,9839,13935,18031, ...} |
| $[2^{*12}\mathbf{X}\text{-}2393_{0}]$ | $[3^{7}\mathbf{X}\text{-}1277_{19}]$ | $[1909_{D}\mathbf{X}\text{-}1116_{D}]$ | 39 | 12 | *{**1703**,5799,9895,13991,18087, ...} |
| $[2^{*12}\mathbf{X}\text{-}2309_{0}]$ | $[3^{7}\mathbf{X}\text{-}1232_{19}]$ | $[1909_{D}\mathbf{X}\text{-}1077_{D}]$ | 40 | 13 | *{**1787**,5883,9979,14075,18171, ...} |
| $[2^{*12}\mathbf{X}\text{-}2273_{0}]$ | $[3^{7}\mathbf{X}\text{-}1213_{19}]$ | $[1909_{D}\mathbf{X}\text{-}1060_{D}]$ | 41 | 14 | *{**1823**,5919,10015,14111,18207, ...} |
| $[2^{*12}\mathbf{X}\text{-}2241_{0}]$ | $[3^{7}\mathbf{X}\text{-}1196_{19}]$ | $[1909_{D}\mathbf{X}\text{-}1045_{D}]$ | 42 | 15 | *{**1855**,5951,10047,14143,18239, ...} |
| $[2^{*12}\mathbf{X}\text{-}2065_{0}]$ | $[3^{7}\mathbf{X}\text{-}1102_{19}]$ | $[1909_{D}\mathbf{X}\text{-}963_{D}]$ | 43 | 16 | *{**2031**,6127,10223,14319,18415, ...} |
| $[2^{*12}\mathbf{X}\text{-}1893_{0}]$ | $[3^{7}\mathbf{X}\text{-}1010_{19}]$ | $[1909_{D}\mathbf{X}\text{-}883_{D}]$ | 44 | 17 | *{**2203**,6299,10395,14491,18587, ...} |
| $[2^{*12}\mathbf{X}\text{-}1857_{0}]$ | $[3^{7}\mathbf{X}\text{-}991_{19}]$ | $[1909_{D}\mathbf{X}\text{-}866_{D}]$ | 45 | 18 | *{**2239**,6335,10431,14527,18623, ...} |
| $[2^{*12}\mathbf{X}\text{-}1745_{0}]$ | $[3^{7}\mathbf{X}\text{-}931_{19}]$ | $[1909_{D}\mathbf{X}\text{-}814_{D}]$ | 46 | 19 | *{**2351**,6447,10543,14639,18735, ...} |
| $[2^{*12}\mathbf{X}\text{-}1509_{0}]$ | $[3^{7}\mathbf{X}\text{-}805_{19}]$ | $[1909_{D}\mathbf{X}\text{-}704_{D}]$ | 47 | 20 | *{**2587**,6683,10779,14875,18971, ...} |
| $[2^{*12}\mathbf{X}\text{-}1505_{0}]$ | $[3^{7}\mathbf{X}\text{-}803_{19}]$ | $[1909_{D}\mathbf{X}\text{-}702_{D}]$ | 48 | 21 | *{**2591**,6687,10783,14879,18975, ...} |
| $[2^{*12}\mathbf{X}\text{-}1189_{0}]$ | $[3^{7}\mathbf{X}\text{-}634_{19}]$ | $[1909_{D}\mathbf{X}\text{-}555_{D}]$ | 49 | 22 | *{**2907**,7003,11099,15195,19291, ...} |
| $[2^{*12}\mathbf{X}\text{-}1121_{0}]$ | $[3^{7}\mathbf{X}\text{-}598_{19}]$ | $[1909_{D}\mathbf{X}\text{-}523_{D}]$ | 50 | 23 | *{**2975**,7071,11167,15263,19359, ...} |
| $[2^{*12}\mathbf{X}\text{-}977_{0}]$ | $[3^{7}\mathbf{X}\text{-}521_{19}]$ | $[1909_{D}\mathbf{X}\text{-}456_{D}]$ | 51 | 24 | *{**3119**,7215,11311,15407,19503, ...} |
| $[2^{*12}\mathbf{X}\text{-}953_{0}]$ | $[3^{7}\mathbf{X}\text{-}508_{19}]$ | $[1909_{D}\mathbf{X}\text{-}445_{D}]$ | 52 | 25 | *{**3143**,7239,11335,15431,19527, ...} |
| $[2^{*12}\mathbf{X}\text{-}801_{0}]$ | $[3^{7}\mathbf{X}\text{-}427_{19}]$ | $[1909_{D}\mathbf{X}\text{-}374_{D}]$ | 53 | 26 | *{**3295**,7391,11487,15583,19679, ...} |
| $[2^{*12}\mathbf{X}\text{-}537_{0}]$ | $[3^{7}\mathbf{X}\text{-}286_{19}]$ | $[1909_{D}\mathbf{X}\text{-}251_{D}]$ | 54 | 27 | *{**3559**,7655,11751,15847,19943, ...} |
| $[2^{*12}\mathbf{X}\text{-}421_{0}]$ | $[3^{7}\mathbf{X}\text{-}224_{19}]$ | $[1909_{D}\mathbf{X}\text{-}197_{D}]$ | 55 | 28 | *{**3675**,7771,11867,15963,20059, ...} |
| $[2^{*12}\mathbf{X}\text{-}185_{0}]$ | $[3^{7}\mathbf{X}\text{-}98_{19}]$ | $[1909_{D}\mathbf{X}\text{-}87_{D}]$ | 56 | 29 | *{**3911**,8007,12103,16199,20295, ...} |
| $[2^{*12}\mathbf{X}\text{-}33_{0}]$ | $[3^{7}\mathbf{X}\text{-}17_{19}]$ | $[1909_{D}\mathbf{X}\text{-}16_{D}]$ | 57 | **30** | *{**4063**,8159,12255,16351,20447, ...} |

**F(7,*12) = 30**/$2^{*12}$ = 0.00732421875

**S(7,*12)** = 0.94482421875

| **IV**-*class* | **TV**-*class* | **D**-*class* | **u** | **q** | $N_X \equiv$ **P** (mod $2^{*13}$) |
|---|---|---|---|---|---|
| [$2^{*13}$**X**-8001$_0$] | [$3^8$**X**-6407$_{21}$] | [1631$_D$**X**-1594$_D$] | 58 | 1 | *{**191**,8383,16575,24767, ...} |
| [$2^{*13}$**X**-7985$_0$] | [$3^8$**X**-6394$_{21}$] | [1631$_D$**X**-1591$_D$] | 59 | 2 | *{**207**,8399,16591,24783, ...} |
| [$2^{*13}$**X**-7937$_0$] | [$3^8$**X**-6356$_{21}$] | [1631$_D$**X**-1581$_D$] | 60 | 3 | *{**255**,8447,16639,24831, ...} |
| [$2^{*13}$**X**-7889$_0$] | [$3^8$**X**-6317$_{21}$] | [1631$_D$**X**-1572$_D$] | 61 | 4 | *{**303**,8495,16687,24879, ...} |
| [$2^{*13}$**X**-7653$_0$] | [$3^8$**X**-6128$_{21}$] | [1631$_D$**X**-1525$_D$] | 62 | 5 | *{**539**,8731,16923,25115, ...} |
| [$2^{*13}$**X**-7649$_0$] | [$3^8$**X**-6125$_{21}$] | [1631$_D$**X**-1524$_D$] | 63 | 6 | *{**543**,8735,16927,25119, ...} |
| [$2^{*13}$**X**-7569$_0$] | [$3^8$**X**-6061$_{21}$] | [1631$_D$**X**-1508$_D$] | 64 | 7 | *{**623**,8815,17007,25199, ...} |
| [$2^{*13}$**X**-7513$_0$] | [$3^8$**X**-6016$_{21}$] | [1631$_D$**X**-1497$_D$] | 65 | 8 | *{**679**,8871,17063,25255, ...} |
| [$2^{*13}$**X**-7473$_0$] | [$3^8$**X**-5984$_{21}$] | [1631$_D$**X**-1489$_D$] | 66 | 9 | *{**719**,8911,17103,25295, ...} |
| [$2^{*13}$**X**-7393$_0$] | [$3^8$**X**-5920$_{21}$] | [1631$_D$**X**-1473$_D$] | 67 | 10 | *{**799**,8991,17183,25375, ...} |
| [$2^{*13}$**X**-7121$_0$] | [$3^8$**X**-5702$_{21}$] | [1631$_D$**X**-1419$_D$] | 68 | 11 | *{**1071**,9263,17455,25647, ...} |
| [$2^{*13}$**X**-7057$_0$] | [$3^8$**X**-5651$_{21}$] | [1631$_D$**X**-1406$_D$] | 69 | 12 | *{**1135**,9327,17519,25711, ...} |
| [$2^{*13}$**X**-7001$_0$] | [$3^8$**X**-5606$_{21}$] | [1631$_D$**X**-1395$_D$] | 70 | 13 | *{**1191**,9383,17575,25767, ...} |
| [$2^{*13}$**X**-6977$_0$] | [$3^8$**X**-5587$_{21}$] | [1631$_D$**X**-1390$_D$] | 71 | 14 | *{**1215**,9407,17599,25791, ...} |
| [$2^{*13}$**X**-6945$_0$] | [$3^8$**X**-5561$_{21}$] | [1631$_D$**X**-1384$_D$] | 72 | 15 | *{**1247**,9439,17631,25823, ...} |
| [$2^{*13}$**X**-6865$_0$] | [$3^8$**X**-5497$_{21}$] | [1631$_D$**X**-1368$_D$] | 73 | 16 | *{**1327**,9519,17711,25903, ...} |
| [$2^{*13}$**X**-6629$_0$] | [$3^8$**X**-5308$_{21}$] | [1631$_D$**X**-1321$_D$] | 74 | 17 | *{**1563**,9755,17947,26139, ...} |
| [$2^{*13}$**X**-6625$_0$] | [$3^8$**X**-5305$_{21}$] | [1631$_D$**X**-1320$_D$] | 75 | 18 | *{**1567**,9759,17951,26143, ...} |
| [$2^{*13}$**X**-6465$_0$] | [$3^8$**X**-5177$_{21}$] | [1631$_D$**X**-1288$_D$] | 76 | 19 | *{**1727**,9919,18111,26303, ...} |
| [$2^{*13}$**X**-6209$_0$] | [$3^8$**X**-4972$_{21}$] | [1631$_D$**X**-1237$_D$] | 77 | 20 | *{**1983**,10175,18367,26559, ...} |
| [$2^{*13}$**X**-6177$_0$] | [$3^8$**X**-4946$_{21}$] | [1631$_D$**X**-1231$_D$] | 78 | 21 | *{**2015**,10207,18399,26591, ...} |
| [$2^{*13}$**X**-6117$_0$] | [$3^8$**X**-4898$_{21}$] | [1631$_D$**X**-1219$_D$] | 79 | 22 | *{**2075**,10267,18459,26651, ...} |
| [$2^{*13}$**X**-6113$_0$] | [$3^8$**X**-4895$_{21}$] | [1631$_D$**X**-1218$_D$] | 80 | 23 | *{**2079**,10271,18463,26655, ...} |
| [$2^{*13}$**X**-6097$_0$] | [$3^8$**X**-4882$_{21}$] | [1631$_D$**X**-1215$_D$] | 81 | 24 | *{**2095**,10287,18479,26671, ...} |
| [$2^{*13}$**X**-5921$_0$] | [$3^8$**X**-4741$_{21}$] | [1631$_D$**X**-1180$_D$] | 82 | 25 | *{**2271**,10463,18655,26847, ...} |
| [$2^{*13}$**X**-5861$_0$] | [$3^8$**X**-4693$_{21}$] | [1631$_D$**X**-1168$_D$] | 83 | 26 | *{**2331**,10523,18715,26907, ...} |
| [$2^{*13}$**X**-5761$_0$] | [$3^8$**X**-4613$_{21}$] | [1631$_D$**X**-1148$_D$] | 84 | 27 | *{**2431**,10623,18815,27007, ...} |
| [$2^{*13}$**X**-5585$_0$] | [$3^8$**X**-4472$_{21}$] | [1631$_D$**X**-1113$_D$] | 85 | 28 | *{**2607**,10799,18991,27183, ...} |
| [$2^{*13}$**X**-5529$_0$] | [$3^8$**X**-4427$_{21}$] | [1631$_D$**X**-1102$_D$] | 86 | 29 | *{**2663**,10855,19047,27239, ...} |
| [$2^{*13}$**X**-5153$_0$] | [$3^8$**X**-4126$_{21}$] | [1631$_D$**X**-1027$_D$] | 87 | 30 | *{**3039**,11231,19423,27615, ...} |

| **IV**-*class* | **TV**-*class* | **D**-*class* | **u** | **q** | $N_X \equiv$ **P** (mod $2^{*13}$) |
|---|---|---|---|---|---|
| [$2^{*13}$**X**-$5125_0$] | [$3^8$**X**-$4103_{21}$] | [$1631_D$**X**-$1022_D$] | 88 | 31 | *{**3067**,11259,19451,27643, ...} |
| [$2^{*13}$**X**-$5057_0$] | [$3^8$**X**-$4049_{21}$] | [$1631_D$**X**-$1008_D$] | 89 | 32 | *{**3135**,11327,19519,27711, ...} |
| [$2^{*13}$**X**-$4737_0$] | [$3^8$**X**-$3793_{21}$] | [$1631_D$**X**-$944_D$] | 90 | 33 | *{**3455**,11647,19839,28031, ...} |
| [$2^{*13}$**X**-$4709_0$] | [$3^8$**X**-$3770_{21}$] | [$1631_D$**X**-$939_D$] | 91 | 34 | *{**3483**,11675,19867,28059, ...} |
| [$2^{*13}$**X**-$4641_0$] | [$3^8$**X**-$3716_{21}$] | [$1631_D$**X**-$925_D$] | 92 | 35 | *{**3551**,11743,19935,28127, ...} |
| [$2^{*13}$**X**-$4505_0$] | [$3^8$**X**-$3607_{21}$] | [$1631_D$**X**-$898_D$] | 93 | 36 | *{**3687**,11879,20071,28263, ...} |
| [$2^{*13}$**X**-$4357_0$] | [$3^8$**X**-$3488_{21}$] | [$1631_D$**X**-$869_D$] | 94 | 37 | *{**3835**,12027,20219,28411, ...} |
| [$2^{*13}$**X**-$4289_0$] | [$3^8$**X**-$3434_{21}$] | [$1631_D$**X**-$855_D$] | 95 | 38 | *{**3903**,12095,20287,28479, ...} |
| [$2^{*13}$**X**-$4225_0$] | [$3^8$**X**-$3383_{21}$] | [$1631_D$**X**-$842_D$] | 96 | 39 | *{**3967**,12159,20351,28543, ...} |
| [$2^{*13}$**X**-$4113_0$] | [$3^8$**X**-$3293_{21}$] | [$1631_D$**X**-$820_D$] | 97 | 40 | *{**4079**,12271,20463,28655, ...} |
| [$2^{*13}$**X**-$4101_0$] | [$3^8$**X**-$3283_{21}$] | [$1631_D$**X**-$818_D$] | 98 | 41 | *{**4091**,12283,20475,28667, ...} |
| [$2^{*13}$**X**-$4033_0$] | [$3^8$**X**-$3229_{21}$] | [$1631_D$**X**-$804_D$] | 99 | 42 | *{**4159**,12351,20543,28735, ...} |
| [$2^{*13}$**X**-$3993_0$] | [$3^8$**X**-$3197_{21}$] | [$1631_D$**X**-$796_D$] | 100 | 43 | *{**4199**,12391,20583,28775, ...} |
| [$2^{*13}$**X**-$3969_0$] | [$3^8$**X**-$3178_{21}$] | [$1631_D$**X**-$791_D$] | 101 | 44 | *{**4223**,12415,20607,28799, ...} |
| [$2^{*13}$**X**-$3941_0$] | [$3^8$**X**-$3155_{21}$] | [$1631_D$**X**-$786_D$] | 102 | 45 | *{**4251**,12443,20635,28827, ...} |
| [$2^{*13}$**X**-$3737_0$] | [$3^8$**X**-$2992_{21}$] | [$1631_D$**X**-$745_D$] | 103 | 46 | *{**4455**,12647,20839,29031, ...} |
| [$2^{*13}$**X**-$3685_0$] | [$3^8$**X**-$2950_{21}$] | [$1631_D$**X**-$735_D$] | 104 | 47 | *{**4507**,12699,20891,29083, ...} |
| [$2^{*13}$**X**-$3333_0$] | [$3^8$**X**-$2668_{21}$] | [$1631_D$**X**-$665_D$] | 105 | 48 | *{**4859**,13051,21243,29435, ...} |
| [$2^{*13}$**X**-$3265_0$] | [$3^8$**X**-$2614_{21}$] | [$1631_D$**X**-$651_D$] | 106 | 49 | *{**4927**,13119,21311,29503, ...} |
| [$2^{*13}$**X**-$3237_0$] | [$3^8$**X**-$2591_{21}$] | [$1631_D$**X**-$646_D$] | 107 | 50 | *{**4955**,13147,21339,29531, ...} |
| [$2^{*13}$**X**-$3169_0$] | [$3^8$**X**-$2537_{21}$] | [$1631_D$**X**-$632_D$] | 108 | 51 | *{**5023**,13215,21407,29599, ...} |
| [$2^{*13}$**X**-$3089_0$] | [$3^8$**X**-$2473_{21}$] | [$1631_D$**X**-$616_D$] | 109 | 52 | *{**5103**,13295,21487,29679, ...} |
| [$2^{*13}$**X**-$3001_0$] | [$3^8$**X**-$2402_{21}$] | [$1631_D$**X**-$599_D$] | 110 | 53 | *{**5191**,13383,21575,29767, ...} |
| [$2^{*13}$**X**-$2917_0$] | [$3^8$**X**-$2335_{21}$] | [$1631_D$**X**-$582_D$] | 111 | 54 | *{**5275**,13467,21659,29851, ...} |
| [$2^{*13}$**X**-$2821_0$] | [$3^8$**X**-$2258_{21}$] | [$1631_D$**X**-$563_D$] | 112 | 55 | *{**5371**,13563,21755,29947, ...} |
| [$2^{*13}$**X**-$2753_0$] | [$3^8$**X**-$2204_{21}$] | [$1631_D$**X**-$549_D$] | 113 | 56 | *{**5439**,13631,21823,30015, ...} |
| [$2^{*13}$**X**-$2585_0$] | [$3^8$**X**-$2069_{21}$] | [$1631_D$**X**-$516_D$] | 114 | 57 | *{**5607**,13799,21991,30183, ...} |
| [$2^{*13}$**X**-$2577_0$] | [$3^8$**X**-$2063_{21}$] | [$1631_D$**X**-$514_D$] | 115 | 58 | *{**5615**,13807,21999,30191, ...} |
| [$2^{*13}$**X**-$2469_0$] | [$3^8$**X**-$1976_{21}$] | [$1631_D$**X**-$493_D$] | 116 | 59 | *{**5723**,13915,22107,30299, ...} |
| [$2^{*13}$**X**-$2405_0$] | [$3^8$**X**-$1925_{21}$] | [$1631_D$**X**-$480_D$] | 117 | 60 | *{**5787**,13979,22171,30363, ...} |

| **IV**-*class* | **TV**-*class* | **D**-*class* | **u** | **q** | $N_X \equiv$ **P** (mod $2^{*13}$) |
|---|---|---|---|---|---|
| [$2^{*13}$**X**-$2321_0$] | [$3^8$**X**-$1858_{21}$] | [$1631_D$**X**-$463_D$] | 118 | 61 | *{**5871**,14063,22255,30447, ...} |
| [$2^{*13}$**X**-$2233_0$] | [$3^8$**X**-$1787_{21}$] | [$1631_D$**X**-$446_D$] | 119 | 62 | *{**5959**,14151,22343,30535, ...} |
| [$2^{*13}$**X**-$2213_0$] | [$3^8$**X**-$1771_{21}$] | [$1631_D$**X**-$442_D$] | 120 | 63 | *{**5979**,14171,22363,30555, ...} |
| [$2^{*13}$**X**-$2145_0$] | [$3^8$**X**-$1717_{21}$] | [$1631_D$**X**-$428_D$] | 121 | 64 | *{**6047**,14239,22431,30623, ...} |
| [$2^{*13}$**X**-$1977_0$] | [$3^8$**X**-$1582_{21}$] | [$1631_D$**X**-$395_D$] | 122 | 65 | *{**6215**,14407,22599,30791, ...} |
| [$2^{*13}$**X**-$1817_0$] | [$3^8$**X**-$1454_{21}$] | [$1631_D$**X**-$363_D$] | 123 | 66 | *{**6375**,14567,22759,30951, ...} |
| [$2^{*13}$**X**-$1633_0$] | [$3^8$**X**-$1307_{21}$] | [$1631_D$**X**-$326_D$] | 124 | 67 | *{**6559**,14751,22943,31135, ...} |
| [$2^{*13}$**X**-$1585_0$] | [$3^8$**X**-$1268_{21}$] | [$1631_D$**X**-$317_D$] | 125 | 68 | *{**6607**,14799,22991,31183, ...} |
| [$2^{*13}$**X**-$1561_0$] | [$3^8$**X**-$1249_{21}$] | [$1631_D$**X**-$312_D$] | 126 | 69 | *{**6631**,14823,23015,31207, ...} |
| [$2^{*13}$**X**-$1445_0$] | [$3^8$**X**-$1156_{21}$] | [$1631_D$**X**-$289_D$] | 127 | 70 | *{**6747**,14939,23131,31323, ...} |
| [$2^{*13}$**X**-$1377_0$] | [$3^8$**X**-$1102_{21}$] | [$1631_D$**X**-$275_D$] | 128 | 71 | *{**6815**,15007,23199,31391, ...} |
| [$2^{*13}$**X**-$1209_0$] | [$3^8$**X**-$967_{21}$] | [$1631_D$**X**-$242_D$] | 129 | 72 | *{**6983**,15175,23367,31559, ...} |
| [$2^{*13}$**X**-$1169_0$] | [$3^8$**X**-$935_{21}$] | [$1631_D$**X**-$234_D$] | 130 | 73 | *{**7023**,15215,23407,31599, ...} |
| [$2^{*13}$**X**-$1113_0$] | [$3^8$**X**-$890_{21}$] | [$1631_D$**X**-$223_D$] | 131 | 74 | *{**7079**,15271,23463,31655, ...} |
| [$2^{*13}$**X**-$933_0$] | [$3^8$**X**-$746_{21}$] | [$1631_D$**X**-$187_D$] | 132 | 75 | *{**7259**,15451,23643,31835, ...} |
| [$2^{*13}$**X**-$817_0$] | [$3^8$**X**-$653_{21}$] | [$1631_D$**X**-$164_D$] | 133 | 76 | *{**7375**,15567,23759,31951, ...} |
| [$2^{*13}$**X**-$793_0$] | [$3^8$**X**-$634_{21}$] | [$1631_D$**X**-$159_D$] | 134 | 77 | *{**7399**,15591,23783,31975, ...} |
| [$2^{*13}$**X**-$697_0$] | [$3^8$**X**-$557_{21}$] | [$1631_D$**X**-$140_D$] | 135 | 78 | *{**7495**,15687,23879,32071, ...} |
| [$2^{*13}$**X**-$561_0$] | [$3^8$**X**-$448_{21}$] | [$1631_D$**X**-$113_D$] | 136 | 79 | *{**7631**,15823,24015,32207, ...} |
| [$2^{*13}$**X**-$401_0$] | [$3^8$**X**-$320_{21}$] | [$1631_D$**X**-$81_D$] | 137 | 80 | *{**7791**,15983,24175,32367, ...} |
| [$2^{*13}$**X**-$345_0$] | [$3^8$**X**-$275_{21}$] | [$1631_D$**X**-$70_D$] | 138 | 81 | *{**7847**,16039,24231,32423, ...} |
| [$2^{*13}$**X**-$281_0$] | [$3^8$**X**-$224_{21}$] | [$1631_D$**X**-$57_D$] | 139 | 82 | *{**7911**,16103,24295,32487, ...} |
| [$2^{*13}$**X**-$225_0$] | [$3^8$**X**-$179_{21}$] | [$1631_D$**X**-$46_D$] | 140 | 83 | *{**7967**,16159,24351,32543, ...} |
| [$2^{*13}$**X**-$145_0$] | [$3^8$**X**-$115_{21}$] | [$1631_D$**X**-$30_D$] | 141 | 84 | *{**8047**,16239,24431,32623, ...} |
| [$2^{*13}$**X**-$89_0$] | [$3^8$**X**-$70_{21}$] | [$1631_D$**X**-$19_D$] | 142 | **85** | *{**8103**,16295,24487,32679, ...} |

**F(8,*13) = 85**/$2^{*13}$ = 0.0103759765625

**S(8,*13) = 0.9**552001953125

# Appendix Beta

The six Flowcharts showing the connections "**UP**" in Type 1 and Type 2

The first pair **of Charts** for *triple* **modulo** 3A = $2^1 3^2$ = 18

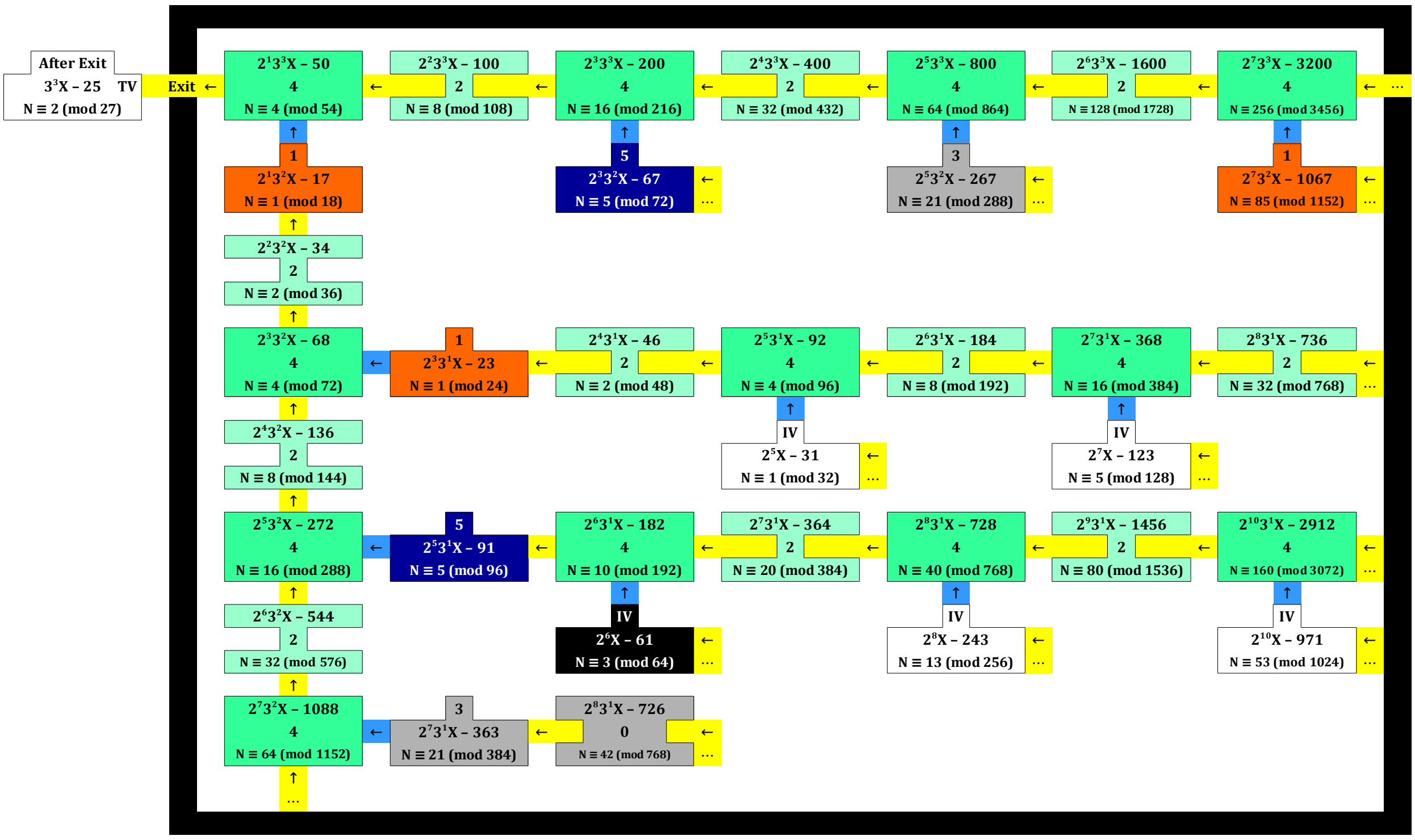

**Flowchart 3a.** Extrapolated from N ≡ 1 (mod $2^1 3^2$)
**OBS**! Relevant for **X** > 1 due to ***The Loop***.

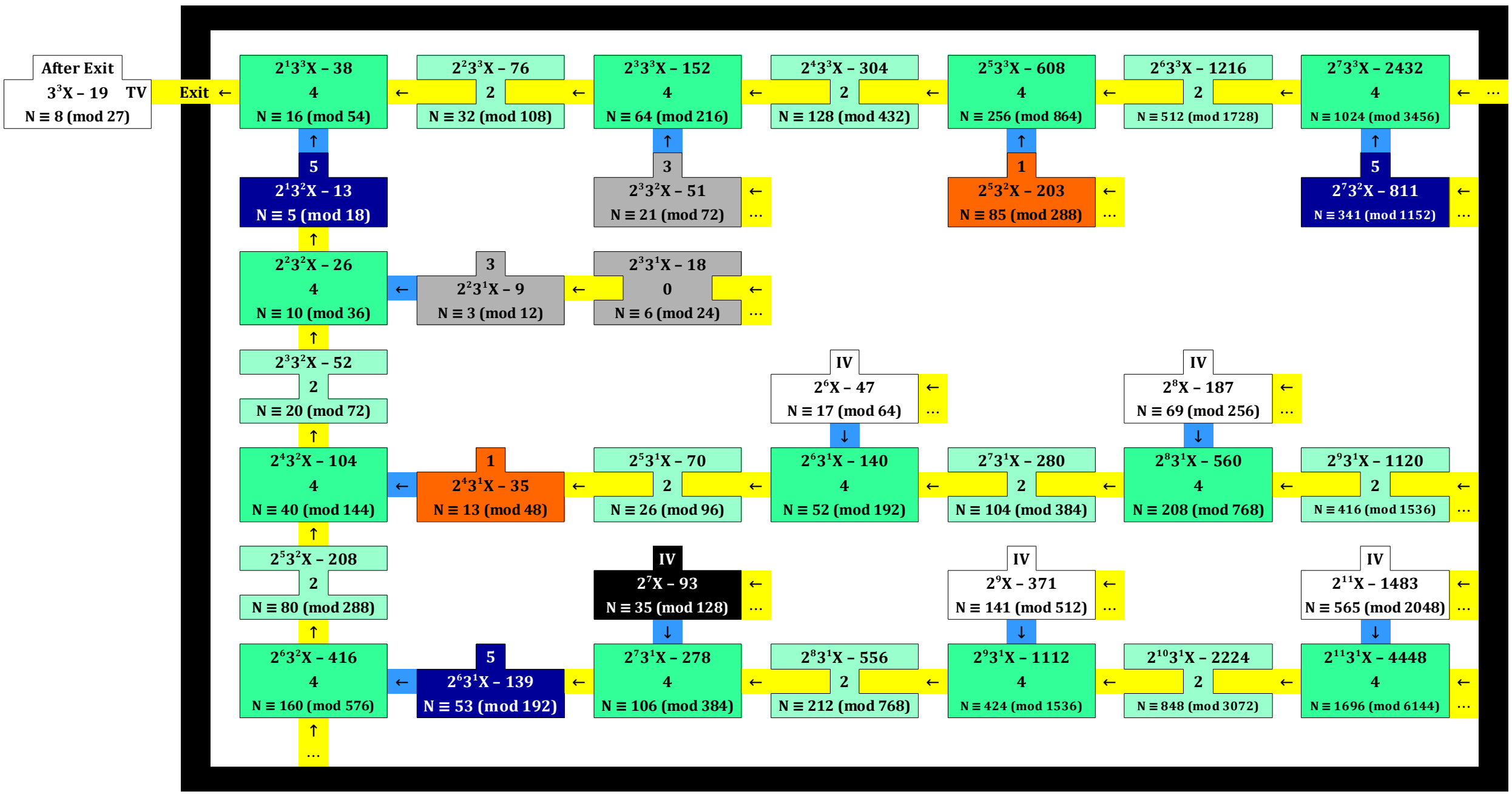

**Flowchart 3b.** Extrapolated from N ≡ 5 (mod $2^1 3^2$)

The second pair **of Charts** for *triple* **modulo** 3A = $2^1 3^2$ = 18

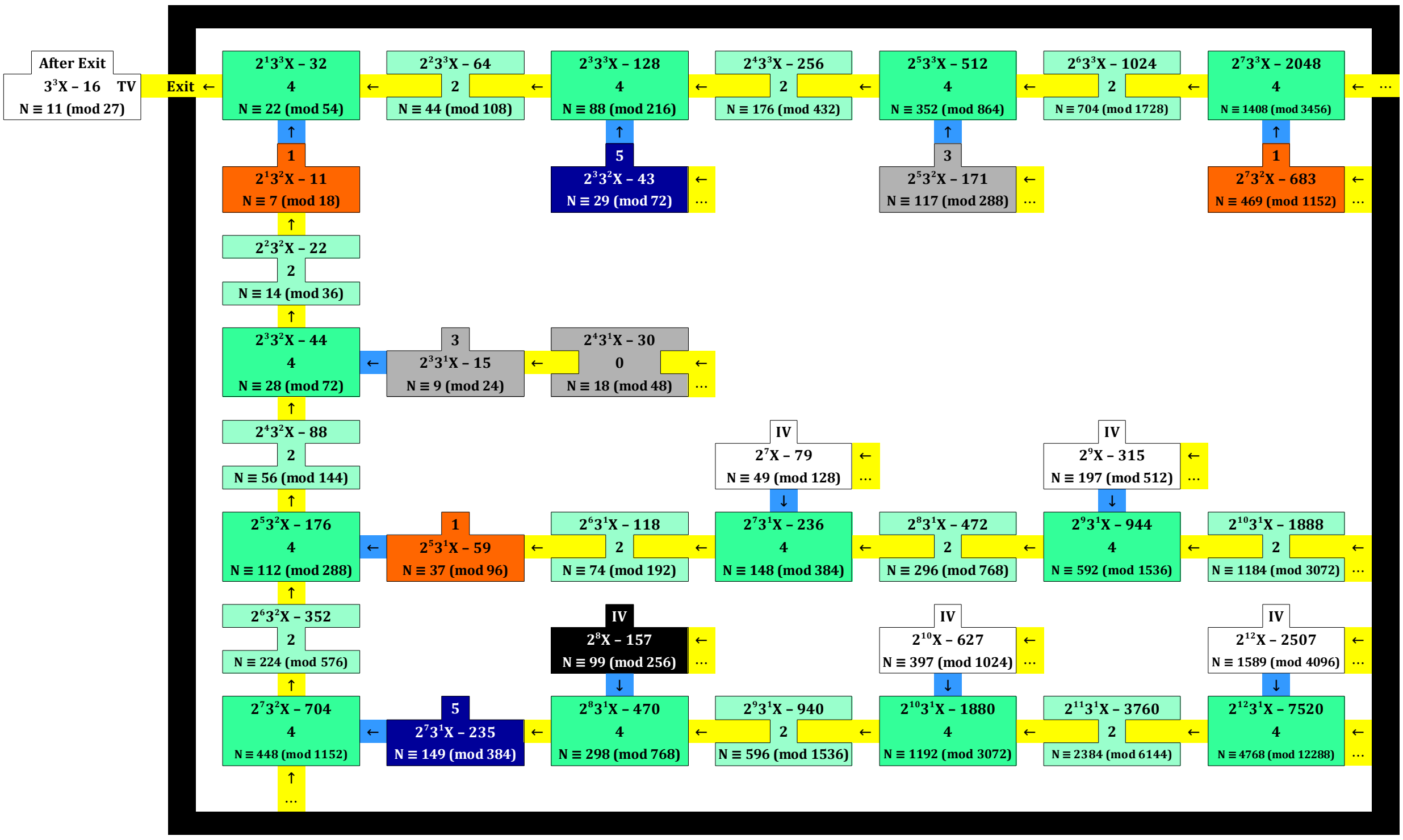

**Flowchart 3c.** Extrapolated from N ≡ 7 (mod $2^1 3^2$)

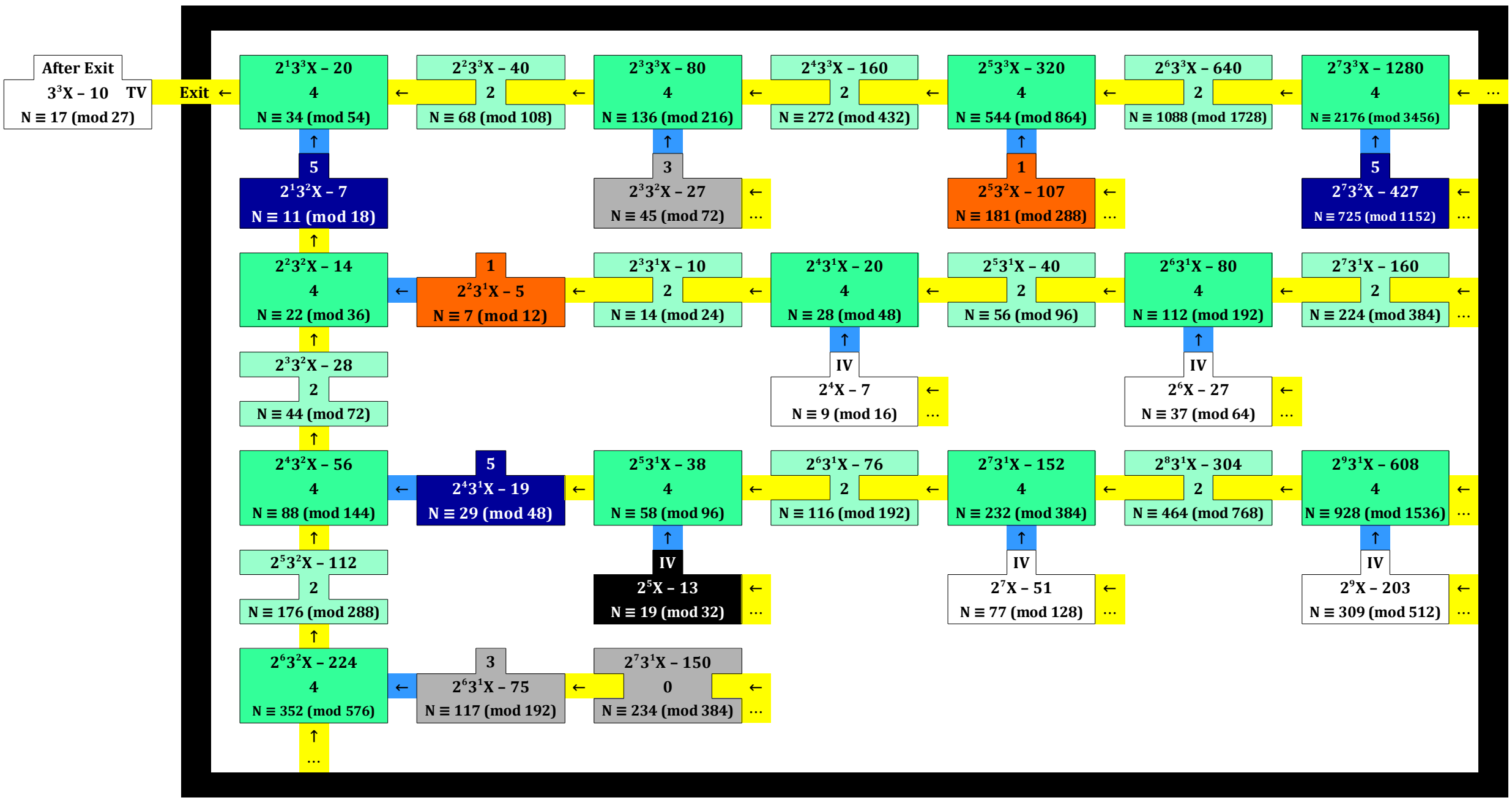

**Flowchart 3d.** Extrapolated from N ≡ 11 (mod $2^1 3^2$)

The third pair **of Charts** for *triple* **modulo** $3A = 2^1 3^2 = 18$

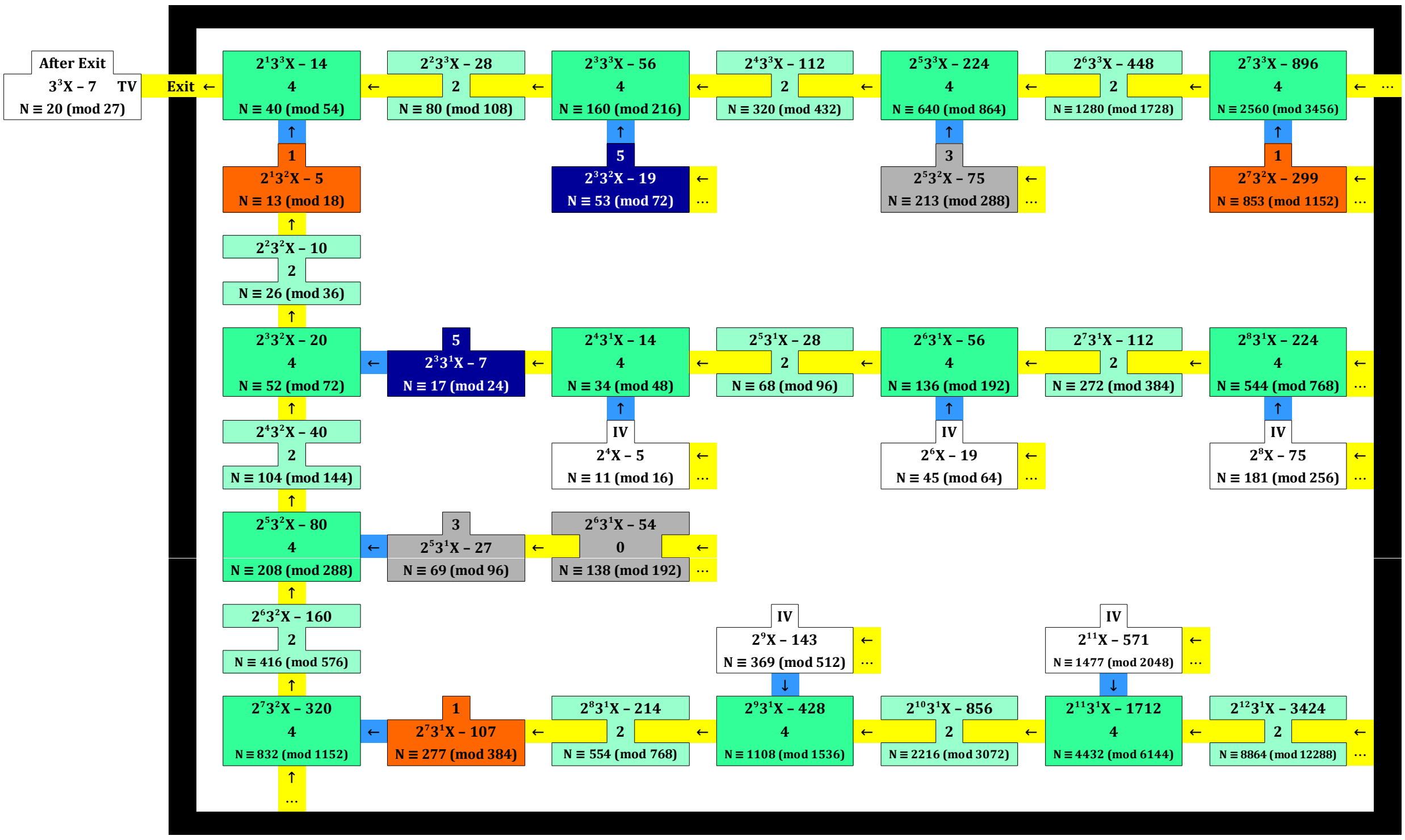

**Flowchart 3e**. Extrapolated from $N \equiv 13 \pmod{2^1 3^2}$

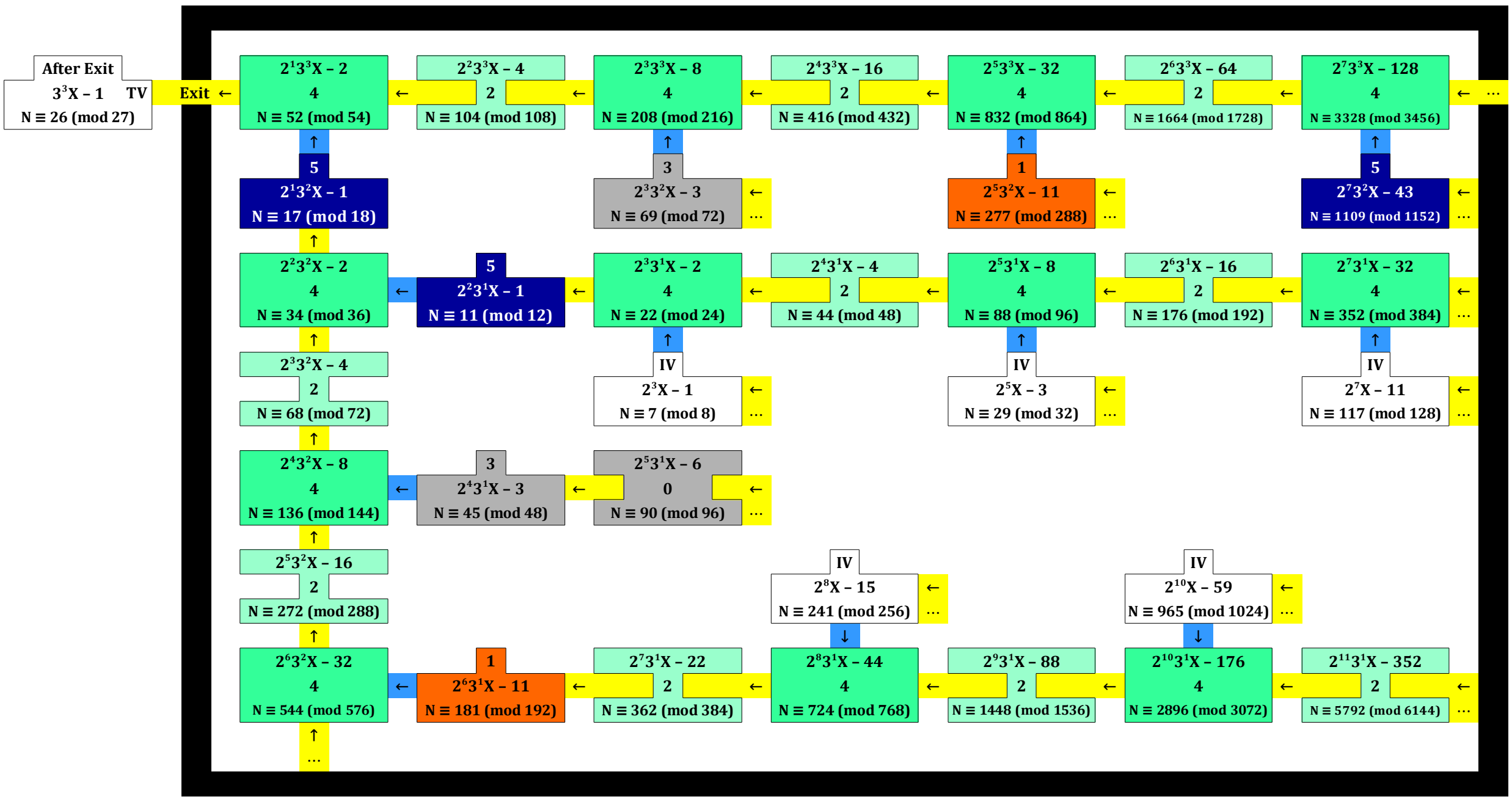

**Flowchart 3f.** Extrapolated from $N \equiv 17 \pmod{2^1 3^2}$

# Appendix Delta

The two (of four) **Flowcharts** extended with a Type Even operation

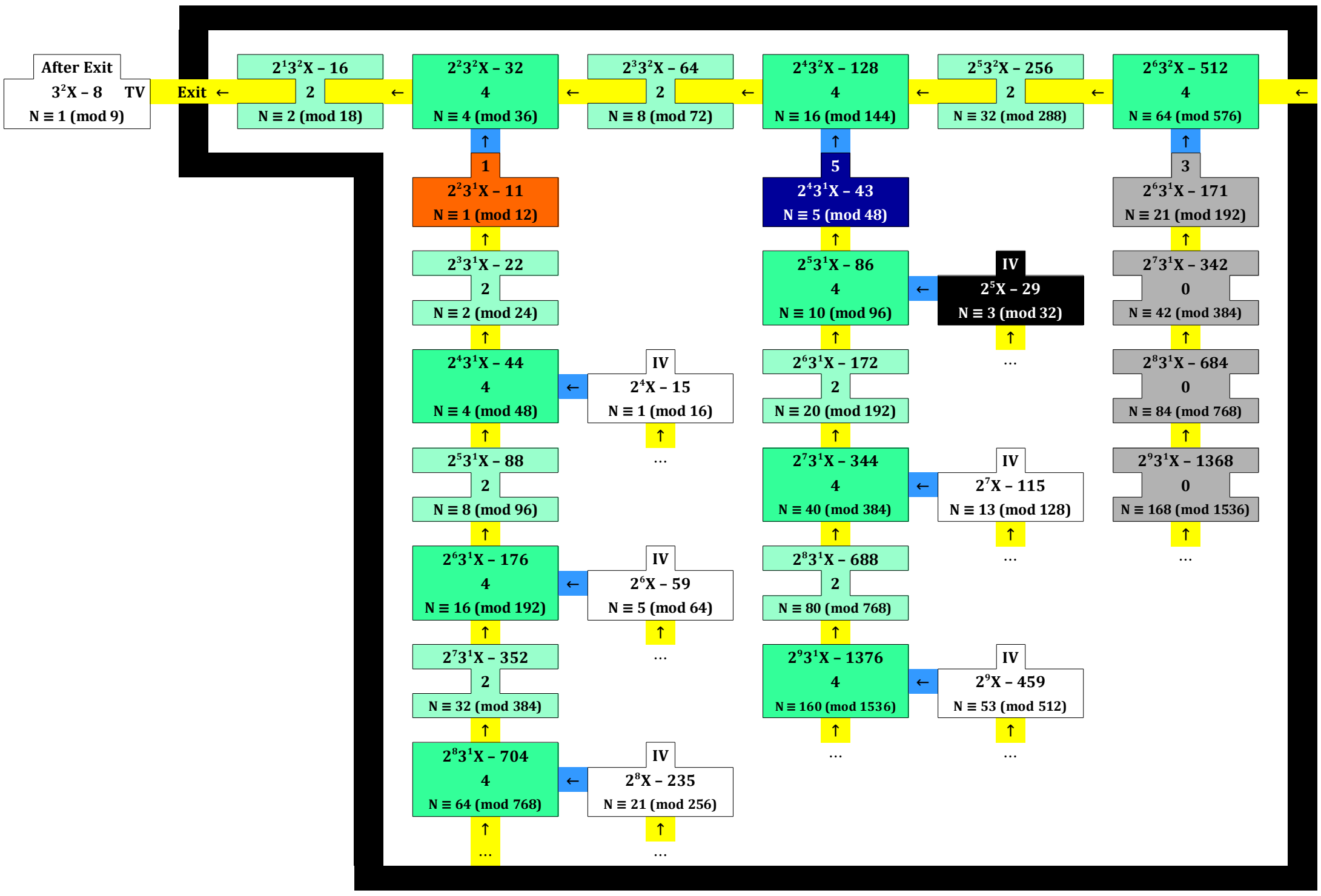

**Flowchart 4a.** Extrapolated from N ≡ 1 (mod $2^23^1$)

**OBS**! Relevant for **X** > 1 due to ***The Loop***.

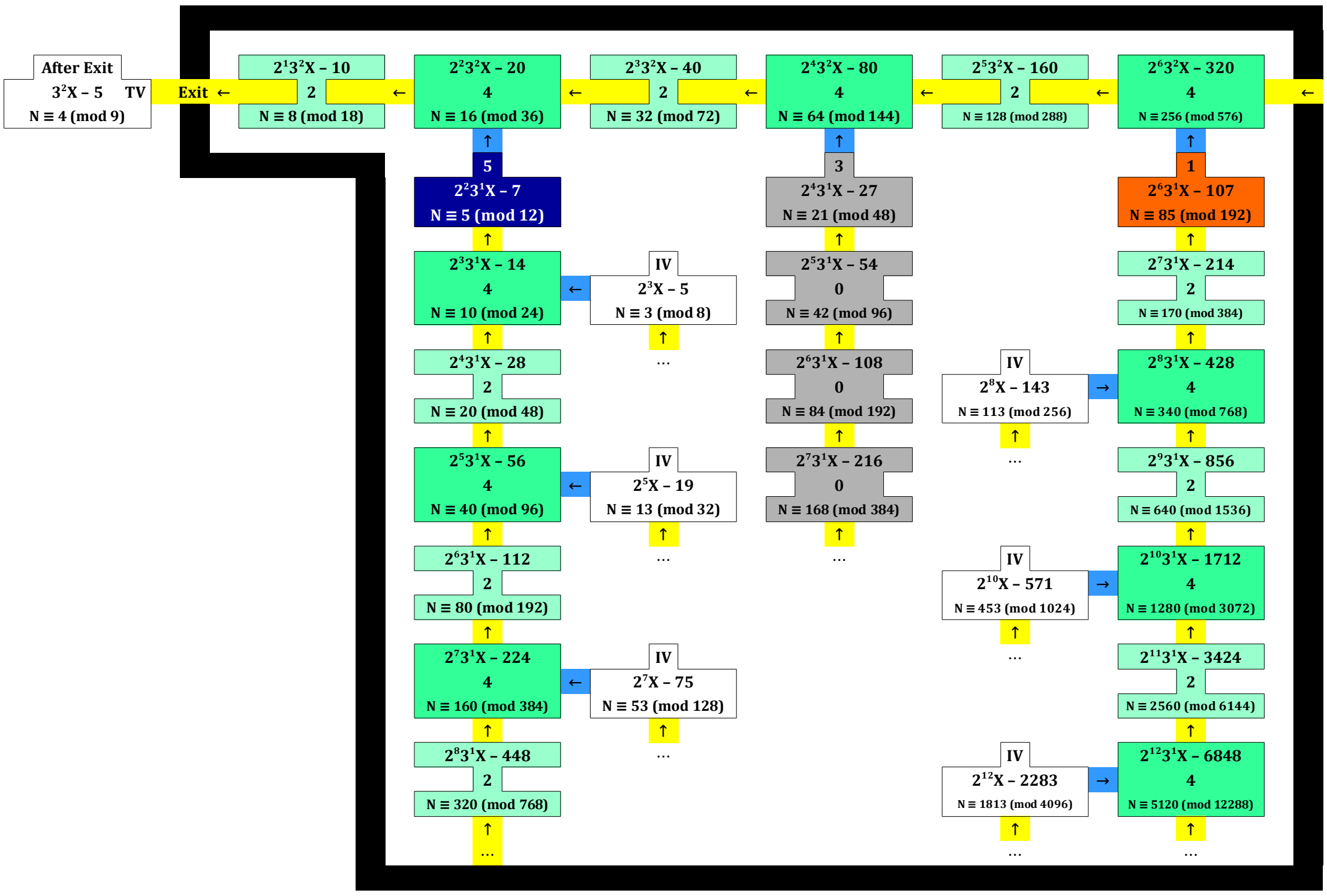

**Flowchart 4b.** Extrapolated from N ≡ 5 (mod $2^23^1$)

## The two (of four) **Flowcharts** extended with a Type Odd+Even operation

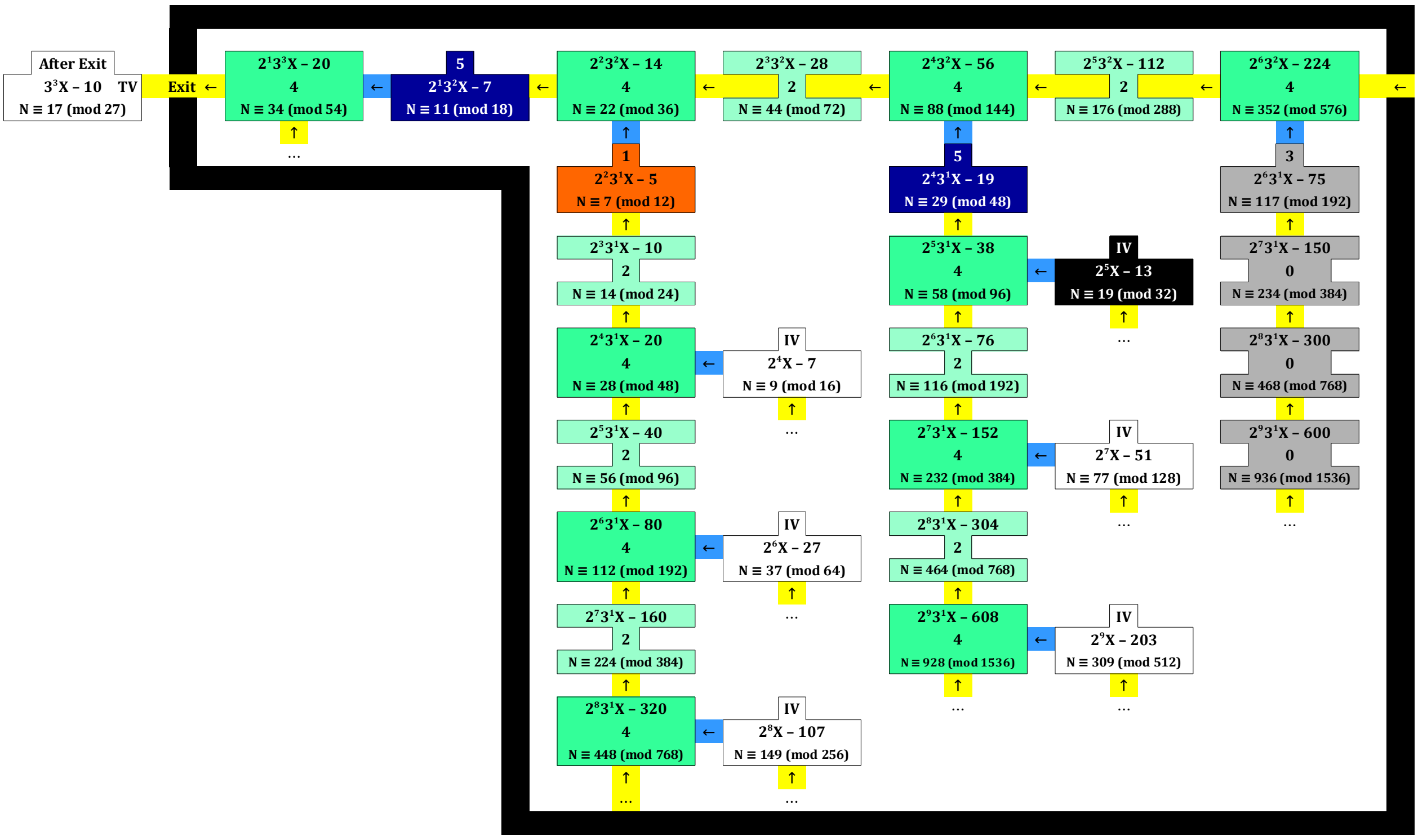

**Flowchart 4c.** Extrapolated from N ≡ 7 (mod $2^2 3^1$)

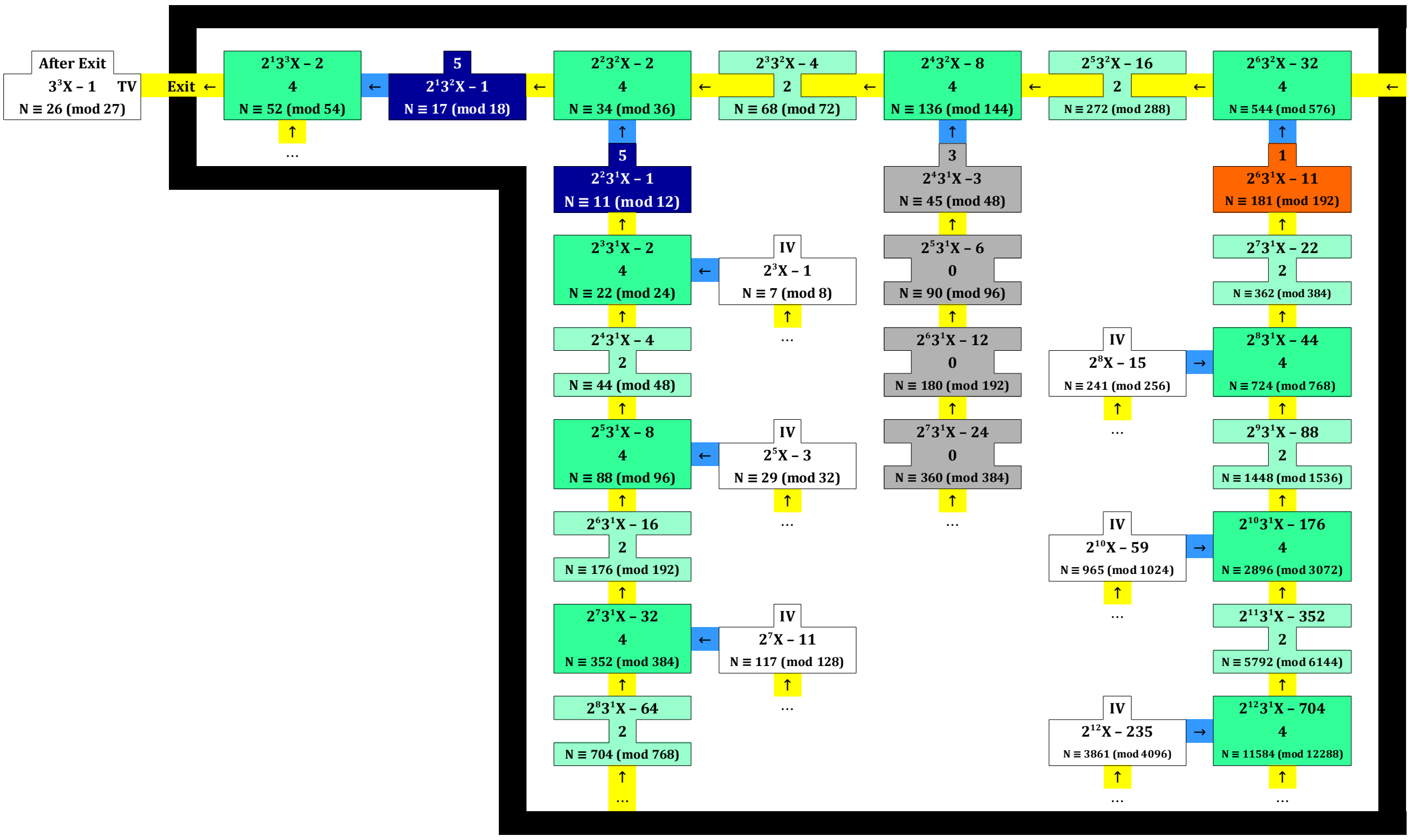

**Flowchart 4d.** Extrapolated from N ≡ 11 (mod $2^2 3^1$)

# Appendix Gamma

The two (of eight) **Flowcharts** extended with two Type Even operations

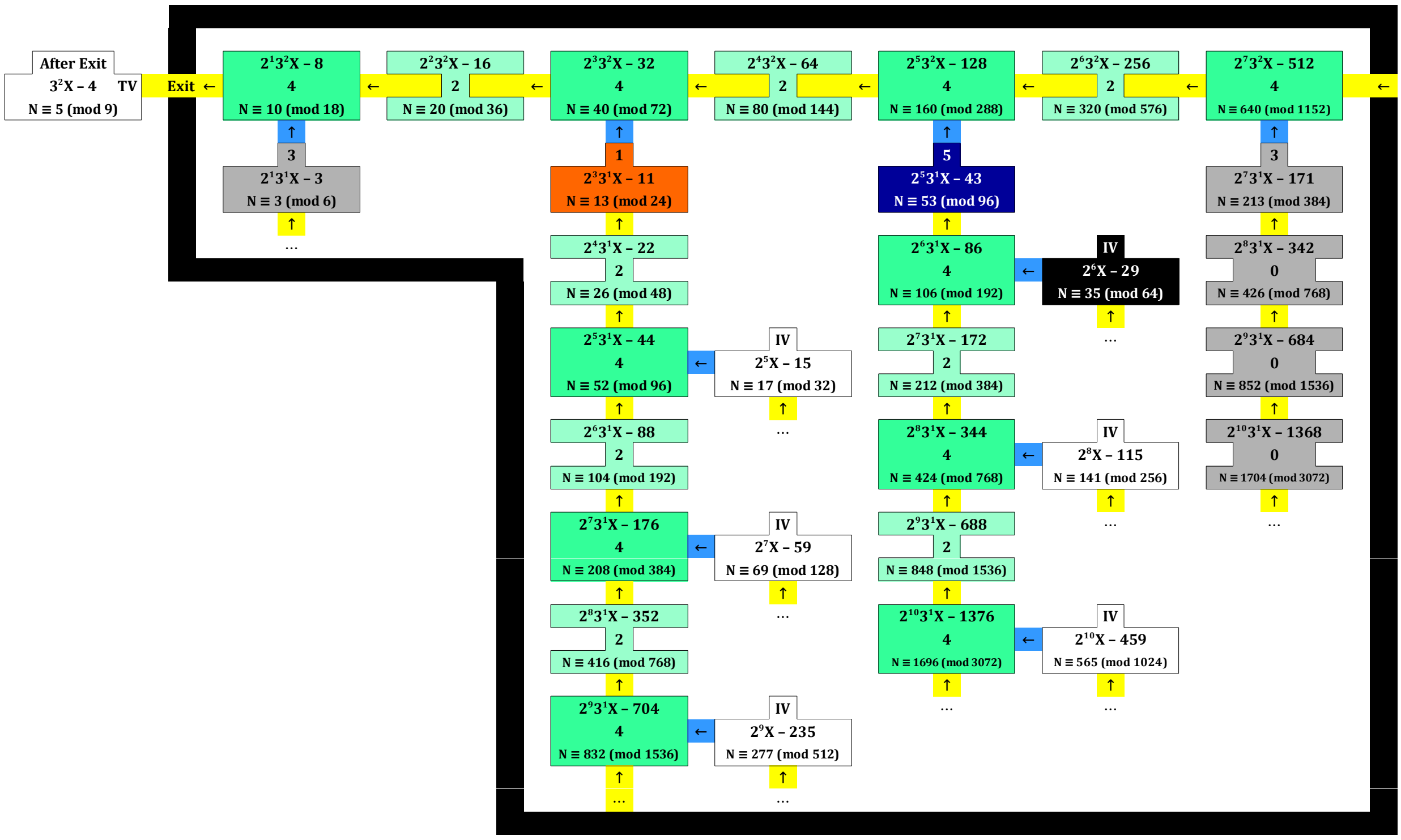

**Flowchart 5a**. Extrapolated from N ≡ 13 (mod $2^3 3^1$)

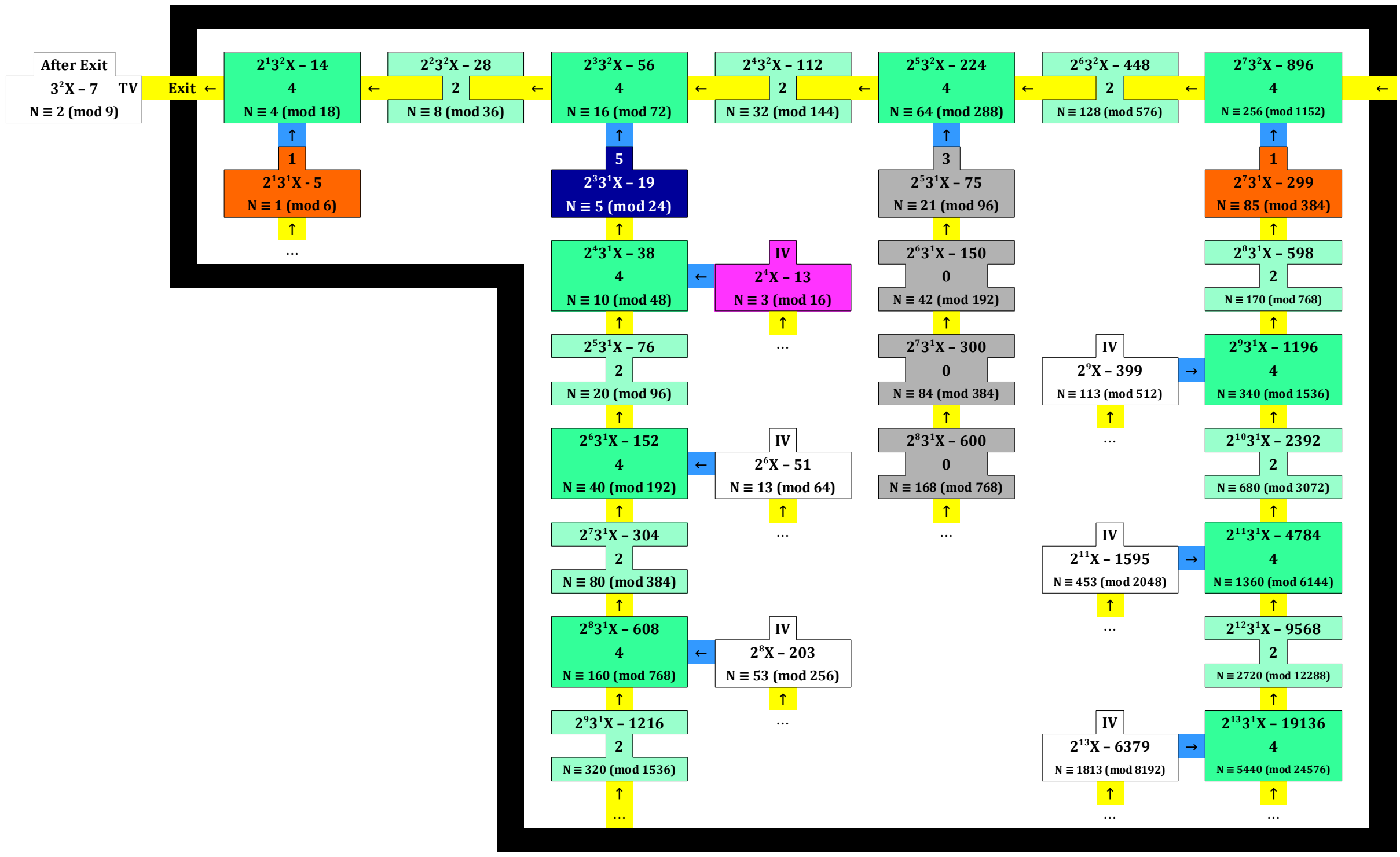

**Flowchart 5b**. Extrapolated from N ≡ 5 (mod $2^3 3^1$)

The two (of eight) **Flowcharts** extended with a Type Even and a Odd+Even operation

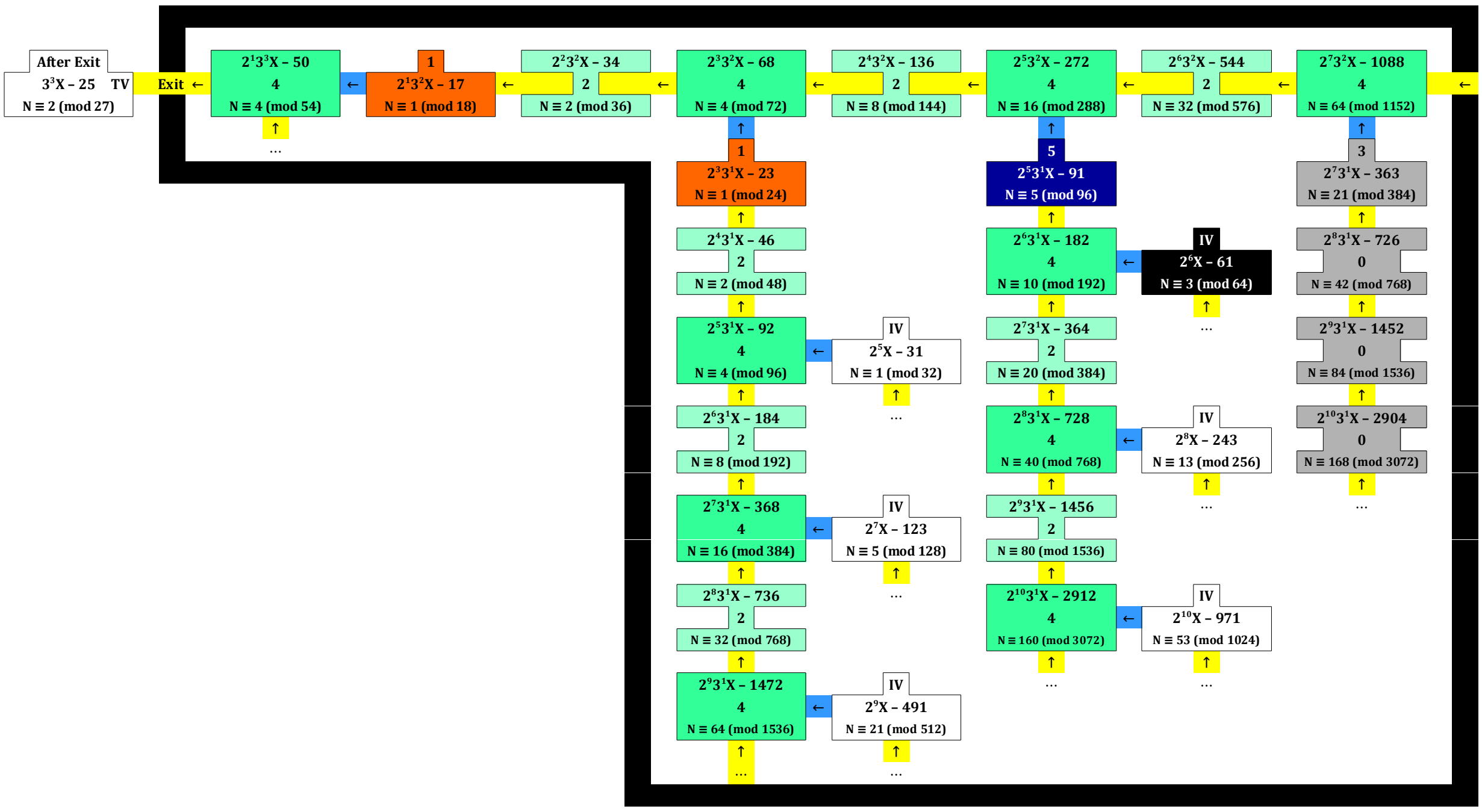

**Flowchart 5c**. Extrapolated from $N \equiv 1 \pmod{2^3 3^1}$
**OBS**! Relevant for **X** > 1 due to ***The Loop***.

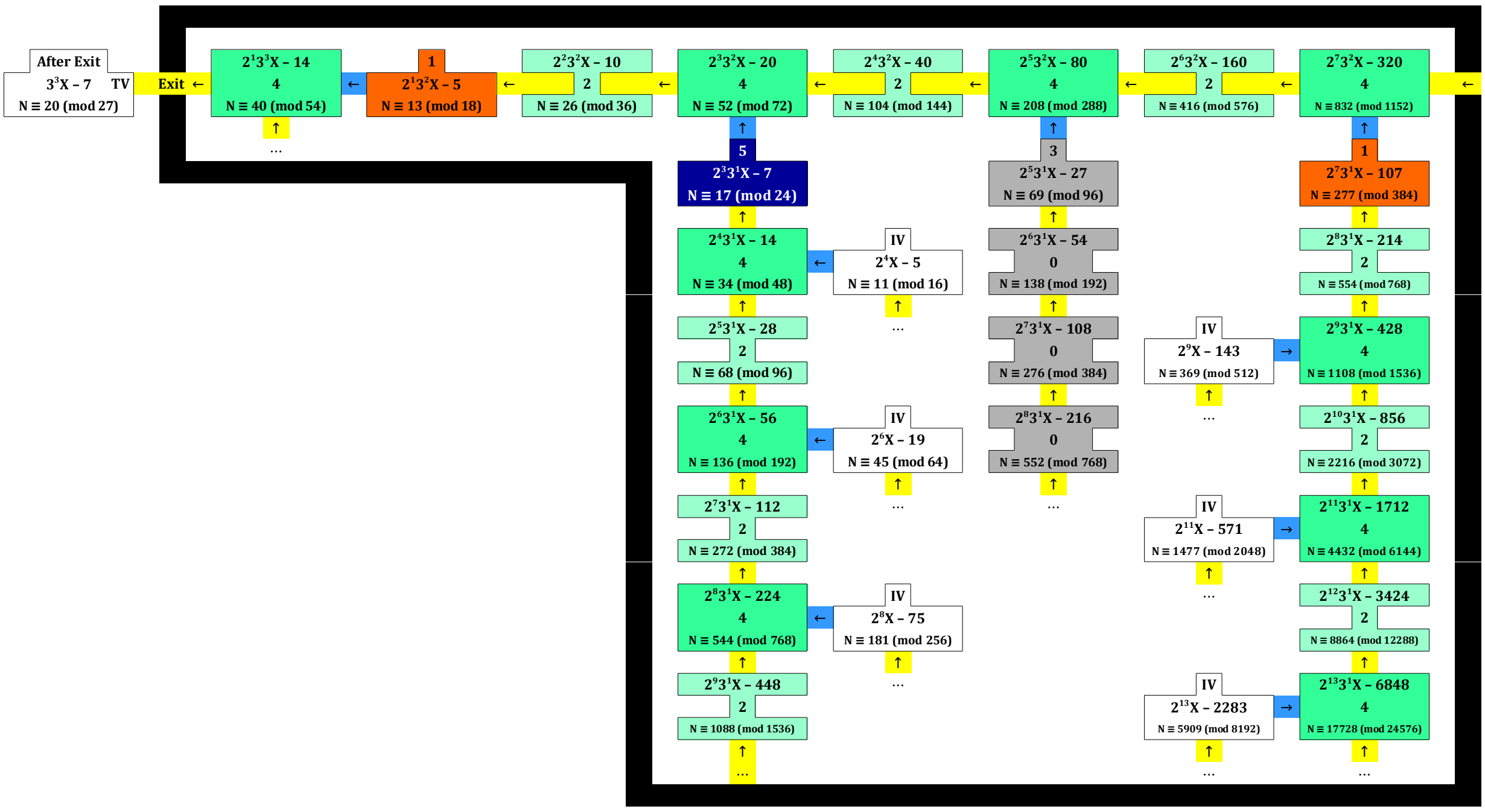

**Flowchart 5d**. Extrapolated from $N \equiv 17 \pmod{2^3 3^1}$

The two (of eight) **Flowcharts** extended with a Odd+Even operation and a Type Even

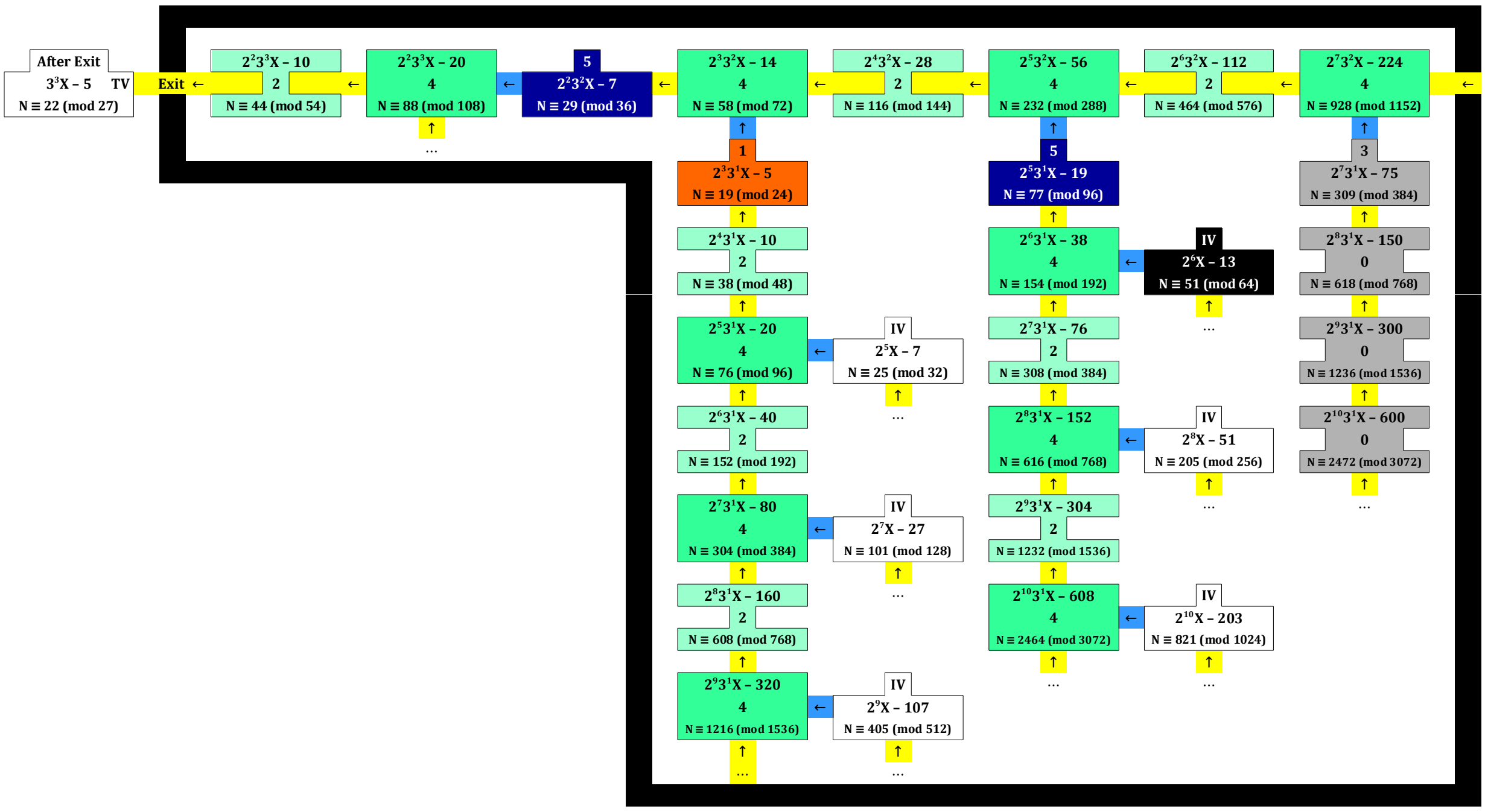

**Flowchart 5e.** Extrapolated from $N \equiv 19 \pmod{2^3 3^1}$

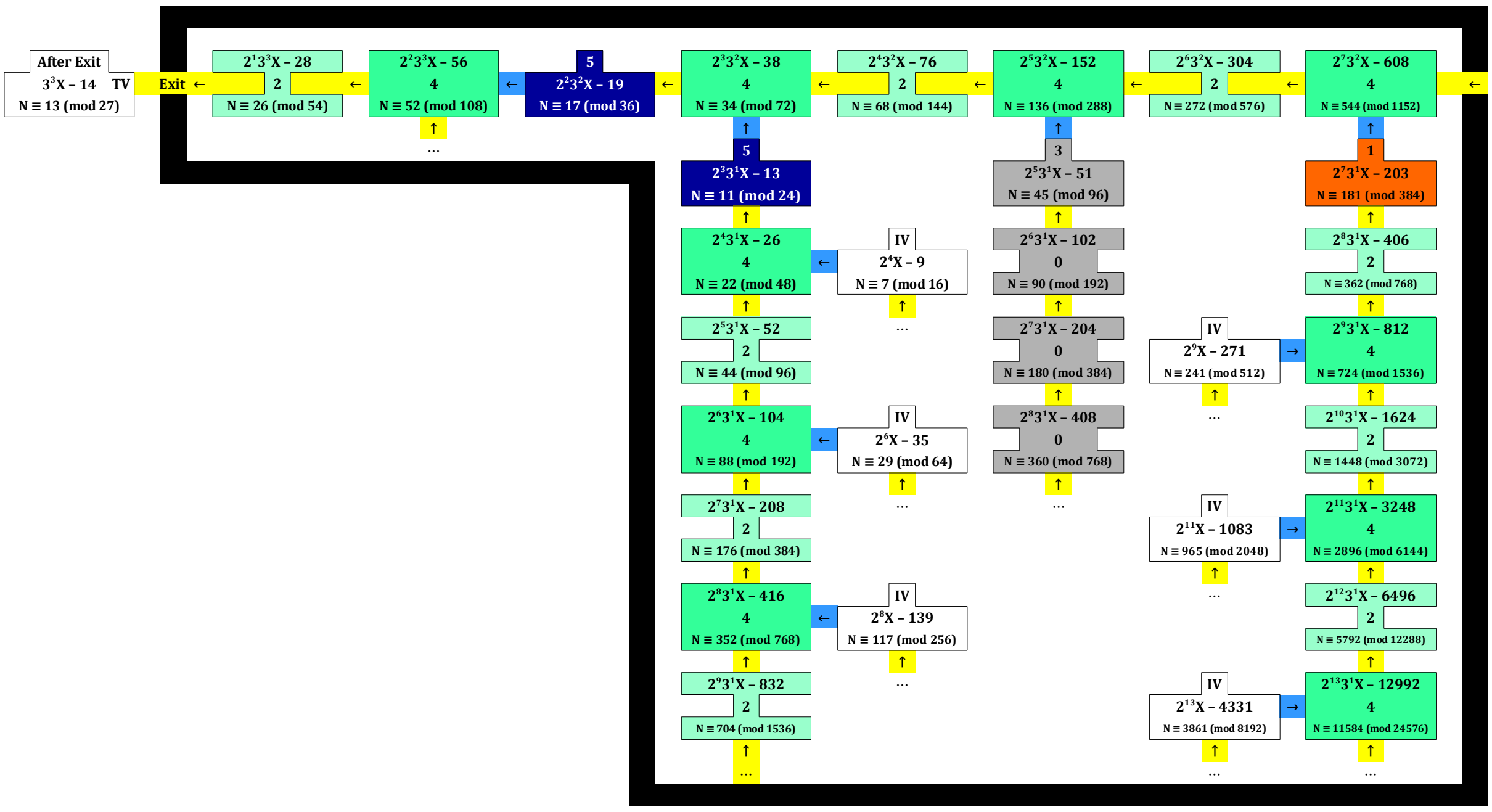

**Flowchart 5f.** Extrapolated from $N \equiv 11 \pmod{2^3 3^1}$

## The two (of eight) **Flowcharts** extended with **two** times Odd+Even operations

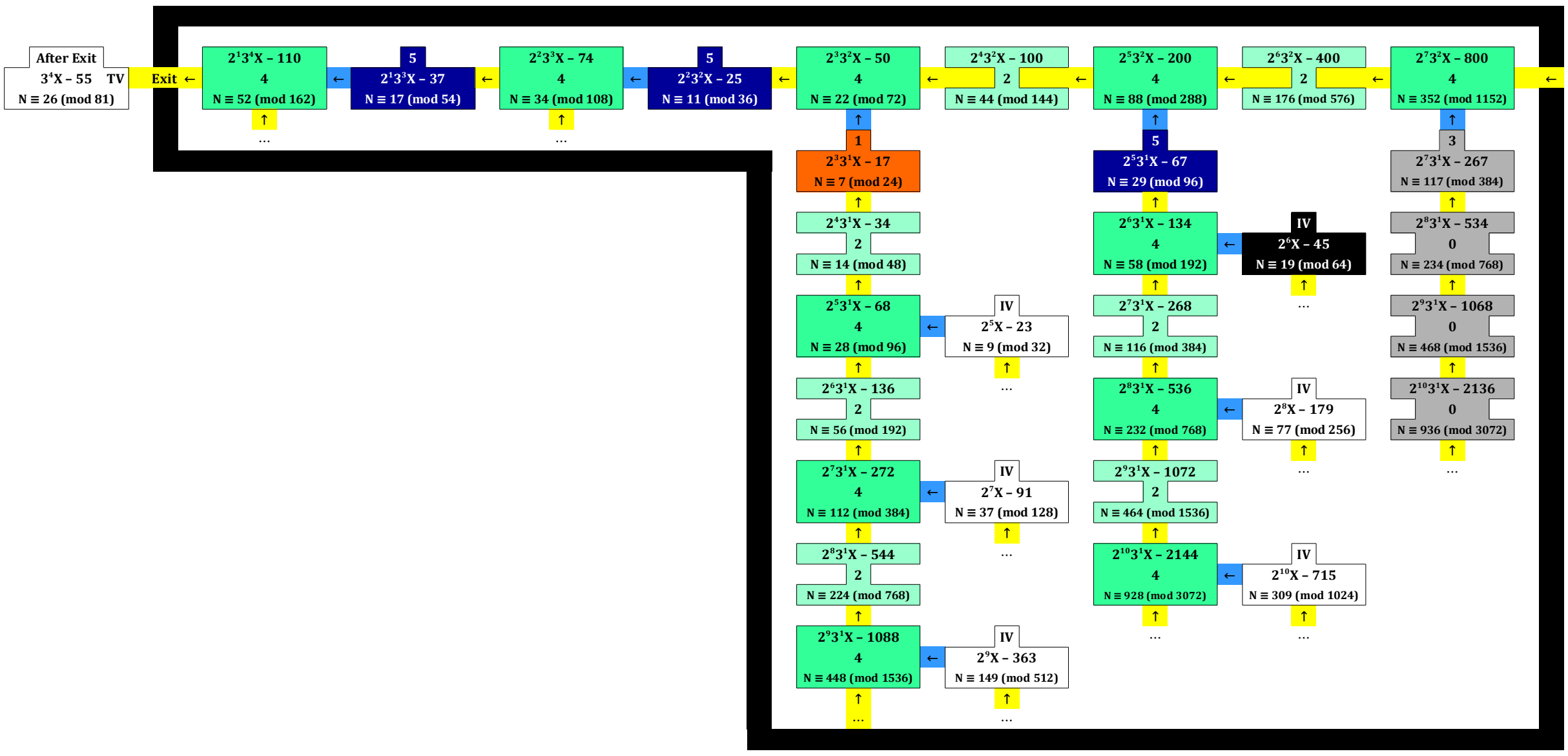

**Flowchart 5g**. Extrapolated from N ≡ 7 (mod $2^3 3^1$)

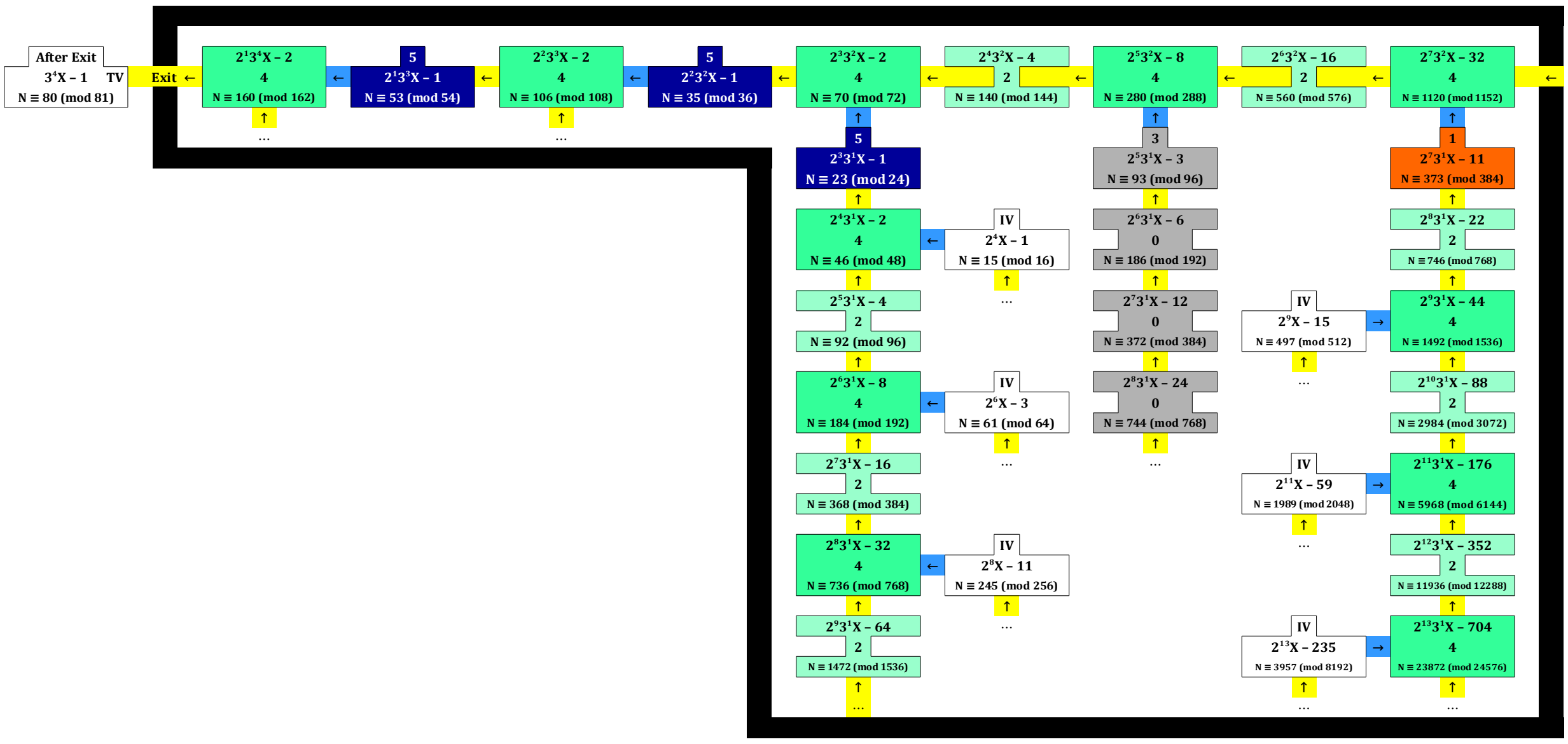

**Flowchart 5h**. Extrapolated from N ≡ 23 (mod $2^3 3^1$)

# Appendix Omega. *The Algorithm* with examples.

**Example for $N_0$ = 15 as Initial Value** (**NB**! The screen-shots are edited a bit due to lay-out):

Please input an integer Test-value N: 15

| | | | | | |
|---|---|---|---|---|---|
| R | | | | Value in Series: 15 | In binary: 1111 |
| O | #Odd : 1 | #Even : 0 | #Total : 1 | Value in Series: 46 | In binary: 101110 |
| E | #Odd : 1 | #Even : 1 | #Total : 2 | Value in Series: 23 | In binary: 10111 |
| O | #Odd : 2 | #Even : 1 | #Total : 3 | Value in Series: 70 | In binary: 1000110 |
| E | #Odd : 2 | #Even : 2 | #Total : 4 | Value in Series: 35 | In binary: 100011 |
| O | #Odd : 3 | #Even : 2 | #Total : 5 | Value in Series: 106 | In binary: 1101010 |
| E | #Odd : 3 | #Even : 3 | #Total : 6 | Value in Series: 53 | In binary: 110101 |
| O | #Odd : 4 | #Even : 3 | #Total : 7 | Value in Series: 160 | In binary: 10100000 |
| E | #Odd : 4 | #Even : 4 | #Total : 8 | Value in Series: 80 | In binary: 1010000 |
| E | #Odd : 4 | #Even : 5 | #Total : 9 | Value in Series: 40 | In binary: 101000 |
| E | #Odd : 4 | #Even : 6 | #Total : 10 | Value in Series: 20 | In binary: 10100 |
| E | #Odd : 4 | #Even : 7 | #Total : 11 | Value in Series: 10 | In binary: 1010 |
| T | **s = 4** | **r = 7** | **s+r = 11** | **Reduced Value: 10** | |

The Test-value 15 N in binary: 1111
is an Element in the *reducing **IV**-class *[*AX-*B] where:
*A = 2^r = 128 *A in binary: 10000000
*B = *A-*P = 113 *B in binary: 1110001
*P = 2^r-*B = 15 *P in binary: 0001111
The Initial Value-class [2^ 7 X - 113 ]

The Test-series for *[*AX-*B] Terminates at the **TV**-class [AX-B] where:
A = 3^s = 81 A in Binary: 1010001
B = A- P = 71 B in Binary: 1000111
P = 3^s- B = 10 P in Binary: 0001010
The **T**erminal**V**alue-class [3^ 4 X - 71 ]

The Test-value 15 represents the X'th value in
the *reducing IV-class *[2^ 7 X - 113 ]
where X = (N+*B)/2^r = 1
and the *Reduced Value 10 represents the same X in
the TV-class [3^ 4 X - 71 ]

The *reducing IV-class *[2^rX-*B], *P = (2^r-*B),
meaning 15 is congruent to 15 (mod 2^ 7 )
share the *Reducing Combination (s,r) = ( 4 , 7 )
i.e. the *reducing parity-sequence *ROEOE...ET is a combination of
s = 4 OE-tiles and (r-s) = 3 E-tiles

| | | | | | |
|---|---|---|---|---|---|
| E | #Odd : 4 | #Even : 8 | #Total : 12 | Value in Series: 5 | In binary: 101 |
| O | #Odd : 5 | #Even : 8 | #Total : 13 | Value in Series: 16 | In binary: 10000 |
| E | #Odd : 5 | #Even : 9 | #Total : 14 | Value in Series: 8 | In binary: 1000 |
| E | #Odd : 5 | #Even : 10 | #Total : 15 | Value in Series: 4 | In binary: 100 |
| E | #Odd : 5 | #Even : 11 | #Total : 16 | Value in Series: 2 | In binary: 10 |
| E | #Odd : 5 | #Even : 12 | #Total : 17 | Value in Series: 1 | In binary: 1 |

The Test-value 15 reaches *The End-value* N = 1 in
#Odd + #Even = 5 + 12 = 17 Total-operations,
and to reach **The Loop** the parity-sequence then have a
UNIQUE combination of 5 OE-tiles & 7 E-tiles
15 belongs to the **END** IV-class [2^dX - (2^d-N)] = [2^ 12 X - 4081 ]

The above example shows the full capacity of ***The Algorithm***. The *next* example do **not** show the *iterationseries* from the ****Reduced Value*** to the **End Value** N = **1** (print-lines deactivated).

Please notice that we find the **same** **reducing* **IV**-class as in the above example and notice in particular for the *iterationseries* that the ***ending*** of the binary strings are *identical*, as we find "enough middle zeros" in the strings, to the effect that the "head" and the "tail" in the strings **can not *interact*** before the *iterationseries* reaches the ****Reduced Value*** after 11 **Operations**.

**Example for $N_0 = (2^7)^2+15 = 2^7(2^7+1)-113_{Con} = 16399$:**

**OBS**! The "**entire** *head*" of the **string** ***left*** of the marked **000** is observed to ***begin*** at $A_{IV} = 2^7$ and ***end*** at $A_{TV} = 3^4$

Please input an integer Test-value N: 16399

| | | | | | | |
|---|---|---|---|---|---|---|
| R | | | | Value in Series: 16399 | In binary: | 100000000001111 |
| O | #Odd : 1 | #Even : 0 | #Total : 1 | Value in Series: 49198 | In binary: | 1 100000000101110 |
| E | #Odd : 1 | #Even : 1 | #Total : 2 | Value in Series: 24599 | In binary: | 1 10000000010111 |
| O | #Odd : 2 | #Even : 1 | #Total : 3 | Value in Series: 73798 | In binary: | 100 10000001000110 |
| E | #Odd : 2 | #Even : 2 | #Total : 4 | Value in Series: 36899 | In binary: | 100 1000000100011 |
| O | #Odd : 3 | #Even : 2 | #Total : 5 | Value in Series: 110698 | In binary: | 1101 1000001101010 |
| E | #Odd : 3 | #Even : 3 | #Total : 6 | Value in Series: 55349 | In binary: | 1101 100000110101 |
| O | #Odd : 4 | #Even : 3 | #Total : 7 | Value in Series: 166048 | In binary: | 101000 100010100000 |
| E | #Odd : 4 | #Even : 4 | #Total : 8 | Value in Series: 83024 | In binary: | 101000 10001010000 |
| E | #Odd : 4 | #Even : 5 | #Total : 9 | Value in Series: 41512 | In binary: | 101000 1000101000 |
| E | #Odd : 4 | #Even : 6 | #Total : 10 | Value in Series: 20756 | In binary: | 101000 100010100 |
| E | #Odd : 4 | #Even : 7 | #Total : 11 | Value in Series: 10378 | In binary: | 101000 10001010 |
| T | **s = 4** | **r = 7** | **s+r = 11** | **Reduced Value: 10378** | | |

The Test-value 16399 N in binary: 100000000001111
is an Element in the *reducing IV-class *[*AX-*B] where:
*A = 2^r = 128 *A in binary: 10000000
*B = *A-*P = 113 *B in binary: 1110001
*P = 2^r-*B = 15 *P in binary: 0001111
The InitialValue-class [2^ 7 X - 113 ]

The Test-series for *[*AX-*B] Terminates at the TV-class [AX-B] where:
A = 3^s = 81 A in Binary: 1010001
B = A- P = 71 B in Binary: 1000111
P = 3^s- B = 10 P in Binary: 0001010
The TerminalValue-class [3^ 4 X - 71 ]

The Test-value 16399 represents the X'th value in
the *reducing IV-class *[2^ 7 X - 113 ]
where X = (N+*B)/2^r = 129
and the *Reduced Value 10378 represents the same X in
the TV-class [3^ 4 X - 71 ]

The *reducing IV-class *[2^rX-*B], *P = (2^r-*B),
meaning 16399 is congruent to 15 (mod 2^ 7 )
share the *Reducing Combination (s,r) = ( 4 , 7 )
i.e. the *reducing parity-sequence *ROEOE...ET is a combination of
s = 4 OE-tiles and (r-s) = 3 E-tiles

The Test-value 16399 reaches *The End-value* N = 1 in
#Odd + #Even = 39 + 76 = 115 Total-operations,
and to reach **The Loop** the parity-sequence then have a
UNIQUE combination of 39 OE-tiles & 37 E-tiles
16399 belongs to the **END** IV-class [2^dX - (2^d-N)] = [2^ 76 X - 75557863725914323402737 ]

## Examples of results from the two Algorithms

The below table contains all values from $^{\Phi}[2^5\mathbf{X}\text{-}1_0]_{IV} \cup {}^{\Phi}[2^5\mathbf{X}\text{-}5_0]_{IV} \cup {}^{\Phi}[2^5\mathbf{X}\text{-}17_0]_{IV} \cup {}^{\Phi}[2^5\mathbf{X}\text{-}25_0]_{IV}$

For $N_X$ is shown the **Reducing Combination* (**s**,***r**$_s$), the *$B_{Con}$ from the **Converging* **IV**-class $^*{}_u[2^{*r_s}\mathbf{X}\text{-}^*B_{Con}]_q$ & the **End Combination** (**s**$_d$,**d**) from the **End IV**-class $^{\dagger}[2^d\mathbf{X}\text{-}(2^d\text{-}N_N)_{End}]_{IV}$, $N_X = N_N$

For $B_0$ is shown the **Reducing Combination* (**s**,***r**$_s$), **P** = ($2^{*r_s}$-*$B_{Con}$) from the **Converging* **IV**-class $^*{}_u[2^{*r_s}\mathbf{X}\text{-}^*B_{Con}]_q$ & the relevant **loop** Combination (**s**,**r**) = (# **O** ,# **E** ) for $B_{j=s+r}$ ∈ {1,5,17}.

| | | | | | | | | | | | | | | | | |
|---|---|---|---|---|---|---|---|---|---|---|---|---|---|---|---|---|
| (s_d,d) | (39_d,67) | (39_d,68) | (38_d,67) | (15_d,31) | (18_d,36) | (14_d,30) | (24_d,46) | (15_d,32) | (13_d,29) | (18_d,37) | (28_d,53) | (14_d,31) | (48_d,85) | (34_d,63) | (17_d,36) | (20_d,41) |
| B = 2^*r - **P** | 2^56-31 | 2^54-63 | 161 | 32641 | 2^21-159 | 8001 | 2^31-223 | 7937 | 737 | 2^21-319 | | 3713 | 32353 | 2^40-447 | 65057 | 2^18-511 |
| (s,*r) | (35,*56) | (34,*54) | (5,*8) | (9,*15) | (13,*21) | (8,*13) | (19,*31) | (8,*13) | (6,*10) | (13,*21) | (5,*8) | (7,*12) | (9,*15) | (25,*40) | (10,*16) | (11,*18) |
| N ≡ 31 (mod 32) | 31 | 63 | 95 | 127 | 159 | 191 | 223 | 255 | 287 | 319 | 351 | 383 | 415 | 447 | 479 | 511 |
| **X** | **1** | **2** | **3** | **4** | **5** | **6** | **7** | **8** | **9** | **10** | **11** | **12** | **13** | **14** | **15** | **16** |
| B ≡ 1 (mod 32) | 1 | 33 | 65 | 97 | 129 | 161 | 193 | 225 | 257 | 289 | 321 | 353 | 385 | 417 | 449 | 481 |
| For B | B-loop | (7,*12) | (11,*18) | (10,*16) | (10,*16) | (5,*8) | (9,*15) | (8,*13) | (23,*37) | (6,*10) | (31,*50) | (10,*16) | (22,*35) | (5,*8) | (6,*10) | (29,*46) |
| **P** = 2^*r - *B_Con | N. A. | 4063 | 262079 | 65439 | 65407 | 95 | 32575 | 7967 | 2^37-257 | 735 | 2^50-321 | 65183 | 2^35-385 | | 575 | 2^46-481 |
| B-loop | | 17 (7,12) | 1 (12,25) | 1 (11,24) | 1 (12,26) | 5 (7,16) | 1 (11,25) | 17 (16,29) | 17 (31,53) | 1 (10,24) | 17 (39,66) | 5 (24,44) | 17 (30,52) | 17 (11,22) | 17 (16,30) | 17 (38,65) |
| | | | | | | | | | | | | | | | | |
| (s_d,d) | (41_d,70) | (10_d,22) | (33_d,59) | (15_d,31) | (30_d,55) | (14_d,30) | (17_d,35) | (22_d,43) | (20_d,40) | (11_d,26) | (45_d,80) | (19_d,39) | (48_d,85) | (17_d,36) | (7_d,20) | (10_d,25) |
| B = 2^*r - **P** | 2^59-27 | 69 | 2^45-91 | 133 | 2^40-155 | | 37 | 2^27-251 | 2^24-283 | | 677 | | 32357 | | | 517 |
| (s,*r) | (37,*59) | (4,*7) | (28,*45) | (5,*8) | (25,*40) | (4,*7) | (5,*8) | (17,*27) | (15,*24) | (4,*7) | (6,*10) | (5,*8) | (9,*15) | (4,*7) | (5,*8) | (6,*10) |
| N ≡ 27 (mod 32) | 27 | 59 | 91 | 123 | 155 | 187 | 219 | 251 | 283 | 315 | 347 | 379 | 411 | 443 | 475 | 507 |
| **X** | **1** | **2** | **3** | **4** | **5** | **6** | **7** | **8** | **9** | **10** | **11** | **12** | **13** | **14** | **15** | **16** |
| B ≡ 5 (mod 32) | 5 | 37 | 69 | 101 | 133 | 165 | 197 | 229 | 261 | 293 | 325 | 357 | 389 | 421 | 453 | 485 |
| For B | B-loop | (5,*8) | (4,*7) | (6,*10) | (5,*8) | (9,*15) | (4,*7) | (25,*40) | (10,*16) | (5,*8) | (4,*7) | (10,*16) | (5,*8) | (7,*12) | (4,*7) | (13,*21) |
| **P** = 2^*r - *B_Con | N. A. | 219 | 59 | 923 | 123 | 32603 | | 2^40-229 | 65275 | | | 65179 | | 3675 | | 2^21-485 |
| B-loop | | 17 (5,9) | 1 (5,14) | 1 (11,24) | 5 (11,22) | 17 (15,27) | 17 (6,13) | 17 (33,56) | 17 (19,34) | 17 (17,31) | 1 (8,21) | 5 (10,22) | 17 (6,14) | 5 (8,19) | 5 (6,16) | 1 (19,39) |
| | | | | | | | | | | | | | | | | |
| (s_d,d) | (5_d,12) | (38_d,66) | (11_d,24) | (24_d,45) | (37_d,66) | (28_d,52) | (31_d,57) | (17_d,35) | (13_d,29) | (13_d,29) | (23_d,45) | (14_d,31) | (43_d,77) | (12_d,28) | (46_d,82) | (34_d,63) |
| B = 2^*r - **P** | 113 | 2^54-47 | 177 | 2^31-111 | | 81 | 7985 | 2^20-239 | | 7889 | | 657 | | | 3633 | 2^27-495 |
| (s,*r) | (4,*7) | (34,*54) | (5,*8) | (19,*31) | (4,*7) | (5,*8) | (8,*13) | (12,*20) | (4,*7) | (8,*13) | (5,*8) | (6,*10) | (4,*7) | (5,*8) | (7,*12) | (17,*27) |
| N ≡ 15 (mod 32) | 15 | 47 | 79 | 111 | 143 | 175 | 207 | 239 | 271 | 303 | 335 | 367 | 399 | 431 | 463 | 495 |
| **X** | **1** | **2** | **3** | **4** | **5** | **6** | **7** | **8** | **9** | **10** | **11** | **12** | **13** | **14** | **15** | **16** |
| B ≡ 17 (mod 32) | 17 | 49 | 81 | 113 | 145 | 177 | 209 | 241 | 273 | 305 | 337 | 369 | 401 | 433 | 465 | 497 |
| For B | B-loop | (6,*10) | (5,*8) | (4,*7) | (8,*13) | (5,*8) | (6,*10) | (4,*7) | (13,*21) | (26,*42) | (5,*8) | (4,*7) | (8,*13) | (5,*8) | (11,*18) | (4,*7) |
| **P** = 2^*r - *B_Con | N. A. | 975 | 175 | 15 | 8047 | 79 | 815 | | 2^21-273 | 2^42-305 | | | 7791 | | 261679 | |
| B-loop | | 17 (6,11) | 5 (7,15) | 1 (9,21) | 1 (10,23) | 5 (12,24) | 17 (11,21) | 5 (6,15) | 1 (17,35) | 17 (34,58) | 17 (15,28) | 17 (13,25) | 5 (8,19) | 1 (9,23) | 17 (21,38) | 1 (7,20) |
| | | | | | | | | | | | | | | | | |
| (s_d,d) | (5_d,11) | (11_d,23) | (37_d,65) | (31_d,56) | (13_d,28) | (23_d,44) | (43_d,76) | (46_d,81) | (27_d,51) | (18_d,37) | (52_d,91) | (16_d,34) | (43_d,77) | (9_d,23) | (12_d,28) | (51_d,90) |
| B = 2^*r - **P** | 121 | 217 | 2^51-71 | 2^42-103 | | 2^29-167 | 57 | 3865 | | | 2^21-327 | 65177 | | 601 | | 2^20-487 |
| (s,*r) | (4,*7) | (5,*8) | (32,*51) | (26,*42) | (4,*7) | (18,*29) | (5,*8) | (7,*12) | (4,*7) | (5,*8) | (13,*21) | (10,*16) | (4,*7) | (6,*10) | (5,*8) | (12,*20) |
| N ≡ 7 (mod 32) | 7 | 39 | 71 | 103 | 135 | 167 | 199 | 231 | 263 | 295 | 327 | 359 | 391 | 423 | 455 | 487 |
| **X** | **1** | **2** | **3** | **4** | **5** | **6** | **7** | **8** | **9** | **10** | **11** | **12** | **13** | **14** | **15** | **16** |
| B ≡ 25 (mod 32) | 25 | 57 | 89 | 121 | 153 | 185 | 217 | 249 | 281 | 313 | 345 | 377 | 409 | 441 | 473 | 505 |
| For B | (6,*10) | (5,*8) | (8,*13) | (4,*7) | (27,*43) | (7,*12) | (5,*8) | (4,*7) | (8,*13) | (5,*8) | (8,*13) | (4,*7) | (10,*16) | (6,*10) | (5,*8) | (4,*7) |
| **P** = 2^*r - *B_Con | 999 | 199 | 8103 | 7 | 2^43-153 | 3911 | 39 | | 7911 | | 7847 | | 35127 | 583 | | |
| B-loop | 17 (6,10) | 1 (9,20) | 5 (12,23) | 5 (6,14) | 17 (34,57) | 17 (13,24) | 1 (9,22) | 1 (7,19) | 5 (9,20) | 17 (1020) | 1 (13,29) | 5 (17,33) | 1 (16,34) | 17 (16,30) | 5 (6,16) | 17 (14,27) |
| (s_d,d) | (49_d,87) | (13_d,30) | (13_d,30) | (47_d,84) | (33_d,62) | (62_d,108) | (14_d,32) | (19_d,40) | (24_d,48) | (48_d,86) | (12_d,29) | (34_d,64) | (41_d,75) | (46_d,83) | (34_d,64) | (20_d,42) |
| B = 2^*r - **P** | 7649 | 449 | | 2^23-639 | 2^39-671 | 2^81-703 | 289 | 64769 | 7393 | 31937 | | 2^24-895 | 2^37-927 | 2^21-959 | 2^26-991 | 2^18-1023 |
| (s,*r) | (8,*13) | (6,*10) | (5,*8) | (14,*23) | (24,*39) | (51,*81) | (6,*10) | (10,*16) | (8,*13) | (9,*15) | (5,*8) | (15,*24) | (23,*37) | (13,*21) | (16,*26) | (11,*18) |
| N ≡ 31 (mod 32) | 543 | 575 | 607 | 639 | 671 | 703 | 735 | 767 | 799 | 831 | 863 | 895 | 927 | 959 | 991 | 1023 |
| **X** | **17** | **18** | **19** | **20** | **21** | **22** | **23** | **24** | **25** | **26** | **27** | **28** | **29** | **30** | **31** | **32** |
| B ≡ 1 (mod 32) | 513 | 545 | 577 | 609 | 641 | 673 | 705 | 737 | 769 | 801 | 833 | 865 | 897 | 929 | 961 | 993 |
| For B | (22,*35) | (11,*18) | (20,*32) | (12,*20) | (30,*48) | (5,*8) | (10,*16) | (6,*10) | (21,*34) | (7,*12) | (11,*18) | (19,*31) | (12,*20) | (5,*8) | (29,*46) | (11,*18) |
| **P** = 2^*r - *B_Con | 2^35-513 | 2^18-545 | 2^32-577 | 2^20-609 | 2^48-641 | | 2^16-705 | 287 | 2^34-769 | 3295 | 261311 | 2^31-865 | 2^20-897 | | 2^46-961 | 2^18-993 |
| B-loop | 17 (31,54) | 1 (17,36) | 17 (29,51) | 1 (15,33) | 17 (39,67) | 17 (15,29) | 5 (24,45) | 17 (13,26) | 17 (30,53) | 5 (8,20) | 5 (20,39) | 17 (28,50) | 1 (21,43) | 17 (21,39) | 17 (38,66) | 5 (18,36) |
| | | | | | | | | | | | | | | | | |
| (s_d,d) | (15_d,33) | (8_d,22) | (23_d,46) | (8_d,22) | (52_d,92) | (21_d,43) | (50_d,89) | (55_d,97) | (36_d,67) | (12_d,29) | (53_d,94) | (14_d,32) | (22_d,45) | (7_d,21) | (10_d,26) | (20_d,42) |
| B = 2^*r - **P** | 7653 | | 64933 | | 2^24-667 | | | 2^20-763 | 2^27-795 | | 64677 | | 101 | | | 3077 |
| (s,*r) | (8,*13) | (4,*7) | (10,*16) | (5,*8) | (15,*24) | (4,*7) | (5,*8) | (12,*20) | (17,*27) | (4,*7) | (10,*16) | (5,*8) | (6,*10) | (4,*7) | (5,*8) | (7,*12) |
| N ≡ 27 (mod 32) | 539 | 571 | 603 | 635 | 667 | 699 | 731 | 763 | 795 | 827 | 859 | 891 | 923 | 955 | 987 | 1019 |
| **X** | **17** | **18** | **19** | **20** | **21** | **22** | **23** | **24** | **25** | **26** | **27** | **28** | **29** | **30** | **31** | **32** |
| B ≡ 5 (mod 32) | 517 | 549 | 581 | 613 | 645 | 677 | 709 | 741 | 773 | 805 | 837 | 869 | 901 | 933 | 965 | 997 |
| For B | (6,*10) | (5,*8) | (4,*7) | (9,*15) | (5,*8) | (6,*10) | (4,*7) | (13,*21) | (11,*18) | (5,*8) | (4,*7) | (9,*15) | (5,*8) | (8,*13) | (4,*7) | (9,*15) |
| **P** = 2^*r - *B_Con | 507 | | | 32155 | | 347 | | 2^21-741 | 261371 | | | 31899 | | 7259 | | 31771 |
| B-loop | 1 (12,28) | 1 (17,36) | 1 (10,25) | 1 (15,33) | 17 (39,67) | 17 (44,75) | 5 (5,15) | 5 (22,42) | 17 (30,53) | 5 (8,20) | 5 (27,50) | 17 (28,50) | 5 (6,17) | 1 (14,32) | 17 (38,66) | 5 (11,25) |
| | | | | | | | | | | | | | | | | |
| (s_d,d) | (27_d,52) | (30_d,57) | (18_d,38) | (47_d,84) | (52_d,92) | (11_d,27) | (16_d,35) | (50_d,89) | (43_d,78) | (48_d,86) | (9_d,24) | (53_d,94) | (12_d,29) | (10_d,26) | (51_d,91) | (32_d,61) |
| B = 2^*r - **P** | | 64977 | | 7569 | | | 7473 | 2^21-751 | | 209 | | 3217 | | | 49 | 2^21-1007 |
| (s,*r) | (4,*7) | (10,*16) | (5,*8) | (8,*13) | (4,*7) | (5,*8) | (8,*13) | (13,*21) | (4,*7) | (6,*10) | (5,*8) | (7,*12) | (4,*7) | (5,*8) | (6,*10) | (13,*21) |
| N ≡ 15 (mod 32) | 527 | 559 | 591 | 623 | 655 | 687 | 719 | 751 | 783 | 815 | 847 | 879 | 911 | 943 | 975 | 1007 |
| **X** | **17** | **18** | **19** | **20** | **21** | **22** | **23** | **24** | **25** | **26** | **27** | **28** | **29** | **30** | **31** | **32** |
| B ≡ 17 (mod 32) | 529 | 561 | 593 | 625 | 657 | 689 | 721 | 753 | 785 | 817 | 849 | 881 | 913 | 945 | 977 | 1009 |
| For B | (9,*15) | (8,*13) | (5,*8) | (4,*7) | (6,*10) | (5,*8) | (28,*45) | (4,*7) | (15,*24) | (8,*13) | (5,*8) | (4,*7) | (11,*18) | (5,*8) | (7,*12) | (4,*7) |
| **P** = 2^*r - *B_Con | 32239 | 7631 | | | 367 | | 2^45-721 | | 2^24-785 | 7375 | | | 261231 | | 3119 | |
| B-loop | 5 (23,43) | 5 (9,21) | 5 (14,29) | 17 (10,21) | 17 (15,29) | 1 (13,30) | 17 (37,64) | 5 (17,34) | 5 (22,42) | 1 (16,35) | 5 (13,28) | 17 (16,31) | 1 (14,32) | 5 (6,17) | 17 (26,47) | 17 (14,28) |
| | | | | | | | | | | | | | | | | |
| (s_d,d) | (20_d,41) | (13_d,30) | (8_d,22) | (23_d,46) | (11_d,27) | (21_d,43) | (21_d,43) | (33_d,62) | (55_d,97) | (24_d,48) | (29_d,56) | (65_d,113) | (41_d,75) | (46_d,83) | (51_d,91) | (15_d,34) |
| B = 2^*r - **P** | | | 441 | 3481 | | 7513 | | 2^24-743 | | | 31929 | 2^35-871 | | 3161 | | 25 |
| (s,*r) | (4,*7) | (5,*8) | (6,*10) | (7,*12) | (4,*7) | (8,*13) | (5,*8) | (15,*24) | (4,*7) | (5,*8) | (9,*15) | (22,*35) | (4,*7) | (7,*12) | (5,*8) | (6,*10) |
| N ≡ 7 (mod 32) | 519 | 551 | 583 | 615 | 647 | 679 | 711 | 743 | 775 | 807 | 839 | 871 | 903 | 935 | 967 | 999 |
| **X** | **17** | **18** | **19** | **20** | **21** | **22** | **23** | **24** | **25** | **26** | **27** | **28** | **29** | **30** | **31** | **32** |
| B ≡ 25 (mod 32) | 537 | 569 | 601 | 633 | 665 | 697 | 729 | 761 | 793 | 825 | 857 | 889 | 921 | 953 | 985 | 1017 |
| For B | (7,*12) | (5,*8) | (6,*10) | (4,*7) | (10,*16) | (8,*13) | (5,*8) | (4,*7) | (8,*13) | (5,*8) | (11,*18) | (4,*7) | (20,*32) | (7,*12) | (5,*8) | (4,*7) |
| **P** = 2^*r - *B_Con | 3559 | | 423 | | 64871 | 7495 | | | 7399 | | 261287 | | 2^32-921 | 3143 | | |
| B-loop | 5 (9,21) | 1 (17,36) | 5 (7,18) | 17 (10,21) | 5 (12,26) | 17 (20,37) | 1 (18,38) | 17 (42,72) | 5 (22,42) | 17 (23,42) | 1 (21,43) | 5 (13,28) | 1 (26,51) | 1 (14,32) | 17 (14,28) | 5 (11,25) |

The entries in the above table matches **pairs** *$B_{Con}$ & **P**, as it is possible to find the **Converging* **IV**-classes $^*{}_u[2^{*r_s}\mathbf{X}\text{-}^*B_{Con}]_q$ by analysing **either** the N-values **or** the B-values, B ∉ {1,5,17}.

The infinity of **Converging* **IV**-classes for all **u** ∈ ℕ ***exist*** and the Union of $^*{}_u[2^{*r_s}\mathbf{X}\text{-}^*B_{Con}]_q$ contains **all** N ∈ ℕ. The Union of $^{\Delta}B_{Redundant} \equiv {}^*B_{Con} \pmod{2^{*r_s}}$ contains all B ∈ $\mathbb{N}_0$ ***except*** {1,5,17}.

## Examples of *reducing IV-classes & End-IV-classes for certain Threshold-values $N_T$:

We *disregard* $N_X \in *_{u=1}[2^{*1_0}\mathbf{X}\text{-}*0_{Con}]_{q=1} = \{\mathbf{2}_1,4,6,8, ...\}$ & $N_X \in *_{u=2}[2^{*2_1}\mathbf{X}\text{-}*3_{Con}]_{q=1} = \{\mathbf{1}_1,5,9,13, ...\}$
We ***focus*** on (*reducing*) **subsets** in ***The Domain-class*** $N_X \in {}^{\Phi}[2^2\mathbf{X}\text{-}1_0] = \{\mathbf{3}_1,7,11,15,19,23, ...\}$

| | |
|---|---|
| The Test-value 3 N in binary: 11 | The Test-value 3 reaches The End-value N = 1 in |
| is an Element in the *reducing IV-class *[*AX-*B] where: | #Odd + #Even = 2 + 5 = 7 Total-operations, |
| *A = 2^r = 16 *A in binary: 10000 | and to reach The Loop the parity-sequence then have a |
| *B = *A-*P = 13 *B in binary: 1101 | UNIQUE combination of 2 OE-tiles & 3 E-tiles |
| *P = 2^r-*B = 3 *P in binary: 0011 | 3 belongs to the END IV-class [2^dX - (2^d-N)] = |
| The Test-series [2^ 4 X - 13 ] — [3^ 2 X - 7 ] | [2^ 5 X - 29 ] |

Data for the Threshold-**value** $N_{j=0} = \mathbf{3}_T = \mathbf{P}$ in the **IV**-class $*_{u=3}[2^{*4_2}\mathbf{X}\text{-}13_{Con}]_{q=1} = \{\mathbf{3}_1, 19_2, 35_3, ...\}$
**NB**! *Brute force analysis* show, that the **lowest value** <u>not</u> in a $*Con_u(*\mathbf{r}_s)_q$ for $*\mathbf{r}_s \leq *4_2$ is $N_T = \mathbf{7}$

| | |
|---|---|
| The Test-value 7 N in binary: 111 | The Test-value 7 reaches The End-value N = 1 in |
| is an Element in the *reducing IV-class *[*AX-*B] where: | #Odd + #Even = 5 + 11 = 16 Total-operations, |
| *A = 2^r = 128 *A in binary: 10000000 | and to reach The Loop the parity-sequence then have a |
| *B = *A-*P = 121 *B in binary: 1111001 | UNIQUE combination of 5 OE-tiles & 6 E-tiles |
| *P = 2^r-*B = 7 *P in binary: 0000111 | 7 belongs to the END IV-class [2^dX - (2^d-N)] = |
| The Test-series [2^ 7 X - 121 ] — [3^ 4 X - 76 ] | [2^ 11 X - 2041 ] |

Data for the Threshold-**value** $N_{j=0} = \mathbf{7}_T = \mathbf{P}$ in the **IV**-class $*_{u=6}[2^{*7_4}\mathbf{X}\text{-}121_{Con}]_{q=1} = \{\mathbf{7}_1, 135_2, ...\}$
***Observe*** *that* $\mathbf{11}_1 \in *_{u=4}[2^{*5_3}\mathbf{X}\text{-}21_{Con}]_{q=1}$ & $\mathbf{23}_1 \in *_{u=5}[2^{*5_3}\mathbf{X}\text{-}9_{Con}]_{q=2}$ & $\mathbf{15}_1 \in *_{u=7}[2^{*7_4}\mathbf{X}\text{-}113_{Con}]_{q=2}$
**NB!** *Brute force analysis* show, that the **lowest value** <u>not</u> in a $*Con_u(*\mathbf{r}_s)_q$ for $*\mathbf{r}_s \leq *7_4$ is $N_T = \mathbf{27}$

| | |
|---|---|
| The Test-value 27 N in binary: 11011 | The Test-value 27 reaches The End-value N = 1 in |
| is an Element in the *reducing IV-class *[*AX-*B] where: | #Odd + #Even = 41 + 70 = 111 Total-operations, |
| *A = 2^r = 576460752303423488 | and to reach The Loop the parity-sequence then have a |
| *B = *A-*P = 576460752303423461 | UNIQUE combination of 41 OE-tiles & 29 E-tiles |
| *P = 2^r-*B = 27 | 27 belongs to the END IV-class [2^dX - (2^d-N)] = |
| [2^ 59 X-576460752303423461]---[3^ 37 X-450283905890997340] | [2^ 70 X - 1180591620717411303397] |

The Threshold-**value** $N_{j=0} = \mathbf{27}_T = \mathbf{P}$ in the **IV**-class $*_{u=?}[2^{*59_{37}}\mathbf{X}\text{-}(2^{*59}\text{-}\mathbf{27})_{Con}]_{q=1} = \{\mathbf{27}_1, (27+2^{*59})_2, ...\}$
**NB**! *Brute force analysis* show the **lowest value** <u>not</u> in a $*Con_u(*\mathbf{r}_s)_q$ for $*\mathbf{r}_s \leq *59_{37}$ is $N_T = \mathbf{703}$

| | |
|---|---|
| The Test-value 703 N in binary: 1010111111 | The Test-value 703 reaches The End-value N = 1 in |
| is an Element in the *reducing IV-class *[*AX-*B] where: | #Odd + #Even = 62 + 108 = 170 Total-operations, |
| *A = 2^r = 2417851639229258349412352 | and to reach The Loop the parity-sequence then have a |
| *B = *A-*P = 2417851639229258349411649 | UNIQUE combination of 62 OE-tiles & 46 E-tiles |
| *P = 2^r-*B = 703 | 703 belongs to the END IV-class [2^dX - (2^d-N)] = |
| [2^81X-2417851639229258349411649]--[3^51X-2153693963075557766310119] | [2^ 108 X - 324518553658426726783156020575553] |

Threshold-**value** $N_0 = \mathbf{703}_T = \mathbf{P}$ in the **IV**-class $*_{u=?}[2^{*81_{51}}\mathbf{X}\text{-}(2^{*81}\text{-}\mathbf{703})_{Con}]_{q=1} = \{\mathbf{703}_1, (703+2^{*81})_2, ...\}$
**NB**! *Brute force analysis* shows the **lowest** $N_T$ <u>not</u> in a $*Con_u(*\mathbf{r}_s)_q$ for $*\mathbf{r}_s \leq *81_{51}$ is $N_T = \mathbf{10087}$

| | |
|---|---|
| The Test-value 10087 N in binary: 10011101100111 | The Test-value 10087 reaches The End-value N = 1 in |
| is an Element in the *reducing IV-class *[*AX-*B] where: | #Odd + #Even = 81 + 142 = 223 Total-operations, |
| *A = 2^r = 40564819207303340847894502572032 | and to reach The Loop the parity-sequence then have a |
| *B = *A-*P = 40564819207303340847894502561945 | UNIQUE combination of 81 OE-tiles & 61 E-tiles |
| *P = 2^r-*B = 10087 | 10087 belongs to the END IV-class [2^dX - (2^d-N)] = |
| [2^105X-40564819207303340847894502561945]-[3^66X-30903154382632612361920641795841] | [2^ 142 X - 5575186299632655785383929568162090376485017] |

The **value** $N_0 = \mathbf{10087}_T = \mathbf{P}$ in the class $*_{u=?}[2^{*105_{66}}\mathbf{X}\text{-}(2^{*105}\text{-}\mathbf{10087})_{Con}]_{q=1} = \{\mathbf{10087}_1, (10087+2^{*105})_2, ...\}$
**NB**! *Brute force analysis shows, that <u>all</u>* $N_X < 2^{14}$ *belongs to a **specific*** $*_u[2^{*r_s}\mathbf{X}\text{-}B_{Con}]_q$, $*\mathbf{r}_s \leq *\mathbf{105}_{66}$

The following code/**Algorithm** can be used in *Thonny-Python* (Freeware):

**NB**! Be aware of line indentations if **copying** from pdf.

```
# This Algorithm uses the Rules in The Collatz Conjecture to test N-values
# Notice that in the print "*" is an Index and NOT "multiplication"
# All N congruent to 1 mod 4 are reduced in three operations, so
# the "interesting" Test-values are the N congruent to 3 mod 4

Reduce : int = 0
Stop : int = 0
j : int = 0
s : int = 0
r : int = 0
X : int = 0

Test = int(input("Please input an integer Test-value N:  "))
print(" ")
Start : int = Test


print(" R                                          Value in Series: ",Test,"    In binary:", str(bin(Test))[2:],)

if Test == 1 :
    print("")
    print("The Algorithm stops at N = 1")
    print("")
    print("Please try again for a Test-value N > 1")
    Reduce = 1
    Stop = 1

if Test == 2 :
    print("")
    print("N = 2 reaches N = 1 in one type Even operation")
    print("")
    print("Please try again for a Test-value N > 2")
    Reduce = 1
    Stop = 1

if Test % 2 == 0 :
    print("")
    print("The Test-value ", Start)
    print("is an Element in the *reducing IV-class *[2X-0]")
    print("")
    Reduce = 1

while Reduce == 0 :

    if Test % 2 == 1 :                  # Rule for Odd Values
        Test = int(3*Test + 1)          # If the two print-lines below are de-activated
        s = s + 1                        # the *reducing iteration-series is not printed
        j = j + 1                       # and only the results are shown
        print(" O  #Odd : ", s , "    #Even : ", r ,"    #Total : ", j , "    Value in Series: ", Test, "    In binary: " , str(bin(Test))[2:])

    if Test % 2 == 0 :                  # Rule for Even Values
        Test = int(Test / 2)
        r = r + 1
        j = j + 1
        print(" E  #Odd : ", s , "    #Even : ", r ,"    #Total : ", j , "    Value in Series: ", Test, "    In binary: " , str(bin(Test))[2:])
```

```
    if Test < Start :                      # Reduced-criterion
        print(" T      s = ",s,"          r = ", r,"         s+r = ", j, "       Reduced Value: ",Test )
        print("")
        print("The Test-value ", Start , "      N in binary: ", str(bin(Start))[2:])
        print("is an Element in the *reducing IV-class *[*AX-*B] where:")
        print(" *A = 2^r      = ", 2**r , "   *A in binary: " , str(bin((2**r)))[2:])
        print(" *B =   *A-*P = ", 2**r - (Start % (2**r)) , "   *B in binary: " , str(bin(2**r - (Start % (2**r))))[2:] )
        print(" *P = 2^r-*B = ", Start % (2**r), "    *P in binary: " , str(bin(Start % (2**r)))[2:])
        print("  The InitialValue-class [2^",r,"X -",2**r - (Start % (2**r)),"]")
        print("")
        print("  The Test-series for *[*AX-*B] Terminates at the TV-class [AX-B] where:")
        print("  A = 3^s      = ", 3**s, "    A in Binary: ", str(bin((3**s)))[2:] )
        print("  B =    A- P = ", 3**s - (Test % (3**s)), "     B in Binary: ", str(bin(3**s - (Test % (3**s))))[2:])
        print("  P = 3^s- B = ", Test % (3**s), "      P in Binary: ", str(bin(Test % (3**s)))[2:])
        print("  The TerminalValue-class [3^",s,"X -",3**s - (Test % (3**s)),"]")
        print("")
        X = int(X + int(Start + 2**r - (Start % (2**r)))/2**r)
        print("The Test-value  ", Start ,"  represents the X'th value in")
        print("the *reducing IV-class *[2^",r,"X -",2**r - (Start % (2**r)),"]" )
        print("where X = (N+*B)/2^r = ", X )
        print("and the *Reduced Value ",Test," represents the same X in")
        print("the TV-class [3^",s,"X -",3**s - (Test % (3**s)),"]")
        print("")
        print("The *reducing IV-class *[2^rX-*B], *P = (2^r-*B), ")
        print("meaning",Start,"is congruent to", Start % (2**r), "(mod 2^",r,")  ")
        print("share the *Reducing Combination  (s,r) = (",s,",",r,")")
        print("i.e. the *reducing parity-sequence *ROEOE...ET is a combination of ")
        print("s =   ",s," OE-tiles  and  (r-s) =   ", r-s," E-tiles")
        print("")
        Reduce = 1

while Stop == 0 :   # The Algorithm continues until N = 1

    if Test % 2 == 1 :                     # Rule for Odd Values
        Test = int(3*Test + 1)             # If the de-activated print-lines below are re-activated
        s = s + 1                           # the iteration-series from the *Reduced Value is printed
        j = j + 1
        #print(" O   #Odd : ", s , "     #Even : ", r ,"     #Total : ", j , "     Value in Series: ", Test, "    In binary: " , str(bin(Test))[2:])

    if Test % 2 == 0 :                     # Rule for Even Values
        Test = int(Test / 2)
        r = r + 1
        j = j + 1
        #print(" E   #Odd : ", s , "     #Even : ", r ,"     #Total : ", j , "     Value in Series: ", Test, "    In binary: " , str(bin(Test))[2:])

    if Test == 1 :                      # Stop-criterion
      # For The Loop
        print("")
        print("The Test-value  ", Start ,"  reaches The End-value N = 1 in ")
        print("#Odd + #Even = ",s,"+",r,"=",j,"  Total-operations,")
        print("and to reach The Loop the parity-sequence then have a ")
        print("UNIQUE combination of ",s," OE-tiles  & ",r-s," E-tiles")
        print(Start,"belongs to the END IV-class [2^dX - (2^d-N)] = [2^",r,"X -",2**r-Start,"]")
        Stop = 1
```

**NB**! As ***The Collatz Conjecture*** is **True** this is actually a "true algorithm" in the sense that for **ANY** input N ∈ ℕ the **Algorithm** <u>will</u> **terminate** and deliver a **result** (*but of course we need "a perfect computer" to analyse high* N *as a* normal *computer have limits for accuracy/precision*).

# The Algorithm to test B-values

```
# This Algorithm uses the Rules in The Collatz Conjecture to test B-values
# Even B are reduced in one operation, so the "interesting"
# Test-values are the odd B except {1,5,17} as these B results in B-loop

Reduce : int = 0
Stop : int = 0
j : int = 0
s : int = 0
r : int = 0
X : int = 0

Test = int(input("Please input an odd integer Test-value B: "))
print(" ")
Start : int = Test
print(" R                                    Value in Series: ",Test)

if Start == 1 :
   print("")
   print("The B-value 1 results in B-loop")
   print("")
   print("Please try again for an odd Test-value B > 1")
   Reduce = 1
   Stop = 1
if Start == 5 :
   print("")
   print("The B-value 5 results in B-loop")
   print("")
   print("Please try again for an odd Test-value B")
   Reduce = 1
   Stop = 1
if Start == 17 :
   print("")
   print("The B-value 17 results in B-loop")
   print("")
   print("Please try again for an odd Test-value B")
   Reduce = 1
   Stop = 1
if Start % 2 == 0 :
   print("")
   print("The Test-B ", Start)
   print("is an Even B belonging to the *reducing IV-class *[2X-0]")
   print("")
   Reduce = 1
   Stop = 1

while Reduce == 0 :

  if Test % 2 == 1 :          # Rule for Odd Values
    Test = int(3*Test - 1)    # If the two print-lines below are de-activated
    s = s + 1                 # the series is not printed
    j = j + 1                 # and only the results are shown
    print(" O #Odd : ", s , " #Even : ", r ," #Total : ", j , " Value in Series: ", Test)

  if Test % 2 == 0 :          # Rule for Even Values
    Test = int(Test / 2)
    r = r + 1
    j = j + 1
    print(" E #Odd : ", s , " #Even : ", r ," #Total : ", j , " Value in Series: ", Test)
```

```
    if Test < Start :      # Reduced-criterion
      print(" T s = ",s," r = ", r," s+r = ", j, " Reduced Value: ",Test )
      print("")
      print("The Test-B ", Start)
      print("is associated with the *reducing IV-class *[*AX-*B] where:")
      print(" *A = 2^r = ", 2**r )
      print(" *P = *A-*B = ", 2**r - (Start % (2**r))  )
      print(" *B = 2^r-*P = ", Start % (2**r))
      print(" The InitialValue-class [2^",r,"X -",Start % (2**r),"]")
      print("")
      print(" The Test-series for *[*AX-*B] Terminates at the TV-class [AX-B] where:")
      print(" A = 3^s = ", 3**s )
      print(" P = A- B = ", 3**s - (Test % (3**s)))
      print(" B = 3^s- P = ", Test % (3**s))
      print(" The TerminalValue-class [3^",s,"X -",Test % (3**s),"]")
      print("")
      X = int(X + int(Start + 2**r - (Start % (2**r)))/2**r)
      print("The Test-B ", Start ," represents a B associated with")
      print("the *reducing IV-class *[2^",r,"X -", Start % (2**r),"]" )
      print("where it is observed: X = (*P+Test-B)/2^r = ", X, "  (if X = 1 then Test-B = *B)" )
      print("")
      print("The *reducing IV-class *[2^rX-*B], *B = (2^r-*P) ")
      print("meaning",Start,"is congruent to", Start % (2**r), "(mod 2^",r,") ")
      print("share the *Reducing Combination (s,r) = (",s,",",r,") i.e. the Total s+r =",j)
      print("and the *reducing parity-sequence *ROEOE...ET is a combination of ")
      print("s = ",s," OE-tiles and (r-s) = ", r-s," E-tiles")
      print("")
      Reduce = 1

  while Stop == 0 :

    if Test == 1 :
      print("The Test-B",Start," ends in the B-loop for B = 1")
      print("")
      print("Number of Odd s =",s," Number of Even r =",r," Total s+r =",j)
      Stop = 1

    if Test == 5 :
      print("The Test-B",Start," ends in the B-loop for B = 5")
      print("")
      print("Number of Odd s =",s," Number of Even r =",r," Total s+r =",j)
      Stop = 1

    if Test == 17 :
      print("The Test-B",Start," ends in the B-loop for B = 17")
      print("")
      print("Number of Odd s =",s," Number of Even r =",r," Total s+r =",j)
      Stop = 1

    if Test % 2 == 1 :        # Rule for Odd Values
      Test = int(3*Test - 1)
      s = s + 1
      j = j + 1

    if Test % 2 == 0 :        # Rule for Even Values
      Test = int(Test / 2)
      r = r + 1
      j = j + 1
```